\documentclass[thmsa,a4paper,oneside,final,notitlepage,onecolumn,german]{article}
\usepackage{amssymb}

\usepackage{sw20jart}
\usepackage{epsfig}

\input{tcilatex}
\evensidemargin\oddsidemargin
\input tcilatex
\begin{document}

\author{\textbf{Dimitrios I. Dais and Martin Henk}}
\title{\textbf{On a series of Gorenstein cyclic quotient singularities admitting a
unique projective crepant resolution}}
\date{}
\maketitle

\noindent \textsc{Abstract.} Let $G$ be a finite subgroup of SL$\left( r,%
\Bbb{C}\right) $. In dimensions $r=2$ and $r=3$, McKay correspondence
provides a natural bijection between the set of irreducible representations
of $G$ and a cohomology-ring basis of the overlying space of a projective,
crepant desingularization of $\Bbb{C}^r/G$. For $r=2$ this desingularization
is unique and is known to be determined by the Hilbert scheme of the $G$%
-orbits. Similar statements (including a method of distinguishing just 
\textit{one} among all possible smooth minimal models of $\Bbb{C}^3/G$), are
very probably true for all $G$'s $\subset $ SL$\left( 3,\Bbb{C}\right) $
too, and recent Hilbert-scheme-techniques due to Ito, Nakamura and Reid, are
expected to lead to a new fascinating uniform theory. For dimensions $r\geq 4
$, however, to apply analogous techniques one needs extra modifications. In
addition, minimal models of $\Bbb{C}^r/G$ are smooth only under special
circumstances. $\Bbb{C}^4/\left( \text{involution}\right) $, for instance,
cannot have any smooth minimal model. On the other hand, all abelian
quotient spaces which are c.i.'s can always be fully resolved by
torus-equivariant, crepant, projective morphisms. Hence, from the very
beginning, the question whether a given Gorenstein quotient space $\Bbb{C}%
^r/G$, $r\geq 4$, admits special desingularizations of this kind, seems to
be absolutely crucial.\smallskip

\noindent In the present paper, after a brief introduction to the
existence-problem of such desingularizations (for abelian $G$'s) from the
point of view of toric geometry, we prove that the Gorenstein cyclic
quotient singularities of type 
\[
\frac 1l\,\left( 1,\ldots ,1,l-\left( r-1\right) \right) 
\]
with $l\geq r\geq 2$, have a \textit{unique }torus-equivariant projective,
crepant, partial resolution, which is ``full'' iff either $l\equiv 0$ mod $%
\left( r-1\right) $ or $l\equiv 1$ mod $\left( r-1\right) $. As it turns
out, if one of these two conditions is fulfilled, then the exceptional locus
of the full desingularization consists of $\QOVERD\lfloor \rfloor {l}{r-1}$
prime divisors, $\QOVERD\lfloor \rfloor {l}{r-1}-1$ of which are isomorphic
to the total spaces of $\Bbb{P}_{\Bbb{C}}^1$-bundles over $\Bbb{P}_{\Bbb{C}%
}^{r-2}$. Moreover, it is shown that intersection numbers are computable
explicitly and that the resolution morphism can be viewed as a composite of
successive (normalized) blow-ups. Obviously, the monoparametrized
singularity-series of the above type contains (as its ``first member'') the
well-known Gorenstein singularity defined by the origin of the affine cone
which lies over the $r$-tuple Veronese embedding of $\Bbb{P}_{\Bbb{C}}^{r-1}$%
.\newpage

\section{Introduction\label{Intro}}

\noindent \textsf{(a) }Let $f:Y\rightarrow X$ be a birational morphism
between two normal, $\mathbb{Q}$-Gorenstein complex varieties $X$ and $Y$ of
index $j$. Denote by $\omega _X=\mathcal{O}\left( K_X\right) $ and $\omega
_Y=\mathcal{O}\left( K_Y\right) $ the dualizing sheaves, and by $K_X$ and $%
K_Y$ representatives of canonical divisors of $X$ and $Y$, respectively. $f$
is called \textit{crepant }if $\omega _X^{\left[ j\right] }\cong f_{*}\left(
\omega _Y^{\otimes j}\right) $, or, in other words, if the \textit{%
discrepancy} $jK_Y-f^{*}\left( jK_X\right) $ vanishes. The ``prototype'' for
a crepant morphism is the proper birational map which desingularizes the
usual double-point-locus 
\[
X=\left\{ \left( z_1,z_2,z_3\right) \in \mathbb{C}^3\ \left| \
z_1^2+z_2^2+z_3^2=0\right. \right\} 
\]
by blowing up $\mathbf{0\in \,}X\subset \mathbb{C}^3$. Crepant birational
morphisms were mainly used in the past two decades in algebraic geometry to
reduce the singularities of complex $3$-folds (and, sometimes, $n$-folds) to
terminal (or even $\mathbb{Q}$-factorial terminal) singularities, and to
treat of \textit{minimal models} in high dimensions. For $X$ being the
underlying space of a Gorenstein quotient singularity, they are ``by
definition'' related to McKay-type correspondences.\medskip

\noindent \textsf{(b) }Let $r$ be an integer $\geq 2$, $G$ a finite subgroup
of GL$\left( r,\mathbb{C}\right) $ containing no pseudoreflections and
acting linearly on $\mathbb{C}^{r}$, and $p:\mathbb{C}^{r}\rightarrow 
\mathbb{C}^{r}/G$ the corresponding quotient map. The underlying space $%
\mathbb{C}^{r}/G$ of the (germ of the) quotient singularity $\left( \mathbb{C%
}^{r}/G,\left[ \mathbf{0}\right] \right) $, with $\left[ \mathbf{0}\right]
:=p\left( \mathbf{0}\right) $, is canonically equipped with the structure of
a normal, Cohen-Macaulay complex variety (or complex-analytic
space).\smallskip \newline
$\bullet $ The singular locus Sing$\left( \mathbb{C}^{r}/G\right) $ of $%
\mathbb{C}^{r}/G$ itself contains always $\left[ \mathbf{0}\right] $, but
for $r\geq 3$, it is possible to possess also other strata of $\mathbb{C}^{r}
$ of codimension $\geq 2$ passing through $\left[ \mathbf{0}\right] $ (cf. 
\ref{SLOC} below).\smallskip \newline
$\bullet $ As it was proved by Watanabe \cite{Wat1}, $\mathbb{C}^{r}/G$ is
Gorenstein iff $G\subset $ SL$\left( r,\mathbb{C}\right) $.\smallskip 
\newline
$\bullet $ If $r=2$, $G\subset $ GL$\left( 2,\mathbb{C}\right) $, the
quotient space $\mathbb{C}^{2}/G$ admits a unique \textit{minimal }%
desingularization 
\begin{equation}
f:\widehat{X}\rightarrow X=\mathbb{C}^{2}/G  \label{MIND}
\end{equation}
(``minimal'' in the sence that the exceptional locus of $f$ does not contain
any curve with self-intersection number $-1$, or equivalently, that there
exists, up to isomorphism, a unique morphism $h:\widetilde{X}\rightarrow 
\widehat{X}$ with $g=f\circ h$, for any desingularization $g:\widetilde{X}%
\rightarrow X$ of $X)$. The description of the prime divisors (rational
curves) consisting the exceptional locus of the above $f$, as well as that
of the way of how these divisors intersect each other (tree configurations),
is due to Hirzebruch \cite{Hirz2} (for cyclic acting groups) and Brieskorn 
\cite{Br2} (for all the other finite subgroups $G$ of GL$\left( 2,\mathbb{C}%
\right) $).\smallskip \newline
$\bullet $ The minimal desingularization (\ref{MIND}) is crepant if and only
if $G\subset $ SL$\left( 2,\mathbb{C}\right) $. In this special case, the
(Gorenstein) quotient spaces $\mathbb{C}^{2}/G$ are embeddable as $A$-$D$-$E$
hypersurfaces in $\mathbb{C}^{3}$ (Klein \cite{Kl}, Du Val \cite{DuVal1}, 
\cite{DuVal2}) and are nothing but the \textit{rational double points}
treated in the classical theory of the simple hypersurface singularities.
Table 1 contains all possible finite subgroups of SL$\left( 2,\mathbb{C}%
\right) $.\medskip 

\begin{center}
\begin{tabular}{|c||c||c|c|c|}
\hline
Nr. & \textbf{groups }$G$ & $
\begin{array}{c}
\text{Dynkin's} \\ 
\text{notation}
\end{array}
$ & $\left| G\right| $ & $\#\ \left\{ 
\begin{array}{c}
\text{conjugacy} \\ 
\text{classes of }G
\end{array}
\right\} $ \\ \hline\hline
1. & $
\begin{array}{c}
\text{cyclic groups }\mathbf{C}_{k} \\ 
\text{of order }k\geq 2
\end{array}
$ & $A_{k-1}$ & $k$ & $k$ \\ \hline
2. & $
\begin{array}{c}
\text{binary dihedral } \\ 
\text{groups \thinspace }\mathbf{D}_{k-2} \\ 
\text{with \ }k\geq 4
\end{array}
$ & $D_{k}$ & $4\left( k-2\right) $ & $k+1$ \\ \hline
3. & $
\begin{array}{c}
\text{binary tetrahedral} \\ 
\text{group\ \ }\mathbf{T}
\end{array}
$ & $E_{6}$ & $24$ & $7$ \\ \hline
4. & $
\begin{array}{c}
\text{binary octahedral} \\ 
\text{group\ \ }\mathbf{O}
\end{array}
$ & $E_{7}$ & $48$ & $8$ \\ \hline
5. & $
\begin{array}{c}
\text{binary icosahedral} \\ 
\text{group\ \ }\mathbf{I}
\end{array}
$ & $E_{8}$ & $120$ & $9$ \\ \hline
\end{tabular}
\medskip

\textbf{Table 1.\bigskip }
\end{center}

\noindent \noindent More precisely, taking into account the above group
classification, one obtains for the quotient spaces $\mathbb{C}^{2}/G,$ $%
G\subset $ SL$\left( 2,\mathbb{C}\right) $, and (\ref{MIND}):

\begin{theorem}
The quotient spaces $\mathbb{C}^{2}/G=$ \emph{Max-S}$\emph{pec}\left( 
\mathbb{C}\left[ \frak{x}_{1},\frak{x}_{2}\right] ^{G}\right) ,$ for $G$ a
finite subgroup of \emph{SL}$\left( 2,\mathbb{C}\right) $, are minimally
embedded as hypersurfaces 
\[
\left\{ \left( z_{1},z_{2},z_{3}\right) \in \mathbb{C}^{3}\ \left| \ h\left(
z_{1},z_{2},z_{3}\right) =0\right. \right\} 
\]
in $\mathbb{C}^{3}$, i.e. $\mathbb{C}\left[ \frak{x}_{1},\frak{x}_{2}\right]
^{G}\cong \mathbb{C}\left[ z_{1},z_{2},z_{3}\right] \ /\ \left( h\left(
z_{1},z_{2},z_{3}\right) \right) $. \emph{(The normal form for the ideal
generator in each individual case is mentioned in the third column of table
2.) }
\end{theorem}

\begin{theorem}
Every quotient singularity $\left( \mathbb{C}^{2}/G,\left[ \mathbf{0}\right]
\right) $, with $G$ a finite subgroup of \emph{SL}$\left( 2,\mathbb{C}%
\right) $, admits a unique minimal \emph{(}$=$ crepant\emph{)} resolution 
\[
f:\left( \widehat{\mathbb{C}^{2}/G},\mathbf{E}\right) \rightarrow \left( 
\mathbb{C}^{2}/G,\left[ \mathbf{0}\right] \right) 
\]
\emph{(}up to isomorphism\emph{)} with exceptional divisor $\mathbf{E}$
consisting of a configuration of rational smooth curves with
self-intersection number $-2$. The intersection form of $\mathbf{E}$ is
negative definite and therefore the Dynkin diagrams of type $A$-$D$-$E$ are
the dual graphs of the irreducible components of $\mathbf{E}$. \emph{(See
table 2; each ``}$\bullet $\emph{''} \emph{intimates a smooth rational curve
and each edge a transversal intersection at one point. For details, see e.g.
Lamotke \cite{Lamotke} and Slodowy \cite{Slodowy}).\medskip }
\end{theorem}

\begin{center}
$
\begin{tabular}{|c|c|c|c|}
\hline
Nr. & \textbf{\ }$G$ & $\ h\left( z_{1},z_{2},z_{3}\right) $ & Dynkin
diagram of the minimal resolution \\ \hline\hline
1. & $\mathbf{C}_{k}\ $ & $
\begin{array}{c}
\ \  \\ 
z_{1}^{2}+z_{2}^{2}+z_{3}^{k} \\ 
\ \ 
\end{array}
$ & $\underbrace{\bullet \!\!-\!\!\bullet \!\!-\!\!\bullet \!\!-\cdots -\!\!\bullet }_{\left( k-1\right) \text{ vertices}}$ \\ \hline
2. & $\mathbf{D}_{k-2}$ & $z_{1}^{2}+z_{2}^{2}\ z_{3}+z_{3}^{k-1}$ & $
\begin{array}{ccc}
\bullet  &  &  \\ 
\mid  &  &  \\ 
\bullet  & \!\!\!-\!\!\!\!\!\!\!\!- & \!\!\!\!\,\bullet \!\!-\!\!\bullet
\!\!-\cdots -\!\!\bullet  \\ 
\mid  &  &  \\ 
\bullet  &  & {}^{\text{(}k\text{ vertices)}}{\ }
\end{array}
$ \\ \hline
3. & $\mathbf{T}$ & $z_{1}^{2}+z_{2}^{3}+z_{3}^{4}$ & $
\begin{array}{ccccccccc}
\bullet  & \!\!-\!\!\!\,-\!\! & \bullet  & \!\!-\!\!\!\,-\!\! & \bullet  & 
-\!\!\!\,- & \bullet  & -\!\!\!\,- & \bullet  \\ 
&  &  &  & \mid  &  &  &  &  \\ 
&  &  &  & \bullet  &  &  &  & 
\end{array}
$ \\ \hline
4. & $\mathbf{O}$ & $z_{1}^{2}+z_{2}^{2}+z_{2}\ z_{3}^{3}$ & $
\begin{array}{ccccccccccc}
\bullet  & \!\!-\!\!\!\,-\!\! & \bullet  & \!\!-\!\!\!\,-\!\! & \bullet  & 
-\!\!\!\,- & \bullet  & -\!\!\!\,- & \bullet  & -\!\!\!\,- & \bullet  \\ 
&  &  &  & \mid  &  &  &  &  &  &  \\ 
&  &  &  & \bullet  &  &  &  &  &  & 
\end{array}
$ \\ \hline
5. & $\mathbf{I}$ & $z_{1}^{2}+z_{2}^{3}+z_{3}^{5}$ & $
\begin{array}{ccccccc}
\bullet  & \!\!-\!\!\!\,-\!\! & \bullet  & \!\!-\!\!\!\,-\!\! & \bullet  & 
-\!\!\!\,\cdots \,\cdots - & \bullet  \\ 
&  &  &  & \mid  & 
\begin{array}{c}
\text{{\scriptsize (by dots are meant here}} \\ 
\text{{\scriptsize 3 additional vertices)}}
\end{array}
&  \\ 
&  &  &  & \bullet  &  & 
\end{array}
$ \\ \hline
\end{tabular}
\medskip $

\textbf{Table 2.} \bigskip
\end{center}

\noindent $\bullet $ Denoting by $\frak{D}$ a simply-laced Dynkin diagram
belonging to the fourth column of Table 2, there is an ``extended'' Dynkin
diagram $\frak{D}^{\text{ext}}$ obtained by adding one vertex $\mathbf{v}_{0}
$ to $\frak{D}$ and by connecting the vertices, say $\mathbf{v}_{1},\ldots ,%
\mathbf{v}_{\nu }$, of $\frak{D}$ with $\mathbf{v}_{0}$ in the following way:
Let $\mathsf{R}$ be the root system corresponding to $\frak{D}$. The
vertices $\mathbf{v}_{1},\ldots ,\mathbf{v}_{\nu }$ give rise to simple
roots $\frak{r}_{1},\ldots ,\frak{r}_{\nu }$ with respect to the system 
\textsf{R}$_{+}$ of the positive roots. Furthermore, $\mathbf{v}_{0}$
corresponds to $\frak{r}_{0}=-\theta $, where $\theta $ is the ``longest''
root (i.e., $\theta $ is positive and $\theta +\frak{r}_{i}$ non-positive
for all $i$, $i\in \left\{ 1,..,\nu \right\} $). In $\frak{D}^{\text{ext}}$
one connects $\mathbf{v}_{0}$ and $\mathbf{v}_{i}$ by an \textit{edge} if
and only if $\left( \frak{r}_{0},\frak{r}_{i}\right) \neq 0$. In the late
seventies, McKay \cite{McK} observed a remarkable connection between the
representation theory of the finite subgroups of SL$\left( 2,\mathbb{C}%
\right) $ and the extended Dynkin diagrams $\frak{D}^{\text{ext}}$. To
formulate it explicitly, let us point out that each vertex $\mathbf{v}_{i}$
of $\frak{D}^{\text{ext}}$ is accompanied by a label $q_{i}\in \mathbb{N}$,
with  $q_{0}=1$ and $\sum_{i=0}^{\nu }\frak{r}_{i}\,q_{i}=0$, or
equivalently, $\theta =$ $\sum_{i=1}^{\nu }\frak{r}_{i}\,q_{i}$. 

\begin{theorem}[Classical McKay Correspondence ]
Let $G$ be a finite subgroup of \emph{SL}$\left( 2,\mathbb{C}\right) $ and $%
\frak{D}$ its Dynkin diagram. Then there is an one-to-one correspondence
\[
\left\{ \text{\emph{vertices of} \thinspace }\frak{D}^{\text{\emph{ext}}%
}\right\} \ni \mathbf{v}_{i}\longleftrightarrow \mathbf{\rho }_{i}\in
\left\{ 
\begin{array}{c}
\text{\emph{equivalence classes of}} \\ 
\text{\emph{irreducible representations of} \thinspace }G
\end{array}
\right\} 
\]
so that $\mathbf{\rho }_{i}$ has dimension $q_{i}$. Moreover, for any two
dimensional representation $\mathbf{\rho }^{\prime }$ of $G$ in $\mathbb{C}%
^{2}$, there exist isomorphisms for all $i$, $i\in \left\{ 0,1,...,\nu
\right\} ,$%
\[
\mathbf{\rho }_{i}\otimes \mathbf{\rho }^{\prime }\stackrel{\cong }{%
\longrightarrow }\left( \bigoplus\limits_{j\text{ \emph{incident to} }i}\,%
\mathbf{\rho }_{j}\right) \ .
\]
\end{theorem}

\noindent$\bullet $ Having this theorem as starting-point,
Gonzalez-Sprinberg, Verdier \cite{GSV}, and later Kn\"{o}rrer \cite{Kn},
constructed a purely geometric, direct correspondence ``of McKay-type''
between the set of irreducible representations of $G$, and the \textit{%
cohomology ring} of $\widehat{X}$ =$\,\widehat{\mathbb{C}^{2}/G}$ via
``tautological sheaves''. Recently, Ito, Nakamura \cite{Ito-Nak3}, and Reid 
\cite{Reid5}, introduced new techniques for the study of McKay
correspondence involving Hilbert schemes of $G$-orbits. \smallskip \newline
$\bullet $ In particular, for $r=2$, the main result of Ito and Nakamura 
\cite{Ito-Nak1}, \cite{Ito-Nak2}, \cite{Ito-Nak3}, \cite{Nak1}, \cite{Nak2}
can be roughly stated as follows:

\begin{theorem}
Let $G$ be a finite subgroup of \emph{SL}$\left( 2,\mathbb{C}\right) $ with $%
l=\left| G\right| $. Then there is a unique irreducible component $\mathcal{H%
}^{G}\left( \mathbb{C}^{2}\right) $of the $G$-fixed point set $\mathcal{H}%
^{l}\left[ \mathbb{C}^{2}\right] ^{G}$ of the Hilbert scheme $\mathcal{H}^{l}%
\left[ \mathbb{C}^{2}\right] $ parametrizing all clusters of length $l$ on $%
\mathbb{C}^{2}$, such that the induced proper birational morphism 
\[
\mathcal{H}^{G}\left( \mathbb{C}^{2}\right) \longrightarrow \mathbb{C}%
^{2}/G=X
\]
gives again the minimal resolution \emph{(\ref{MIND})} of $X$ \emph{(}up to
isomorphism\emph{).} Moreover, the original correspondence of \emph{\cite
{McK}, \cite{GSV}, \cite{Kn},} between the non-trivial representations of $G$
and the exceptional prime divisors of \emph{(\ref{MIND}) }can be
reinterpreted exclusively in terms of suitable ideals of $\mathcal{H}%
^{G}\left( \mathbb{C}^{2}\right) $.
\end{theorem}

\noindent $\bullet $ More generally, for arbitrary $r$, if $f:\widehat{X}%
\rightarrow X$ denotes a projective crepant (``full'') desingularization of $%
X=\mathbb{C}^r/G$, then the expected bijections are those of the following
box:\smallskip

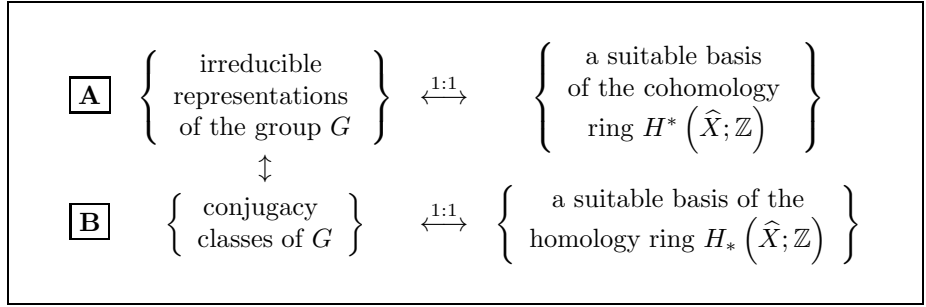
\[
\fbox{$
\begin{array}{ccc}
&  &  \\ 
& 
\begin{array}{cccc}
\fbox{\textbf{A}} & \left\{ 
\begin{array}{c}
\text{irreducible } \\ 
\text{representations } \\ 
\text{of the group }G
\end{array}
\right\}  & \stackrel{1:1}{\longleftrightarrow } & \left\{ 
\begin{array}{c}
\text{a suitable basis } \\ 
\text{of the cohomology } \\ 
\text{ring }H^{\ast }\left( \widehat{X};\mathbb{Z}\right) 
\end{array}
\right\}  \\ 
& \updownarrow  &  &  \\ 
\fbox{\textbf{B}} & \left\{ 
\begin{array}{c}
\text{conjugacy } \\ 
\text{classes of }G
\end{array}
\right\}  & \stackrel{1:1}{\longleftrightarrow } & \left\{ 
\begin{array}{c}
\text{a suitable basis of the} \\ 
\text{homology ring }H_{\ast }\left( \widehat{X};\mathbb{Z}\right) 
\end{array}
\right\} 
\end{array}
&  \\ 
&  & 
\end{array}
$}\smallskip 
\]
which is now known as \textit{Reid's slogan} : 
\[
\text{representation theory of \ }G\ \ \ \text{``}=\text{''\ }\ \text{\
homology of\ \ }\widehat{X}.
\]
This conjecture is accompanied by the remark that the above bijections
probably satisfy certain ``compatibilities'' (as for $r=2$) with respect to
the behaviour of the cup product, the image of the character table of $G$,
the duality interrelation etc.\smallskip \newline
$\bullet $ Although on the level of ``counting'' dimensions of \textit{%
rational }cohomology groups, or even on that of providing formal
correspondences between the left and the right hand side, the required
techniques are meanwhile well-understood (by toric methods \cite{BD} for $G$
abelian, and by recent results of Batyrev \cite{Bat2} involving
non-archimedian integrals, for arbitrary $G$'s), there are still lots of
open questions of how one might work with (co)homology groups the
coefficients of which are taken from $\mathbb{Z}$.\smallskip \newline
$\bullet $ Reid's approach to this generalized McKay-type-conjecture is
two-fold. The first idea (concerning correspondence \textbf{B}) relies on
the application of the following \textit{Ito-Reid theorem }\cite{Ito-Reid}
in order to construct a suitable collection of loci within $\widehat{X}$
(i.e., a collection of centers of monomial valuations on $\mathbb{C}\left(
X\right) $) generating $H_{\ast }\left( \widehat{X};\mathbb{Z}\right) $.

\begin{theorem}
Let $G$ be a finite subgroup of \emph{SL}$\left( r,\Bbb{C}\right) $ acting
linearly on $\Bbb{C}^r,\,r\geq 2$, and $X=\Bbb{C}^r/G$. Then there is a
canonical one-to-one correspondence between the junior conjugacy classes in $%
G$ and the crepant discrete valuations of $X$.
\end{theorem}

\noindent $\bullet $ The second idea (w.r.t. bijection \textbf{A}) is to
consider the tautological sheaves $\frak{F}_\rho $ assigned to the
irreducible representations $\rho $ of $G$. Reid \cite{Reid5} conjectures
(and proves for several examples) that appropriate $\Bbb{Z}$-linear
combinations of the Chern classes of $\frak{F}_\rho $'s lead to a canonical $%
\Bbb{Z}$-basis of the cohomology ring $H^{*}\left( \widehat{X};\Bbb{Z}%
\right) $. Moreover, if $\widehat{X}$ happens to be isomorphic to $\mathcal{H%
}^G\left( \Bbb{C}^r\right) $, then these sheaves enjoy very good
algebraic-geometric properties (they are generated by their global sections,
are vector bundles, their first Chern classes induce nef linear systems
etc., cf. \cite{Reid5}, 5.5). In particular, in this case $\mathcal{H}%
^G\left( \Bbb{C}^r\right) $ is birationally distinguished among all the
other projective crepant resolutions of $X$.\medskip

\noindent $\bullet $ As both \textbf{A} and \textbf{B} are clarified in
dimension $2$, let us recall what is  known for $r=3$. Next theorem is due
to Markushevich \cite{Mark}, \cite{Markush}, Ito \cite{Ito1}, \cite{Ito2}, 
\cite{Ito3}, and Roan \cite{Roan1}, \cite{Roan2}, \cite{Roan3}, \cite{Roan4}.

\begin{theorem}
The underlying spaces of all $3$-dimensional Gorenstein quotient
singularities possess crepant resolutions.
\end{theorem}

\begin{conjecture}
\label{3Hilb}Let $G$ be any finite subgroup of \emph{SL}$\left( 3,\Bbb{C}%
\right) $. Then $\mathcal{H}^G\left( \Bbb{C}^3\right) $ is a crepant
resolution of $\Bbb{C}^3/G$.
\end{conjecture}

\noindent For abelian $G$'s Conjecture \ref{3Hilb} was proved by Nakamura 
\cite{Nak3} (for an outline of the proof see also \cite{Reid5}, \S\ 7); work
on the non-abelian case is in progress. The complete verification of \ref
{3Hilb} would mean that for all quotients $\Bbb{C}^3/G$ there is always a
distinguished\footnote{\noindent In contrast to dimension $2$, in dimension $%
3$ minimal models are unique only up to isomorphisms in codimension $1$.
Moreover, there exist lots of examples of acting groups $G$, for which $%
\Bbb{C}^3/G$ has crepant, full, \textit{non-projective} resolutions (see
below \ref{NONPROJ}).} smooth minimal model available, satisfying all the
above mentioned peculiar properties.\medskip

\noindent $\bullet $ In dimensions $r\geq 4$, however, there are certain
additional troubles already from the very beginning and Reid's complementary
question (\cite{Reid4}, \cite{Ito-Reid}, \S\ 4.5, \cite{Reid5}, 5.4) still
remains an unanswered enigma:\smallskip

\noindent \underline{\textbf{Reid's question}}\textbf{\ : Under which
conditions on the acting groups} 
\[
G\subset \mathbf{SL}\left( r,\Bbb{C}\right) ,\ r\geq 4, 
\]
\textbf{do the quotient spaces }$\Bbb{C}^r/G$ \textbf{have projective
crepant desingularizations?\smallskip }\newline
$\bullet $ Note that the existence of \textit{terminal} Gorenstein
singularities implies automatically that \textit{not all} Gorenstein
quotient spaces $\Bbb{C}^r/G$, $r\geq 4$, can have such desingularizations
(cf. Morrison-Stevens \cite{Mo-Ste}).\smallskip \newline
$\bullet $ Moreover, in contrast to what is valid in the ``low'' dimensions $%
2$ and $3$, the Hilbert scheme $\mathcal{H}^G\left( \Bbb{C}^r\right) $ for $%
r\geq 4$ might be \textit{singular}, even if the quotient $\Bbb{C}^r/G$
being under consideration is known to possess projective, crepant
resolutions. (We are indebted to I.Nakamura and M.Reid for this information%
\footnote{%
Nakamura's typical counterexample is the so-called $\left( 4;2\right) $%
-hypersurface-singularity (in the terminology of \cite{DHZ1}), with \textit{%
non-smooth }$\mathcal{H}^G$. This singularity has projective crepant
resolutions (cf. \cite{DHZ1}, cor. 6.3, or \cite{Roan2}, \S\ 5).}%
).\smallskip \smallskip

\noindent $\bullet $ On the other hand, as it was proved in \cite{DHZ1} by
making use of Watanabe's classification \cite{Wat2} of all abelian quotient
singularities $\left( \Bbb{C}^r/G,\left[ \mathbf{0}\right] \right) $, $%
G\subset $ SL$\left( r,\Bbb{C}\right) $, (up to analytic isomorphism) whose
underlying spaces are embeddable as complete intersections (``c.i.'s'') of
hypersurfaces into an affine complex space, and methods of toric and
discrete geometry,

\begin{theorem}
The underlying spaces of all abelian quotient c.i.-singularities admit of
torus-equivariant projective, crepant resolutions \emph{(}and therefore 
\emph{smooth minimal models) }in all dimensions.
\end{theorem}

\noindent In particular, taking into account the specific structure of these
singularities depending on the \textit{free parameters }of the so-called
``Watanabe forests'', this theorem guarantees the existence of \textit{%
infinitely many }(isomorphism classes of) Gorenstein quotient singularities 
\textit{in each dimension} having resolutions with the required properties.
Nevertheless, these c.i.-singularities are of special nature, and they form
a relatively ``sparse'' subclass of the class of all Gorenstein abelian
quotient singularities. (For instance, all Gorenstein \textit{cyclic}
quotient msc-singularities in dimensions $\geq 3$ are \textit{not }c.i.'s!).
Thus, as ``next step'' it is natural to ask what happens with respect to the 
\textit{non-c.i.}'s. Various necessary existence-criteria working quite well
in the framework of this most general consideration, partially sufficient
conditions, and certain theoretical and algorithmic difficulties which arise
from LP-feasibility problems (cf. rem. \ref{FIRLA} below), as well as
further families of non-c.i. abelian quotient singularities for which it is
possible to apply a direct, constructive method to obtain the desired
resolutions, will be discussed in detail in \cite{DHZ2}. In the present
paper, we shall study another special (but, again, infinite) series of
Gorenstein, non-c.i. (for $r\geq 3$), cyclic quotient singularities
admitting a \textit{uniquely determined }torus-equivariant projective,
crepant resolution under very simple and absolutely well-controllable
(necessary \textit{and} sufficient) number-theoretic conditions. This
resolution will be defined as an immediate generalization of the most
well-known example in the literature, namely of the ``single blow-up'' of
the affine cone over the $l$-tuple Veronese embedding of $\Bbb{P}_{\Bbb{C}%
}^{r-1}$ at the origin. Let us first formulate it explicitly.

\begin{proposition}
\label{VER}Let $G$ be the finite cyclic group of analytic automorphisms of $%
\Bbb{C}^r$, $r\geq 2$, of order $l\geq 2$, generated by 
\begin{equation}
g:\Bbb{C}^r\ni \left( z_1,\ldots ,z_r\right) \mapsto \left( e^{\frac{2\pi 
\sqrt{-1}}l}\ z_1,\ldots ,\,e^{\frac{2\pi \sqrt{-1}}l}\ z_r\right) \in \Bbb{C%
}^r\smallskip \,  \label{ALL1}
\end{equation}
\emph{(i)} $\left( \Bbb{C}^r/G,\left[ \mathbf{0}\right] \right) $ is an
isolated singularity, and its underlying space is embedded into $\Bbb{C}^{%
\binom{r+l-1}l}$ as the zero-set\smallskip 
\[
\left\{ t_{\mathbf{u}}\,t_{\mathbf{v}}\,\cdots \,t_{\mathbf{w}}-t_{\mathbf{u}%
^{\prime }}\,t_{\mathbf{v}^{\prime }}\,\cdots \,t_{\mathbf{w}^{\prime }}=0\
\left| \ \text{\emph{sort}}\left( \mathbf{u\,v\,}\cdots \,\mathbf{w}\right) =%
\text{\emph{sort}}\left( \mathbf{u}^{\prime }\mathbf{\,v}^{\prime }\mathbf{\,%
}\cdots \,\mathbf{w}^{\prime }\right) \right. \right\} \ ,\smallskip 
\]
where $\mathbf{u},\mathbf{v},\ldots $ are defined as multiple index sets 
\[
\mathbf{u}=u_1\,u_{2\,}\,\ldots \,u_l=\,\stackunder{i_1\text{\emph{-times}}}{%
\underbrace{1\,1\,\ldots 1\,}}\,\ \stackunder{i_2\text{\emph{-times}}}{%
\underbrace{2\,2\,\ldots \,2}}\,\,\stackunder{i_3\text{\emph{-times}}}{\,%
\underbrace{3\,3\,\ldots \,3}}\,\,\ldots \,\,\,\stackunder{i_r\text{\emph{%
-times}}}{\underbrace{r\,r\,\ldots \,r}}
\]
with 
\[
i_1+i_2+\cdots +i_r=l\ \ \ \text{\emph{\& \ }}\ 0\leq i_j\leq l,\ \forall
j,\ 1\leq j\leq r,
\]
and \emph{sort}$\left( \,\cdot \,\right) $ denotes the sorting of any string
of the alphabet $\left\{ 1,2,3,\ldots ,r\right\} $ into weakly increasing
order.\smallskip \newline
\emph{(ii) }The singularity $\left( \Bbb{C}^r/G,\left[ \mathbf{0}\right]
\right) $ is 
\[
\text{ }\left\{ 
\begin{array}{ccc}
\text{\emph{terminal}} & \Longleftrightarrow  & r>l \\ 
&  &  \\ 
\text{\emph{canonical}} & \Longleftrightarrow  & r\geq l \\ 
&  &  \\ 
\text{\emph{Gorenstein}} & \Longleftrightarrow  & l\ \left| \ r\right. 
\end{array}
\right. 
\]
\smallskip \smallskip \newline
\emph{(iii) } If $\mathbf{Bl}_{\mathbf{0}}\left( \Bbb{C}^r\right) $ denotes
the \emph{(}usual\emph{) }blow-up of $\Bbb{C}^r$ at the origin, then the
action of $G$ on $\Bbb{C}^r$ can be extended onto $\mathbf{Bl}_{\mathbf{0}%
}\left( \Bbb{C}^r\right) $, 
\begin{equation}
\mathbf{Bl}_{\mathbf{0}}\left( \Bbb{C}^r\right) /G\rightarrow \Bbb{C}^r/G
\label{BLVER}
\end{equation}
is a \emph{(}full\emph{) }resolution of $\left[ \mathbf{0}\right] \in \Bbb{C}%
^r/G$ , the exceptional prime divisor $D$ is isomorphic to $\Bbb{P}_{\Bbb{C}%
}^{r-1}$, and the corresponding relative canonical divisor equals $\left(
\frac rl-1\right) D$.\smallskip \newline
\emph{(iv) (\ref{BLVER}) }is \emph{(the unique)} crepant resolution of $\Bbb{%
C}^r/G$ if and only if $r=l$ .$\Bbb{\smallskip }$
\end{proposition}

\noindent \textit{Proof. }(i) is an easy exercise (one has just to compute
the generators of the ring of invariants and their relations, cf. \cite{KKMS},
p. 40, and Sturmfels \cite{Sturm}, p. 141, for the explanation); (ii)
follows from the general theorems of Reid \cite{Reid1} (cf. (1.5), p. 277,
(3.1), p. 292) or directly from (iii). The construction of (\ref{BLVER}) is
due to Ueno \cite{Ueno2}, (see also \cite{Ueno3}, pp. 199-211), who called
it the ``canonical resolution'' of $\Bbb{C}^r/G$ and used it to obtain
generalized Kummer manifolds by resolving special quotient spaces defined by
discrete groups acting on complex tori. (iv) was pointed out by Hirzebruch
\& H\"{o}fer in \cite{HH}, p. 257, and follows from (iii). For the
uniqueness (up to isomorphism), see e.g. Roan's comments in \cite{Roan4},
Ex. 1, p. 135. $_{\Box }\medskip \smallskip $\newline
\textsf{(c) }For \textit{fixed} dimension $r\geq 2$, a ``$2$-parameter''
series of cyclic quotient singularities containing (\ref{ALL1}) of prop. \ref
{VER} as its ``first member'' ($\mu =1$) is that of \textit{type} \newline
\begin{equation}
\dfrac 1l\ \left( \stackunder{\left( r-1\right) \text{-times}}{\underbrace{%
1,1,\ldots ,1,1}},\ \mu \right) ,\,\,\mu \geq 1  \label{2PAR}
\end{equation}
(see \ref{type} for the definition and notation). A unique ``canonical
desingularization'' for each of its members was given by Ueno \cite{Ueno1}, 
\S\ 4, pp. 53-63, and Fujiki \cite{Fujiki}, pp. 316-318, in the case in
which the dimension equals $r=3$, and turned out to be very useful for the
characterization of exceptional fibers belonging to smooth threefolds
fibered over a curve and having normally polarized abelian surfaces as
generic fibers. (This was, in fact, a direct generalization of Hirzebruch's
``continued fraction algorithm'' \cite{Hirz2}, pp. 15-20, to the next coming
dimension). Since we are mainly interested in Gorenstein singularities and
in the existence (or non-existence) of projective, crepant, full resolutions
in dimensions $r\geq 4$, we shall consider for (\ref{2PAR}) only the case
where $\mu =l-\left( r-1\right) $ (cf. \ref{Gor-prop}, \ref{wnor} below) and
consequently $l$ as the only available parameter. Our main motivation to
work out this series was a recent remark of Reid in \cite{Reid5}, 5.4;
namely, that the four-dimensional cyclic singularities of type $\frac
1l\left( 1,1,1,l-3\right) $ with gcd$\left( l,3\right) =1$ are to be
resolved by a crepant morphism if and only if $l\equiv 1$ mod 3$.$ We give a
generalization in all dimensions, even without assuming the isolatedness of
the corresponding singularity, show the uniqueness and projectivity of the
crepant morphism, describe the exceptional prime divisors and their
intersection numbers, and compute the cohomology dimensions of the
desingularized space.\medskip \newline
Exactly as in \cite{DHZ1}, we shall exclusively work with the machinery of
toric geometry. More precisely, the paper has the following structuring : In
section $2$ we give the toric glossary which will be used in the sequel.
(The reader who is familiar with this matter may skip it). Sections $3$ and $%
4$ are complementary. (In fact, the reason for adding in \S\ $3$ some
lengthy explanations is that there is a potential for confusion between the
usual blow-up of a toric subvariety $V\left( \tau \right) $ of an $X\left(
N,\Delta \right) $ and the starring subdivision w.r.t. $\tau $. These are
identical only for smooth $X\left( N,\Delta \right) $'s! On the other hand,
to lend an algebraic-geometric characterization to even very simple
combinatorially motivated cone subdivisions, it is absolutely natural to
blow-up also not necessarily reduced subschemes). In \S\ $4$ we deal with a
high dimensional analogue of the so-called Hirzebruch-surfaces and make
certain remarks concerning its embeddings and intersection theory. (It turns
out that all but one exceptional prime divisors which will arise later on in
our desingularizations are of this sort). Sections $5$-$6$ outline a first
systematic approach to the general problem of the existence or non-existence
of crepant (preferably projective) resolutions of Gorenstein abelian
quotient singularities of dimension $\geq 4$ by basic (and coherent)
triangulations of the junior simplex. In \S\ $7$ we take a closer look at
the low dimensions. Sections $8$ and $9$ contain our main results. Though
the singularity-series which we study is rather special, we hope at least
that it will become clear how one may apply our techniques to more demanding
singularities. In particular, in section\ $9$, the factorization of the
desingularizing morphism is reduced to a ``game'' with the available
simplices. Finally, in \S\ $10$ we give a foretaste of what may be done for
the series generalizing $\frac 17\left( 1,2,4\right) $ and state the
GPSS-conjecture.\bigskip \newline
\textit{Terminology and general notation}. By a \textit{complex variety} is
meant an integral, separated algebraic scheme over $\Bbb{C}$. A complex
variety is therefore an irreducible, reduced ringed space $\left( X,\mathcal{%
O}_X\right) $ with structure sheaf $\mathcal{O}_X$ which is locally
determined by the canonical structure sheaf of the spectrum of an affine
complex coordinate ring. Sing$\left( X\right) $ denotes the \textit{singular
locus} of $X$, i.e., the set of all points $x\in X$ with $\mathcal{O}_{X,x}$
a non-regular local ring. Analogously, $x\in X$ is normal, Cohen-Macaulay,
Gorenstein etc., if $\mathcal{O}_{X,x}$ is of this type. A \textit{subvariety%
} $Y$ of $X$ is a closed integral subscheme of $X$. If codim$_X\left(
Y\right) =1$, then $Y$ is especially called a \textit{prime divisor}. By CDiv%
$\left( X\right) $, WDiv$\left( X\right) $, Pic$\left( X\right) $ and $%
A_{\bullet }\left( X\right) =\oplus _{k\geq 0}A_k\left( X\right) $ we denote
the groups of Cartier and Weil divisors, the Picard group, and the graded
Chow ring of $X$, respectively. (For $X$ smooth, $A^{\bullet }\left(
X\right) =\oplus _{k\geq 0}A^k\left( X\right) $, with $A^k\left( X\right)
=A_{\text{dim}X-k}\left( X\right) $). Just as in \cite{BD}, \cite{DHZ1}, by
a \textit{desingularization} (or \textit{resolution of singularities}) $f:%
\widehat{X}\rightarrow X$ of a non-smooth $X$, we mean a ``full'' or
``overall'' desingularization (if not mentioned), i.e., Sing$\left( \widehat{%
X}\right) =\varnothing $. When we deal with \textit{partial}
desingularizations, we mention it explicitly. A birational morphism $%
f:X^{\prime }\rightarrow X$ is \textit{projective} if $X^{\prime }$ admits
an $f$-ample Cartier divisor. The \textit{intersection numbers} of Cartier
divisors are defined as in \cite{Ful} (see below \S\ \ref{TORIC} \textsf{(i)}%
).

\section{ Preliminaries from toric geometry\label{TORIC}}

\noindent We recall some basic facts from the theory of toric varieties and
fix the notation which will be used in the sequel. For details the reader is
referred to the standard textbooks of Oda \cite{Oda}, Fulton \cite{Fulton},
and Ewald \cite{Ewald}, and to the lecture notes \cite{KKMS}. \medskip

\noindent \textsf{(a)} The \textit{linear hull, }the\textit{\ affine hull},
the \textit{positive hull} and \textit{the convex hull} of a set $B$ of
vectors of $\Bbb{R}^r$, $r\geq 1,$ will be denoted by lin$\left( B\right) $,
aff$\left( B\right) $, pos$\left( B\right) $ (or $\Bbb{R}_{\geq 0}\,B$) and
conv$\left( B\right) $ respectively. The \textit{dimension} dim$\left(
B\right) $ of a $B\subset \Bbb{R}^r$ is defined to be the dimension of its
affine hull. \newline
\newline
\textsf{(b) }\textit{Polyhedral cones. }Let $N\cong \Bbb{Z}^r$ be a free $%
\Bbb{Z-}$module of rank $r\geq 1$. $N$ can be regarded as a \textit{lattice }%
in $N_{\Bbb{R}}:=N\otimes _{\Bbb{Z}}\Bbb{R}\cong \Bbb{R}^r$. (For fixed
identification, we shall represent the elements of $N_{\Bbb{R}}$ by
column-vectors in $\Bbb{R}^r$). If $\left\{ n_1,\ldots ,n_r\right\} $ is a $%
\Bbb{Z}$-basis of $N$, then 
\[
\text{det}\left( N\right) :=\left| \text{det}\left( n_1,\ldots ,n_r\right)
\right| 
\]
is the \textit{lattice determinant}. An $n\in N$ is called \textit{primitive}
if conv$\left( \left\{ \mathbf{0},n\right\} \right) \cap N$ contains no
other points except $\mathbf{0}$ and $n$.\smallskip

Let $N\cong \Bbb{Z}^r$ be as above, $M:=$ Hom$_{\Bbb{Z}}\left( N,\Bbb{Z}%
\right) $ its dual lattice, $N_{\Bbb{R}},M_{\Bbb{R}}$ their real scalar
extensions, and $\left\langle .,.\right\rangle :N_{\Bbb{R}}\times M_{\Bbb{R}%
}\rightarrow \Bbb{R}$ the natural $\Bbb{R}$-bilinear pairing. (For fixed
identification $M_{\Bbb{R}}\cong \Bbb{R}^r$, we analogously represent the
elements of $M_{\Bbb{R}}$ by row-vectors in $\Bbb{R}^r$). A subset $\sigma $
of $N_{\Bbb{R}}$ is called \textit{strongly convex polyhedral cone }(\textit{%
s.c.p.c.}, for short), if there exist $n_1,\ldots ,n_k\in N_{\Bbb{R}}$, such
that $\sigma =$ pos$\left( \left\{ n_1,\ldots ,n_k\right\} \right) $ and $%
\sigma \cap \left( -\sigma \right) =\left\{ \mathbf{0}\right\} $. Its 
\textit{relative interior }int$\left( \sigma \right) $ (resp. its \textit{%
relative boundary }$\partial \sigma $) is the usual topological interior
(resp. the usual topological boundary) of it, considered as subset of lin$%
\left( \sigma \right) $. The \textit{dual cone} of $\sigma $ is defined by 
\[
\sigma ^{\vee }:=\left\{ \mathbf{x}\in M_{\Bbb{R}}\ \left| \ \left\langle 
\mathbf{x},\mathbf{y}\right\rangle \geq 0,\ \forall \mathbf{y},\ \mathbf{y}%
\in \sigma \right. \right\} 
\]
and satisfies: $\sigma ^{\vee }+\left( -\sigma ^{\vee }\right) =M_{\Bbb{R}}$
and dim$\left( \sigma ^{\vee }\right) =r$. \ A subset $\tau $ of a s.c.p.
cone $\sigma $ is called a \textit{face} of $\sigma $ (notation: $\tau \prec
\sigma $), if $\tau =\left\{ \mathbf{y}\in \sigma \ \left| \ \left\langle
m_0,\mathbf{y}\right\rangle =0\right. \right\} $, for some $m_0\in \sigma
^{\vee }$. A s.c.p.c $\sigma =$ pos$\left( \left\{ n_1,\ldots ,n_k\right\}
\right) $ is called \textit{simplicial} (resp. \textit{rational}) if $%
n_1,\ldots ,n_k$ are $\Bbb{R}$-linearly independent (resp. if $n_1,\ldots
,n_k\in N_{\Bbb{Q}}$, where $N_{\Bbb{Q}}:=N\otimes _{\Bbb{Z}}\Bbb{Q}$).%
\newline
\newline
\textsf{(c) }\textit{Monoids.} If $\sigma \subset N_{\Bbb{R}}$ is a rational
s.c.p. cone, then $\sigma $ has $\mathbf{0}$ as its apex and the
subsemigroup $\sigma \cap N$ of $N$ is a monoid. The following two
propositions describe the fundamental properties of this monoid $\sigma \cap
N$ and their proofs go essentially back to Gordan \cite{Gordan}, Hilbert 
\cite{Hilbert} and van der Corput \cite{vanderCorput1}, \cite{vanderCorput2}.

\begin{proposition}[Gordan's lemma]
\label{Gorlem}$\sigma \cap N$ is finitely generated as additive semigroup,
i.e. there exist 
\[
\ n_1,n_2,\ldots ,n_\nu \in \sigma \cap N\ ,\text{ \ such that\ \ \ }\sigma
\cap N=\Bbb{Z}_{\geq 0}\ n_1+\Bbb{Z}_{\geq 0}\ n_2+\cdots +\Bbb{Z}_{\geq 0}\
n_\nu \ .
\]
\end{proposition}

\begin{proposition}[Minimal generating system]
\label{MINGS}Among all systems of generators of $\sigma \cap N$, there
is a system $\mathbf{Hlb}_N\left( \sigma \right) $ of \emph{minimal
cardinality}, which is uniquely determined \emph{(}up to the ordering of its
elements\emph{)} by the following characterization \emph{:\smallskip } 
\begin{equation}
\fbox{$\mathbf{Hlb}_N\left( \sigma \right) =\left\{ n\in \sigma \cap \left(
N\smallsetminus \left\{ \mathbf{0}\right\} \right) \ \left| \ 
\begin{array}{c}
n\ \text{\emph{cannot be expressed }} \\ 
\text{\emph{as the sum of two other }} \\ 
\text{\emph{\ vectors belonging }} \\ 
\text{\emph{\ to\ } }\sigma \cap \left( N\smallsetminus \left\{ \mathbf{0}%
\right\} \right) 
\end{array}
\right. \right\} $}  \label{Hilbbasis}
\end{equation}
\end{proposition}

\noindent \textit{Proof. }See e.g. Schrijver \cite{Schrijver}, p. 233. $%
_{\Box }\medskip $

\begin{definition}
\emph{\ }$\mathbf{Hlb}_N\left( \sigma \right) $ \emph{is often called }%
\textit{the Hilbert basis of }$\sigma $ w.r.t. $N.$
\end{definition}

\noindent About \textit{algorithms} for the determination of Hibert bases of
pointed rational cones, we refer to Pottier \cite{Pottier1}, \cite{Pottier2}%
, Sturmfels \cite{Sturm} (13.2, p. 128), and Henk-Weismantel \cite{HEW}, and
to the other references therein.\smallskip \medskip

\noindent \textsf{(d) }\textit{Algebraic tori defined via }$N$. Let $\Bbb{C}%
^{*}$ be the multiplicative group of non-zero complex numbers. For $N\cong 
\Bbb{Z}^r$ we define an $r$-dimensional \textit{algebraic torus }$T_N\cong
\left( \Bbb{C}^{*}\right) ^r$ by : 
\[
T_N:=\text{Hom}_{\Bbb{Z}}\left( M,\Bbb{C}^{*}\right) =N\otimes _{\Bbb{Z}}%
\Bbb{C}^{*}. 
\]
Every $m\in M$ assigns a character $\mathbf{e}\left( m\right)
:T_N\rightarrow \Bbb{C}^{*}$ with 
\[
\mathbf{e}\left( m\right) \left( t\right) :=t\left( m\right) ,\ \ \forall
t,\ \ t\in T_N. 
\]
We have : 
\[
\mathbf{e}\left( m+m^{\prime }\right) =\mathbf{e}\left( m\right) \cdot 
\mathbf{e}\left( m^{\prime }\right) \ ,\ \ \text{for\ \ }m,m^{\prime }\in M%
\text{ , \ and \ }\mathbf{e}\left( \mathbf{0}\right) =1\text{ .} 
\]
Moreover, for each $n\in N$, we define an $1$-parameter subgroup 
\[
\gamma _n:\Bbb{C}^{*}\rightarrow T_N\ \ \ \text{with\ \ \ }\gamma _n\left(
\lambda \right) \left( m\right) :=\lambda ^{\left\langle m,n\right\rangle }%
\text{, \ \ for\ \ \ }\lambda \in \Bbb{C}^{*},\ m\in M,\ \ 
\]
\ ($\gamma _{n+n^{\prime }}=\gamma _n\circ $\ $\gamma _{n^{\prime }}$, \ \
for\ \ \ $n,n^{\prime }\in N$). We can therefore identify $M$ with the
character group of $T_N$ and $N$ with the group of $1$-parameter subgroups
of $T_N$. If $\left\{ n_1,\ldots ,n_r\right\} $ is a $\Bbb{Z}$-basis of $N$
and $\left\{ m_1,\ldots ,m_r\right\} $ the dual basis (of $M$) and if we set 
$u_j:=\mathbf{e}\left( m_j\right) $, $\forall j$, $1\leq j\leq r$, then
there exists an isomorphism 
\[
T_N\ni t\stackrel{\cong }{\longmapsto }\left( u_1\left( t\right) ,\ldots
,u_r\left( t\right) \right) \in \left( \Bbb{C}^{*}\right) ^r\ , 
\]
and $\left\{ u_1,\ldots ,u_r\right\} $ plays the role of a coordinate system
of $T_N$. Hence, to an 
\[
m=\sum_{j=1}^ra_jm_j\in M\ \ \ \ \text{( resp. to an }n=\sum_{j=1}^rb_jn_j%
\text{ )} 
\]
we associate the character (``Laurent monomial'') $\ \ \mathbf{e}\left(
m\right) =u_1^{a_1}\cdot u_1^{a_2}\cdots \ u_r^{a_r}\smallskip $ 
\[
\ \ \ \ \text{( resp. the 1-parameter subgroup }\gamma _n:\Bbb{C}^{*}\ni
\lambda \mapsto \left( \lambda ^{b_1},\ldots ,\lambda ^{b_r}\right) \in
\left( \Bbb{C}^{*}\right) ^r\text{ )\ .\smallskip } 
\]
On the other hand, for a rational s.c.p.c. $\sigma $ with 
\[
M\cap \sigma ^{\vee }=\Bbb{Z}_{\geq 0}\ m_1+\Bbb{Z}_{\geq 0}\ m_2+\cdots +%
\Bbb{Z}_{\geq 0}\ m_k, 
\]
we associate to the finitely generated, normal, monoidal $\Bbb{C}$%
-subalgebra $\Bbb{C}\left[ M\cap \sigma ^{\vee }\right] $ of $\Bbb{C}\left[
M\right] $ an affine complex variety 
\[
U_\sigma :=\text{Max-Spec}\left( \Bbb{C}\left[ M\cap \sigma ^{\vee }\right]
\right) , 
\]
which can be identified with the set of semigroup homomorphisms : 
\[
U_\sigma =\left\{ u:M\cap \sigma ^{\vee }\ \rightarrow \Bbb{C\ }\left| 
\begin{array}{c}
\ u\left( \mathbf{0}\right) =1,\ u\left( m+m^{\prime }\right) =u\left(
m\right) \cdot u\left( m^{\prime }\right) ,\smallskip \  \\ 
\text{for all \ \ }m,m^{\prime }\in M\cap \sigma ^{\vee }
\end{array}
\right. \right\} \ , 
\]
where $\mathbf{e}\left( m\right) \left( u\right) :=u\left( m\right) ,\
\forall m,\ m\in M\cap \sigma ^{\vee }\ $ and\ $\forall u,\ u\in U_\sigma $ .

\begin{proposition}[Embedding by binomials]
\label{EMB}In the analytic category, $U_\sigma $, identified with its image
under the injective map $\left( \mathbf{e}\left( m_1\right) ,\ldots ,\mathbf{%
e}\left( m_k\right) \right) :U_\sigma \hookrightarrow \Bbb{C}^k$, can be
regarded as an analytic set determined by a system of equations of the form%
\emph{:} \emph{(monomial) = (monomial).} This analytic structure induced on $%
U_\sigma $ is independent of the semigroup generators $\left\{ m_1,\ldots
,m_k\right\} $ and each map $\mathbf{e}\left( m\right) $ on $U_\sigma $ is
holomorphic w.r.t. it. In particular, for $\tau \prec \sigma $, $U_\tau $ is
an open subset of $U_\sigma $. Moreover, if $\#\left( \mathbf{Hlb}_M\left(
\sigma ^{\vee }\right) \right) =d\ \left( \leq k\right) $, then $d$ is
nothing but the \emph{embedding dimension} of $U_\sigma $, i.e. the \emph{%
minimal} number of generators of the maximal ideal of the local $\Bbb{C}$%
-algebra $\mathcal{O}_{U_\sigma ,\ \left( \mathbf{0}\in \Bbb{C}^d\right) }$.
\end{proposition}

\noindent \textit{Proof. }See Oda \cite{Oda} prop. 1.2 and 1.3., pp. 4-7. $%
_{\Box }\bigskip $

\noindent \textsf{(e) }\textit{Fans.} A \textit{fan }w.r.t.\textit{\ }$%
N\cong \Bbb{Z}^r$ is a finite collection $\Delta $ of rational s.c.p. cones
in $N_{\Bbb{R}}$, such that :\smallskip \newline
(i) any face $\tau $ of $\sigma \in \Delta $ belongs to $\Delta $,
and\smallskip \newline
(ii) for $\sigma _1,\sigma _2\in \Delta $, the intersection $\sigma _1\cap
\sigma _2$ is a face of both $\sigma _1$ and $\sigma _2.\smallskip $\newline
The union $\left| \Delta \right| :=\cup \left\{ \sigma \ \left| \ \sigma \in
\Delta \right. \right\} $ is called the \textit{support }of $\Delta $.
Furthermore, we define 
\[
\Delta \left( i\right) :=\left\{ \sigma \in \Delta \ \left| \ \text{dim}%
\left( \sigma \right) =i\right. \right\} \ ,\ \text{for \ \ }0\leq i\leq r\
. 
\]
If $\varrho \in \Delta \left( 1\right) $, then there exists a unique
primitive vector $n\left( \varrho \right) \in N\cap \varrho $ with $\varrho =%
\Bbb{R}_{\geq 0}\ n\left( \varrho \right) $ and each cone $\sigma \in \Delta 
$ can be therefore written as 
\[
\sigma =\sum\Sb \varrho \in \Delta \left( 1\right)  \\ \varrho \prec \sigma 
\endSb \ \Bbb{R}_{\geq 0}\ n\left( \varrho \right) \ \ . 
\]
The set Gen$\left( \sigma \right) :=\left\{ n\left( \varrho \right) \ \left|
\ \varrho \in \Delta \left( 1\right) ,\varrho \prec \sigma \right. \right\} $
is called the\textit{\ set of minimal generators }(within the pure first
skeleton) of $\sigma $. For $\Delta $ itself one defines analogously 
\[
\text{Gen}\left( \Delta \right) :=\bigcup_{\sigma \in \Delta }\text{ Gen}%
\left( \sigma \right) \ .\medskip 
\]
\textsf{(f) }\textit{Toric varieties, orbits and stars.} The \textit{toric
variety X}$\left( N,\Delta \right) $ associated to a fan\textit{\ }$\Delta $
w.r.t. the lattice\textit{\ }$N$ is by definition the identification space 
\[
X\left( N,\Delta \right) :=\left( \left( \bigcup_{\sigma \in \Delta }\
U_\sigma \right) \ /\ \sim \right) 
\]
with $U_{\sigma _1}\ni u_1\sim u_2\in U_{\sigma _2}$ if and only if there is
a $\tau \in \Delta ,$ such that $\tau \prec \sigma _1\cap \sigma _2$ and $%
u_1=u_2$ within $U_\tau $ (cf. lemma \ref{GLSCH} below). As complex variety, 
$X\left( N,\Delta \right) $ turns out to be irreducible, normal,
Cohen-Macaulay and to have at most rational singularities (cf. \cite{Fulton}%
, p. 76, and \cite{Oda}, thm. 1.4, p. 7, and cor. 3.9, p. 125). $X\left(
N,\Delta \right) $ is called \textit{simplicial }if all cones of $\Delta $
are simplicial.\medskip \newline
$\bullet $ $X\left( N,\Delta \right) $ admits a canonical $T_N$-action which
extends the group multiplication of $T_N=U_{\left\{ \mathbf{0}\right\} }$ : 
\begin{equation}
T_N\times X\left( N,\Delta \right) \ni \left( t,u\right) \longmapsto t\cdot
u\in X\left( N,\Delta \right)  \label{torus action}
\end{equation}
where, for $u\in U_\sigma $, $\left( t\cdot u\right) \left( m\right)
:=t\left( m\right) \cdot u\left( m\right) ,\ \forall m,\ m\in M\cap \sigma
^{\vee }$ . The orbits w.r.t. the action (\ref{torus action}) are
parametrized by the set of all the cones belonging to $\Delta $. For a $\tau
\in \Delta $, we denote by orb$\left( \tau \right) $ (resp. by $V\left( \tau
\right) $) the orbit (resp. the closure of the orbit) which is associated to 
$\tau $. The spaces orb$\left( \tau \right) $ and $V\left( \tau \right) $
have the following properties (cf. \cite{Fulton}, pp. 52-55, \cite{Oda}, \S\
1.3) :\smallskip \newline
(i) For $\tau \in \Delta $, it is 
\[
V\left( \tau \right) =\coprod \left\{ \text{orb}\left( \sigma \right) \
\left| \ \sigma \in \Delta ,\tau \prec \sigma \right. \right\} 
\]
and 
\[
\text{orb}\left( \tau \right) =V\left( \tau \right) \smallsetminus \bigcup \
\left\{ V\left( \sigma \right) \ \left| \ \tau \precneqq \sigma \right.
\right\} \ .\smallskip 
\]
\newline
(ii) If $\tau \in \Delta $, then $V\left( \tau \right) =X\left( N\left( \tau
\right) ,\text{ Star}\left( \tau ;\Delta \right) \right) $ is itself a toric
variety w.r.t. 
\[
N\left( \tau \right) :=N/N_\tau \ ,\ \ \ N_\tau :=N\cap \text{lin}\left( 
\Bbb{\tau }\right) \Bbb{\ },\ \ \ \text{Star}\left( \tau ;\Delta \right)
:=\left\{ \overline{\sigma }\ \left| \ \sigma \in \Delta ,\ \tau \prec
\sigma \right. \right\} \ , 
\]
where $\overline{\sigma }\ =\left( \sigma +\left( N_\tau \right) _{\Bbb{R}%
}\right) /\left( N_\tau \right) _{\Bbb{R}}$ denotes the image of $\sigma $
in $N\left( \tau \right) _{\Bbb{R}}=N_{\Bbb{R}}/\left( N_\tau \right) _{\Bbb{%
R}}$.\smallskip \smallskip \newline
(iii) For $\tau \in \Delta $, the closure $V\left( \tau \right) $ is
equipped with an affine open covering 
\[
\left\{ U_\sigma \left( \tau \right) \ \left| \ \tau \prec \sigma \right.
\right\} 
\]
consisting of ``intermediate'' subvarieties 
\[
U_\tau \left( \tau \right) =\text{ orb}\left( \tau \right) \hookrightarrow
U_\sigma \left( \tau \right) \hookrightarrow U_\sigma 
\]
being defined by : $U_\sigma \left( \tau \right) :=$ Max-Spec$\left( \Bbb{C}%
\left[ \overline{\sigma }^{\vee }\cap \ M\left( \tau \right) \right] \right) 
$, with $M\left( \tau \right) $ denoting the dual of $N\left( \tau \right) $%
.\bigskip \newline
\textsf{(g)} \textit{Smoothness and compactness criterion. }Let $N\cong \Bbb{%
Z}^r$ be a lattice of rank $r$ and $\sigma \subset N_{\Bbb{R}}$ a
simplicial, rational s.c.p.c. of dimension $k\leq r$. $\sigma $ can be
obviously written as $\sigma =\varrho _1+\cdots +\varrho _k$, for distinct $%
1 $-dimensional cones $\varrho _1,\ldots ,\varrho _k$. We denote by 
\[
\mathbf{Par}\left( \sigma \right) :=\left\{ \mathbf{y}\in \left( N_\sigma
\right) _{\Bbb{R}}\ \left| \ \mathbf{y}=\sum_{j=1}^k\ \varepsilon _j\
n\left( \varrho _j\right) ,\ \text{with\ \ }0\leq \varepsilon _j<1,\ \forall
j,\ 1\leq j\leq k\right. \right\} 
\]
the \textit{fundamental }(\textit{half-open})\textit{\ parallelotope }which
is associated to\textit{\ }$\sigma $. The \textit{multiplicity} mult$\left(
\sigma ;N\right) $ of $\sigma $ with respect to $N$ is defined as 
\[
\text{mult}\left( \sigma ;N\right) :=\#\left( \mathbf{Par}\left( \sigma
\right) \cap N_\sigma \right) =\text{Vol}\left( \mathbf{Par}\left( \sigma
\right) ;N_\sigma \right) \ , 
\]
where Vol$\left( \mathbf{Par}\left( \sigma \right) \right) $ denotes the
usual volume of $\mathbf{Par}\left( \sigma \right) $ and 
\[
\text{Vol}\left( \mathbf{Par}\left( \sigma \right) ;N_\sigma \right) :=\frac{%
\text{Vol}\left( \mathbf{Par}\left( \sigma \right) \right) }{\text{det}%
\left( N_\sigma \right) } 
\]
its the relative volume w.r.t. $N_\sigma $.

\begin{proposition}
\label{SMCR}The affine toric variety $U_\sigma $ is smooth iff \emph{mult}$%
\left( \sigma ;N\right) =1$. \emph{(}Correspondingly, an arbitrary toric
variety $X\left( N,\Delta \right) $ is smooth if and only if it is
simplicial and each s.c.p. cone $\sigma \in \Delta $ satisfies this
condition.\emph{)}
\end{proposition}

\noindent \textit{Proof. }It follows from \cite{Oda}, thm. 1.10, p. 15. $%
_{\Box }\medskip $

\noindent $\bullet $ For the systematic study of toric singularities it is
useful to introduce the notion of the ``splitting codimension'' of the
closed point orb$\left( \sigma \right) $ of an $U_\sigma $. For the germ $%
\left( U_\sigma ,\text{orb}\left( \sigma \right) \right) $ of an affine $r$%
-dimensional toric variety w.r.t. a \textit{singular point }orb$\left(
\sigma \right) $, the \textit{splitting codimension} splcod$\left( \text{orb}%
\left( \sigma \right) ;U_\sigma \right) $ of orb$\left( \sigma \right) $ in $%
U_\sigma $ is defined as : 
\[
\text{splcod}\left( \text{orb}\left( \sigma \right) ;U_\sigma \right) :=%
\text{max\emph{\ }}\left\{ \varkappa \in \left\{ 2,\ldots ,r\right\} \
\left| \ 
\begin{array}{c}
U_\sigma \cong U_{\sigma ^{\prime }}\times \Bbb{C}^{r-\varkappa },\ \text{%
s.t.} \\ 
\text{dim}\left( \sigma ^{\prime }\right) =\varkappa \text{ \ and} \\ 
\text{Sing}\left( U_{\sigma ^{\prime }}\right) \neq \varnothing
\end{array}
\right. \right\} 
\]
If splcod$\left( \text{orb}\left( \sigma \right) ;U_\sigma \right) =r$, then
orb$\left( \sigma \right) $\emph{\ }will be called an \textit{msc-singularity%
}, i.e. a singularity having the maximum splitting codimension.\smallskip 
\newline
\noindent $\bullet $ Next theorem gives a necessary and sufficient condition
for $X\left( N,\Delta \right) $ to be compact.

\begin{theorem}
\label{COMPACT}A toric variety $X\left( N,\Delta \right) $ is compact if and
only if $\Delta $ is a \emph{complete fan}, i.e., $\left| \Delta \right| =N_{%
\Bbb{R}}$.
\end{theorem}

\noindent \textit{Proof. }See Oda \cite{Oda}, thm. 1.11, p. 16. $_{\Box
}\bigskip $

\noindent \textsf{(h) }\textit{Order functions, support functions and
divisors. }If $X\left( N,\Delta \right) $ is a (not necessarily compact)
toric variety associated to a fan $\Delta $ w.r.t. a lattice $N$, $M$ its
dual, and $\mathbf{\iota }:T_N\hookrightarrow X\left( N,\Delta \right) $ the
canonical inclusion, then $\mathbf{\iota }_{*}\left( \mathcal{O}%
_{T_N}\right) $ is a $T_N$-invariant quasi-coherent sheaf of $\mathcal{O}%
_{X\left( N,\Delta \right) }$-modules canonically embedded into the constant
sheaf $\Bbb{C}\left( X\left( N,\Delta \right) \right) $ of rational
functions of $X\left( N,\Delta \right) $. Let $\mathcal{F}\neq \,$\underline{%
$0$} be a $T_N$-invariant coherent sheaf of fractional ideals over $X\left(
N,\Delta \right) $ contained in $\mathbf{\iota }_{*}\left( \mathcal{O}%
_{T_N}\right) $. Fix a s.c.p. cone $\sigma \in \Delta $, $n\in N\cap \sigma
, $ and consider the corresponding $1$-parameter group\footnote{%
Here we drop the prefix Max- because we do not work only with closed points
but also with closed subsets, and the $\Bbb{C}$-scheme structure is
essential for the arguments.}: 
\[
\gamma _n:\text{ Spec}\left( \Bbb{C}\left[ w,w^{-1}\right] \right) =\Bbb{C}%
^{*}\longrightarrow T_N=\text{Spec}\left( \Bbb{C}\left[ M\right] \right) 
\]
Since $\exists $ $\stackunder{\lambda \rightarrow 0}{\text{lim}}\,\gamma
_n\left( \lambda \right) \in U_\sigma $ (cf. \cite{Oda}, 1.6.(v), p. 10), $%
\gamma _n$ is extendable to a map 
\[
\overline{\gamma _n}:\text{ Spec}\left( \Bbb{C}\left[ w\right] \right) =\Bbb{%
A}_{\Bbb{C}}^1\longrightarrow U_\sigma =\text{Spec}\left( \Bbb{C}\left[
M\cap \sigma ^{\vee }\right] \right) \ . 
\]
The coherence of $\mathcal{F}$ implies that $\mathcal{F}\left| _{U_\sigma
}\right. $ is of type $J_\sigma ^{\,\sim }$ (cf. \cite{Hart}, pp. 110-111),
with $J_\sigma $ being an $M$-graded complex vector-subspace of 
\[
H^0\left( T_N,\mathcal{O}_{T_N}\right) =\Bbb{C}\left[ M\right]
=\bigoplus_{m\in M}\,\Bbb{C\,}\mathbf{e}\left( m\right) 
\]
on the one hand, and a finitely generated $\Bbb{C}\left[ M\cap \sigma ^{\vee
}\right] $-module 
\[
J_\sigma =\sum_{j=1}^{q_\sigma }\,\Bbb{C}\left[ M\cap \sigma ^{\vee }\right]
\cdot \mathbf{e}\left( m_j^{\left( \sigma \right) }\right) ,\,\text{ for
some \ }m_1^{\left( \sigma \right) },\ldots ,m_{q_\sigma }^{\left( \sigma
\right) }\in M\text{, and\ \ \thinspace }q_\sigma \in \Bbb{N\,}\text{,} 
\]
on the other. The pullback $\overline{\gamma _n}^{\,*}\mathcal{F}\left|
_{U_\sigma }\right. $ is realized via the finitely generated $\Bbb{C}\left[
w\right] $-module\smallskip 
\[
\sum_{j=1}^{q_\sigma }\,\Bbb{C}\left[ w\right] \cdot w^{\left\langle
m_j^{\left( \sigma \right) },\,n\right\rangle }\ \subset \Bbb{C}\left(
w\right) . 
\]
Define the \textit{order-function }w.r.t. $\mathcal{F}$ by 
\[
\text{ord}_{\mathcal{F}}\left( n\right) :=\text{ inf\thinspace }\left\{
\left\langle m_j^{\left( \sigma \right) },n\right\rangle \ \left| \ 1\leq
j\leq q_\sigma \right. \right\} \in \Bbb{Z\ }.\smallskip 
\]
This ord$_{\mathcal{F}}\left( n\right) $ is exactly the image ord$_0\left( 
\overline{\gamma _n}^{\,*}\mathcal{F}\left| _{U_\sigma }\right. \right) $
under the usual order function 
\[
\text{ord}_0:\Bbb{C}\left( w\right) ^{*}\rightarrow \Bbb{Z} 
\]
of the discrete valuation ring $\mathcal{O}_{\Bbb{A}_{\Bbb{C}}^1,0}$ with $w$
as uniformizing parameter. Since the above definition depends only on $%
\sigma $, one extends ord$_{\mathcal{F}}$ to the entire $\left| \Delta
\right| $ by setting\smallskip 
\[
\text{ord}_{\mathcal{F}}\left( \mathbf{y}\right) :=\text{ inf\thinspace }%
\left\{ \left\langle m_j^{\left( \sigma \right) },\mathbf{y}\right\rangle \
\left| \ 1\leq j\leq q_\sigma \right. \right\} ,\ \forall \mathbf{y},\ \ 
\mathbf{y}\in \sigma ,\ \ \sigma \in \Delta \ .\smallskip 
\]
The order function ord$_{\mathcal{F}}$ is $\Bbb{R}$-valued and has the
following characteristic properties:\smallskip \newline
(i) \ \ it is \textit{positively homogeneous}, i.e. ord$_{\mathcal{F}}\left(
c\,\mathbf{y}\right) =c$ ord$_{\mathcal{F}}\left( \mathbf{y}\right) $, for
all $c\in \Bbb{R}_{\geq 0}\,$,\smallskip \newline
(ii) \ ord$_{\mathcal{F}}\left| _\sigma \right. $ is piecewise linear on
each $\sigma \in \Delta $,\smallskip \newline
(iii) ord$_{\mathcal{F}}\left( N\cap \left| \Delta \right| \right) \subset 
\Bbb{Z}$, and\smallskip \newline
(iv) for all $\sigma \in \Delta $, ord$_{\mathcal{F}}\left| _\sigma \right. $
is \textit{upper convex}, i.e., 
\[
\text{ord}_{\mathcal{F}}\left| _\sigma \right. \left( \mathbf{y+y}^{\prime
}\right) \geq \text{ord}_{\mathcal{F}}\left| _\sigma \right. \left( \mathbf{y%
}\right) +\text{ord}_{\mathcal{F}}\left| _\sigma \right. \left( \mathbf{y}%
^{\prime }\right) ,\ \ \text{for any pair\ \ }\mathbf{y},\mathbf{y}^{\prime
}\in \sigma . 
\]

\begin{definition}
\emph{Let }$X\left( N,\Delta \right) $\emph{\ be a toric variety. A function 
}$\psi :\left| \Delta \right| \rightarrow \Bbb{R}$\emph{\ is called} \textit{%
integral PL-support function }\emph{if it satisfies the above properties
(i)-(iv). We define } 
\[
\begin{array}{c}
\text{\emph{PL-SF}}\left( N,\Delta \right) :=\left\{ 
\begin{array}{c}
\text{\emph{all integral PL-support\smallskip }} \\ 
\text{\emph{functions defined on} }\left| \Delta \right| 
\end{array}
\right\} \smallskip 
\end{array}
\]
\emph{and the sets of }\textit{integral \emph{(}}$\Delta $-\textit{linear%
\emph{) }support functions}\emph{\footnote{%
An integral $\Delta $-linear support function $\psi $ is called \textit{%
strictly upper convex }if it is upper convex on $\left| \Delta \right| $ and
if for any two distinct maximal-dimensional cones $\sigma $ and $\sigma
^{\prime }$, the linear functions $m_\sigma $, $m_{\sigma ^{\prime }}\in M=$
Hom$_{\Bbb{Z}}\left( N,\Bbb{Z}\right) \subset $ Hom$_{\Bbb{R}}\left( N_\Bbb{R},\Bbb{R%
}\right) =M_{\Bbb{R}}$ defining $\psi \left| _\sigma \right. =\left\langle
m_\sigma ,\bullet \right\rangle $ and $\psi \left| _{\sigma ^{\prime
}}\right. =\left\langle m_{\sigma ^{\prime }},\bullet \right\rangle $ are
different. }:\smallskip } 
\[
\begin{array}{c}
\text{\emph{SF}}\left( N,\Delta \right) :=\left\{ \psi \in \text{\emph{PL-SF}%
}\left( N,\Delta \right) \ \left| \ \psi \left| _\sigma \right. \ \text{%
\emph{linear on each} \ }\sigma \in \Delta \text{\ }\right. \right\}
\smallskip  \\ 
\bigcup \smallskip  \\ 
\text{\emph{UCSF}}\left( N,\Delta \right) :=\left\{ \psi \in \text{\emph{SF}}%
\left( N,\Delta \right) \ \left| 
\begin{array}{c}
\ \psi \text{ \ \emph{upper convex\smallskip }} \\ 
\text{\emph{on the whole }}\left| \Delta \right| 
\end{array}
\right. \right\} \smallskip  \\ 
\bigcup \smallskip  \\ 
\text{\emph{SUCSF}}\left( N,\Delta \right) :=\left\{ \psi \in \text{\emph{%
UCSF}}\left( N,\Delta \right) \ \left| 
\begin{array}{c}
\ \psi \text{ \ \emph{strictly upper }} \\ 
\text{\emph{convex on }}\left| \Delta \right| 
\end{array}
\right. \right\} .\smallskip 
\end{array}
\]
\emph{Analogously, one defines the sets} \emph{SF}$_{\Bbb{Q}}\left( N,\Delta
\right) $\emph{, UCSF}$_{\Bbb{Q}}\left( N,\Delta \right) $\emph{, SUCSF}$_{%
\Bbb{Q}}\left( N,\Delta \right) $ \emph{of }\textit{rational }\emph{support
functions by modifying property (iii) into : }$\psi \left( N_{\Bbb{Q}}\cap
\left| \Delta \right| \right) \subset \Bbb{Q}$\emph{. All the above sets are
equipped with the usual additive group structure.}
\end{definition}

\begin{theorem}
For a function $\psi \in $ \emph{PL-SF}$\left( N,\Delta \right) $, and an
arbitrary cone $\sigma \in \Delta $ define 
\[
\left( J_\psi \right) _\sigma :=\bigoplus_{m\in M}\,\left\{ \Bbb{C\,}\mathbf{%
e}\left( m\right) \ \left| \ \left\langle m,\mathbf{y}\right\rangle \geq
\psi \left( \mathbf{y}\right) ,\ \forall \mathbf{y},\ \mathbf{y}\in \sigma
\right. \right\} .
\]
The family of $T_N$-invariant sheaves $\left\{ \left( J_\psi \right) _\sigma
^{\sim }\ \left| \ \sigma \in \Delta \right. \right\} $ being associated to
the family of ideals $\left\{ \left( J_\psi \right) _\sigma \ \left| \
\sigma \in \Delta \right. \right\} $ can be glued together \emph{(cf. \cite
{Hart}, Ex.II.1.22, p. 69) }to construct a coherent sheaf $\mathcal{F}_\psi $
of $T_N$-invariant fractional ideals over $X\left( N,\Delta \right) $
contained in $\mathbf{\iota }_{*}\left( \mathcal{O}_{T_N}\right) $.
Moreover,\smallskip \ \newline
\emph{(i) ord}$_{\mathcal{F}_\psi }=\psi ,\smallskip $\newline
\emph{(ii) }$\mathcal{F}_{\emph{ord}_{\mathcal{F}}}$ is the \emph{completion}
of $\mathcal{F}$ \emph{(in Zariski's sence);\smallskip }\newline
\emph{(iii) }mapping\emph{\ \ }$\mathcal{F}\longmapsto $ \emph{ord}$_{%
\mathcal{F}}$ \ and $\ \psi \longmapsto \mathcal{F}_\psi $ \ one obtains a
bijection\smallskip 
\[
\text{\emph{PL-SF}}\left( N,\Delta \right) \stackunder{\text{\emph{1:1}}}{%
\longleftrightarrow }\left\{ 
\begin{array}{c}
\text{\emph{coherent sheaves of} }T_N\text{\emph{-invariant} } \\ 
\text{\emph{complete fractional ideals over} \thinspace }X\left( N,\Delta
\right) 
\end{array}
\right\} \ ,\smallskip 
\]
\emph{(iv) }$\mathcal{F}\subset \mathcal{F}_\psi \Leftrightarrow $ \emph{ord}%
$_{\mathcal{F}}\geq \psi ,\ \ $\emph{ord}$_{\mathcal{F}_1\cdot \mathcal{F}%
_2}=$ \emph{ord}$_{\mathcal{F}_1}+$ \emph{ord}$_{\mathcal{F}_2}$,
and\smallskip \newline
\emph{(v) }$\mathcal{F}_{\psi _1}\cong \mathcal{F}_{\psi _2}$ \emph{(}as $%
\mathcal{O}_{X\left( N,\Delta \right) }$-module sheaves\emph{)} $%
\Leftrightarrow \psi _1-\psi _2$ is linear.
\end{theorem}

\noindent \textit{Proof. }See Saint-Donat \cite{KKMS}, ch. I, \S\ 2, thm. 9,
pp. 28-31. $_{\Box }$

\begin{definition}
\emph{Let }$X\left( N,\Delta \right) $\emph{\ be a toric variety and }$%
\sigma \in \Delta $\emph{. The} \textit{convex interpolation} $\psi
_\vartheta $ \emph{of a function} $\vartheta :$ \emph{Gen}$\left( \sigma
\right) \rightarrow \Bbb{Z}$ \emph{is defined by\smallskip } 
\[
\sigma \ni \mathbf{y\longmapsto \,}\psi _\vartheta \left( \mathbf{y}\right)
:=\text{\emph{inf\thinspace }}\left\{ \left\langle m,\mathbf{y}\right\rangle
\ \left| 
\begin{array}{c}
\ m\in M\text{ \ \emph{and} \ }\left\langle m,\mathbf{y}\right\rangle \geq
\vartheta \left( n\left( \varrho \right) \right) ,\smallskip \  \\ 
\forall \varrho ,\ \ \varrho \in \Delta \left( 1\right) ,\text{ \ \emph{with}
}\varrho \prec \sigma 
\end{array}
\right. \right\} \in \Bbb{R}\,.\smallskip 
\]
\emph{Correspondingly, by the} \textit{convex interpolation} \emph{of a real
function} $\vartheta :$ \emph{Gen}$\left( \Delta \right) \rightarrow \Bbb{Z}$
\emph{is meant a function} $\psi _\vartheta :\left| \Delta \right|
\rightarrow \Bbb{R}$\emph{, such that} $\psi _\vartheta \left| _\sigma
\right. $ \emph{is the convex interpolation of }$\vartheta \left| _\sigma
\right. $\emph{\ (in the above sence) for all }$\sigma \in \Delta $\emph{.
Obviously, such a }$\psi _\vartheta $\emph{\ belongs to PL-SF}$\left(
N,\Delta \right) $\emph{, and conversely each integral PL-support function
has this form.\smallskip }
\end{definition}

\noindent $\bullet $ Let $X$ be any $r$-dimensional \textit{normal} complex
variety and $D\in $ WDiv$\left( X\right) $. The correspondence\footnote{%
\textit{Reflexive} coherent sheaves $\mathcal{F}$ are those which are
isomorphic to their biduals $\mathcal{F}^{\vee \,\vee }$ (where $\mathcal{F}%
^{\vee \,}:=$ \textit{Hom}$_{\mathcal{O}_X}\left( \mathcal{F},\mathcal{O}%
_X\right) $). For $\mathcal{F}\subset \Bbb{C}\left( X\right) $ of rank one
they are also called \textit{divisorial.}}:\smallskip 
\[
A_{r-1}\left( X\right) \ni \left\{ D\right\} \stackrel{\mathbf{\delta }}{%
\longmapsto }\left\{ \mathcal{O}_X\left( D\right) \right\} \in \left\{ 
\begin{array}{c}
\text{\textit{reflexive} coherent } \\ 
\text{(i.p.~torsion-free) sheaves } \\ 
\text{of fractional ideals } \\ 
\text{over }X\text{ having rank one}
\end{array}
\right\} \,/\,H^0\left( X,\mathcal{O}_X^{*}\right) \smallskip 
\]
with $\mathcal{O}_X\left( D\right) $ defined by sending every non-empty open
subset $U$ of $X$ onto 
\[
U\longmapsto \mathcal{O}_X\left( D\right) \left( U\right) :=\left\{ \varphi
\in \Bbb{C}\left( X\right) ^{*}\,\left| \,\left( \text{div}\left( \varphi
\right) +D\right) \left| _U\right. \geq 0\right. \right\} , 
\]
induces a $\Bbb{Z}$-module isomorphism (cf. Reid \cite{Reid1}, App. to \S\
1); in fact, to avoid torsion, one defines this $\Bbb{Z}$-module structure
by setting: 
\[
\mathbf{\delta }\left( D_1+D_2\right) :=\left( \mathcal{O}_X\left(
D_1\right) \otimes \mathcal{O}_X\left( D_2\right) \right) ^{\vee \,\vee },\ 
\mathbf{\delta }\left( \kappa \,D\right) :=\mathcal{O}_X\left( D\right)
^{\left[ \kappa \right] }=\mathcal{O}_X\left( \kappa \,D\right) ^{\vee
\,\vee }, 
\]
for any $D,D_1,D_2\in $ WDiv$\left( X\right) $ and $\kappa \in \Bbb{Z}$%
.\medskip

\noindent $\bullet $ For $X=X\left( N,\Delta \right) $ any $r$-dimensional
toric variety let now 
\[
T_N\text{-WDiv}\left( X\right) \ \ \text{ and \ \ }T_N\text{-CDiv}\left(
X\right) =T_N\text{-WDiv}\left( X\right) \cap \text{CDiv}\left( X\right) 
\]
denote the groups of $T_N$-invariant Weil and Cartier divisors,
respectively. If $\Delta $ is not contained in any proper subspace of $N_{%
\Bbb{R}}$, then\footnote{%
In particular, for $X$ smooth and compact, Pic$\left( X\right) $ is torsion
free and the Picard number equals $\#\Delta \left( 1\right) -r.$}\smallskip 
\[
\begin{array}{ccc}
0 &  & 0 \\ 
\uparrow &  & \uparrow \\ 
\text{Pic}\left( X\right) & \hookrightarrow & A_{r-1}\left( X\right) \\ 
\uparrow &  & \uparrow \\ 
T_N\text{-CDiv}\left( X\right) & \hookrightarrow & T_N\text{-WDiv}\left(
X\right) \\ 
\uparrow &  & \uparrow \\ 
M & =\!=\!= & M \\ 
\uparrow &  & \uparrow \\ 
0 &  & 0
\end{array}
\smallskip 
\]
is a commutative diagram with exact columns. $T_N$-WDiv$\left( X\right) $
has as $\Bbb{Z}$-basis: 
\[
T_N\text{-WDiv}\left( X\right) =\bigoplus \Bbb{Z}\,\left\{ V\left( \varrho
\right) \ \left| \ \varrho \in \Delta \left( 1\right) \right. \right\} \ . 
\]

\begin{theorem}[Divisors and support functions]
\label{INVERT}There exist one-to-one correspondences\emph{:} 
\[
\begin{array}{ccc}
\emph{PL}\text{\emph{-}}\emph{SF}\left( N,\Delta \right)  & 
\longleftrightarrow  & \left\{ 
\begin{array}{c}
\text{\emph{evaluation functions}} \\ 
\vartheta :\text{\emph{Gen}}\left( \Delta \right) \rightarrow \Bbb{Z}
\end{array}
\right\}  \\ 
\updownarrow  &  & \updownarrow  \\ 
T_N\text{\emph{-WDiv}}\left( X\right)  & \longleftrightarrow  & \left\{ 
\begin{array}{c}
\text{\emph{reflexive coherent (i.p.~torsion-free)}} \\ 
\text{\emph{sheaves of rank }}1\ \text{\emph{of} } \\ 
T_N\text{\emph{-invariant,} \emph{complete }} \\ 
\text{\emph{fractional ideals}} \\ 
\text{\emph{\ over} \thinspace }X=X\left( N,\Delta \right) 
\end{array}
\right\} 
\end{array}
\]
\emph{(}in fact, $\Bbb{Z}$-module isomorphisms\emph{), }induced by
mapping\smallskip 
\[
\vartheta \longmapsto \psi _\vartheta ,\ \ \psi =\psi _\vartheta \longmapsto
D=D_{\psi _{}}\,,
\]
and 
\[
D=D_\psi \stackrel{\mathbf{\delta }}{\longmapsto }\mathcal{F}_\psi =\mathcal{%
O}_X\left( D\right) \,\longmapsto \text{\emph{ord}}_{\mathcal{F}_\psi }\,,
\]
with $\mathbf{\delta }$ as above and 
\[
D_{\psi _\vartheta }:=-\sum_{\varrho \in \Delta \left( 1\right) }\,\vartheta
\left( n\left( \varrho \right) \right) \,V\left( \varrho \right) \ .
\]
Moreover, 
\[
\left\{ 
\begin{array}{c}
D=D_\psi \in T_N\text{\emph{-CDiv}}\left( X\right) \smallskip  \\ 
\text{\emph{(i.e.,} }\mathcal{F}_\psi =\mathcal{O}_X\left( D\right) \ \text{%
\emph{is invertible)}}
\end{array}
\right\} \Longleftrightarrow \text{ }\psi =\text{\emph{ord}}_{\mathcal{F}%
_\psi }\in \text{\emph{SF}}\left( N,\Delta \right) \,.
\]
\end{theorem}

\noindent \textit{Proof. }See Saint-Donat \cite{KKMS}, ch. I, \S\ 2, thm. 9,
pp. 28-31. $_{\Box }$

\begin{theorem}[$H^0$-generated]
\label{UCON}If $X=X\left( N,\Delta \right) $ is a compact toric variety and $%
\psi \in $ \emph{SF}$\left( N,\Delta \right) $, then\smallskip 
\[
\mathcal{O}_X\left( D_\psi \right) \text{\emph{is generated by its global
sections }}\Longleftrightarrow \psi \in \text{\emph{UCSF}}\left( N,\Delta
\right) .
\]
\end{theorem}

\noindent \textit{Proof. }See Oda \cite{Oda}, thm. 2.7, p. 76. $_{\Box }$

\begin{theorem}[Ampelness]
Let $X=X\left( N,\Delta \right) $ be a \emph{(}not necessarily compact\emph{%
) }toric variety\emph{.} Then a divisor 
\[
D=D_\psi \in T_N\text{\emph{-CDiv}}\left( X\right) \text{ \ \ \ \emph{(resp.}
}D=D_\psi \in T_N\text{\emph{-CDiv}}\left( X\right) \otimes _{\Bbb{Z}}\Bbb{Q}%
\text{\emph{)}}
\]
is \emph{ample} if and and only if 
\[
\psi \in \text{\emph{SUCSF}}\left( N,\Delta \right) \ \ \ \ \emph{\ }\text{%
\emph{(resp.} }\psi \in \text{\emph{SUCSF}}_{\Bbb{Q}}\left( N,\Delta \right) 
\text{\emph{)}}
\]
\end{theorem}

\noindent \textit{Proof. }See Kempf \cite{KKMS}, ch. I, \S\ 3, thm. 13, p.
48. $_{\Box }$

\begin{corollary}[Quasiprojectivity]
\label{quasiproj}Let $X\left( N,\Delta \right) $ be a toric variety \emph{(}%
resp. a compact toric variety\emph{).} Then $X\left( N,\Delta \right) $ is
quasiprojective \emph{(}resp. projective\emph{) }if and only if 
\[
\text{\emph{SUCSF}}\left( N,\Delta \right) \neq \varnothing \text{ \ \ \ 
\emph{(or equivalently SUCSF}}_{\Bbb{Q}}\left( N,\Delta \right) \neq
\varnothing \text{\emph{) .\smallskip }}
\]
\end{corollary}

\noindent \textsf{(i) }\textit{Intersection numbers of Cartier divisors on
toric varieties. }If\textit{\ }$X$ is a normal complex variety of dimension $%
r$, and $D_1,\ldots ,D_r$ Cartier divisors on $X$, such that $%
W:=\bigcap_{i=1}^r\left( \text{supp}\left( D_i\right) \right) $ is compact,
then their \textit{intersection number }is defined to be the degree 
\[
\left( D_1\,\cdots \,D_r\right) :=\text{ deg}_W\left( \left\{ D_1\,\cdots
\,D_r\right\} \right) \in \Bbb{Z} 
\]
of the zero-cycle 
\[
\left\{ D_1\cdot D_2\,\cdots \,D_r\right\} :=\left\{ D_1\cdot \left(
D_2\,\cdots \,D_r\right) \right\} \in A_0\left( W\right) \smallskip 
\]
determined inductively as usual (i.e. probably after passing to the
corresponding pseudodivisors). See Fulton \cite{Ful}, Ch. 1-2; in
particular, I.1.4, p. 13, and pp. 38-39. For $X=X\left( N,\Delta \right) $ a 
\textit{smooth }toric variety, $D_i=V\left( \varrho _i\right) $, for all $i$%
, $1\leq i\leq r$, and $\varrho _i$'s pairwise distinct rays, we
have\smallskip 
\begin{equation}
\left( D_1\,\cdots \,D_r\right) =\left\{ 
\begin{array}{lll}
1 & , & \text{if \thinspace \ \ }\varrho _1+\cdots +\varrho _r\in \Delta
\smallskip \\ 
0 & , & \text{otherwise}
\end{array}
\right.  \label{IN10}
\end{equation}
Moreover, one obtains by general techniques:

\begin{lemma}
Let $X\left( N,\Delta \right) $ be an $r$-dimensional smooth toric variety
and 
\[
D_1=V\left( \varrho _1\right) ,\ldots ,D_r=V\left( \varrho _r\right) 
\]
divisors on $X$. Suppose that either $X$ itself or at least $W$ is compact,
and that $\varrho _1=\varrho _2,$ while all the other rays are distinct.
Then \emph{:\smallskip }\newline
\emph{(i)} \ If $\tau :=\varrho _2+\cdots +\varrho _r\notin \Delta $, the
intersection number $\left( D_1\,\cdots \,D_r\right) $ vanishes.\smallskip 
\newline
\emph{(ii)} If $\tau \in \Delta $, there exist rays $\varrho ^{\prime
},\varrho ^{\prime \prime }\in \Delta \left( 1\right) $ and integers $\kappa
_2,\ldots ,\kappa _r$, such that 
\[
n\left( \varrho ^{\prime }\right) +n\left( \varrho ^{\prime \prime }\right)
+\sum_{j=2}^r\kappa _j\,n\left( \varrho _j\right) =0,\ \ \ \ \varrho
^{\prime }+\tau \in \Delta \left( r\right) ,\ \ \ \ \varrho ^{\prime \prime
}+\tau \in \Delta \left( r\right) .
\]
In this case 
\begin{equation}
\left( D_2^2\cdot D_3\,\cdots \,D_r\right) =\kappa _2  \label{INTD}
\end{equation}
\end{lemma}

\noindent \textit{Proof. }See Oda \cite{Oda}, p. 81. $_{\Box }\medskip $%
\newline
There are also easy generalizations for $\Delta $ simplicial but we shall
not use them because we shall work exclusively with intersection numbers of
divisors on smooth $X\left( N,\Delta \right) $'s. Another method for the
evaluation of intersection numbers is based on mixed polytope volumes.

\begin{proposition}
Let $X\left( N,\Delta \right) $ be an $r$-dimensional compact toric variety.
If $\psi \in $ \emph{UCSF}$\left( N,\Delta \right) $, then the
self-intersection number of $D_\psi $ equals 
\begin{equation}
D_\psi ^r=\left( r!\right) \cdot \text{\emph{Vol}}\left( P_\psi \right) 
\label{SELF}
\end{equation}
where 
\[
P_\psi :=P_{D_\psi }=\left\{ \mathbf{x}\in M_{\Bbb{R}}\ \left| \
\left\langle \mathbf{x},n\left( \varrho \right) \right\rangle \geq \psi
\left( n\left( \varrho \right) \right) \right. \right\} 
\]
is the lattice polytope associated to the divisor $D=D_\psi $ defined in 
\emph{\ref{INVERT}}. More generally, for $r$ upper convex functions $\psi _1$%
,$\ldots ,\psi _r$, one has 
\begin{equation}
\left( D_{\psi _1}\,\cdots \,D_{\psi _r}\right) =\left( r!\right) \cdot 
\text{\emph{Vol}}\left( P_{\psi _1},\ldots ,\,P_{\psi _r}\right) 
\label{MIXED}
\end{equation}
where \emph{Vol}$\left( P_{\psi _1},\ldots ,\,P_{\psi _r}\right) $ denotes
the mixed volume of the polytopes $P_{\psi _1},\ldots ,\,P_{\psi _r}$.
\end{proposition}

\noindent \textit{Proof. }See Oda \cite{Oda}, prop. 2.10, p. 79. $_{\Box
}\bigskip $\newline
\textsf{(j) }\textit{Euler-Poincar\'{e} characteristic. }The topological
Euler-Poincar\'{e} characteristic of a (not necessarily compact) toric
variety can be easily read off from the maximal cones of the defining fan.

\begin{proposition}
Let $X\left( N,\Delta \right) $ be an $r$-dimensional toric variety
associated to $\Delta $. Then the topological Euler-Poincar\'{e}
characteristic 
\[
\chi \left( X\left( N,\Delta \right) \right) =\sum_{i=0}^{2r}\left(
-1\right) ^i\emph{dim}_{\Bbb{Q}}\ H^i\left( X\left( N,\Delta \right) ;\Bbb{Q}%
\right) 
\]
of $X\left( N,\Delta \right) $ is equal to the number of $r$-dimensional
cones, i.e. 
\begin{equation}
\chi \left( X\left( N,\Delta \right) \right) =\#\left( \Delta \left(
r\right) \right)   \label{Eul-Poi}
\end{equation}
\end{proposition}

\noindent \textit{Proof. }See Fulton \cite{Fulton}, p. 59. $_{\Box
}\smallskip \medskip $

\noindent \textsf{(k) }\textit{Maps of fans}. A map of fans\textit{\ }$%
\varpi :\left( N^{\prime },\Delta ^{\prime }\right) \rightarrow \left(
N,\Delta \right) $ is a $\Bbb{Z}$-linear homomorphism $\varpi :N^{\prime
}\rightarrow N$ whose scalar extension $\varpi :N_{\Bbb{R}}^{\prime
}\rightarrow N_{\Bbb{R}}$ satisfies the property: 
\[
\forall \sigma ^{\prime },\ \sigma ^{\prime }\in \Delta ^{\prime }\ \ \text{ 
}\exists \ \sigma ,\ \sigma \in \Delta \ \ \text{ with\ \ }\varpi \left(
\sigma ^{\prime }\right) \subset \sigma \ . 
\]
$\varpi \otimes _{\Bbb{Z}}$id$_{\Bbb{C}^{*}}:T_{N^{\prime }}=N^{\prime
}\otimes _{\Bbb{Z}}\Bbb{C}^{*}\rightarrow T_N=N\otimes _{\Bbb{Z}}\Bbb{C}^{*}$
is a homomorphism from $T_{N^{\prime }}$ to $T_N$ and the scalar extension $%
\varpi ^{\vee }:M_{\Bbb{R}}\rightarrow M_{\Bbb{R}}^{\prime }$ of the dual $%
\Bbb{Z}$-linear map $\varpi ^{\vee }:M\rightarrow M^{\prime }$ induces an 
\textit{equivariant holomorphic map }$\varpi _{*}:X\left( N^{\prime },\Delta
^{\prime }\right) \rightarrow X\left( N,\Delta \right) $ as follows: If $%
\varpi \left( \sigma ^{\prime }\right) \subset \sigma $ \ for $\ \sigma \in
\Delta $, $\sigma ^{\prime }\in \Delta ^{\prime }$, then obviously $\varpi
^{\vee }\left( M\cap \sigma ^{\vee }\right) \subset M^{\prime }\cap \left(
\sigma ^{\prime }\right) ^{\vee }$, and the holomorphic map 
\[
\varpi _{*}:U_{\sigma ^{\prime }}^{\prime }\rightarrow U_\sigma \ \ \ \text{%
with\ \ }\ \varpi _{*}\left( u^{\prime }\right) \left( m\right) :=u^{\prime
}\left( \varpi ^{\vee }\left( m\right) \right) ,\ \forall m,\ m\in M\cap
\sigma ^{\vee }, 
\]
\ is equivariant because 
\[
\begin{array}{l}
\varpi _{*}\left( t^{\prime }\cdot u^{\prime }\right) \left( m\right)
=\left( t^{\prime }\cdot u^{\prime }\right) \left( \varpi ^{\vee }\left(
m\right) \right) =\smallskip \\ 
=t^{\prime }\left( \varpi ^{\vee }\left( m\right) \right) \cdot u^{\prime
}\left( \varpi ^{\vee }\left( m\right) \right) =\varpi _{*}\left( t^{\prime
}\right) \left( m\right) \cdot \varpi _{*}\left( u^{\prime }\right) \left(
m\right) ,
\end{array}
\]
for all $t^{\prime }\in T_{N^{\prime }},\ m\in M\cap \sigma ^{\vee }.$ After
gluing together the affine charts of $\Delta $ and of $\Delta ^{\prime }$ we
determine a well-defined map $\varpi _{*}:X\left( N^{\prime },\Delta
^{\prime }\right) \rightarrow X\left( N,\Delta \right) .$

\begin{theorem}[Properness]
\label{PRBIR}If $\varpi :\left( N^{\prime },\Delta ^{\prime }\right)
\rightarrow \left( N,\Delta \right) $ is a map of fans, $\varpi _{*}$ is 
\emph{proper} if and only if $\varpi ^{-1}\left( \left| \Delta \right|
\right) =\left| \Delta ^{\prime }\right| .$ In particular, if $N=N^{\prime }$
and $\Delta ^{\prime }$ is a refinement of $\Delta $, i.e. if each cone of $%
\Delta $ is a union of cones of $\Delta ^{\prime }$, then the holomorphic
map \emph{id}$_{*}:X\left( N,\Delta ^{\prime }\right) \rightarrow X\left(
N,\Delta \right) $ is \emph{proper }and \emph{birational.}
\end{theorem}

\noindent \textit{Proof. }See Oda \cite{Oda}, thm. 1.15, p. 20, and cor.
1.18, p. 23. $_{\Box }$

\section{Blow-ups and resolutions of toric varieties\label{TRESOL}}

\noindent One of the most fundamental cornerstones of various significant
constructions of birational morphisms between complex varieties is the
blowing up along subvarieties or -more general- along closed
subschemes.\medskip

\noindent \textsf{(a) }\textit{Local construction.} Let $U=$ Max-Spec$\left(
R\right) $ be (the closed point set of) an affine noetherian scheme, $I$ an
ideal of $R$, $Z=$ Max-Spec$\left( R/I\right) $ and $S\left( R,I\right)
:=\bigoplus_{d\geq 0}I^d$. The homogeneous spectrum Proj$\left( S\left(
R,I\right) \right) $ of $S\left( R,I\right) $, together with the structure
morphism 
\[
\text{ Bl}_I\left( R\right) :=\text{ Bl}_Z^I\left( U\right) :=\text{ Proj}%
\left( S\left( R,I\right) \right) \stackrel{\beta }{\longrightarrow }U 
\]
is called the blow-up of $U$ w.r.t. $I$ or \textit{the blow-up} of $U$ 
\textit{along} $Z$ (or of $U$ with\textit{\ center} $Z$). If $\left\{
h_0,\ldots ,h_\mu \right\} $ is a set of generators of $I$, then 
\[
\text{Bl}_Z^I\left( U\right) =\text{ }\bigcup_{i=0}^\mu \ \text{Max-Spec}%
\left( R\left[ \frac{h_0}{h_i},\ldots ,\frac{h_\mu }{h_i}\right] \right) 
\]
with $R\left[ \frac{h_0}{h_i},\ldots ,\frac{h_\mu }{h_i}\right] $ viewed as
an $R$-subalgebra of $R_{h_i}$. The exceptional locus of $\beta $ is Exc$%
\left( \beta \right) :=\beta ^{-1}\left( Z\right) $ and its contraction
locus $=Z$. Moreover, $\beta ^{-1}\left( I^{\,\sim }\right) $ is invertible
(with $I^{\,\sim }$ denoting here the sheaf being associated to the ideal $I$%
, and with $I$ regarded as an $R$-module, cf. \cite{Hart}, p. 110), and $%
\beta ^{-1}\left( I^{\,\sim }\right) \cong S\left( R,I\right) \left(
1\right) ^{\,\sim }.\bigskip $\newline
\textsf{(b) }\textit{Globalization by gluing lemma.} Let $\left\{ \frak{X}%
_j\ \left| \ j\in J\right. \right\} $ be a family of schemes.

\begin{lemma}[Gluing schemes]
\label{GLSCH}If there exists a collection $\left\{ \frak{X}_{j,k}\ \left| \
j,k\in J\right. \right\} $ of open sets $\frak{X}_{j,k}$ of $\frak{X}_j$ and
isomorphisms of schemes\smallskip 
\[
\eta _{k,j}:\left( \frak{X}_{j,k},\mathcal{O}_{\frak{X}_j}\left| _{\frak{X}%
_{j,k}}\right. \right) \stackrel{\cong }{\longrightarrow }\left( \frak{X}%
_{k,j},\mathcal{O}_{\frak{X}_j}\left| _{\frak{X}_{k,j}}\right. \right) 
\]
satisfying the conditions \emph{:\smallskip }\newline
$\bullet $ $\frak{X}_j=\frak{X}_{j,j}$ and $\eta _{j,j}=$ \emph{id}$_{\frak{X%
}_j},\smallskip $\newline
$\bullet $ $\eta _{j,k}\circ \eta _{k,j}=$ \emph{id}$_{\frak{X}_{j,k}},$
and\smallskip \newline
$\bullet $ $\left( \eta _{j,k}\circ \eta _{j,\xi }\right) \left| _{\frak{X}%
_{j,k,\xi }}\right. =\eta _{j,\xi }\left| _{\frak{X}_{j,k,\xi }}\right. $,
where $\frak{X}_{j,k,\xi }:=\frak{X}_{j,k}\cap \eta _{j,\xi }^{-1}\left( 
\frak{X}_{k,j}\right) ,\smallskip $\newline
then there exists a scheme $W$, an open cover $\left\{ W_j\ \left| \ j\in
J\right. \right\} $ of $W$, and a collection of isomorphisms $\left\{ f_j:W_j%
\stackrel{\cong }{\longrightarrow }\frak{X}_j\ \left| \ j\in J\right.
\right\} $, such that $\eta _{k,j}=f_k\circ f_j^{-1}\left| _{\frak{X}%
_{j,k}}\right. $ for all $i,j\in J$.\smallskip \newline
\emph{The scheme }$W$\emph{\ obtained by gluing the members of the family }$%
\left\{ \text{ }\frak{X}_j\ \left| \ j\in J\right. \right\} $\emph{\ via the
above isomorphisms will be denoted by} 
\[
W=\coprod_{\eta _{k,j}}\text{ }\frak{X}_j\text{\ \ \emph{or simply by }\ \ }%
\,W=\coprod^{\blacktriangledown }\text{ }\frak{X}_j
\]
\emph{if the gluing isomorphisms are self-evident from the context.}
\end{lemma}

\noindent \textit{Proof.} See e.g. \cite{Hart}, Ex. II.2.12, p. 80. $_{\Box
}\medskip $\newline
Let now $\mathcal{I}$ be a coherent (non-zero) sheaf of $\mathcal{O}_X$%
-ideals over a complex variety $X$ and 
\[
\left\{ U_j=\text{ Max-Spec}\left( R_j\right) \ \left| \ j\in J\right.
\right\} 
\]
an affine cover of $X$. For every $\ j\in J$, we have $\mathcal{I}\left|
_{U_j}\right. \cong I_j^{\,\sim }$ for some ideal $I_j$ of $R_j$ (cf. \cite
{Hart}, II.5.4). Considering 
\[
\left\{ \beta _j:\text{ Bl}_{Z_j}^{I_j}\left( U_j\right) \longrightarrow
U_j\subset X\ \left| \ j\in J\right. \right\} 
\]
as in \textsf{(a)}, where $Z_j=$ Max-Spec$\left( R_j/I_j\right) $, and $Z$
the closed subscheme of $X$ defined by 
\[
Z=\text{ supp}\left( \mathcal{O}_X/\mathcal{I}\right) :=\left\{ \mathbf{x}%
\in X\ \left| \ \left( \mathcal{O}_X/\mathcal{I}\right) _{\mathbf{x}}\neq
0\right. \right\} \ , 
\]
we determine natural isomorphisms 
\[
\beta _j^{-1}\left( U_j\cap U_k\right) \stackrel{\cong }{\stackunder{\eta
_{k,j}}{\longrightarrow }}\beta _k^{-1}\left( U_j\cap U_k\right) \ . 
\]
Applying lemma \ref{GLSCH} to the family $\left\{ \beta _j^{-1}\left(
U_j\right) \ \left| \ j\in J\right. \right\} $ we construct a birational,
proper, surjective morphism $\pi =\pi _{\mathcal{I}}\smallskip $%
\[
\fbox{$
\begin{array}{ccc}
&  &  \\ 
& \mathbf{Bl}_Z^{\mathcal{I}}\left( X\right) :=\dcoprod\limits_{\eta
_{k,j}}\beta _j^{-1}\left( U_j\right) =\mathbf{Proj}\left(
\bigoplus\limits_{d\geq 0}\ \mathcal{I}^{\,d}\right) \stackrel{\pi }{%
\longrightarrow }X &  \\ 
&  & 
\end{array}
$}\smallskip 
\]
as the natural projection induced by $\mathcal{O}_X\stackrel{\cong }{%
\rightarrow }\mathcal{I}^{\,0}\hookrightarrow \bigoplus_{d\geq 0}\ \mathcal{I%
}^{\,d}$, with $\mathbf{Proj}$ denoting the \textit{global} homogeneous
spectrum (as in \cite{Hart}, p. 160). $\left( \mathbf{Bl}_Z^{\mathcal{I}%
}\left( X\right) ,\pi \right) $ is \textit{the blow-up} of $X$ \textit{along}
$\mathcal{I}$ or \textit{the monoidal transformation} w.r.t. $\mathcal{I}$
(with \textit{center} $Z$). Let us recall its main properties :\newline
$\bullet $ $\pi $ induces an isomorphism $\mathbf{Bl}_Z^{\mathcal{I}}\left(
X\right) \smallsetminus \pi ^{-1}\left( Z\right) \stackrel{\cong }{%
\rightarrow }X\smallsetminus Z$, i.e. Exc$\left( \pi \right) =\pi
^{-1}\left( Z\right) $.\smallskip \newline
$\bullet $ The algebraic scheme $\mathbf{Bl}_Z^{\mathcal{I}}\left( X\right) $
is a complex variety (cf. \cite{Hart}, II. 7.16 (a), p. 166).\smallskip 
\newline
$\bullet $ If both $X$ and $Z$ are smooth, then $\mathbf{Bl}_Z^{\mathcal{I}%
}\left( X\right) $ is smooth too.\smallskip \newline
$\bullet $ The preimage sheaf of $\mathcal{I}$, $\pi ^{-1}\mathcal{I}\cdot 
\mathcal{O}_X:=$ Im$\left( \pi ^{*}\mathcal{I}\rightarrow \mathcal{O}_{%
\mathbf{Bl}_Z^{\mathcal{I}}\left( X\right) }\right) $, is invertible and
determines Exc$\left( \pi \right) .$ Hence, Exc$\left( \pi \right) $
represents a Cartier (not necessarily prime) divisor of $\mathbf{Bl}_Z^{%
\mathcal{I}}\left( X\right) $ which is isomorphic to the projectivization $%
\Bbb{P}\left( \text{NC}_Z^{\mathcal{I}}\left( X\right) \right) $ of \textit{%
the normal cone} 
\[
\text{NC}_Z^{\mathcal{I}}\left( X\right) :=\text{ Spec}\left(
\bigoplus_{d\geq 0}\,\mathcal{I}^d\,/\,\mathcal{I}^{d+1}\right) 
\]
of $X$ \textit{along} $Z$. The relation between $\mathbf{Bl}_Z^{\mathcal{I}%
}\left( X\right) $ and Exc$\left( \pi \right) $ is described by the
isomorphisms 
\[
\mathcal{N}_{\text{Exc}\left( \pi \right) \,\left| \,\mathbf{Bl}_Z^{\mathcal{%
I}}\left( X\right) \right. }\cong \mathcal{O}_{\mathbf{Bl}_Z^{\mathcal{I}%
}\left( X\right) }\left( \text{Exc}\left( \pi \right) \right) \left| _{\,%
\text{Exc}\left( \pi \right) }\right. \cong \mathcal{O}_{\Bbb{P}\left( \text{%
NC}_Z^{\mathcal{I}}\left( X\right) \right) }\left( -1\right) \,. 
\]
In particular, if $Z$ is a local complete intersection in $X$, the canonical
epimorphism from the $d$-th symmetrizer sheaf of $\mathcal{I}\,/\,\mathcal{I}%
^2$ onto the $d$-th part of the normal cone graded algebra 
\[
\text{Sym}^d\,\left( \mathcal{I}\,/\,\mathcal{I}^2\right) \twoheadrightarrow 
\mathcal{I}^d\,/\,\mathcal{I}^{d+1} 
\]
becomes an isomorphism and therefore 
\[
\text{Exc}\left( \pi \right) \cong \mathbf{Proj}\left( \bigoplus_{d\geq 0}\,%
\text{Sym}^d\,\left( \mathcal{I}\,/\,\mathcal{I}^2\right) \right) \cong \Bbb{%
P}\left( \mathcal{I}\,/\,\mathcal{I}^2\right) =\Bbb{P}\left( \mathcal{N}%
_{Z\,\left| \,X\right. }^{\,\,\vee }\right) \,. 
\]
(Here $\mathcal{N}_{\ldots }$, $\mathcal{N}_{\ldots }^{\,\,\vee }$ denote
the corresponding normal and conormal sheaves.)\bigskip \newline
\textsf{(c) }\textit{Universal property of blowing up}. If $g:X^{\prime
}\rightarrow X$ is any morphism, and $Z$ (resp. $Z^{\prime }$) is defined by
the ideal sheaf $\mathcal{I}$ (resp. by $\mathcal{I}^{\prime }:=\pi ^{-1}%
\mathcal{I}\cdot \mathcal{O}_{X^{\prime }}$), where $\pi $ denotes the
blow-up-morphism of $X$ along $Z$, then composing the morphism $\mathbf{Bl}%
_{Z^{\prime }}^{\mathcal{I}^{\prime }}\left( X^{\prime }\right) \rightarrow 
\mathbf{Bl}_Z^{\mathcal{I}}\left( X\right) \times _XX^{\prime }$ with the
projection to $X^{\prime }$, we get the commutative diagram: 
\[
\begin{array}{ccccc}
&  & \mathbf{Bl}_Z^{\mathcal{I}}\left( X\right) & \stackrel{\pi }{%
\longrightarrow } & X \\ 
& \nearrow & \uparrow & \nwarrow ^h & \ \uparrow \ g \\ 
\mathbf{Bl}_{Z^{\prime }}^{\mathcal{I}^{\prime }}\left( X^{\prime }\right) & 
\rightarrow & \mathbf{Bl}_Z^{\mathcal{I}}\left( X\right) \times _XX^{\prime }
& \longrightarrow & X^{\prime }
\end{array}
\]
If $g^{-1}\mathcal{I}\cdot \mathcal{O}_{X^{\prime }}$ is invertible, then it
is easy to deduce an isomorhism $\mathbf{Bl}_{Z^{\prime }}^{\mathcal{I}%
^{\prime }}\left( X^{\prime }\right) \cong X^{\prime }$. Hence, there is a
unique morphism $h$ factorizing $g$. This means that $\left( \mathbf{Bl}_Z^{%
\mathcal{I}}\left( X\right) ,\pi \right) $ is universal among all pairs $%
\left( X^{\prime },g\right) $ having invertible ideal sheaves $g^{-1}%
\mathcal{I}\cdot \mathcal{O}_{X^{\prime }}$ .\bigskip \newline
\textsf{(d) }\textit{Blowing up intermediate subschemes}. Let $Z\subsetneqq
X $, $W\subsetneqq X$ be two closed subschemes of $X$ which are defined by
the ideal sheaves $\mathcal{I}$ and $\mathcal{J}$ respectively, such that $%
Z\cap W=$ supp$\left( \mathcal{O}_X\,/\,\mathcal{I}+\mathcal{J}\right) \neq
\varnothing $ and $Z\cap W$ is nowhere dense in $W$, and let 
\[
\pi =\pi _{\mathcal{I}}:\mathbf{Bl}_Z^{\mathcal{I}}\left( X\right)
\longrightarrow X,\ \ \ \pi \left| _{\text{restr.}}\right. :\mathbf{Bl}%
_{Z\cap W}^{\mathcal{I}+\mathcal{J}}\left( W\right) \longrightarrow W 
\]
denote the corresponding blow-ups. Then, by the above mentioned universality
of the morphism $\pi $, one verifies easily the following isomorphism: 
\[
\begin{array}{ccccc}
&  & \mathbf{Bl}_Z^{\mathcal{I}}\left( X\right) \smallskip & \stackrel{\pi }{%
\longrightarrow } & X \\ 
&  & \bigcup &  & \bigcup \\ 
\left\{ 
\begin{array}{c}
\text{scheme-theoretic closure of the preimage} \\ 
\text{ }\pi ^{-1}\left( W\smallsetminus \left( Z\cap W\right) \right) \ 
\text{ \ within \ \ }\mathbf{Bl}_Z^{\mathcal{I}}\left( X\right)
\end{array}
\right\} & \cong & \mathbf{Bl}_{Z\cap W}^{\mathcal{I}+\mathcal{J}}\left(
W\right) & \stackrel{\pi \left| _{\text{restr.}}\right. }{\longrightarrow }
& W
\end{array}
\]
This closure is called \textit{the strict transform} STR$\left( W,\mathcal{J}%
,\pi \right) $ of $\left( W,\pi \right) $ \textit{under} $\pi $. If we
assume that $W\nsubseteq Z$ and $Z\nsubseteq W$, then STR$\left( Z,\mathcal{I%
},\pi _{\mathcal{I}}\right) \cap $ STR$\left( W,\mathcal{J},\pi _{\mathcal{J}%
}\right) =\varnothing $ (cf. \cite{Hart}, Ex. 7.12, p. 171).\bigskip \newline
\textsf{(e) }\textit{Normalization process}. Even if the complex variety $X$
itself is normal and $Z=$ supp$\left( \mathcal{O}_X\,/\,\mathcal{I}\right) $
smooth, with $Z\cap $ Sing$\left( X\right) \neq \varnothing $, $\mathbf{Bl}%
_Z^{\mathcal{I}}\left( X\right) $ is \textit{not} necessarily normal. Using
an affine cover of $X$ 
\[
\left\{ U_j=\text{ Max-Spec}\left( R_j\right) \ \left| \ j\in J\right.
\right\} \,\,\,\ \,\text{and \thinspace \thinspace \thinspace \thinspace
\thinspace \ }Z_j=\text{Max-Spec}\left( R_j/I_j\right) 
\]
(as in \textsf{(b)}) and the (finite) normalization morphisms onto Bl$%
_{Z_j}^{I_j}\left( U_j\right) $ : 
\[
\mathbf{\nu }_j:\text{Norm\thinspace }\left[ \text{Bl}_{Z_j}^{I_j}\left(
U_j\right) \right] \longrightarrow \text{Bl}_{Z_j}^{I_j}\left( U_j\right)
\,, 
\]
we define \textit{the normalized blow-up} $\left( \text{Norm\thinspace }%
\left[ \mathbf{Bl}_Z^{\mathcal{I}}\left( X\right) \right] ,\,\pi _{\mathcal{I%
}}\circ \mathbf{\nu }_{\mathcal{I}}\right) $ of $X$ w.r.t. $Z$ (or with 
\textit{center }$Z$) by patching the affine pieces together :\smallskip 
\[
\fbox{$
\begin{array}{ccc}
&  &  \\ 
& \text{Norm\thinspace }\left[ \mathbf{Bl}_Z^{\mathcal{I}}\left( X\right)
\right] :=\dcoprod\limits^{\blacktriangledown }\,\text{Norm\thinspace }%
\left[ \text{Bl}_{Z_j}^{I_j}\left( U_j\right) \right] \stackrel{\mathbf{\nu }%
_{\mathcal{I}}}{\longrightarrow }\mathbf{Bl}_Z^{\mathcal{I}}\left( X\right) 
\stackrel{\pi _{\mathcal{I}}}{\longrightarrow }X &  \\ 
&  & 
\end{array}
$}\smallskip 
\]
The combination of the universal property of the normalization morphism $%
\mathbf{\nu }_{\mathcal{I}}$ (see \cite{Hart}, Ex. 3.8, p. 91) with that of $%
\pi _{\mathcal{I}}$ (cf. \textsf{(c)}) leads to the universal property of
normalized blowing up: If $g:X^{\prime }\rightarrow X$ is any proper
morphism, $X^{\prime }$ normal and $g^{-1}\mathcal{I}\cdot \mathcal{O}%
_{X^{\prime }}$ is invertible, then there exists a unique morphism $h$
factorizing $g$ : 
\[
\begin{array}{ccc}
\text{Norm\thinspace }\left[ \mathbf{Bl}_Z^{\mathcal{I}}\left( X\right)
\right] & \stackrel{\pi _{\mathcal{I}}\circ \mathbf{\nu }_{\mathcal{I}}}{%
\longrightarrow } & X \\ 
& \nwarrow ^h & \,\uparrow \,g \\ 
&  & X^{\prime }
\end{array}
\]
\textsf{(f) }\textit{Some warnings. }To blow up arbitrary subschemes $Z$ of $%
X$ requires great care.\smallskip \newline
$\bullet $ Sometimes, for different coherent ideal $\mathcal{O}_X$-sheaves $%
\mathcal{I}$, $\mathcal{I}^{\prime }$ over $X$ and 
\[
Z:=\text{supp}\left( \mathcal{O}_X\,/\,\mathcal{I}\right) ,\ \ \ Z^{\prime
}:=\text{supp}\left( \mathcal{O}_X\,/\,\mathcal{I}^{\prime }\right) , 
\]
it is possible to have an isomorphism $\mathbf{Bl}_Z^{\mathcal{I}}\left(
X\right) \cong \mathbf{Bl}_{Z^{\prime }}^{\mathcal{I}^{\prime }}\left(
X\right) $. For instance, for $\mathcal{I}$ an arbitrary coherent sheaf of
ideals and $\mathcal{I}^{\prime }$ an invertible sheaf of ideals, $\mathcal{I%
}$ and $\mathcal{I}\cdot \mathcal{I}^{\prime }$ give isomorphic blow-ups.
The same remains true for $\mathcal{I}$ and $\mathcal{I}^{\,d}$ with $d\geq
2 $ (see \cite{Hart}, Ex. 7.11, p. 171).\smallskip \newline
$\bullet $ On the other hand, if $Z$ is endowed with two different scheme
structures 
\[
Z=\text{supp}\left( \mathcal{O}_X\,/\,\mathcal{I}\right) =\text{supp}\left( 
\mathcal{O}_X\,/\,\mathcal{J}\right) ,\ \ \mathcal{I}\neq \mathcal{J}, 
\]
then $\mathbf{Bl}_Z^{\mathcal{I}}\left( X\right) $ and $\mathbf{Bl}_Z^{%
\mathcal{J}}\left( X\right) $ are in general too far from being isomorphic
to each other.

\begin{example}
\label{2BLUP}\emph{Perhaps the simplest example is to consider }$X=\Bbb{A}_{%
\Bbb{C}}^2$\emph{\ to be the affine complex plane with coordinate system }$%
\left\{ x,y\right\} $\emph{, }$Z=\left\{ \left( 0,0\right) \right\} $\emph{\
the zero point, }$\mathcal{I}=\left( x,y\right) \cdot \mathcal{O}_{\Bbb{A}_{%
\Bbb{C}}^2}$\emph{\ the maximal ideal of }$\mathcal{O}_{\Bbb{A}_{\Bbb{C}%
}^2,\left( 0,0\right) }$\emph{\ and }$\mathcal{J}=\left( x,y^2\right) \cdot 
\mathcal{O}_{\Bbb{A}_{\Bbb{C}}^2}$\emph{. Then } 
\[
\mathbf{Bl}_{\left\{ \left( 0,0\right) \right\} }^{\mathcal{I}}\left( \Bbb{A}%
_{\Bbb{C}}^2\right) \stackrel{\pi _{\mathcal{I}}}{\longrightarrow }\Bbb{A}_{%
\Bbb{C}}^2
\]
\emph{\ is the usual blow-up with exceptional set} 
\[
\pi _{\mathcal{I}}^{-1}\left( \left\{ \left( 0,0\right) \right\} \right)
\cong \Bbb{P}_{\Bbb{C}}^1\ \ \emph{\ and\ \ \ Sing}\left( \mathbf{Bl}%
_{\left\{ \left( 0,0\right) \right\} }^{\mathcal{I}}\left( \Bbb{A}_{\Bbb{C}%
}^2\right) \right) =\varnothing \emph{.}
\]
\emph{\ In contrast to the usual blow-up case, and although the blow-up } 
\[
\mathbf{Bl}_{\left\{ \left( 0,0\right) \right\} }^{\mathcal{J}}\left( \Bbb{A}%
_{\Bbb{C}}^2\right) \stackrel{\pi _{\mathcal{J}}}{\longrightarrow }\Bbb{A}_{%
\Bbb{C}}^2
\]
\emph{w.r.t. the prime ideal }$\mathcal{J}$\emph{\ has also a smooth
rational curve as exceptional set, the singular locus of }$\mathbf{Bl}%
_{\left\{ \left( 0,0\right) \right\} }^{\mathcal{J}}\left( \Bbb{A}_{\Bbb{C}%
}^2\right) $\emph{\ is non-empty. As one may easily verify, it consists of a
single ordinary double point which lies on this rational curve. The real
parts }$\pi _{\mathcal{I}}^{-1}\left( \frak{D}\right) $\emph{\ and }$\pi _{%
\mathcal{J}}^{-1}\left( \frak{D}\right) $\emph{\ over a small disk }$\frak{D}%
\subset \Bbb{A}_{\Bbb{R}}^2\subset \Bbb{A}_{\Bbb{C}}^2$\emph{\ centered at }$%
\left( 0,0\right) $\emph{\ are illustrated in figures \textbf{1} (a) and 
\textbf{2}.  The real part of }$\pi _{\mathcal{I}}^{-1}\left( \frak{D}\right) 
$\emph{\ is used to be viewed as a} \textit{spiral staircase} \emph{(whose
stairs extend in both directions). Away from the origin we get an
isomorphism, and the points of Exc}$\left( \pi _{\mathcal{I}}\right) $\emph{%
\ are in 1-1 correspondence with the set of straight lines passing through }$%
\left( 0,0\right) $\emph{\ (see fig.\textbf{1} (a)). One may, of course,
think of it just as an} \textit{enlargement} \emph{(}\textit{Aufblasung}%
\emph{) of the origin spreading out the ``tangent directions'' (see fig.%
\textbf{1} (b)). In fact, it can be shown that the topological space of the
real part of }$\pi _{\mathcal{I}}^{-1}\left( \frak{D}\right) $ \emph{is
homeomorphic to a M\"{o}bius strip, while, in the analogous setting, }$\pi _{%
\mathcal{J}}^{-1}\left( \frak{D}\right) $ \emph{has a marked ``twisted''
point which corresponds to the occuring singularity (see fig.\textbf{2}). In
both cases the dotted meridian line indicates the exceptional set.}
\end{example}

\def\makespace{\par\vspace{1.4\baselineskip}\par}%
\begin{figure*}[htbp]
\begin{center}
\input{fig1a.pstex_t}
\makespace
 {\textbf{Fig.} \textbf{1 }(a)}\par
\end{center}
\end{figure*}
\begin{figure*}[htbp]
\begin{center}
\input{fig1b.pstex_t}
\makespace
 {\textbf{Fig.} \textbf{1 }(b)}\par
\end{center}
\end{figure*}
\begin{figure*}[htbp]
\begin{center}
\input{fig2.pstex_t}
\makespace
 {\textbf{Fig.} \textbf{2 }}
\end{center}
\end{figure*}

\vspace{4cm}
\noindent $\bullet $ The closed subscheme $Z$ of $X$ can be in general
equipped with lots of scheme structures $\mathcal{O}_Z$, with $\mathcal{O}_Z=%
\mathcal{O}_X\,/\,\mathcal{I}$ for different $\mathcal{I}$'s. There is,
however, a unique scheme structure above them, which is the ``smallest
one''\thinspace ; namely, the reduced induced structure $\mathcal{O}_Z\,/\,%
\frak{R}_Z$ (supporting the same underlying topological space) where $\frak{R%
}_Z$ is the $\mathcal{O}_Z$-ideal sheaf defined by 
\[
\left\{ \text{open sets of\ \thinspace }Z\right\} \ni U\longmapsto \frak{R}%
_Z\left( U\right) :=\left\{ \text{the nilradical of\ \thinspace }\mathcal{O}%
_Z\left( U\right) \right\} \, 
\]
(see \cite{Hart}, II.3.2.6, p. 86). Let $Z_{\mathbf{red}}=$ supp$\left( 
\mathcal{O}_Z\,/\,\frak{R}_Z\right) $ denote the reduced subscheme
associated to $Z$. The singular locus\footnote{%
Since we allow $Z$ to carry a non-reduced structure (in contrast to $X$
which is assumed to be a complex variety, and therefore always reduced), Sing%
$\left( Z\right) $ might become very ``big''. In particular, in the above
example \ref{2BLUP}, for $Z=\left\{ \left( 0,0\right) \right\} =$ supp$%
\left( \mathcal{O}_{\Bbb{A}_{\Bbb{C}}^2}\,/\,\mathcal{J}\right) $, we obtain 
$Z=$ Sing$\left( Z\right) \,!$} of $Z$ can be written as 
\[
\text{Sing}\left( Z\right) =\text{Sing}\left( Z_{\mathbf{red}}\right) \cup
\left\{ \text{non-reduced points of }Z\right\} \,.\smallskip 
\]
From now on, the uniquely determined blow-up $\mathbf{Bl}_{Z_{\mathbf{red}%
}}^{\mathcal{I}\left| _Z\right. /\frak{R}_Z}\left( X\right) $ of $X$ along
this reduced subscheme $Z_{\mathbf{red}}$ will be denoted simply by $\mathbf{%
Bl}_Z^{\text{red}}\left( X\right) $ and will be called the \textit{usual }%
blow-up of $X$ along $Z$.\smallskip \newline
$\bullet $ Another notable property of a blow-up is the projectivity of its
defining morphism $\pi _{\mathcal{I}}$ (and consequently the projectivity of 
$\pi _{\mathcal{I}}\circ \mathbf{\nu }_{\mathcal{I}}$, because $\mathbf{\nu }%
_{\mathcal{I}}$ is finite).

\begin{proposition}
If $X$ is quasiprojective \emph{(}resp. projective\emph{)}, then both $%
\mathbf{Bl}_Z^{\mathcal{I}}\left( X\right) $ and \emph{Norm}\thinspace $%
\left[ \mathbf{Bl}_Z^{\mathcal{I}}\left( X\right) \right] $ will be
quasiprojective \emph{(}resp. projective\emph{) }as well.
\end{proposition}

\noindent \textit{Proof. }See Hartshorne \cite{Hart}, prop. II.7.16, p. 166. 
$_{\Box }$

\begin{corollary}
\label{MBLU}Let $X$ be a quasiprojective \emph{(}resp. projective\emph{) }%
complex variety and $Z=$ \emph{supp}$\left( \mathcal{O}_X\,/\,\mathcal{I}%
\right) $ a closed subscheme of $X$. If $X_0:=X$, $Z_0:=Z$, $\mathcal{I}_0:=%
\mathcal{I}$ \thinspace \thinspace and 
\[
\left\{ \pi _{\mathcal{I}_{j-1}}:X_j=\mathbf{Bl}_{Z_{j-1}}^{\mathcal{I}%
_{j-1}}\left( X_{j-1}\right) \longrightarrow X_{j-1}\ \left| \ 1\leq j\leq
l\right. \right\} 
\]
is a finite sequence of blow-ups with $Z_j=$ \emph{supp}$\left( \mathcal{O}%
_{X_j}\,/\,\mathcal{I}_j\right) $, then $X_l$ is also quasiprojective \emph{(%
}resp. projective\emph{).\smallskip }
\end{corollary}

\noindent $\bullet $ The troubles with the projectivity begin whenever one
wishes to blow up (once or more times) ``overlapping one another pieces'' $%
U_1,\ldots ,U_k$ of a fixed closed (possibly reducible, non-reduced or
singular) subscheme $Z$ of a quasiprojective complex variety with $%
Z=\bigcup_{i=1}^kU_i$, and try to glue the blown up (non-disjoint) overlying
parts, say $W_1,\ldots ,W_k$, together. Even if this gluing procedure is
absolutely natural, the resulting morphism 
\[
f:X^{\prime }=\dcoprod\limits^{\blacktriangledown }\,W_i\longrightarrow X 
\]
(with contraction locus $Z$) is \textit{not }always projective. The two
classical examples of this kind, having as starting-point two smooth curves
meeting trasversally at two points and a curve with an ordinary double point
on a smooth complete threefold, respectively, are due to Hironaka. See
Hartshorne \cite{Hart}, pp. 443-445, for the details of the construction.
The non-projectivity of such an $f$ implies automatically not only the
non-quasiprojectivity of $X^{\prime }$ but also the fact that $f$ cannot be
expressed as a composite of finitely many blow-ups (by cor. \ref{MBLU}). For
another simple example, see rem. \ref{NONPROJ} below.\bigskip \newline
\textsf{(g) }\textit{Desingularization by successive blow-ups with smooth
centers. }We just formulate here the famous theorem of Hironaka which
guarantees the \textit{existence} of (full) desingularizations by performing
a finite number of monoidal transformations with smooth centers. (Of course,
we should again stress, that \textit{not every }desingularization can be
composed of finitely many blow-ups.)

\begin{theorem}
\label{HIRTH}\emph{(\textbf{Hironaka's Theorem}\footnote{%
The original result of Hironaka is more general. It is valid for any
equicharacteristic zero excellent scheme $X$. The centers $Z_j$ are normally
flat in the ambient space.}, \cite{Hiro}). }Let $X$ be any complex variety.
Then there exists always a finite sequence of blow-ups 
\[
\left\{ \pi _{\mathcal{I}_{j-1}}:X_j=\mathbf{Bl}_{Z_{j-1}}^{\mathcal{I}%
_{j-1}}\left( X_{j-1}\right) \longrightarrow X_{j-1}\ \left| \ 1\leq j\leq
l\right. \right\} \ ,
\]
\newline
such that $X_0=X$ and \emph{Sing}$\left( X_l\right) =\varnothing $, where%
\emph{\ }$Z_j$ is a \emph{smooth} subvariety of \emph{Sing}$\left(
X_j\right) ,\forall j$, $0\leq j\leq l-1$,
\end{theorem}

\noindent (There are meanwhile considerable simplifications of the original
proof of thm. \ref{HIRTH}, like those due to Bierstone \& Milman \cite{BM1}, 
\cite{BM2}, who introduced an appropriate discrete local invariant for
points of $X$ whose maximum locus determines a center of blow-up, leading to
constructive desingularization.)\bigskip

\noindent \textsf{(h) }\textit{Toric blow-ups. }Working within the category
of toric varieties, both blow-up and resolution are much more easier as they
can be translated into purely combinatorial operations involving specific
cones of our fans.

\begin{theorem}[Toric normalized blow-up]
\label{TORNBL}Let $X\left( N,\Delta \right) $ be a toric variety, $\mathcal{I%
}\neq \,$\underline{$0$} a $T_N$-invariant coherent ideal sheaf defining a
subscheme $Z$ of $X\left( N,\Delta \right) $ with $\mathcal{F}=\left( 
\mathcal{O}_{X\left( N,\Delta \right) }\,/\,\mathcal{I}\right) \left|
_Z\right. $ contained in $\mathbf{\iota }_{*}\left( \mathcal{O}_{T_N}\right) 
$ \emph{(}where $\mathbf{\iota }:T_N\hookrightarrow X\left( N,\Delta \right) 
$ denotes the canonical injection\emph{),} and 
\[
\emph{Norm}\,\left[ \mathbf{Bl}_Z^{\mathcal{I}}\left( X\left( N,\Delta
\right) \right) \right] \stackrel{\pi _{\mathcal{I}}\circ \mathbf{\nu }_{%
\mathcal{I}}}{\longrightarrow }X\left( N,\Delta \right) 
\]
the normalized blow-up of $X\left( N,\Delta \right) $ along $Z$. Then $\pi _{%
\mathcal{I}}\circ \mathbf{\nu }_{\mathcal{I}}$ is the $T_N$-equivariant
holomorphic map\emph{\ }and the overlying space is identified with the toric
variety 
\[
\fbox{$
\begin{array}{ccc}
&  &  \\ 
& \emph{Norm}\,\left[ \mathbf{Bl}_Z^{\mathcal{I}}\left( X\left( N,\Delta
\right) \right) \right] =X\left( N,\Delta _{\mathbf{bl}}\left( \mathcal{I}%
;Z\right) \right)  &  \\ 
&  & 
\end{array}
$}\,
\]
where $\Delta _{\mathbf{bl}}\left( \mathcal{I};Z\right) $ denotes the
following refinement of $\Delta $ \emph{:\smallskip } 
\[
\Delta _{\mathbf{bl}}\left( \mathcal{I};Z\right) :=\left\{ 
\begin{array}{c}
\text{\emph{the fan in} }N_{\Bbb{R}}\ \text{\emph{defined by rational s.c.p.
cones}} \\ 
\text{\emph{which constitute the maximal subdivision }} \\ 
\text{\emph{(w.r.t. usual inclusion) of the cones of} }\Delta \text{\emph{,
\ }} \\ 
\text{\emph{so that the restriction of \ ord}}\left( \mathcal{F}\right) \ 
\\ 
\text{\emph{on each of them becomes an integral }} \\ 
\text{\textit{linear}\emph{\ support function}}
\end{array}
\right\} 
\]
\end{theorem}

\noindent \textit{Proof. }Since $\Delta _{\mathbf{bl}}\left( \mathcal{I}%
;Z\right) $ is a refinement of the initial fan $\Delta $, the toric map 
\[
\text{id}_{*}:X\left( N,\Delta _{\mathbf{bl}}\left( \mathcal{I};Z\right)
\right) \longrightarrow X\left( N,\Delta \right) 
\]
is a proper birational morphism by thm. \ref{PRBIR}. Since $\mathcal{F}%
=\left( \mathcal{O}_{X\left( N,\Delta \right) }\,/\,\mathcal{I}\right)
\left| _Z\right. $ is $T_N$-invariant, $T_N$ acts on $\mathbf{Bl}_Z^{%
\mathcal{I}}\left( X\left( N,\Delta \right) \right) $ too. $\mathcal{F}$ is
invertible over the open set $T_N$ of $X\left( N,\Delta \right) $. Hence,
the normalized blow-up is an isomorphism over $T_N,$ 
\[
T_N\hookrightarrow \mathbf{Bl}_Z^{\mathcal{I}}\left( X\left( N,\Delta
\right) \right) 
\]
is an equivariant immersion and $\pi _{\mathcal{I}}\circ \mathbf{\nu }_{%
\mathcal{I}}$ a torus-equivariant map. Combining the universal property of
the normalized blow-up (see \textsf{(e)}) with the fact, that for any $%
\sigma \in \Delta $, the inverse image of $\mathcal{I}$ restricted onto $%
U_\sigma $ is an invertible sheaf if and only if ord$_{\mathcal{F}}\left|
_\sigma \right. $ is linear (by thm. \ref{INVERT}), we deduce that id$%
_{*}=\pi _{\mathcal{I}}\circ \mathbf{\nu }_{\mathcal{I}}$ (up to isomorphism
over $X\left( N,\Delta \right) $). $_{\Box }$

\begin{remark}
\label{TBLPATH}\emph{(i) If }$X\left( N,\Delta \right) $\emph{\ is} \textit{%
smooth}, 
\[
\Delta \left( 1\right) =\left\{ \varrho _1,\ldots ,\varrho _k\right\} \emph{%
,\ \ \ Gen}\left( \Delta \right) =\left\{ n_1,\ldots ,n_k\right\} \emph{,}
\]
\emph{\ with }$n_i=n\left( \varrho _i\right) $\emph{, for }$1\leq i\leq k$%
\emph{,} 
\[
n_0:=n_1+\cdots +n_k
\]
\emph{and} 
\[
\tau _j:=\left\{ 
\begin{array}{lll}
\text{\emph{pos}}\left( n_0,n_1,\ldots ,n_{j-1},n_{j+1},\ldots ,n_k\right) 
& ,\ \ \ \text{\emph{for}} & 1\leq j\leq k-1 \\ 
&  &  \\ 
\text{\emph{pos}}\left( n_0,n_1,\ldots ,n_{k-1}\right)  & ,\ \ \ \text{\emph{%
for}} & \ j=k
\end{array}
\right. \smallskip 
\]
\emph{then, fixing a cone }$\tau \in \Delta $\emph{,} \emph{\ every cone }$%
\sigma \in \Delta $\emph{\ with }$\tau \prec \sigma $\emph{\ can be written
as } 
\[
\sigma =\tau +\sigma ^{\prime }\text{\emph{, \ \ for some cone} }\sigma
^{\prime }\in \Delta \text{\emph{, }\ \ \emph{with} }\sigma ^{\prime }\cap
\tau =\left\{ \mathbf{0}\right\} \ .
\]
\emph{Setting }$\sigma _j:=\tau _j+\sigma ^{\prime }$\emph{, for all }$j$%
\emph{, }$1\leq j\leq k$\emph{, and} 
\[
\Delta ^{*}\left( \tau \right) :=\left( \Delta \smallsetminus \left\{ \sigma
\in \Delta \ \left| \ \emph{\ }\tau \prec \sigma \right. \right\} \right)
\cup \left\{ \text{\emph{faces of } }\sigma _j\ \left| 
\begin{array}{l}
\ \sigma \in \Delta ,\ \tau \prec \sigma  \\ 
\text{\emph{for all}\ }j,1\leq j\leq k
\end{array}
\right. \right\} 
\]
\emph{we get} 
\[
X\left( N,\Delta ^{*}\left( \tau \right) \right) =\mathbf{Bl}_{V\left( \tau
\right) }^{\emph{red}}\left( X\left( N,\Delta \right) \right) 
\]
\emph{with }$\Delta ^{*}\left( \tau \right) =\Delta _{\mathbf{bl}}\left( 
\mathcal{I};V\left( \tau \right) \right) $ \emph{the fan ``starring'' }$\tau 
$\emph{, where }$\mathcal{I}$\emph{\ is now the usual ideal (with no
nilpotent elements) defining the closure }$V\left( \tau \right) $ \emph{as
subvariety of }$X\left( N,\Delta \right) $ \emph{(cf. Oda \cite{Oda}, prop.
1.26, pp. 38-39; see also Ewald \cite{Ewald}, \S\ VI.7, for an equivalent
combinatorial characterization in terms of ``stellar
subdivisions''.)\smallskip }\newline
\emph{(ii) If }$X\left( N,\Delta \right) $\emph{\ is} \textit{not} \emph{%
smooth, then performing the above starring subdivision w.r.t. }$\tau $\emph{%
, we get a normalized blow-up of }$X\left( N,\Delta \right) $ \emph{along }$%
V\left( \tau \right) $\emph{, with }$V\left( \tau \right) $\emph{\ being
equipped with a} \textit{not} \emph{necessarily reduced scheme structure!
For a simple example, see \ref{DIFFBL} below.\smallskip }\newline
\emph{(iii) Even if} $X\left( N,\Delta \right) $\emph{\ }\textit{is} \emph{%
smooth, applying theorem \ref{TORNBL} for a} \textit{non-reduced }\emph{%
subscheme }$Z$ \emph{as center, the resulting normal complex} \emph{variety} 
$\emph{Norm}\left[ \mathbf{Bl}_Z^{\mathcal{I}}\left( X\left( N,\Delta
\right) \right) \right] $ \emph{is} \textit{not} \emph{necessarily smooth.
The simplest example is to take }$N=\Bbb{Z}e_1\oplus \Bbb{Z}e_2$\emph{,} 
\emph{i.e., the standard lattice with }$e_1$\emph{, }$e_2$\emph{\ the unit
vectors, $M=\Bbb{Z}e_1^{\vee}\oplus\Bbb{Z}e_2^{\vee}$ its dual,} 
\[
\sigma =\emph{pos}\left( e_1,e_2\right) \emph{,\ \ }U_\sigma \cong \Bbb{C}^2%
\emph{,\ \ }Z=\left( 0,0\right) =\emph{orb}\left( \sigma \right) \emph{,\ \
and\ \ \ }\mathcal{I}=\left(\hbox{\bf e}(e_1^{\vee}),\hbox{\bf e}(2 e_2^{\vee})\right) .
\]
\emph{Then the blown up space corresponds to the fan consisting of the cones}
\[
\sigma _1=\emph{\ pos}\left( \left( 2,1\right) ^{\intercal },e_2\right) 
\emph{,\ \ \ }\sigma _2=\emph{\ pos}\left( e_1,\left( 2,1\right) ^{\intercal
}\right) \emph{,}
\]
\emph{together with their faces. Obviously, mult}$\left( \sigma _1;N\right)
=2$\emph{, and we rediscover the second example of \ref{2BLUP}. (In fact,
this is nothing but the so-called }\textit{weighted blow-up} \emph{of the
origin w.r.t. }$\left( 2,1\right) ^{\intercal }$\emph{\ in Reid's
terminology; see \cite{Reid1}, p. 297).\smallskip } \emph{\ } \newline
\emph{(iv) On the other hand, the normalization of the }\textit{usual}\emph{%
\ blow-up of a }\textit{not}\emph{\ necessarily smooth, affine toric variety}
$U_\sigma =$ \emph{Max-Spec}$\left( \Bbb{C}\left[ M\cap \sigma ^{\vee
}\right] \right) $ \emph{along a  subvariety can be described
intrinsically by making use of thm. \ref{TORNBL} and arguments coming from
an embedding. Next proposition treats of the case in which one blows up }$%
U_\sigma $ \emph{at orb}$\left( \sigma \right) $.
\end{remark}

\begin{proposition}[Usual normalized blow-up at the closed point]
\label{USNBL} \ \smallskip \newline
Let $N$ be a lattice of rank $r$, $\sigma $ an $r$-dimensional rational
s.c.p. cone in $N_{\Bbb{R}}$ and 
\[
U_\sigma =\emph{Max-Spec}\left( \Bbb{C}\left[ M\cap \sigma ^{\vee }\right]
\right) 
\]
the associated affine toric variety with $M=$ \emph{Hom}$_{\Bbb{Z}}\left( N,%
\Bbb{Z}\right) $. Moreover, let 
\[
\mathbf{Hlb}_M\left( \sigma ^{\vee }\right) =\left\{ m_1,\ldots ,m_d\right\}
\smallskip ,\ \ r\leq d,
\]
denote a fixed enumeration of the members of the Hilbert basis of $\sigma
^{\vee }$ w.r.t. $M$ \emph{(cf. prop. \ref{MINGS}). }Then\smallskip 
\[
\text{\emph{Norm}}\,\left[ \mathbf{Bl}_{\text{\emph{orb}}\left( \sigma
\right) }^{\emph{red}}\left( U_\sigma \right) \right] =X\left( N,\Delta _{%
\mathbf{bl}}\left[ \emph{orb}\left( \sigma \right) \right] \right) \,,
\]
where 
\[
\Delta _{\mathbf{bl}}\left[ \emph{orb}\left( \sigma \right) \right] :=\Delta
_{\mathbf{bl}}\left( \left( 
\begin{array}{c}
\text{\emph{maximal }} \\ 
\text{\emph{ideal}}
\end{array}
\right) ;\ \emph{orb}\left( \sigma \right) \right) =\left\{ 
\begin{array}{c}
\left\{ \sigma _j\,\left| \,1\leq j\leq d\right. \right\} \text{\emph{%
together}} \\ 
\text{\emph{with all their faces}}
\end{array}
\right\},
\]
with 
\[
\sigma _j:=\left\{ \mathbf{y}\in \sigma \ \left| \ \left\langle m_i-m_j,%
\mathbf{y}\right\rangle \geq 0,\ \forall i,\ i\in \left\{ 1,\ldots
,j-1,j+1,\ldots ,d\right\} \right. \right\} ,\smallskip \text{ }
\]
for all $j,\ 1\leq j\leq d$. \emph{(Warning. Though the union of the above }$%
\sigma _j$\emph{'s forms always a fan, it might happen that }$\sigma
_j=\sigma _{j^{\prime }}$ \emph{or that }$\sigma _j$\emph{\ is a proper face
of }$\sigma _{j^{\prime }}$\emph{, for certain distinct indices }$%
j,j^{\prime }$\emph{\ belonging to }$\left\{ 1,\ldots ,d\right\} $\emph{. In
this case, we just ignore the superfluous cones and introduce a new
index-enumeration for the rest, preferably by considering only the
maximal-dimensional ones.)}
\end{proposition}

\noindent \textit{Proof. }Under the above assumptions, we may use prop. \ref
{EMB} and embed $U_\sigma $ by $\left( \mathbf{e}\left( m_1\right) ,\ldots ,%
\mathbf{e}\left( m_d\right) \right) $ ``minimally'' into the affine complex
space 
\[
\Bbb{C}^d=\text{Max-Spec}\left( \Bbb{C}\left[ \widetilde{M}\cap \widetilde{%
\sigma }^{\vee }\right] \right) , 
\]
where 
\[
\widetilde{\sigma }^{\vee }:=\text{ pos}\left( \left\{ m_1,\ldots
,m_d\right\} \right) \subset \left( \widetilde{M}\right) _{\Bbb{R}}\text{, \
\ \ }\widetilde{M}:=\Bbb{Z\,}m_1\oplus \cdots \oplus \Bbb{Z\,}m_d\ . 
\]
orb$\left( \sigma \right) $ is mapped onto $\mathbf{0=}$ orb$\left( 
\widetilde{\sigma }\right) \in \Bbb{C}^d$, with 
\[
\widetilde{\sigma }=\widetilde{\sigma }^{\vee \,\vee }\subset \left( 
\widetilde{N}\right) _{\Bbb{R}},\ \ \widetilde{N}=\text{Hom}_{\Bbb{Z}}\left( 
\widetilde{M},\Bbb{Z}\right) \cong \Bbb{Z}^d. 
\]
Using the embedding $N\hookrightarrow \widetilde{N}$, it is possible to
describe $\mathbf{Bl}_{\text{orb}\left( \sigma \right) }^{\text{red}}\left(
U_\sigma \right) $ as the strict transform of $U_\sigma $ under the usual
blow-up morphism of $\Bbb{C}^d$ at the origin (i.e. just by applying what we
mentioned in \textsf{(d) }for $Z=$ orb$\left( \sigma \right) $, $W=U_\sigma $
and $X=\Bbb{C}^d$). Hence, we obtain the following commutative diagrams of
torus-equivariant holomorphic maps: 
\[
\begin{array}{lllllll}
&  & \text{Norm}\,\left[ \mathbf{Bl}_{\text{orb}\left( \sigma \right) }^{%
\text{red}}\left( U_\sigma \right) \right] &  &  &  &  \\ 
&  & \downarrow & \searrow &  &  &  \\ 
\text{Exc}\left( \pi \left| _{\text{restr.}}\right. \right) & \in & \mathbf{%
Bl}_{\text{orb}\left( \sigma \right) }^{\text{red}}\left( U_\sigma \right) & 
\hookrightarrow & \mathbf{Bl}_{\mathbf{0}}^{\text{red}}\left( \Bbb{C}%
^d\right) & \ni & \text{Exc}\left( \pi \right) \\ 
&  & \downarrow \ \pi \left| _{\text{restr.}}\right. &  & \downarrow \ \pi & 
&  \\ 
\text{orb}\left( \sigma \right) & \in & U_\sigma & \hookrightarrow & \Bbb{C}%
^d & \ni & \mathbf{0}
\end{array}
\]
Since $\mathbf{Bl}_{\mathbf{0}}^{\text{red}}\left( \Bbb{C}^d\right) $ is
realized as toric variety by the $d$ cones of the barycentric subdivision of 
$\widetilde{\sigma }$, 
\[
\tau _j:=\left\{ \mathbf{y}\in \widetilde{\sigma }\ \left| \ \left\langle
m_i-m_j,\mathbf{y}\right\rangle \geq 0,\ \forall i,\ i\in \left\{ 1,\ldots
,j-1,j+1,\ldots ,d\right\} \right. \right\} ,\ 1\leq j\leq d, 
\]
(cf. \ref{TBLPATH} (i)), the above defined $\sigma _j$'s are exactly the
restrictions of $\tau _j$'s on $\sigma $ and determine obviously $\Delta _{%
\mathbf{bl}}\left[ \text{orb}\left( \sigma \right) \right] $ as it was given
initially in thm. \ref{TORNBL}. $_{\Box }\bigskip $\newline
\textsf{(i) }\textit{Resolutions of toric singularities. }To resolve toric
singularities is actually equivalent to subdividing simplicial cones into
others of smaller multiplicity.

\begin{theorem}[Resolution of toric singularities. Weak version]
\label{RESW} \ \smallskip \newline
If $X\left( N,\Delta \right) $ is an arbitrary toric variety, then there
exists a refinement $\widetilde{\Delta }$ of $\Delta $, such that 
\[
f=\emph{id}_{*}:X\left( N,\widetilde{\Delta }\right) \longrightarrow X\left(
N,\Delta \right) 
\]
is a \emph{(}full\emph{) }desingularization of $X\left( N,\Delta \right) $.
\end{theorem}

\noindent \textit{Sketch of proof. }Considering the multiplicity of a
simplicial cone $\sigma $\ as a volume and using the well-behaved volume
properties under subdivisions w.r.t. lattice points of $\mathbf{Par}\left(
\sigma \right) \cap N_\sigma $, one can easily desingularize equivariantly
any toric variety $X\left( N,\Delta \right) $:\smallskip \newline
$\bullet $ At first we refine the cones of $\Delta $\ in order to make it
simplicial. That this is always possible follows basically from
Carath\'{e}odory's theorem concerning convex polyhedral cones (cf. \cite
{Schrijver}, p. 94; for a simple proof see also G.Ewald \cite{Ewald}, III
2.6, p. 75, and V 4.2, p. 158).\smallskip \newline
$\bullet $ In the second step this new simplicial $\Delta $ will be
subdivided further into subcones of strictly smaller multiplicities than
those of the cones of the starting-point. After finitely many subdivisions
of this kind one can construct a refinement $\Delta ^{\prime }$ of $\Delta $%
, so that $f=$ id$_{*}$ gives rise to a resolution of singularities of $%
X\left( N,\Delta \right) $ (by thm. \ref{PRBIR} and prop. \ref{SMCR}). $%
_{\Box }$ \medskip

\noindent If fact, for toric varieties, a \textit{single} normalized blow-up
of a suitable ideal sheaf is able to provide a full resolution.

\begin{theorem}[Resolution of toric singularities. Strong version]
\label{RESSTR} \ \newline
If $X\left( N,\Delta \right) $ is an arbitrary toric variety, then there
exists a $T_N$-invariant coherent sheaf $\mathcal{I}$ of $\mathcal{O}%
_{X\left( N,\Delta \right) }$-ideals with 
\[
\emph{Sing}\left( X\left( N,\Delta \right) \right) =\emph{supp}\left( 
\mathcal{O}_{X\left( N,\Delta \right) }\,/\,\mathcal{I}\right) ,
\]
such that 
\[
\emph{Norm}\,\left[ \mathbf{Bl}_{\text{\emph{Sing}}\left( X\left( N,\Delta
\right) \right) }^{\mathcal{I}}\left( X\left( N,\Delta \right) \right)
\right] \stackrel{\pi _{\mathcal{I}}\circ \mathbf{\nu }_{\mathcal{I}}=\emph{%
id}_{*}}{\longrightarrow }X\left( N,\Delta \right) 
\]
forms a \emph{(}full, projective\emph{) }desingularization of $X\left(
N,\Delta \right) $. \emph{(}$\mathcal{I}$ is in general not uniquely
determined by this property and might contain nilpotent elements\emph{). }
\end{theorem}

\noindent \textit{Proof. }See Saint-Donat \cite{KKMS}, thm. 11, pp. 32-35,
and Brylinski \cite{Bryl}, pp. 273-279.$_{\Box }$

\section{Toric $\Bbb{P}^k$-bundles over projective spaces\label{Bundles}}

\noindent This section is a brief excursus to a part of the theory of toric
bundles over projective spaces which will be used later on (in \S\ \ref
{MONSE}) for the precise description of the exceptional prime divisors
occuring in our desingularizations.\medskip \newline
\textsf{(a) }An equivariant holomorphic map $\varpi _{*}:X\left( N,\Delta
\right) \rightarrow X\left( N^{\prime },\Delta ^{\prime }\right) $ of toric
varieties, induced by a map of fans $\varpi :\left( N,\Delta \right)
\rightarrow \left( N^{\prime },\Delta ^{\prime }\right) $ (cf. \S \ref{TORIC}
\textsf{(k)}), can be viewed as the projection map of an equivariant fiber
bundle (\textit{toric bundle}) over $X\left( N^{\prime },\Delta ^{\prime
}\right) $ with typical fiber $X\left( N^{\prime \prime },\Delta ^{\prime
\prime }\right) $, $N^{\prime \prime }=$ Ker$\left( \varpi :N\rightarrow
N^{\prime }\right) $, if and only if $\varpi :N\rightarrow N^{\prime }$ is
surjective and the cones of $\Delta $ are representable as ``joins'' of the
cones of a fan $\Delta ^{\prime \prime }$ ($\left| \Delta ^{\prime \prime
}\right| \subset N_{\Bbb{R}}^{\prime \prime }$) with those of another fan $%
\widetilde{\Delta }^{\prime }\subset \Delta $, so that the supports $\left| 
\widetilde{\Delta }^{\prime }\right| $ and $\left| \Delta ^{\prime }\right| $
are homeomorphic to each other. (See Oda \cite{Oda}, prop. 1.33, p. 58, and
Ewald \cite{Ewald}, thm. VI.6.7, p. 246). In the case in which its total
space is assumed to be smooth and compact, this criterion can be
considerably simplified by means of the notion of ``primitive collections''
introduced by Batyrev \cite{Bat}.

\begin{definition}
\emph{If }$X\left( N,\Delta \right) $\emph{\ is an }$r$\emph{-dimensional
smooth, compact toric variety, then a non-empty subset }$\frak{N}=\left\{
n_1,n_2,\ldots ,n_k\right\} $\emph{\ of Gen}$\left( \Delta \right) $\emph{, }%
$k\geq 2,$\emph{\ is defined to be a }\textit{primitive collection}\emph{\
if it satisfies anyone of the following equivalent conditions :\smallskip }%
\newline
\emph{(i) For each }$n_i\in \frak{N}$\emph{, }$1\leq i\leq k,$\emph{\ one
has }$\frak{N}\smallsetminus \left\{ n_i\right\} =$\emph{\ Gen}$\left(
\sigma _i\right) $\emph{, for some }$\sigma _i$ \emph{belonging to }$\Delta
\left( k-1\right) $\emph{, while }$\frak{N}$ \emph{\ itself cannot be the
set of minimal generators of any }$k$\emph{-dimensional cone of }$\Delta $%
\emph{.\smallskip \newline
(ii) For each subset of indices }$\left\{ j_1,\ldots ,j_q\right\} \subset
\left\{ 1,\ldots ,k\right\} ,1\leq q< k,$ \emph{the set }$\left\{
n_{j_1},\ldots ,n_{j_q}\right\} $\emph{\ coincides with the set of minimal
generators of a }$q$\emph{-dimensional cone of }$\Delta $\emph{, while }$%
\frak{N}$ \emph{\ itself cannot be the set of minimal generators of any }$k$%
\emph{-dimensional cone of }$\Delta $\emph{. }
\end{definition}

\begin{proposition}[Characterization of toric $\QTR{Bbb}{P}_{\QTR{Bbb}{C}}^k$-bundles]
\label{chbudl}Let $X\left( N,\Delta \right) $ be a smooth, compact toric
variety of dimension $r$ and $k$ a positive integer $\leq r$. $X\left(
N,\Delta \right) $ is the total space of a toric $\Bbb{P}_{\Bbb{C}}^k$%
-bundle over a smooth $\left( r-k\right) $-dimensional toric variety if and
only if there exists a primitive collection 
\[
\frak{N}=\left\{ n_1,n_2,\ldots ,n_{k+1}\right\} \subset \emph{Gen}\left(
\Delta \right) \ ,
\]
such that\smallskip \newline
\emph{(i)} $n_1+n_2+\cdots +n_{k+1}=\mathbf{0}_N$, and\smallskip \newline
\emph{(ii)} $\frak{N}\cap \frak{N}^{\prime }=\varnothing ,$ for all
primitive collections $\frak{N}^{\prime }\subset $ \emph{Gen}$\left( \Delta
\right) $ with $\frak{N}\neq \frak{N}^{\prime }.$
\end{proposition}

\begin{definition}
\emph{Let X}$\left( N,\Delta \right) $\emph{\ be a smooth, compact toric
variety. The fan }$\Delta $\emph{\ is called a }\textit{splitting fan }\emph{%
if any two different primitive collections within Gen}$\left( \Delta \right) 
$ \emph{have no common elements.}
\end{definition}

\begin{theorem}
\label{splfan}If $\Delta $ is a splitting fan, then $X\left( N,\Delta
\right) $ is a projectivization of a decomposable bundle over a toric
variety being associated to a splitting fan of smaller dimension.
\end{theorem}

\noindent \textit{Proof. }See Batyrev \cite{Bat}, thm. 4.3, p. 577. $_{\Box
}\medskip \smallskip $

\noindent \textsf{(b) }The projectivized decomposable bundles over
projective spaces, having only twisted hyperplane bundles as summands, can
be easily described as toric bundles in terms of splitting fans with exactly
two disjoint primitive collections. In particular, applying \ref{chbudl} and 
\ref{splfan}, we obtain:

\begin{lemma}
\label{YPSILON}Let $k$ and $s$ be two positive integers, $\left( \lambda
_1,\lambda _2,\ldots ,\lambda _k\right) $ a $k$-tuple of non-negative
integers and $N_{\left( r;\lambda _1,\lambda _2,\ldots ,\lambda _k\right) }$ 
\emph{(}resp. $N_s^{\prime }$\emph{) }a lattice of rank $r=k+s$ \emph{(}%
resp. of rank $s$\emph{)} generated by $\frak{N}\cup \frak{N}^{\prime }$ 
\emph{(}resp. by $\frak{N}^{\prime }$\emph{) }where 
\[
\frak{N}:=\left\{ \text{\ }n_1,n_2,\ldots ,n_{k+1}\right\} ,\ \ \frak{N}^{\
\prime }:=\left\{ n_1^{\prime },n_2^{\prime },\ldots ,n_s^{\prime
},n_{s+1}^{\prime }\right\} 
\]
with the two relations 
\[
n_1+n_2+\cdots +n_k+n_{k+1}=\mathbf{0,}\text{ \ \ }n_1^{\prime }+n_2^{\prime
}+\cdots +n_{s+1}^{\prime }=\lambda _1n_1+\cdots +\lambda _kn_k\ .
\]
If we define the $r$-dimensional, smooth, compact toric variety\smallskip 
\[
\fbox{$
\begin{array}{ccc}
&  &  \\ 
& Y\left( r;\lambda _1,\lambda _2,\ldots ,\lambda _k\right) :=X\left(
N_{\left( r;\lambda _1,\lambda _2,\ldots ,\lambda _k\right) },\ \Delta
_{\left( r;\lambda _1,\lambda _2,\ldots ,\lambda _k\right) }\right)  &  \\ 
&  & 
\end{array}
$}
\]
by means of the fan 
\[
\Delta _{\left( r;\lambda _1,\lambda _2,\ldots ,\lambda _k\right) }:=\left\{ 
\begin{array}{l}
\text{\emph{pos}}\left( \left( \frak{N}\cup \frak{N}^{\prime }\right)
\smallsetminus \left\{ n_i,n_j^{\prime }\right\} \right) ,\smallskip  \\ 
\text{\emph{for all }}\left( i,j\right) \in \left\{ 1,\ldots ,k+1\right\}
\times \left\{ 1,\ldots ,s+1\right\} ,\smallskip  \\ 
\text{\emph{together with all their faces}}
\end{array}
\right\} \ ,
\]
then it is isomorphic to the total space of the $\Bbb{P}_{\Bbb{C}}^k$-bundle 
\[
Y\left( r;\lambda _1,\lambda _2,\ldots ,\lambda _k\right) \cong \Bbb{P}%
\left( \mathcal{O}_{\Bbb{P}_{\Bbb{C}}^s}\oplus \mathcal{O}_{\Bbb{P}_{\Bbb{C}%
}^s}\left( \lambda _1\right) \oplus \mathcal{O}_{\Bbb{P}_{\Bbb{C}}^s}\left(
\lambda _2\right) \oplus \cdots \oplus \mathcal{O}_{\Bbb{P}_{\Bbb{C}%
}^s}\left( \lambda _k\right) \right) \rightarrow \Bbb{P}_{\Bbb{C}}^s\text{ }
\]
over $\Bbb{P}_{\Bbb{C}}^s=X\left( N_s^{\prime },\Delta _s^{\prime }\right) ,$
with 
\[
\Delta _s^{\prime }:=\left\{ 
\begin{array}{l}
\text{\emph{pos}}\left( \frak{N}^{\prime }\smallsetminus \left\{ n_j^{\prime
}\right\} \right) \ \text{\emph{in} }\left( N_s^{\prime }\right) _{\Bbb{R}%
},\smallskip  \\ 
\text{\emph{for all }}j\in \left\{ 1,\ldots ,s+1\right\} ,\smallskip  \\ 
\text{\emph{together with all their faces}}
\end{array}
\right\} \ .\text{ }
\]
\end{lemma}

\noindent \textit{Proof. }By construction, $\frak{N}$ and $\frak{N}^{\
\prime }$ are the only primitive collections within Gen$\left( \Delta
_{\left( r;\lambda _1,\ldots ,\lambda _k\right) }\right) $, and $\frak{N}%
\cap \frak{N}^{\ \prime }=\varnothing $. Hence, $\Delta _{\left( r;\lambda
_1,\ldots ,\lambda _k\right) }$ is a splitting fan. Since $\#\left( \text{Gen%
}\left( \Delta _s^{\prime }\right) \right) =\#\left( \Delta _s^{\prime
}\left( 1\right) \right) =s+1$, the base-space of the smooth toric $\Bbb{P}_{%
\Bbb{C}}^k$-bundle $Y\left( r;\lambda _1,\ldots ,\lambda _k\right) $ has to
be isomorphic to $\Bbb{P}_{\Bbb{C}}^s$. Finally, the decomposable bundle
over $X\left( N_s^{\prime },\Delta _s^{\prime }\right) $ (whose
projectivization gives the total space $Y\left( r;\lambda _1,\ldots ,\lambda
_k\right) $) is isomorphic to $\mathcal{O}_{\Bbb{P}_{\Bbb{C}}^s}\oplus 
\mathcal{O}_{\Bbb{P}_{\Bbb{C}}^s}\left( \lambda _1\right) \oplus \cdots
\oplus \mathcal{O}_{\Bbb{P}_{\Bbb{C}}^s}\left( \lambda _k\right) $ because
the fan corresponding to the typical fiber consists of cones which are the
images of the cones of $\Delta _s^{\prime }$ under the linear map 
\[
\left( N_s^{\prime }\right) _{\Bbb{R}}\hookrightarrow \left( N_{\left(
r;\lambda _1,\ldots ,\lambda _k\right) }\right) _{\Bbb{R}} 
\]
sending a $\mathbf{y}^{\prime }$ to $\left( \mathbf{y}^{\prime },-\left(
\lambda _1n_1+\cdots +\lambda _kn_k+n_{k+1}\right) \right) $. $_{\Box }$

\begin{theorem}[Classification theorem of Kleinschmidt]
\label{KLCL} \ \smallskip \newline
Every smooth, compact $r$-dimensional toric variety \emph{(}$r\geq 2$\emph{)}
with Picard number $2$, is isomorphic to one of the varieties $Y\left(
r;\lambda _1,\lambda _2,\ldots ,\lambda _k\right) $.
\end{theorem}

\noindent \textit{Proof. }See \cite{Klein}, \S\ 2, pp. 256-261. $_{\Box }$

\begin{example}[Hirzebruch surfaces]
\emph{Setting} $k=1,\lambda =\lambda _1,$\emph{\ and} $r=2,$ \emph{one gets
the rational scrolls } 
\[
\Bbb{F}_\lambda :=Y\left( 2;\lambda \right) \cong \Bbb{P}\left( \mathcal{O}_{%
\Bbb{P}_{\Bbb{C}}^1}\oplus \mathcal{O}_{\Bbb{P}_{\Bbb{C}}^1}\left( \lambda
\right) \right) 
\]
\emph{over }$\Bbb{P}_{\Bbb{C}}^1$\emph{, i.e. the so-called }\textit{%
Hirzebruch surfaces} \emph{whose topological, analytic and birational
properties were studied in the early fifties in \cite{Hirz1}. (Certain
birational properties of them were already investigated by Segre and Del
Pezzo around the end of the last century in connection with other types of
scrolls and ruled surfaces. See e.g. Segre \cite{Segre}). It is well-known
that all }$\Bbb{F}_\lambda $\emph{'s, }$\lambda \neq 1$\emph{, together with 
}$\Bbb{P}_{\Bbb{C}}^2$\emph{,\ exhaust the class of all minimal, smooth,
rational, projective complex surfaces, and that }$\Bbb{F}_\lambda $ \emph{%
can be considered as the zero-set :} \emph{\ } 
\begin{equation}
\left\{ \left( \left[ z_0:z_1:z_2\right] ,\left[ t_1,t_2\right] \right) \in 
\Bbb{P}_{\Bbb{C}}^2\times \Bbb{P}_{\Bbb{C}}^1\ \left| \ \ z_1\cdot
t_1^\lambda -z_2\cdot t_2^\lambda =0\right. \text{\emph{\ }}\right\} \ 
\label{HIRZEMB}
\end{equation}
$\Bbb{F}_0$\emph{\ is therefore }$\Bbb{P}_{\Bbb{C}}^1\times \Bbb{P}_{\Bbb{C}%
}^1$\emph{\ , }$\Bbb{F}_1$\emph{\ is isomorphic to the usual blow-up of }$%
\Bbb{P}_{\Bbb{C}}^2$\emph{\ at a (}$T_{N_{\left( 2;\lambda \right) }}$\emph{%
-fixed) point, and for }$\lambda \geq 2$ \emph{there exists a natural }$%
\lambda $\emph{-sheeted ramified covering }$\Bbb{F}_\lambda \rightarrow \Bbb{%
F}_1$\emph{\ over }$\Bbb{F}_1$ \emph{(see \cite{Hirz1}, p. 82). The
differential-topological and diffeomorphism classification theory, and the
deformation theory of the }$\Bbb{P}_{\Bbb{C}}^k$\emph{-bundles} $Y\left(
k+1;\lambda _1,\ldots ,\lambda _k\right) $ \emph{over} $\Bbb{P}_{\Bbb{C}}^1$%
\emph{\ were developed in Brieskorn's work \cite{Br1} in the sixties. (For
another, purely geometric approach to the theory of rational scrolls over} $%
\Bbb{P}_{\Bbb{C}}^1$\emph{, see Reid \cite{Reid6}).}
\end{example}

\noindent Generalizing the embedding (\ref{HIRZEMB}) for arbitrary $k$ and $%
s $, one obtains :

\begin{proposition}[Bihomogeneous binomial representation]
The toric variety $Y\left( r;\lambda _1,\ldots ,\lambda _k\right) $ is
embeddable into the projective space $\Bbb{P}_{\Bbb{C}}^{\left( s+1\right)
k}\times \Bbb{P}_{\Bbb{C}}^s$\emph{; }in fact, if 
\[
\left\{ z_0,z_{1,1},z_{2,1},\ldots ,z_{s+1,1},z_{1,2},\ldots
,z_{s+1,2},\ldots ,z_{1,k},\ldots ,z_{s+1,k};t_1,t_2,\ldots ,t_{s+1}\right\} 
\]
denote bihomogeneous coordinates, it is representable as the zero-set 
\[
\left\{ z_{\nu ,i}\cdot t_\mu ^{\lambda _i}-z_{\mu ,i}\cdot t_\nu ^{\lambda
_i}=0\ \left| \ 
\begin{array}{l}
\text{\emph{for all triples\ }}\left( \mu ,\nu ,i\right) \ \text{\emph{with}}
\\ 
1\leq i\leq k,\ 1\leq \mu ,\nu \leq s+1,\ \mu \neq \nu 
\end{array}
\right. \text{\emph{\ }}\right\} 
\]
\end{proposition}

\noindent It is possible to embed these, sometimes called \textit{%
Hirzebruch-Kleinschmidt varieties}, $Y\left( r;\lambda _1,\ldots ,\lambda
_k\right) $, into a single projective space by using the Segre-embedding
\[
\Bbb{P}_{\Bbb{C}}^{\left( s+1\right) k}\times \Bbb{P}_{\Bbb{C}%
}^s\hookrightarrow \Bbb{P}_{\Bbb{C}}^{\left( \left( s+1\right) k+1\right)
\left( s+1\right) -1}
\]
but as it was proved by Ewald and Schmeinck in \cite{Ew-Sch}, this can be
done in a more economical way (w.r.t. the degrees of the defining
homogeneous binomials), namely by only considering quadrics within a
suitably higher dimensional projective space.

\begin{theorem}[Representation by quadrics]
\label{EWS} \ \smallskip \newline
The toric varieties $Y\left( r;\lambda _1,\ldots ,\lambda _k\right) $ are
embeddable into the projective space $\Bbb{P}_{\Bbb{C}}^d$ of dimension 
\[
d=r-k+\sum_{i=1}^k\ \dbinom{\lambda _i+r-k+1}{r-k}
\]
\emph{(}depending on $\lambda _i$'s\emph{) }and their defining ideals \emph{(%
}w.r.t. homogeneous coordinates $\left\{ z_0,z_1,\ldots ,z_d\right\} $
of $\Bbb{P}_{\Bbb{C}}^d$\emph{)} are generated by the quadratic
binomials 
\[
z_{\nu _1}\ z_{\nu _2}-z_{\nu _3}\ z_{\nu _4},\ \ \text{\emph{for all \ \ }}%
\left( \nu _1,\nu _2,\nu _3,\nu _4\right) \in \left( \left\{ 0,1,\ldots
,d\right\} \right) ^4\ .\smallskip 
\]
\end{theorem}

\noindent \textsf{(c) }The intersection theory on the varieties $Y\left(
r;\lambda _1,\ldots ,\lambda _k\right) $ is more complicated compared with
that on $\Bbb{F}_\lambda $'s. For example, although for $s=1$, the Chern
numbers of the total spaces of these $\Bbb{P}_{\Bbb{C}}^k$-bundles over $%
\Bbb{P}_{\Bbb{C}}^1$ for fixed dimension $r=k+1$ are constant (cf. Brieskorn 
\cite{Br1}, Satz 2.4.(ii), p. 348), this is no more true for decomposable $%
\Bbb{P}_{\Bbb{C}}^k$-bundles over $\Bbb{P}_{\Bbb{C}}^s$ with $s\geq 2$,
because there is an obvious dependence on $\lambda _i$'s, i.e. on the given
``twisting numbers''. On the other hand, the isomorphism for the Picard
group 
\[
\text{Pic}\left( Y\left( r;\lambda _1,\ldots ,\lambda _k\right) \right)
\cong \text{ Pic}\left( \Bbb{P}_{\Bbb{C}}^s\right) \times \Bbb{Z\cong
Z\times Z} 
\]
is in general valid (cf. \cite{Hart}, Ex. 7.9, p. 170), and is useful as
long as one makes a specific choice of a $\Bbb{Z}$-basis and expresses each
member of any examined $r$-tuple of divisors as concrete $\Bbb{Z}$-linear
combinations of its two elements. (The most natural choice is to consider
the classes in Pic corresponding to a typical fiber of the bundle map and to
a hyperplane section under the embedding of $Y\left( r;\lambda _1,\ldots
,\lambda _k\right) $ in $\Bbb{P}_{\Bbb{C}}^d$ of thm. \ref{EWS},
respectively). Let us now give, in broad outline, three complementary
practical methods for the computation of intersection numbers.\medskip

\noindent $\bullet $ \textit{First method.} For arbitrary $k$ and $s\geq 1$,
there exist two towers of birational morphisms 
\[
\begin{array}{ll}
Y_0\longleftarrow Y_1\longleftarrow Y_2\longleftarrow \cdots \longleftarrow
Y_{\nu -1}\longleftarrow & Y_\nu =W_\mu \\ 
\parallel & \downarrow \\ 
Y\left( r;\lambda _1,\ldots ,\lambda _k\right) & W_{\mu -1} \\ 
& \downarrow \\ 
& \,\vdots \\ 
& \downarrow \\ 
& W_1 \\ 
& \downarrow \\ 
& W_0=\Bbb{P}_{\Bbb{C}}^r
\end{array}
\]
which are nothing but usual toric blow-ups with toric subvarieties of
codimension $\geq 2$ as centers. (For the precise description of this
procedure in terms of convex geometry, i.e. via barycentric stellar
subdivisions of cones of $\Delta _{\left( r;\lambda _1,\lambda _2,\ldots
,\lambda _k\right) }$, and the algorithmic determination of $\nu $, $\mu
\geq 0$, the reader is referred to Kleinschmidt \cite{Klein}, pp. 264-265).
So the evaluation of intersection numbers of divisors on $Y\left( r;\lambda
_1,\ldots ,\lambda _k\right) $ can be reduced to another one w.r.t. divisors
sitting on $\Bbb{P}_{\Bbb{C}}^r$. The problem here is that one blows up and
down subvarieties of varying dimensions and must therefore control carefully
the intersection behaviour of the proper transforms of divisors in each
step. The simplest example is $Y\left( r;1\right) $ (with $k=1$, $r=s+1$)
which is $\Bbb{P}_{\Bbb{C}}^r$ blown up at a ($T_{N\left( r;1\right) }$%
-fixed) point (with $\nu =0$, $\mu =1$); but, for instance, already $Y\left(
3;2\right) $ is the blow-up of $\Bbb{P}_{\Bbb{C}}^3$ along an entire ($%
T_{N\left( 3;2\right) }$-fixed) curve followed by the contraction of another
($T_{N\left( 3;2\right) }$-fixed) curve (i.e., $\nu =$ $\mu =1$, in this
case).\medskip

\noindent $\bullet $ \textit{Second method. }This method can be applied to
any projective toric variety (or even to any compact toric variety), but
demands familiarity with mixed volumes of ``virtual'' lattice polytopes.
(Virtual polytopes are defined to be finite families of suitably translated
dual cones of cones of a given fan, which are equipped with a $\Bbb{Z}$%
-module structure w.r.t. formal addition and scalar multiplication, though
their intersection might be not a polytope in the usual sence). If all the
line bundles, being associated to the Cartier divisors whose intersection
number is to be computed, are generated by their global sections (in other
words, if the corresponding integral support functions are upper convex, cf.
thm. \ref{UCON}), then we may apply formula (\ref{MIXED}) and evaluate the
normalized mixed volumes of the arising lattice polytopes w.r.t. the dual
lattice. For arbitrary Cartier divisors, however, we need the combinatorial
version of ``moving lemma'', i.e., to write down each of their associated
support functions as the difference of two upper convex (or even strictly
upper convex) support functions, and to calculate afterwards the desired
intersection numbers as the mixed volume of the difference of two virtual
polytopes. (See in Ewald's book \cite{Ewald}, thm. V.\thinspace 5.15, pp.
175-177, and thm. VII.\thinspace 6.3, pp. 292-295. For an intrinsic,
algorithmic characterization of moving lemma of this kind, we refer to
Wessels' thesis \cite{Wessels}, Satz 3.2.4, p. 83, and 3.2.18, p. 89. In
contrast to \cite{Ewald}, he works directly with ``generalized'' mixed
volumes of virtual polytopes the normalizations of which take integer, but
not necessarily only non-negative values). A simple example: the
self-intersection number of the canonical divisor of $Y\left( 3;2\right) $.
Since $Y\left( 3;2\right) $ is a Fano 3-fold, $-K_{Y\left( 3;2\right) }$ is
ample, and it is easy to verify that the lattice $3$-polytope $%
P_{-K_{Y\left( 3;2\right) }}$ ($3$-polytope, in the usual sence), induced by
the anticanonical integral strictly upper convex support function, can be
realized (up to an affine integral transformation) by\smallskip 
\[
P_{-K_{Y\left( 3;2\right) }}=\text{conv\thinspace }\left( \left\{ 
\begin{array}{c}
\left( 
\begin{array}{r}
0 \\ 
-1 \\ 
-1
\end{array}
\right) ,\left( 
\begin{array}{r}
-1 \\ 
0 \\ 
-1
\end{array}
\right) ,\left( 
\begin{array}{r}
-1 \\ 
-1 \\ 
-1
\end{array}
\right) , \\ 
\  \\ 
\left( 
\begin{array}{r}
-1 \\ 
-1 \\ 
1
\end{array}
\right) ,\left( 
\begin{array}{r}
4 \\ 
-1 \\ 
1
\end{array}
\right) ,\left( 
\begin{array}{r}
-1 \\ 
4 \\ 
1
\end{array}
\right)
\end{array}
\right\} \right) \subset \Bbb{R}^3\ ,\smallskip 
\]
i.e., by the polar of a lattice bypiramid over a triangle or, equivalently,
by a \textit{triangular lattice prism}. Fig.~{\bf 3} shows $-P_{-K_{Y\left( 3;2\right) }}$.

\begin{figure*}[htbp]
\begin{center}
\input{fig3.pstex_t}
\makespace
 {\textbf{Fig. 3}}
\end{center}
\end{figure*}

\noindent An immediate calculation gives 
\[
\text{Vol}\left( P_{-K_{Y\left( 3;2\right) }}\right) =\frac{31}%
3\Longrightarrow K_{Y\left( 3;2\right) }^3=-\left( -K_{Y\left( 3;2\right)
}\right) ^3=-62 
\]
(by formula (\ref{SELF})).\medskip \smallskip

\noindent $\bullet $ \textit{Third method. }If $\mathcal{E}$ is a locally
free sheaf with $q=$ rk$\left( \mathcal{E}\right) $ defined over a smooth
projective complex variety $Z$ and 
\[
\pi :\Bbb{P}\left( \mathcal{E}\right) =\mathbf{Proj}\left( \text{Sym}%
^{\bullet }\,\left( \mathcal{E}\right) \right) \longrightarrow Z 
\]
the associated projective bundle, then by Grothendieck's direct construction
of the Chern classes of $\mathcal{E}$, 
\[
c_i\left( \mathcal{E}\right) \in A^i\left( Z\right) ,\ \ \ \ i\in \left\{
0,1,\ldots ,q\right\} , 
\]
i.e., by setting $c_0\left( \mathcal{E}\right) =1$ and 
\begin{equation}
\sum_{i=0}^q\ \left( -1\right) ^i\ \pi ^{*}c_i\left( \mathcal{E}\right)
\cdot \left[ E^{q-i}\right] =0  \label{Groth}
\end{equation}
(cf. Hartshorne \cite{Hart}, p. 429, or Fulton \cite{Ful}, 3.2.4, p. 55),
the Chow ring $A^{\bullet }\left( \Bbb{P}\left( \mathcal{E}\right) \right) $
can be viewed as a free $A^{\bullet }\left( Z\right) $-module generated by
the classes $\left[ 1\right] $, $\left[ E\right] $, $\ldots $, $\left[
E^{q-1}\right] $, where here $E$ denotes the divisor on $\Bbb{P}\left( 
\mathcal{E}\right) $ corresponding to $\mathcal{O}_{\Bbb{P}\left( \mathcal{E}%
\right) }\left( 1\right) $. The relationship between the Chern polynomial of
the tangent bundle of $\Bbb{P}\left( \mathcal{E}\right) $ and that of the
pullback of the tangent bundle of $Z$ follows from the relative tangent
bundle exact sequence (see Fulton \cite{Ful}, 3.2.11, p. 59). In particular,
computing the first Chern class, we deduce \textit{the canonical bundle
formula} : 
\begin{equation}
K_{\Bbb{P}\left( \mathcal{E}\right) }\stackunder{\text{lin}}{\sim }\pi
^{*}\left( K_Z+\text{det}\left( \mathcal{E}\right) \right) -qE  \label{CANLB}
\end{equation}
The equation (\ref{Groth}) for $Z=\Bbb{P}_{\Bbb{C}}^s$ and $\mathcal{E}=%
\mathcal{O}_{\Bbb{P}_{\Bbb{C}}^s}\oplus \mathcal{O}_{\Bbb{P}_{\Bbb{C}%
}^s}\left( \lambda _1\right) \oplus \cdots \oplus \mathcal{O}_{\Bbb{P}_{\Bbb{%
C}}^s}\left( \lambda _k\right) $ turns out to be a quite powerful tool for
manipulating intersection numbers. In the next proposition we compute two
basic self-intersection numbers for $Y\left( s+1;\lambda \right) $ (as $\Bbb{%
P}_{\Bbb{C}}^1$-bundle over $\Bbb{P}_{\Bbb{C}}^s$) which will be used in \S\ 
\ref{MONSE}, and leave to the reader as exercise to examine further (and
more general) examples of various $r$-tuples of divisors by applying the
above mentioned methods.

\begin{proposition}
For $k=1$, $r=s+1$, $\lambda _1=\lambda \neq 0$, the self-intersection
number of the canonical divisor of the total space of the $\Bbb{P}_{\Bbb{C}%
}^1$-bundle $Y\left( s+1;\lambda \right) \rightarrow \Bbb{P}_{\Bbb{C}}^s$ is
given by the formula 
\begin{equation}
K_{Y\left( r;\lambda \right) }^r=\sum_{i=0}^{r-1}\ \binom ri\ \left(
-2\right) ^{r-i}\ \left( \lambda -r\right) ^i\ \lambda ^{r-i-1}\ 
\label{Formula1}
\end{equation}
Moreover, the self-intersection number of the divisor $E=V\left( \text{\emph{%
pos}}\left( \left\{ n_2\right\} \right) \right) $ \emph{(}in the notation of 
\emph{\ref{YPSILON})}\ equals 
\begin{equation}
E^r=\lambda ^{r-1}  \label{Formula2}
\end{equation}
\end{proposition}

\noindent \textit{Proof. }For $Z=\Bbb{P}_{\Bbb{C}}^s$, $\mathcal{E}=\mathcal{%
O}_{\Bbb{P}_{\Bbb{C}}^s}\oplus \mathcal{O}_{\Bbb{P}_{\Bbb{C}}^s}\left(
\lambda \right) $, and $H$ a divisor of $\Bbb{P}\left( \mathcal{E}\right) $
associated to the pullback $\pi ^{*}c_1\left( \mathcal{E}\right) $, we have 
\begin{equation}
E\cdot \left( E-\lambda H\right) =E^2-\lambda \left( H\cdot E\right) =0
\label{EQ1}
\end{equation}
(by (\ref{Groth}) \& $c_2\left( \mathcal{E}\right) =0$) and 
\begin{equation}
H^r=0,\ \ \left( H^{r-1}\cdot E\right) =1  \label{EQ2}
\end{equation}
(by definition). Furthermore, by the canonical bundle formula (\ref{CANLB}), 
\[
K_{Y\left( r;\lambda \right) }\stackunder{\text{lin}}{\sim }\pi ^{*}\left(
K_{\Bbb{P}_{\Bbb{C}}^s}+\text{det}\left( \mathcal{E}\right) \right) -2E%
\stackunder{\text{lin}}{\sim }-\left( s+1\right) H+\lambda H-2E=\left(
\lambda -r\right) H-2E\ .
\]
For the self-intersection number we get 
\begin{equation}
K_{Y\left( r;\lambda \right) }^r=\left( \left( \lambda -r\right) H-2E\right)
^r=\sum_{i=0}^{r-1}\ \binom ri\ \left( -2\right) ^{r-i}\ \left( \lambda
-r\right) ^i\ \left( H^i\cdot E^{r-i}\right)   \label{Formula3}
\end{equation}
(by the first of equations (\ref{EQ2})), and applying successively (\ref{EQ1}%
) and the second of the equations (\ref{EQ2}), 
\begin{eqnarray*}
\left( H^i\cdot E^{r-i}\right)  &=&\left( H^i\cdot \left( \lambda H\right)
\cdot E^{r-i-1}\right) =\left( H^i\cdot \left( \lambda H\right) ^2\cdot
E^{r-i-2}\right) =\cdots \smallskip  \\
\cdots \ \cdots  &=&\left( H^i\cdot \left( \lambda H\right) ^{r-i-1}\cdot
E\right) =\lambda ^{r-i-1}\left( H^{r-1}\cdot E\right) =\lambda ^{r-i-1}\ .
\end{eqnarray*}
Thus, formula (\ref{Formula1}) follows directly from (\ref{Formula3}).
Finally, (\ref{Formula2}) is proved similarly via (\ref{EQ1}) and (\ref{EQ2}%
). $_{\Box }$

\section{Toric description of abelian quotient singularities\label{ABELQ}}

\noindent Abelian quotient singularities can be investigated by means of the
theory of toric varieties in a direct manner. If $G$ is a finite subgroup of
GL$\left( r,\Bbb{C}\right) $, then $\left( \Bbb{C}^{*}\right) ^r/G$ is
automatically an algebraic torus embedded in $\Bbb{C}^r/G$.\medskip

\noindent \textit{Notation. }We shall henceforth use the following extra
notation. For $\nu \in \Bbb{N}$, $\mu \in \Bbb{Z}$, we denote by $\left[
\,\mu \,\right] _\nu $ the (uniquely determined) integer for which 
\[
0\leq \left[ \,\mu \,\right] _\nu <\nu ,\ \ \ \,\mu \,\equiv \left[ \,\mu
\,\right] _\nu \left( \text{mod }\nu \right) . 
\]
If $q\in \Bbb{Q}$, we define $\left\lfloor q\right\rfloor $ to be the
greatest integer number $\leq q$. ``gcd'' will be abbreviation for greatest
common divisor, and diag$\left( \frak{z}_1,\ldots ,\frak{z}_r\right) $ for
complex diagonal $r\times r$ matrices having $\frak{z}_1,\ldots ,\frak{z}_r$
as diagonal elements. Furthermore, for integers $\nu \geq 2$, we denote by $%
\zeta _\nu :=e^{\frac{2\pi \sqrt{-1}}\nu }$ the ``first'' $\nu $-th
primitive root of unity.\medskip

\noindent \textsf{(a) }Let $G$ be a finite subgroup of GL$\left( r,\Bbb{C}%
\right) $ which is \textit{small}, i.e. with no pseudoreflections, acting
linearly on $\Bbb{C}^r$, and $p:$ $\Bbb{C}^r\rightarrow \Bbb{C}^r/G$ the
quotient map. Denote by $\left( \Bbb{C}^r/G,\left[ \mathbf{0}\right] \right) 
$ the (germ of the) corresponding quotient singularity with $\left[ \mathbf{0%
}\right] :=p\left( \mathbf{0}\right) $.

\begin{proposition}[Singular locus]
\label{SLOC}If $G$ is a small finite subgroup of \emph{GL}$\left( r,\Bbb{C}%
\right) $, then 
\[
\text{\emph{Sing}}\left( \Bbb{C}^r/G\right) =p\left( \left\{ \mathbf{z}\in 
\Bbb{C}^r\ \left| \ G_{\mathbf{z}}\neq \left\{ \text{\emph{Id}}\right\}
\right. \right\} \right) 
\]
where $G_{\mathbf{z}}:=\left\{ g\in G\ \left| \ g\cdot \mathbf{z=z}\right.
\right\} $ is the isotropy group of $\mathbf{z=}\left( z_1,\ldots
,z_r\right) \in \Bbb{C}^r$.
\end{proposition}

\begin{theorem}[Prill's group-theoretic isomorphism criterion]
\label{prill-th}Let $G_1$, $G_2$ be two small finite subgroups of \emph{GL}$%
\left( r,\Bbb{C}\right) $. Then there exists an analytic isomorphism 
\[
\left( \Bbb{C}^r/G_1,\left[ \mathbf{0}\right] \right) \cong \left( \Bbb{C}%
^r/G_2,\left[ \mathbf{0}\right] \right) 
\]
if and only if $G_1$ and $G_2$ are conjugate to each other within \emph{GL}$%
\left( r,\Bbb{C}\right) $.
\end{theorem}

\noindent \textit{Proof. }See Prill \cite{Prill}, thm. 2, p. 382, and
Brieskorn \cite{Br2}, Satz 2.3, p. 341. $_{\Box }$ \medskip \smallskip

\noindent \textsf{(b) } Let $G$ be a finite, small, \textit{abelian}
subgroup of GL$\left( r,\Bbb{C}\right) $, $r\geq 2$, having order $l=\left|
G\right| \geq 2$. Define 
\[
\left\{ e_1=\left( 1,0,\ldots ,0,0\right) ^{\intercal },\ldots ,e_r=\left(
0,0,\ldots ,0,1\right) ^{\intercal }\right\} 
\]
to denote the standard basis of $\Bbb{Z}^r$, $N_0:=\sum_{i=1}^r\Bbb{Z}e_i$
the standard lattice, $M_0$ its dual, and\smallskip 
\[
T_{N_0}:=\text{Max-Spec}\left( \Bbb{C}\left[ \frak{x}_1^{\pm 1},\ldots ,%
\frak{x}_r^{\pm 1}\right] \right) =\left( \Bbb{C}^{*}\right) ^r\,\,. 
\]
Clearly, 
\[
T_{N_G}:=\text{Max-Spec}\left( \Bbb{C}\left[ \frak{x}_1^{\pm 1},\ldots ,%
\frak{x}_r^{\pm 1}\right] ^G\right) =\left( \Bbb{C}^{*}\right) ^r/G 
\]
is an $r$-dimensional algebraic torus with $1$-parameter group $N_G$ and
with group of characters $M_G$. Using the exponential map 
\[
\left( N_0\right) _{\Bbb{R}}\ni \left( y_1,\ldots ,y_r\right) ^{\intercal }=%
\mathbf{y\longmapsto }\text{ exp}\left( \mathbf{y}\right) :=\left( e^{\left(
2\pi \sqrt{-1}\right) y_1},\ldots ,e^{\left( 2\pi \sqrt{-1}\right)
y_r}\right) ^{\intercal }\in T_{N_0} 
\]
and the injection $\mathbf{\iota }:T_{N_0}\hookrightarrow $ GL$\left( r,\Bbb{%
C}\right) $ defined by 
\[
T_{N_0}\ni \left( t_1,\ldots ,t_r\right) ^{\intercal }=\mathbf{t}%
\hookrightarrow \text{ }\mathbf{\iota }\left( \mathbf{t}\right) :=\text{diag}%
\left( t_1,\ldots ,t_r\right) \in \text{ GL}\left( r,\Bbb{C}\right) \ , 
\]
we have obviously 
\[
N_G=\left( \mathbf{\iota \circ }\text{ exp}\right) ^{-1}\left( G\right) 
\text{ \ \ \ \ \ \ (\ and determinant \ det}\left( N_G\right) =\dfrac 1l\ 
\text{)} 
\]
(as long as we choose eigencoordinates to diagonalize the action of the
elements of $G$ on $\Bbb{C}^r$) with\smallskip 
\[
M_G=\left\{ m\in M_0\,\left| 
\begin{array}{c}
\frak{x}^m=\,\frak{x}_1^{\mu _1}\,\cdots \,\frak{x}_r^{\mu _r}\text{ \ \ is
a }G\text{-invariant } \\ 
\text{Laurent monomial (}m=\left( \mu _1,\ldots ,\mu _r\right) \text{)}
\end{array}
\right. \right\} \,\,\text{(and det}\left( M_G\right) =l\text{)}. 
\]
\newline
$\bullet $ If we define 
\[
\fbox{$\sigma _0:$=$\ $pos$\left( \left\{ e_1,..,e_r\right\} \right) $} 
\]
to be the $r$-dimensional positive orthant, and $\Delta _G$ to be the
fan\smallskip 
\[
\fbox{$
\begin{array}{ccc}
&  &  \\ 
& \Delta _G:=\left\{ \sigma _0\text{ together with its faces}\right\} &  \\ 
&  & 
\end{array}
$} 
\]
then by the exact sequence 
\[
0\rightarrow G\cong N_G/N_0\rightarrow T_{N_0}\rightarrow T_{N_G}\rightarrow
0 
\]
induced by the canonical duality pairing 
\[
M_0/M_G\times N_G/N_0\rightarrow \Bbb{Q}/\Bbb{Z\hookrightarrow C}%
^{*}\smallskip 
\]
(cf. \cite{Fulton}, p. 34, and \cite{Oda}, pp. 22-23), we get as projection
map :\smallskip 
\[
\Bbb{C}^r=X\left( N_0,\Delta _G\right) \rightarrow X\left( N_G,\Delta
_G\right) \ , 
\]
where 
\[
\fbox{$
\begin{array}{ccc}
&  &  \\ 
& X\left( N_G,\Delta _G\right) =U_{\sigma _0}=\Bbb{C}^r/G=\text{Max-Spec}%
\left( \Bbb{C}\left[ \frak{x}_1,\ldots ,\frak{x}_r\right] ^G\right)
\hookleftarrow T_{N_G} &  \\ 
&  & 
\end{array}
$}\smallskip \smallskip 
\]
$\bullet $ Formally, we identify $\left[ \mathbf{0}\right] $ with orb$\left(
\sigma _0\right) $. Moreover, in these terms, the singular locus of $X\left(
N_G,\Delta _G\right) $ can be written (by \ref{SLOC} and \ref{SMCR}) as the
union 
\[
\text{Sing}\left( X\left( N_G,\Delta _G\right) \right) =\text{ orb}\left(
\sigma _0\right) \cup \left( \bigcup \left\{ U_{\sigma _0}\left( \tau
\right) \,\left| 
\begin{array}{c}
\,\tau \precneqq \sigma _0\text{,\thinspace dim}\left( \tau \right) \geq
2\smallskip \text{ \thinspace } \\ 
\text{and \thinspace mult}\left( \tau ;N_G\right) \geq 2
\end{array}
\right. \right\} \right) \ . 
\]
$\bullet $ In particular, if the acting group $G$ is \textit{cyclic}, then,
fixing diagonalization of the action on $\Bbb{C}^r$, we may assume that $G$
is generated by the element 
\[
\text{diag}\left( \zeta _l^{\alpha _1}\,,\ldots ,\,\zeta _l^{\alpha
_r}\right) 
\]
for $r$ integers $\alpha _1,\ldots ,\alpha _r\in \left\{ 0,1,\ldots
,l-1\right\} $, at least $2$ of which are $\neq 0$. This $r$-tuple $\left(
\alpha _1,\ldots ,\alpha _r\right) $ of \textit{weights} is unique only up
to the usual conjugacy relations (see \ref{ISOCYC} below), and $N_G$ is to
be identified with the so-called \textit{lattice of weights} 
\[
N_G=N_0+\Bbb{Z\ }\left( \frac 1l\left( \alpha _1,\ldots ,\alpha _r\right)
^{\intercal }\ \right) \ \ \ \ 
\]
containing all lattice points representing the elements of 
\[
G=\left\{ \text{diag}\left( \zeta _l^{\left[ \lambda \alpha _1\right]
_l},\ldots ,\zeta _l^{\left[ \lambda \alpha _r\right] _l}\right) \ \left| \
\ \ \lambda \in \Bbb{Z},\ \ 0\leq \lambda \leq l-1\right. \right\} \ . 
\]

\begin{definition}
\label{type}\emph{Under these conditions, we say that the quotient
singularity} $\left( X\left( N_G,\Delta _G\right) ,\text{\emph{orb}}\left(
\sigma _0\right) \right) $\emph{\ is} \textit{of type} 
\begin{equation}
\fbox{$\dfrac 1l\left( \alpha _1,\ldots ,\alpha _r\right) $}  \label{ttype}
\end{equation}
\emph{(This is identical with the definition given in \cite{Reid3}, \S\ 4.2,
p. 370, up to the predeterminated fixing of the primitive root }$\zeta _l$ 
\emph{of unity to be the ``first'' one. In fact, this extra assumption is
not a significant restriction, and by fixing in advance the isomorphism }$%
G\cong N_G/N_0$\emph{\ one just simplifies certain technical arguments. Even
if we would let }$\zeta _l^{}$\emph{\ denote an arbitrary primitive root of
unity, all results would remain the same up to an obvious multiplication of
the exponents of the diagonalized elements of }$G$\emph{\ by a suitable
integer which would be relatively prime to }$l$\emph{; see also the comments
in \cite{Ito-Reid}, at the top of p. 225.)}
\end{definition}

\noindent $\bullet $ The existence of torus-equivariant resolutions of
cyclic quotient singularities was proved by Ehlers \cite{Ehlers}, \S\ I.3 \&
III.1-3, along the same lines as the more general theorem \ref{RESW}, i.e.,
by appropriate subdivisions of $\sigma _0$ into smaller cones of
multiplicity $1$. (Essentially the same result, expressed in the past
language of gluings of affine pieces, is due to Fujiki \cite{Fujiki}, \S\
1.3).\smallskip

\noindent $\bullet $ Note that, since $G$ is small, gcd$\left( l,\alpha
_1,\ldots ,\widehat{\alpha _i},\ldots ,\alpha _r\right) =1$, for all $i$, $%
1\leq i\leq r$. (The symbol $\widehat{\alpha _i}$ means here that $\alpha _i$
is omitted.)

\begin{lemma}
\label{ISOL}\emph{(i)} A cyclic quotient singularity of type $\emph{(\ref
{ttype})}$ has splitting codimension $\varkappa \in \left\{ 2,\ldots
,r-1\right\} $ if and only if there exists an index-subset 
\[
\left\{ \nu _1,\nu _2,\ldots ,\nu _{r-\varkappa }\right\} \subset \left\{
1,\ldots ,r\right\} ,
\]
such that 
\[
\alpha _{\nu _1}=\alpha _{\nu _2}=\cdots =\alpha _{\nu _{r-\varkappa }}=0\ ,
\]
which is, in addition, \emph{maximal} w.r.t. this property.\smallskip 
\newline
\emph{(ii)} A cyclic quotient msc-singularity of type $\emph{(\ref{ttype})}$
is isolated if and only if 
\[
\emph{gcd}\left( \alpha _i,l\right) =1,\ \forall i,\ \ 1\leq i\leq r\ .
\]
\end{lemma}

\noindent \textit{Proof. }It is immediate by the way we let $G$ act on $\Bbb{%
C}^r$. $_{\Box }\medskip \smallskip $

\noindent \textsf{(c) }For two integers $l,\,r\geq 2$, we define 
\[
\Lambda \left( l;r\right) :=\left\{ \left( \alpha _1,..,\alpha _r\right) \in
\left\{ 0,1,2,..,l-1\right\} ^r\left| 
\begin{array}{c}
\ \text{gcd}\left( l,\alpha _1,..,\widehat{\alpha _i},..,\alpha _r\right) =1,
\\ 
\text{for all }i,\ \ 1\leq i\leq r
\end{array}
\right. \right\} 
\]
and for $\left( \left( \alpha _1,\ldots ,\alpha _r\right) ,\left( \alpha
_1^{\prime },\ldots ,\alpha _r^{\prime }\right) \right) \in \Lambda \left(
l;r\right) \times \Lambda \left( l;r\right) $ the relation\smallskip 
\[
\left( \alpha _1,\ldots ,\alpha _r\right) \backsim \left( \alpha _1^{\prime
},\ldots ,\alpha _r^{\prime }\right) :\Longleftrightarrow \left\{ 
\begin{array}{l}
\text{there exists a permutation } \\ 
\theta :\left\{ 1,\ldots ,r\right\} \rightarrow \left\{ 1,\ldots ,r\right\}
\\ 
\text{and an integer }\lambda \text{, }1\leq \lambda \leq l-1\text{, } \\ 
\text{with gcd}\left( \lambda ,l\right) =1\text{, such that} \\ 
\alpha _{\theta \left( i\right) }^{\prime }=\left[ \lambda \cdot \alpha
_i\right] _l\text{, }\forall i,\ 1\leq i\leq r
\end{array}
\right\} \smallskip 
\]
It is easy to see that $\backsim $ is an equivalence relation on $\Lambda
\left( l;r\right) \times \Lambda \left( l;r\right) $.

\begin{corollary}[Isomorphism criterion for cyclic acting groups]
\label{ISOCYC} \ \smallskip \newline
Let $G$, $G^{\prime }$ be two small, cyclic finite subgroups of \emph{GL}$%
\left( r,\Bbb{C}\right) $ acting on $\Bbb{C}^r$, and let the corresponding
quotient singularities be of type $\frac 1l\,\left( \alpha _1,\ldots ,\alpha
_r\right) $ and $\frac 1{l^{\prime }}\,\left( \alpha _1^{\prime },\ldots
,\alpha _r^{\prime }\right) $ respectively. Then there exists an analytic 
\emph{(}torus-equivariant\emph{) }isomorphism 
\[
\left( X\left( N_G,\Delta _G\right) ,\text{\emph{orb}}\left( \sigma
_0\right) \right) \cong \left( X\left( N_{G^{\prime }},\Delta _{G^{\prime
}}\right) ,\text{\emph{orb}}\left( \sigma _0\right) \right) 
\]
if and only if $l=l^{\prime }$ and $\left( \alpha _1,\ldots ,\alpha
_r\right) \backsim \left( \alpha _1^{\prime },\ldots ,\alpha _r^{\prime
}\right) $ within $\Lambda \left( l;r\right) $.
\end{corollary}

\noindent \textit{Proof. }It follows easily from \ref{prill-th} (cf. Fujiki 
\cite{Fujiki}, lemma 2, p. 296). $_{\Box }$

\begin{proposition}[Gorenstein-condition]
\label{Gor-prop} \ \smallskip \newline
Let $\left( \Bbb{C}^r/G,\left[ \mathbf{0}\right] \right) =\left( X\left(
N_G,\Delta _G\right) ,\text{\emph{orb}}\left( \sigma _0\right) \right) $ be
an abelian quotient singularity. Then the following conditions are
equivalent \emph{:\smallskip } \newline
\emph{(i)} $X\left( N_G,\Delta _G\right) =U_{\sigma _0}=\Bbb{C}^r/G$ is
Gorenstein,\smallskip \newline
\emph{(ii) }$G\subset $ \emph{SL}$\left( r,\Bbb{C}\right) $,\smallskip 
\newline
\emph{(iii)} $\left\langle \left( 1,1,\ldots .1,1\right) ,n\right\rangle
\geq 1$, for all $n,\,n\in \sigma _0\cap \left( N_G\smallsetminus \left\{ 
\mathbf{0}\right\} \right) $,\smallskip \newline
\emph{(iv) }$\left( X\left( N_G,\Delta _G\right) ,\text{\emph{orb}}\left(
\sigma _0\right) \right) $ is a canonical singularity of index $1$%
.\smallskip \newline
In particular, if $\left( \Bbb{C}^r/G,\left[ \mathbf{0}\right] \right) $ is
cyclic of type $\frac 1l\,\left( \alpha _1,\ldots ,\alpha _r\right) $, then 
\emph{(i)-(iv)} are equivalent to 
\[
\sum_{j=1}^r\alpha _j\equiv 0\ \text{\emph{(mod} }l\text{\emph{)}}
\]
\end{proposition}

\noindent \textit{Proof. }See Watanabe \cite{Wat1} and Reid \cite{Reid1} . $%
_{\Box }\bigskip $

\noindent $\bullet $ If $X\left( N_G,\Delta _G\right) $ is Gorenstein, then
the cone $\sigma _0=$ pos$\,\left( \frak{s}_G\right) $ is supported by the
so-called \textit{junior lattice simplex\smallskip } 
\[
\fbox{$\frak{s}_G=$ conv$\left( \left\{ e_1,..,e_r\right\} \right) $}%
\smallskip 
\]
(w.r.t. $N_G$, cf. \cite{Ito-Reid}, \cite{BD}). Note that up to $\mathbf{0}$
there is no other lattice point of $\sigma _0\cap N_G$ lying ``under'' the
affine hyperplane of $\Bbb{R}^r$ containing $\frak{s}_G$. Moreover, the
lattice points representing the $l-1$ non-trivial group elements are exactly
those belonging to the intersection of a dilation $\lambda \,\frak{s}_G$ of $%
\frak{s}_G$ with $\mathbf{Par}\left( \sigma _0\right) $, for some integer $%
\lambda $, $1\leq \lambda \leq r-1$.

\section{Lattice triangulations and crepant projective resolutions 
of Gorenstein abelian quotient singularities}\label{LATR}

\noindent In this section we briefly formulate (and only partially prove)
some general theorems concerning the study of projective, crepant
resolutions of Gorenstein abelian quotient singularities in terms of
appropriate lattice triangulations of the junior simplex. (For detailed
expositions we refer to \cite{DHZ1}, \cite{DHZ2}). Diagrams and formulae in
boxes outline actually the only essential prerequisites for the reading of
the rest of the paper.\smallskip \medskip

\noindent \textsf{(a) }By vert$\left( \mathcal{S}\right) $ we denote the set
of vertices of a polyhedral complex $\mathcal{S}$. By a \textit{%
triangulation }$\mathcal{T}$ of a polyhedral complex $\mathcal{S}$ we mean a
geometric simplicial subdivision of $\mathcal{S}$ with vert$\left( \mathcal{S%
}\right) \subset $ vert$\left( \mathcal{T}\right) $. A polytope $P$ will be,
as usual, identified with the polyhedral complex consisting of $P$ itself
together with all its faces.\bigskip

\noindent \textsf{(b) }A triangulation $\mathcal{T}$ of an $r$-dimensional
polyhedral complex $\mathcal{S}$ is called \textit{coherent }(or \textit{%
regular}) if there exists a strictly upper convex $\mathcal{T}$-support
function $\psi :\left| \mathcal{T}\right| \rightarrow \Bbb{R}$, i.e. a
piecewise-linear real function defined on the underlying space $\left| 
\mathcal{T}\right| $ of $\mathcal{T}$, for which 
\[
\psi \left( t\ \mathbf{x}+\left( 1-t\right) \mathbf{\ y}\right) \geq t\ \psi
\left( \mathbf{x}\right) +\left( 1-t\right) \ \psi \left( \mathbf{y}\right) ,%
\text{ for all \ }\mathbf{x},\mathbf{y}\in \left| \mathcal{T}\right| ,\ 
\text{and\ }\emph{\ }t\in \left[ 0,1\right] \ , 
\]
so that for every maximal simplex $\mathbf{s}$ of $\mathcal{T}$, there is a
linear function $\eta _{\mathbf{s}}:\left| \mathbf{s}\right| \mathbf{%
\rightarrow }\Bbb{R}$ satisfying $\psi \left( \mathbf{x}\right) \leq \eta _{%
\mathbf{s}}\left( \mathbf{x}\right) $, for all $\mathbf{x}\in $\emph{\ }$%
\left| \mathcal{T}\right| $, with equality being valid only for those $%
\mathbf{x}$ belonging to $\mathbf{s}$. The set of all strictly upper convex $%
\mathcal{T}$-support functions will be denoted by SUCSF$_{\Bbb{R}}\left( 
\mathcal{T}\right) $. A useful lemma to create a new ``global'' strictly
upper convex support function by gluing together given ``local'' ones, is
the following :

\begin{lemma}[Patching Lemma]
\label{PATCH}Let $P\subset \Bbb{R}^r$ be a $r$-polytope, $\mathcal{T}%
=\left\{ \mathbf{s}_i\ \left| \ i\in I\right. \right\} $ $\ $\emph{(}with $I$
a finite set\emph{) }a coherent triangulation of $P$, and $\mathcal{T}%
_i=\left\{ \mathbf{s}_{i,j}\ \left| \ j\in J_i\right. \right\} \ $\emph{(}$%
J_i$ finite, for all $i\in I$\emph{) }a coherent triangulation of $\mathbf{s}%
_i$, for all $i\in I$. If $\psi _i:\left| \mathcal{T}_i\right| \rightarrow 
\Bbb{R}$ denote strictly upper convex $\mathcal{T}_i$-support functions,
such that 
\[
\fbox{$
\begin{array}{ccc}
&  &  \\ 
& \psi _i\left| _{\mathbf{s}_{i\,}\cap \,\mathbf{s}_{i^{\prime }}}\right.
=\psi _{i^{\prime }}\left| _{\mathbf{s}_{i\,}\cap \,\mathbf{s}_{i^{\prime
}}}\right.  &  \\ 
&  & 
\end{array}
$}
\]
for all $\left( i,i^{\prime }\right) \in I\times I$, then 
\[
\widetilde{\mathcal{T}}:=\left\{ \text{\emph{all the simplices }}\mathbf{s}%
_{i,j}\ \left| \ \ \forall j,\ \ j\in J_i,\ \text{\emph{and \ }}\forall i,\
i\in I\right. \right\} 
\]
forms a \emph{coherent triangulation }of the initial polytope $P$\emph{\ (}%
because the above $\psi _i$'s can be canonically \emph{``}patched together%
\emph{'' }to construct an element $\psi $ of \emph{SUCSF}$_{\Bbb{R}}\left( 
\widetilde{\mathcal{T}}\right) $\emph{).}
\end{lemma}

\noindent \textit{Proof. }See Bruns-Gubeladge-Trung \cite{BGT}, lemma\
2.2.2, pp. 143-145. \medskip \smallskip \smallskip $_{\Box }$

\noindent \textsf{(c) }Let $N$ denote an $r$-dimensional lattice. By a 
\textit{lattice polytope }(w.r.t. $N$)\textit{\ }is meant a polytope in $N_{%
\Bbb{R}}\Bbb{\cong R}^r$ with vertices belonging to $N$. If $\left\{
n_0,n_1,\ldots ,n_k\right\} $ is a set of $k\leq r$ affinely independent
lattice points, $\mathbf{s}$ the lattice\emph{\ }$k$-dimensional simplex $%
\mathbf{s}=$ conv$\left( \left\{ n_0,n_1,n_2,\ldots ,n_k\right\} \right) $,
and $N_{\mathbf{s}}:=$ lin$\left( \left\{ n_1-n_0,\ldots ,n_k-n_0\right\}
\right) \cap N$, then\smallskip \newline
$\bullet $ we say that $\mathbf{s}$ is an \textit{elementary simplex}\emph{\ 
}if 
\[
\left\{ \mathbf{y}-n_0\ \left| \ \mathbf{y}\in \mathbf{s}\right. \right\}
\cap N_{\mathbf{s}}=\left\{ \mathbf{0},\,n_1-n_0,\ldots ,\,n_k-n_0\right\} . 
\]
\newline
$\bullet $ $\mathbf{s}$ is \textit{basic} if it has anyone of the following
equivalent properties:\smallskip \newline
(i) $\left\{ n_1-n_0,n_2-n_0,\ldots ,n_k-n_0\right\} $ is a $\Bbb{Z}$-basis
of $N_{\mathbf{s}}$,\smallskip \newline
(ii) $\mathbf{s\ }$has relative volume\ Vol$\left( \mathbf{s};N_{\mathbf{s}%
}\right) =\dfrac{\text{Vol}\left( \mathbf{s}\right) }{\text{det}\left( N_{%
\mathbf{s}}\right) }=\dfrac 1{k!}\ \ $(w.r.t. $N_{\mathbf{s}}$) .

\begin{lemma}
\label{EL-BA}\emph{(i)} Every basic lattice simplex is elementary.\newline
\emph{(ii)} Elementary lattice simplices of dimension $\leq 2$ are basic.
\end{lemma}

\noindent \textit{Proof. }(i) Let $\mathbf{s=}$ conv$\left( \left\{
n_0,n_1,\ldots ,n_k\right\} \right) $ be a basic lattice simplex. Since 
\[
\left\{ n_1-n_0,n_2-n_0,\ldots ,n_k-n_0\right\} 
\]
is\emph{\ }a $\Bbb{Z}$-basis of $N_{\mathbf{s}}$, if $n^{\prime }\in \mathbf{%
s\,\cap \,}\left( N_{\mathbf{s}}\smallsetminus \left\{ \mathbf{0}\right\}
\right) $, then obviously $n^{\prime }=n_i$ for some index $i,\,1\leq i\leq
k $. Thus, 
\[
\left\{ \mathbf{y}-n_0\ \left| \ \mathbf{y}\in \mathbf{s}\right. \right\}
\cap N_{\mathbf{s}}=\left\{ \mathbf{0},\,n_1-n_0,\ldots ,\,n_k-n_0\right\} . 
\]
\newline
(ii) The relative volume of an elementary lattice $1$-simplex (resp. of an
elementary lattice $2$-simplex) is always equal to $1$ (resp. equal to $1/2$%
). $_{\Box }$

\begin{example}
\label{Counter}\emph{The lattice }$r$\emph{-simplex} 
\[
\mathbf{s}=\text{\emph{\ conv}}\left( \left\{ \mathbf{0},e_1,e_2,\ldots
,e_{r-2},e_{r-1},\left( 1,1,\ldots ,1,1,2\right) ^{\intercal }\right\}
\right) \subset \Bbb{R}^r\text{\emph{,\thinspace } }r\geq 3\text{\emph{,}}
\]
\emph{(w.r.t. }$\Bbb{Z}^r$\emph{)} \emph{serves as example of an elementary
but non-basic simplex because }$\mathbf{s}\cap \Bbb{Z}^r=$ \emph{vert}$%
\left( \mathbf{s}\right) $ \emph{and} 
\[
r!\,\text{\emph{Vol}}\left( \mathbf{s};\Bbb{Z}^r\right) =\left| \text{\emph{%
det}}\left( e_1,\ldots ,e_{r-1},\left( 1,1,\ldots ,1,1,2\right) ^{\intercal
}\right) \right| =2\neq 1\ .
\]
\end{example}

\begin{definition}
\emph{A triangulation} $\mathcal{T}$ \emph{of a lattice polytope} $P\subset
N_{\Bbb{R}}\cong \Bbb{R}^r$ \emph{(w.r.t.} $N$\emph{)} \emph{is called} 
\textit{lattice triangulation }\emph{if vert}$\left( P\right) \subset $ 
\emph{vert}$\left( \mathcal{T}\right) \subset N$\emph{. The set of all
lattice triangulations of a lattice polytope }$P$ \emph{(w.r.t.} $N$\emph{)
will be denoted by} $\mathbf{LTR}_N\left( P\right) $.
\end{definition}

\begin{definition}
\emph{A lattice triangulation} $\mathcal{T}$ \emph{of }$P\subset N_{\Bbb{R}%
}\cong \Bbb{R}^r$ \emph{(w.r.t.} $N$\emph{)} \emph{is called }\textit{%
maximal triangulation }\emph{if vert}$\left( \mathcal{T}\right) =N\cap P$%
\emph{. A lattice triangulation} $\mathcal{T}$ \emph{of }$P$ \emph{is
obviously maximal if and only if each simplex }$\mathbf{s}$ \emph{of }$%
\mathcal{T}$ \emph{is elementary. A lattice triangulation} $\mathcal{T}$ 
\emph{of }$P$ \emph{is said to be }\textit{basic}\emph{\ if }$\mathcal{T}$ 
\emph{consists of exclusively basic simplices. We define :} 
\[
\begin{array}{l}
\mathbf{LTR}_N^{\emph{\scriptsize max}}\,\left( P\right) :=\left\{ \mathcal{T%
}\in \mathbf{LTR}_N\left( P\right) \ \left| \ \mathcal{T}\text{\emph{\ \ is
a maximal triangulation of \ }}P\right. \right\} \smallskip \ ,\smallskip 
\\ 
\mathbf{LTR}_N^{\emph{\scriptsize basic}}\left( P\right) :=\left\{ \mathcal{T%
}\in \mathbf{LTR}_N^{\emph{\scriptsize max}}\,\left( P\right) \ \left| \ 
\mathcal{T}\text{\emph{\ \ is a basic triangulation of \ }}P\right. \right\}
\ .
\end{array}
\]
\emph{(Moreover, adding the prefix} $\mathbf{Coh}$\emph{-} \emph{to anyone
of the above sets, we shall mean the subsets of their elements which are
coherent). The }\textit{hierarchy of lattice triangulations }\emph{of a }$P$%
\emph{\ (as above) is given by the following inclusion-diagram:}
\end{definition}

\[
\fbox{$
\begin{array}{ccccc}
\mathbf{LTR}_N^{\text{{\scriptsize basic}}}\left( P\right) & \subset & 
\mathbf{LTR}_N^{\text{{\scriptsize max}}}\,\left( P\right) & \subset & 
\mathbf{LTR}_N\left( P\right) \smallskip \\ 
\bigcup &  & \bigcup &  & \bigcup \smallskip \\ 
\mathbf{Coh}\text{-}\mathbf{LTR}_N^{\text{{\scriptsize basic}}}\left(
P\right) & \subset & \mathbf{Coh}\text{-}\mathbf{LTR}_N^{\text{{\scriptsize %
max}}}\,\left( P\right) & \subset & \mathbf{Coh}\text{-}\mathbf{LTR}_N\left(
P\right)
\end{array}
$} 
\]

\begin{proposition}
\label{NON-EM}For any lattice polytope $P\subset N_{\Bbb{R}}\cong \Bbb{R}^r$ 
\emph{(}w.r.t. $N$\emph{)} the set of maximal coherent triangulations\emph{\ 
}$
\begin{array}{ccc}
\mathbf{Coh}\text{-}\mathbf{LTR}_N^{\emph{max}}\emph{\,}\left( P\right)  & 
\subset  & \mathbf{Coh}\text{-}\mathbf{LTR}_N\left( P\right) 
\end{array}
$ of $P$ is non-empty.
\end{proposition}

\noindent \textit{Proof. }Consider a s.c.p. cone supported by $P$ in $\Bbb{R}%
^{r+1}$, and then make use of \cite{OP}, cor. 3.8, p. 394. $_{\Box }$

\begin{remark}
\label{Hibi}\emph{(i)} \emph{Already in dimension }$r=2$\emph{, there are
lots of examples of }$P$\emph{'s with } 
\[
\mathbf{LTR}_N^{\text{\emph{basic}}}\left( P\right) \smallsetminus \mathbf{%
Coh}\text{-}\mathbf{LTR}_N^{\text{\emph{basic}}}\left( P\right) \neq
\varnothing \text{ .}
\]
\emph{See, for instance, \ref{NONPROJ} below.\smallskip }\newline
\emph{(ii) In addition, for lattice polytopes }$P$\emph{,} $\mathbf{LTR}_N^{%
\text{\emph{basic}}}\left( P\right) \neq \varnothing $\emph{\ does not imply
necessarily } 
\[
\mathbf{Coh}\text{-}\mathbf{LTR}_N^{\text{\emph{basic}}}\left( P\right) \neq
\varnothing \ \emph{.}
\]
\emph{As it was proved recently by Hibi and Ohsugi \cite{HO},} \emph{there exists a }$9$%
\emph{-dimensional }$0$\emph{/}$1$\emph{-polytope (with }$15$ \emph{%
vertices) possessing basic triangulations but no coherent }\textit{and}\emph{%
\ basic ones. \medskip }
\end{remark}

\noindent \textsf{(d) }To pass from triangulations to desingularizations we
need to introduce some extra notation.\smallskip

\begin{definition}
\emph{Let }$\left( X\left( N_G,\Delta _G\right) ,\text{\emph{orb}}\left(
\sigma _0\right) \right) $ \emph{be an} $r$\emph{-dimensional abelian
Gorenstein quotient singularity} \emph{(}$r\geq 2$\emph{), and }$\frak{s}_G$ 
\emph{the }$\left( r-1\right) $\emph{-dimensional junior simplex. For any
simplex }$\mathbf{s}$ \emph{of a lattice triangulation $\mathcal{T}$ of }$%
\frak{s}_G$ \emph{let }$\sigma _{\mathbf{s}}$ \emph{denote the s.c.p. cone} 
\[
\sigma _{\mathbf{s}}:=\left\{ \lambda \,\mathbf{y}\in \left( N_G\right) _{%
\Bbb{R}}\ \left| \ \lambda \in \Bbb{R}_{\geq 0}\,,\ \mathbf{y}\in \mathbf{s}%
\right. \right\} \ \,\left( \,=\text{\emph{pos}}\left( \mathbf{s}\right) \,\,%
\text{\emph{within} \thinspace }\left( N_G\right) _{\Bbb{R}}\right) \ 
\]
\emph{supporting} $\mathbf{s}$. \emph{We define the fan } 
\[
\widehat{\Delta _G}\left( \mathcal{T}\right) :=\left\{ \sigma _{\mathbf{s}}\
\left| \ \mathbf{s}\in \mathcal{T}\right. \right\} 
\]
\emph{of s.c.p. cones in }$\left( N_G\right) _{\Bbb{R}}\cong \Bbb{R}^r$\emph{%
,} \emph{and } 
\[
\begin{array}{ll}
\mathbf{PCDES}\left( X\left( N_G,\Delta _G\right) \right) := & \left\{ 
\begin{array}{c}
\text{\textit{partial crepant} }T_{N_G}\text{\emph{-equivariant }} \\ 
\text{\emph{desingularizations} \emph{of} \ }X\left( N_G,\Delta _G\right) 
\\ 
\begin{array}{c}
\ \text{\emph{with overlying spaces having }} \\ 
\text{\emph{\ at most (}}\Bbb{Q}\text{\emph{-factorial) }\textit{canonical} }
\\ 
\text{\emph{singularities (of index }}1\text{\emph{)}}
\end{array}
\end{array}
\right\} \ ,\smallskip \medskip  \\ 
\mathbf{PCDES}^{\emph{max}}\left( X\left( N_G,\Delta _G\right) \right) := & 
\left\{ 
\begin{array}{c}
\text{\textit{partial crepant} }T_{N_G}\text{\emph{-equivariant }} \\ 
\text{\emph{desingularizations} \emph{of} \ }X\left( N_G,\Delta _G\right) 
\\ 
\begin{array}{c}
\ \text{\emph{with overlying spaces having }} \\ 
\text{\emph{\ at most (}}\Bbb{Q}\text{\emph{-factorial) }\textit{terminal} }
\\ 
\text{\emph{singularities (of index }}1\text{\emph{)}}
\end{array}
\end{array}
\right\} \smallskip \ ,\medskip  \\ 
\mathbf{CDES}\left( X\left( N_G,\Delta _G\right) \right) := & \left\{ 
\begin{array}{c}
\text{\textit{crepant} }T_{N_G}\text{\emph{-equivariant (full)}} \\ 
\text{\emph{\ desingularizations of} }X\left( N_G,\Delta _G\right) 
\end{array}
\right\} \ .\medskip \smallskip 
\end{array}
\]
\emph{(Whenever we put the prefix} $\mathbf{QP}$\emph{- in the front of
anyone of them, we shall mean the corresponding subsets of them consisting
of those} \emph{desingularizations whose overlying spaces are }\textit{%
quasiprojective}\emph{.) }
\end{definition}

\begin{theorem}[Desingularizing by triangulations]
\label{TRCOR} \ \smallskip \newline
Let $\left( X\left( N_G,\Delta _G\right) ,\text{\emph{orb}}\left( \sigma
_0\right) \right) $ be an $r$-dimensional abelian Gorenstein quotient
singularity \emph{(}$r\geq 2$\emph{).} Then there exist one-to-one
correspondences \emph{:} 
\[
\fbox{$
\begin{array}{ccc}
&  &  \\ 
& 
\begin{array}{ccc}
\mathbf{LTR}_{N_G}^{\emph{basic}}\left( \frak{s}_G\right) \smallskip  & 
\stackrel{\text{\emph{1:1}{\scriptsize \smallskip }}}{\longleftrightarrow }
& \mathbf{CDES}\left( X\left( N_G,\Delta _G\right) \right)  \\ 
\bigcap  &  & \bigcap  \\ 
\mathbf{LTR}_{N_G}^{\emph{max}}\,\left( \frak{s}_G\right) \smallskip  & 
\stackrel{\text{\emph{1:1}{\scriptsize \smallskip }}}{\longleftrightarrow }
& \mathbf{PCDES}^{\emph{max}}\left( X\left( N_G,\Delta _G\right) \right)  \\ 
\bigcap  &  & \bigcap  \\ 
\mathbf{LTR}_{N_G}\left( \frak{s}_G\right)  & \stackrel{\text{\emph{1:1}%
{\scriptsize \smallskip }}}{\longleftrightarrow } & \mathbf{PCDES}\left(
X\left( N_G,\Delta _G\right) \right) 
\end{array}
&  \\ 
&  & 
\end{array}
$}
\]
as well as\smallskip 
\[
\fbox{$
\begin{array}{ccc}
&  &  \\ 
& 
\begin{array}{ccc}
\mathbf{Coh}\text{-}\mathbf{LTR}_{N_G}^{\emph{basic}}\left( \frak{s}%
_G\right) \smallskip  & \stackrel{\text{\emph{1:1}{\scriptsize \smallskip }}%
}{\longleftrightarrow } & \mathbf{QP}\text{-}\mathbf{CDES}\left( X\left(
N_G,\Delta _G\right) \right)  \\ 
\bigcap  &  & \bigcap  \\ 
\mathbf{Coh}\text{-}\mathbf{LTR}_{N_G}^{\emph{max}}\,\left( \frak{s}%
_G\right) \smallskip  & \stackrel{\text{\emph{1:1}{\scriptsize \smallskip }}%
}{\longleftrightarrow } & \mathbf{QP}\text{-}\mathbf{PCDES}^{\emph{max}%
}\left( X\left( N_G,\Delta _G\right) \right)  \\ 
\bigcap  &  & \bigcap  \\ 
\mathbf{Coh}\text{-}\mathbf{LTR}_{N_G}\left( \frak{s}_G\right)  & \stackrel{%
\text{\emph{1:1}{\scriptsize \smallskip }}}{\longleftrightarrow } & \mathbf{%
QP}\text{-}\mathbf{PCDES}\left( X\left( N_G,\Delta _G\right) \right) 
\end{array}
&  \\ 
&  & 
\end{array}
$}\smallskip \smallskip 
\]
which are realized by crepant $T_{N_G}$-equivariant birational morphisms of
the form\smallskip \smallskip 
\begin{equation}
\fbox{$
\begin{array}{ccc}
&  &  \\ 
& f_{\mathcal{T}}=\text{\emph{id}}_{*}:X\left( N_G,\text{ }\widehat{\Delta _G%
}\left( \mathcal{T}\right) \right) \longrightarrow X\left( N_G,\Delta
_G\right)  &  \\ 
&  & 
\end{array}
$}  \label{DESING}
\end{equation}
induced by mapping 
\[
\mathcal{T}\longmapsto \widehat{\Delta _G}\left( \mathcal{T}\right) ,\ \ \ \
\ \widehat{\Delta _G}\left( \mathcal{T}\right) \longmapsto X\left( N_G,\text{
}\widehat{\Delta _G}\left( \mathcal{T}\right) \right) \ .
\]
\end{theorem}

\noindent \textit{Sketch of proof.} $X\left( N_G,\Delta _G\right) $ is
Gorenstein and has at most rational singularities, i.e. canonical
singularities of index $1$. Moreover, its dualizing sheaf is trivial. Let 
\[
f=\text{id}_{*}:X\left( N_G,\text{ }\widetilde{\Delta _G}\right)
\longrightarrow X\left( N_G,\Delta _G\right) 
\]
denote an arbitrary partial desingularization. Studying either the behaviour
of the highest rational differentials on $X\left( N_G,\text{ }\widetilde{%
\Delta _G}\right) $ (see \cite{Reid1}, \S\ 3, \cite{Reid3}, \S\ 4.8, \cite
{Mark}, prop. 3, or \cite{DHZ1}, prop. 4.1), or the support function
associated to $K_{X\left( N_G,\text{ }\widetilde{\Delta _G}\right) }$ (cf. 
\cite{Roan1}, \S\ 2), one proves 
\[
K_{X\left( N_G,\text{ }\widetilde{\Delta _G}\right) }=f^{*}\left( K_{X\left(
N_G,\Delta _G\right) }\right) -\sum_{\varrho \in \widetilde{\Delta _G}\left(
1\right) }\ \left( \left\langle \left( 1,\ldots ,1\right) ,n\left( \varrho
\right) \right\rangle -1\right) \ D_{n\left( \varrho \right) }\ , 
\]
where $D_{n\left( \varrho \right) }:=V\left( \Bbb{\varrho }\right) =V\left( 
\text{pos}\left( \left\{ n\left( \varrho \right) \right\} \right) \right) $.
Obviously, $f$ is crepant if and only if 
\[
\text{Gen}\left( \widetilde{\Delta _G}\right) \subset \left\{ \mathbf{y=}%
\left( y_1,\ldots ,y_r\right) ^{\intercal }\in \left( N_G\right) _{\Bbb{R}}\
\left| \ \sum_{i=1}^r\ y_i=1\right. \right\} \ , 
\]
and since the number of crepant exceptional prime divisors is independent of
the specific choice of $f$, the first and second $1$-$1$ correspondences
(from below) of the first box are obvious by the adjunction-theoretic
definition of terminal (resp. canonical) singularities. In particular, all $%
T_{N_G}$-equivariant partial crepant desingularizations of $X\left(
N_G,\Delta _G\right) $ of the form (\ref{DESING}) have overlying spaces with
at most $\Bbb{Q}$-factorial singularities, and conversely, each partial $%
T_{N_G}$-equivariant crepant desingularization with overlying space
admitting at most $\Bbb{Q}$-factorial singularities, has to be of this form.
($\Bbb{Q}$-factoriality is here equivalent to the consideration only of
triangulations instead of more general polyhedral subdivisions. Furthermore,
by maximal triangulations you exhaust all the crepant exceptional prime
divisors). The top $1$-$1$ correspondence of the first box follows from the
equivalence 
\[
\mathcal{T}\ni \mathbf{s}\text{\ is a basic simplex}\Longleftrightarrow
\left\{ 
\begin{array}{l}
\text{mult}\left( \sigma _{\mathbf{s}};N_G\right) =1\text{ \ for the cone}
\\ 
\, \\ 
\sigma _{\mathbf{s}}\in \widehat{\Delta _G}\left( \mathcal{T}\right) \text{
supporting it}
\end{array}
\right\} \ . 
\]
It remains to prove the $1$-$1$ correspondences of the second box. As it was
explained in \cite{DHZ1}, \S\ 4, for every $\psi \in $ SUCSF$_{\Bbb{Q}%
}\left( N_G,\text{ }\widehat{\Delta _G}\left( \mathcal{T}\right) \right) $,
the restriction $\psi \left| _{\mathcal{T}}\right. $ belongs to SUCSF$_{\Bbb{%
R}}\left( \mathcal{T}\right) $, and conversely, to any support function $%
\psi \in $ SUCSF$_{\Bbb{R}}\left( \mathcal{T}\right) $, one may canonically
assign a rational (or even an integral) strictly upper convex support
function defined on the entire $\left| \text{ }\widehat{\Delta _G}\left( 
\mathcal{T}\right) \right| $. To finish the proof of the theorem we apply
corollary \ref{quasiproj}. $_{\Box }$

\begin{remark}
\emph{(i) }\textit{Flops.} \emph{Every pair of triangulations} $\mathcal{T}_1
$, $\mathcal{T}_2\in \mathbf{Coh}$-$\mathbf{LTR}_{N_G}^{{\rm max}}\,\left( \frak{s}%
_G\right) $ \emph{gives rise to the determination of a birational morphism} 
\[
X\left( N_G,\text{ }\widehat{\Delta _G}\left( \mathcal{T}_1\right) \right)
\longrightarrow X\left( N_G,\text{ }\widehat{\Delta _G}\left( \mathcal{T}%
_2\right) \right) 
\]
\emph{which is composed of a finite number of flops (cf. \cite{OP}, \S\ 3).
Since all triangulations of }$\mathbf{Coh}$-$\mathbf{LTR}_{N_G}\,\left( 
\frak{s}_G\right) $ \emph{are parametrized by the vertices of the so-called }%
\textit{secondary polytope }\emph{of }$\frak{s}_G$\emph{, this
transition-map is induced by performing successively} \textit{bistellar
operations}\emph{, i.e., by passing from the vertex of the secondary
polytope of} $\frak{s}_G$ \emph{representing }$\mathcal{T}_1$\emph{\ to that
one representing }$\mathcal{T}_2$\emph{\ following a (not necessarily
uniquely determined) path which connects these two vertices. (In dimension }$%
3$\emph{\ this is nothing but Danilov's theorem \cite{Danilov}, cf. \ref
{WHAT3} (vi) below). For detailed presentations of the theory of secondary
polytopes we refer to the articles of Billera-Filliman-Sturmfels \cite{BFS}
and Oda-Park \cite{OP}, and to the treatment of Gelfand-Kapranov-Zelevinsky 
\cite{GKZ}, ch. 7.\smallskip }\newline
\emph{(ii)} \textit{Factorization}\emph{. The birational morphisms
corresponding to members of }$\mathbf{QP}$-$\mathbf{PCDES}\left( X\left(
N_G,\Delta _G\right) \right) $\emph{\ can be decomposed into more elementary
toric contractions by Reid's toric version of ``MMP'' (\cite{Reid2},
(0.2)-(0.3)). In several cases these contractions can be directly
interpreted as inversed (normalized) blow-ups. For concrete examples see
rem. \ref{DIFFBL} and \S\ \ref{FACTOR} below.\smallskip }\newline
\emph{(iii) }\textit{Exceptional divisors.} \emph{Each irreducible component
of an exceptional divisor w.r.t. an }$f=f_{\mathcal{T}}$\emph{\ (as above)
is a }$\Bbb{Q}$\emph{-Cartier prime divisor which carries itself the
structure of an }$\left( r-1\right) $\emph{-dimensional toric variety
determined by the corresponding star within } $\widehat{\Delta _G}\left( 
\mathcal{T}\right) $ \emph{(cf. above \S\ \ref{TORIC} \textbf{(f)}, (ii)).
An exceptional prime divisor is compact if and only if the lattice point
representing it belongs to int}$\left( \sigma _0\right) $\emph{\ (see thm. 
\ref{COMPACT}). }
\end{remark}

\begin{theorem}[Number-theoretic version of McKay correspondence]
\label{COHOMOLOGY} \ \smallskip \newline
Let $\left( X\left( N_G,\Delta _G\right) ,\text{\emph{orb}}\left( \sigma
_0\right) \right) $ be an abelian Gorenstein quotient singularity of
dimension $r\geq 2$, 
\[
f=\text{\emph{id}}_{*}:X\left( N_G,\widehat{\Delta _G}\right) \rightarrow
X\left( N_G,\Delta _G\right) 
\]
a $T_{N_G}$-equivariant crepant, full\emph{\ }resolution and $\mathbf{F:}%
=f^{-1}\left( \left[ \mathbf{0}\right] \right) $ the central fiber over the
origin $\left[ \mathbf{0}\right] =$ \emph{orb}$\left( \sigma _0\right) $%
\emph{:}\smallskip 
\[
\ 
\begin{array}{l}
\mathbf{F}=\bigcup_{\varrho \in \widehat{\Delta _G}\left( 1\right) }\
\left\{ D_{n\left( \varrho \right) }\ \left| \ n\left( \varrho \right) \in 
\text{\emph{int}}\left( \sigma _0\right) \cap N_G\right. \right\} \cup  \\ 
\\ 
\cup \left( \bigcup \ \left\{ V\left( \sigma \right) \ \left| \ \sigma \in
\bigcup_{i=2}^{r-1}\ \widehat{\Delta _G}\left( i\right) \ ,\ \text{\emph{int}%
}\left( \sigma \right) \subset \text{\emph{int}}\left( \sigma _0\right)
\right. \right\} \right) \ .
\end{array}
\]
Then $\mathbf{F}$ is a strong deformation retract of $X\left( N_G,\widehat{%
\Delta _G}\right) $, and only the even cohomology groups of $\mathbf{F}$ are
non-trivial. Their dimensions \emph{(}over $\Bbb{Q}$\emph{) }are given by
the formulae \emph{:\smallskip } 
\begin{equation}
\fbox{$\text{\emph{dim}}_{\Bbb{Q}}H^{2\,i}\left( \mathbf{F};\Bbb{Q}\right)
=\left\{ 
\begin{array}{lll}
1 & , & \text{\emph{if \ }}i=0 \\ 
\ \  &  & \ \  \\ 
\#\,\left( \left( i\,\frak{s}_G\right) \cap \mathbf{Par}\left( \sigma
_0\right) \cap N_G\right)  & , & \text{\emph{if\ }}\ 1\leq i\leq r-1 \\ 
\ \  &  & \ \  \\ 
0 & , & \text{\emph{otherwise}}
\end{array}
\right. $}\smallskip   \label{COH-DIM}
\end{equation}
In particular, the topological Euler-Poincar\'{e} characteristic of $X\left(
N_G,\widehat{\Delta _G}\right) $ equals\emph{:} 
\begin{equation}
\fbox{$\chi \left( X\left( N_G,\widehat{\Delta _G}\right) \right) =\chi
\left( \mathbf{F}\right) =l=\left| G\right| $}  \label{EP-DES}
\end{equation}
Obviously, the numbers \emph{(\ref{COH-DIM}), (\ref{EP-DES})} are \emph{%
independent} of the choice of triangulations $\mathcal{T}$ of $\frak{s}_G$
by means of which we construct $\widehat{\Delta _G}\ \left( =\widehat{\Delta
_G}\left( \mathcal{T}\right) \right) $ \emph{(cf. (\ref{DESING})).}
\end{theorem}

\noindent \textit{Proof. }See Batyrev-Dais \cite{BD}, thm. 5.4, p. 910, and 
\cite{DHZ2} for further comments. $_{\Box }$

\begin{corollary}
If $\left( \Bbb{C}^r/G,\left[ \mathbf{0}\right] \right) =\left( X\left(
N_G,\Delta _G\right) ,\text{\emph{orb}}\left( \sigma _0\right) \right) $ is
a Gorenstein cyclic quotient singularity of type $\frac 1l\left( \alpha
_1,\ldots ,\alpha _r\right) $, then, maintaining the above notation and
assumption, we obtain\smallskip 
\begin{equation}
\ \ \fbox{$\emph{dim}_{\Bbb{Q}}H^{2i}\left( \mathbf{F};\Bbb{Q}\right) =\#\
\left\{ \left. \lambda \in \left[ 0,l\right) \cap \Bbb{Z}\ \right| \
\tsum\limits_{j=1}^r\ \left[ \,\lambda \,\alpha _j\,\right] _l=i\cdot
l\right\} $}  \label{CYCCOH}
\end{equation}
\end{corollary}

\begin{remark}
\emph{The numbers of the right hand side of (\ref{COH-DIM}) and (\ref{CYCCOH}%
) make sence, even without assuming the existence of a }$T_{N_G}$\emph{%
-equivariant crepant, full desingularization of }$X\left( N_G,\Delta
_G\right) $\emph{, and were used in \cite{BD} as the ``correctional terms''
for introducing the formal definition of the so-called }\textit{%
string-theoretic Hodge numbers }\emph{of Calabi-Yau varieties (or, more
general, of Gorenstein compact complex varieties) which have at most abelian
quotient singularities.\smallskip }
\end{remark}

\noindent \textsf{(e) }By theorem \ref{TRCOR} it is now clear that Reid's
question (formulated in \S \ref{Intro}), restricted to the category of
torus-equivariant desingularizations of $X\left( N_G,\Delta _G\right) $, can
be restated as follows : \smallskip

\noindent \underline{\textbf{Question}}\textbf{\ : For a} \textbf{Gorenstein
abelian quotient singularity} 
\[
\left( X\left( N_G,\Delta _G\right) ,\text{orb}\left( \sigma _0\right)
\right) 
\]
\textbf{with junior simplex }$\frak{s}_G$\textbf{\ , what kind of conditions
on the acting group }$G$\textbf{\ would guarantee the existence of a basic,
coherent triangulation of }$\frak{s}_G$\textbf{?\smallskip \smallskip }

\noindent Though (as already mentioned in \S \ref{Intro}) this question will
be treated explicitly in \cite{DHZ2}, we at least shall explain here how a
very simple necessary existence-criterion works, and apply it efficiently in
our special singularity-series in \S\ \ref{MONSE}.\smallskip

\begin{lemma}
\label{S1LEMMA}Let $\left( X\left( N_G,\Delta _G\right) ,\emph{orb}\left(
\sigma _0\right) \right) $ be a Gorenstein abelian quotient singularity and $%
\frak{s}_G$ the junior simplex. If $\frak{s}_G$ admits a basic triangulation 
$\mathcal{T}$, then \smallskip for all $n\in \left( N_G\smallsetminus
\left\{ \mathbf{0}\right\} \right) \cap \sigma _0$ there exist $r$ lattice
points $n_1,\ldots ,n_r$ belonging to $\frak{s}_G\cap N_G$ so that 
\[
n\in \Bbb{Z}_{\geq 0}\ n_1+\cdots +\Bbb{Z}_{\geq 0}\ n_r
\]
\end{lemma}

\noindent \textit{Proof}. Since $\mathcal{T}$ is a basic triangulation
inducing a subdivision of $\sigma _0$ into cones of multiplicity $1$, $n$
belongs to a subcone $\sigma ^{\prime }$ of $\sigma _0$ of the form 
\[
\sigma ^{\prime }=\Bbb{R}_{\geq 0}\ n_1+\cdots +\Bbb{R}_{\geq 0}\ n_r,\ \ \
\ \text{mult}\left( \sigma ^{\prime };N_G\right) =1,\ 
\]
and can be therefore written as a linear combination 
\[
n=\mu _1n_1+\cdots +\mu _rn_r,\ \ \ \ \ \left( \mu _i\in \Bbb{R}_{\geq 0},\
\forall i,\ 1\leq i\leq r\right) 
\]
Now $\left\{ n_1,\ldots ,n_r\right\} $ is a $\Bbb{Z}$-basis of $N_G$, which
means that $n$ can be also written as 
\[
n=\mu _1^{\prime }n_1+\cdots +\mu _r^{\prime }n_r,\ \ \ \ \ \left( \mu
_i^{\prime }\in \Bbb{Z},\ \forall i,\ 1\leq i\leq r\right) 
\]
By linear independence we get $\mu _i=\mu _i^{\prime }\in \Bbb{Z}\cap \Bbb{R}%
_{\geq 0}=\Bbb{Z}_{\geq 0}$, $\forall i,\ 1\leq i\leq r$. $_{\Box }$

\begin{theorem}[Necessary Existence-Criterion]
\label{KILLER} \ \smallskip \newline
Let\emph{\ }$\left( X\left( N_G,\Delta _G\right) ,\emph{orb}\left( \sigma
_0\right) \right) $ be a Gorenstein abelian quotient singularity. If $\frak{s%
}_G$ admits a basic triangulation $\mathcal{T}$, then 
\begin{equation}
\fbox{$
\begin{array}{ccc}
&  &  \\ 
& \mathbf{Hlb}_{N_G}\left( \sigma _0\right) =\frak{s}_G\cap N_G &  \\ 
&  & 
\end{array}
$}  \label{HILBCON}
\end{equation}
i.e. all members of the Hilbert basis of $\sigma _0$ have to be either \emph{%
``}junior\emph{''} elements or to belong to $\left\{ e_1,\ldots ,e_r\right\} 
$.
\end{theorem}

\noindent \textit{Proof}. The inclusion ``$\supset $'' is always true
(without any further assumption about the existence or non-existence of such
a triangulation) and is obvious by the definition of Hilbert basis. Now if
there were an element $n\in \mathbf{Hlb}_{N_G}\left( \sigma _0\right)
\smallsetminus \left( \frak{s}_G\cap N_G\right) $, then by lemma \ref
{S1LEMMA} this could be written as a non-negative integer linear combination 
\[
\ n=\mu _1n_1+\cdots +\mu _rn_r 
\]
of $r$ elements of $\frak{s}_G\cap N_G$. Since $\mathbf{0\notin Hlb}%
_{N_G}\left( \sigma _0\right) $, there were at least one index $j=\
j_{\bullet }\in \left\{ 1,\ldots ,r\right\} $, for which $\mu _{\ j_{\bullet
}}\neq 0$. If $\mu _{\ j_{\bullet }}=1$ and $\mu _j=0$ for all $j$, $j\in $ $%
\left\{ 1,\ldots ,r\right\} \smallsetminus \left\{ \ j_{\bullet }\right\} $,
then $n=n_{\ j_{\bullet }}\in \frak{s}_G\cap N_G$ which would contradict
our assumption. But even the cases in which either $\mu _{\ j_{\bullet }}=1$
and some other $\mu _j$'s were $\neq 0$, or $\mu _{\ j_{\bullet }}\geq 2$,
would be exluded as impossible because of the characterization (\ref
{Hilbbasis}) of the Hilbert basis $\mathbf{Hlb}_{N_G}\left( \sigma _0\right) 
$ as the set of additively irreducible vectors of $\sigma _0\cap \left(
N_G\smallsetminus \left\{ \mathbf{0}\right\} \right) $. Hence, $\mathbf{Hlb}%
_{N_G}\left( \sigma _0\right) \subset \frak{s}_G\cap N_G$. $_{\Box }$

\begin{remark}
\emph{\label{FIRLA}(i) For a long time it was expected that condition (\ref
{HILBCON}) might be sufficient too for the existence of a basic
triangulation of }$\frak{s}_G$\emph{, and, as we shall see in \cite{DHZ2},
this} \textit{is} \emph{the case for some special choices of cyclic group
actions on }$\Bbb{C}^r$\emph{. Nevertheless, Firla and Ziegler (\cite{Firla}, \S 4.2 \& \cite{Fir-Zi}) discovered recently the first }\textit{counterexamples}%
. \emph{Among them, the counterexample of the }$4$\emph{-dimensional
Gorenstein cyclic quotient singularity with the smallest possible acting
group-order, fulfilling property (\ref{HILBCON}) and admitting no crepant,
torus-equivariant resolutions, is that of type }$\frac 1{39}\,\left(
1,5,8,25\right) $\emph{.\smallskip }\newline
\emph{(ii) To apply necessary criterion \ref{KILLER} in practice, in order
to exclude ``candidates'' for having} \emph{crepant, torus-equivariant
resolutions, one has first to determine all the elements of the Hilbert
basis }$\mathbf{Hlb}_{N_G}\left( \sigma _0\right) $ \emph{and then to test
if at least one of them breaks away from the junior simplex or not. From the
point of view of complexity theory of algorithms, however, this procedure
is in general ``NP-hard'' (cf. Henk-Weismantel \cite{HEW}, \S\ 3). }
\end{remark}

\begin{exercise}
\emph{\ For the singularity of type }$\frac 19\left( 1,2,3,3\right) $\emph{\
determine explicitly the Hilbert basis }$\mathbf{Hlb}_{N_G}\left( \sigma
_0\right) $ \emph{and show that it does not possess any crepant
torus-equivariant resolution because} \emph{condition (\ref{HILBCON}) is
violated. [Hint. Verify that }$\frac 19\left( 5,1,6,6\right) ^{\intercal
}\in \mathbf{Hlb}_{N_G}\left( \sigma _0\right) $\emph{].}
\end{exercise}

\section{Peculiarities of dimensions $2$ and $3$}

\noindent In the low dimensions $r\in \left\{ 2,3\right\} $, we have always 
\[
\begin{array}{ccc}
\mathbf{CDES}\left( X\left( N_G,\Delta _G\right) \right) & = & \mathbf{PCDES}%
^{\text{{\small max}}}\left( X\left( N_G,\Delta _G\right) \right)
\end{array}
\]
by lemma \ref{EL-BA} (ii) and thm. \ref{TRCOR}. This is exactly the lemma
which fails in general if $r\geq 4$, and makes the high dimensions so
exciting (cf. \ref{Counter}). However, low dimensional Gorenstein abelian
quotient singularities are still valuable as testing ground for lots of
interesting, related problems and pave the way to systematic
generalizations.\bigskip \newline
\textsf{(a) }In dimension $2$ we meet only the ``classical'' $A_{l-1}$%
-singularities (i.e. cyclic quotient singularities of type $1/l\left(
1,l-1\right) $).

\begin{lemma}
\label{DIM2}All $2$-dimensional Gorenstein abelian quotient singularities 
\[
\left( X\left( N_G,\Delta _G\right) ,\text{\emph{orb}}\left( \sigma
_0\right) \right) 
\]
are cyclic of type $\frac 1l\left( 1,l-1\right) $ \emph{(with} $l\geq 2$%
\emph{). }They\emph{\ }admit a unique projective crepant \emph{(}$=$ minimal%
\emph{) }resolution \emph{(\ref{DESING})} induced by the triangulation 
\[
\mathcal{T}=\left\{ \,\text{\emph{conv}}\left( \left\{ \frac 1l\,\left(
j-1,l-\left( j-1\right) \right) ^{\intercal },\,\frac 1l\,\left(
j,l-j\right) ^{\intercal }\right\} \right) \ \left| \ 1\leq j\leq l\right.
\right\} \ \smallskip 
\]
and the $l-1$ exceptional prime divisors 
\[
D_j:=V\left( \text{\emph{pos}}\left( \left\{ \frac 1l\,\left( j,l-j\right)
^{\intercal }\right\} \right) \right) ,\ \ 1\leq j\leq l-1,
\]
are smooth rational curves having intersection numbers\smallskip 
\[
\left( D_i\cdot D_j\right) =\left\{ 
\begin{array}{rrr}
1 & ,\ \text{\emph{if}} & \left| i-j\right| =1\smallskip  \\ 
-2 & ,\ \text{\emph{if}} & \text{\emph{\ }}i=j\smallskip  \\ 
0 & ,\ \text{\emph{if}} & \left| i-j\right| \geq 2
\end{array}
\right. 
\]
for all $i,j,\ 1\leq i\leq j\leq l-1.$
\end{lemma}

\begin{remark}
\label{DIFFBL}\emph{Obviously, the subdivision of }$\sigma _0$\emph{\ into }$%
l$\emph{\ cones of multiplicity }$1$\emph{\ induced by the above }$\mathcal{T%
}$\emph{\ can be done successively in }$k$\emph{\ ``steps'', with }$1\leq
k\leq l-1$\emph{, i.e., by drawing} \emph{in }$k$\emph{\ steps} \emph{the }$%
l-1$\emph{\ required rays in any order you would wish. (The possibility of
drawing more than one rays in a step is not excluded). Such a procedure
gives rise to decomposing the full desingularization into a series of
partial ones. It should nevertheless be stressed that starting-points of
different series of choices correspond to blow-ups of }\textit{different} 
\emph{ideal sheaves }$\mathcal{I}$\emph{\ with orb}$\left( \sigma _0\right) =
$\emph{\ supp}$\left( \mathcal{O}_{U_{\sigma _0}}\,/\,\mathcal{I}\right) $.%
\emph{\ (Note that in dimension }$2$ \emph{we do not need extra
normalizations). Consequently, there are lots of factorizations of the
crepant resolution morphism }$f=f_{\mathcal{T}}$. \emph{Let us illustrate it
by considering the example of the }$A_4$\emph{-singularity }$\frac 15\left(
1,4\right) $\emph{. The morphism }$f$\emph{\ admits of two different natural
factorizations} \emph{\ } 
\[
f=g_2\circ g_1=h_4\circ h_3\circ h_2\circ h_1
\]
\emph{which are} \emph{depicted in figures \textbf{4} and \textbf{5},
respectively.\smallskip }\newline
\emph{(i) Since } 
\[
\mathbf{Hlb}_{M_G}\left( \sigma _0^{\vee }\right) =\left\{ m_1=\left(
5,0\right) ,\ m_2=\left( 1,1\right) ,\ m_3=\left( 0,5\right) \right\}
\,,\smallskip 
\]
\emph{using the induced embedding }$\iota :U_{\sigma _0}\hookrightarrow \Bbb{%
C}^3$ \emph{of }$U_{\sigma _0}=X\left( N_G,\Delta _G\right) $\emph{%
,\smallskip } 
\[
U_{\sigma _0}=\left\{ \left( z_1,z_2,z_3\right) \in \Bbb{C}^3\ \left| \
z_2^5-z_1\,z_3=0\right. \right\} ,\ \ z_i=\mathbf{e}\left( m_i\right) ,\
\forall i,\ 1\leq i\leq 3,\smallskip 
\]
\emph{and proposition \ref{USNBL}, we obtain }$g_1$\emph{\ and }$g_2$\emph{\
coming from blowing up only maximal }$0$\emph{-dimensional ideals:\smallskip 
} 
\[
X\left( N_G,\text{ }\widehat{\Delta _G}\left( \mathcal{T}\right) \right) =%
\mathbf{Bl}_{\text{\emph{orb}}\left( \sigma _2\right) }^{\text{\emph{red}}%
}\left( \mathbf{Bl}_{\text{\emph{orb}}\left( \sigma _0\right) }^{\text{\emph{%
red}}}\left( U_{\sigma _0}\right) \right) \stackrel{g_2}{\longrightarrow }%
\mathbf{Bl}_{\text{\emph{orb}}\left( \sigma _0\right) }^{\text{\emph{red}}%
}\left( U_{\sigma _0}\right) \stackrel{g_1}{\longrightarrow }U_{\sigma
_0}\,.\smallskip 
\]
\emph{More precisely, we first subdivide }$\sigma _0$\emph{\ into the three
subcones }$\sigma _1$\emph{, }$\sigma _2$\emph{, }$\sigma _3$\emph{, with} 
\[
\begin{array}{l}
\sigma _1=\left\{ \mathbf{y}\in \sigma _0\ \left| \ \left\langle m_2-m_1,%
\mathbf{y}\right\rangle \geq 0\ \text{\emph{\&} }\left\langle m_3-m_1,%
\mathbf{y}\right\rangle \geq 0\right. \right\} =\text{\emph{pos}}\left(
\left\{ \left( \frac 15,\frac 45\right) ^{\intercal },e_2\right\} \right) ,
\\ 
\\ 
\sigma _2=\left\{ \mathbf{y}\in \sigma _0\ \left| \ \left\langle m_1-m_2,%
\mathbf{y}\right\rangle \geq 0\ \text{\emph{\&} }\left\langle m_3-m_2,%
\mathbf{y}\right\rangle \geq 0\right. \right\} =\text{\emph{pos}}\left( \left( \frac 45,\frac
15\right) ^{\intercal },\left( \frac 15,\frac 45\right) ^{\intercal }\right),
\\ 
\\ 
\sigma _3=\left\{ \mathbf{y}\in \sigma _0\ \left| \ \left\langle m_1-m_3,%
\mathbf{y}\right\rangle \geq 0,\left\langle m_2-m_3,\mathbf{y}\right\rangle
\geq 0\right. \right\} =\text{\emph{pos}}\left(
\left\{ e_1,\left( \frac 45,\frac 15\right) ^{\intercal }\right\} \right)
.
\end{array}
\]
\emph{After that we perform again a usual blow-up, but now for the cone }$%
\sigma _2$\emph{\ (instead of }$\sigma _0$\emph{) which contains the
remaining two lattice points.}
\end{remark}
\begin{figure*}[hb]
\begin{center}
\input{fig4.pstex_t}
\makespace
 {\textbf{Fig. 4}\emph{\ }}\smallskip
\end{center}
\end{figure*}

\noindent (ii) In the second factorization, $h_1$ gives rise to a
``starring'' subdivision of $\sigma _0$ into only two cones 
\[
\sigma _1=\text{pos}\left( \left\{ \left( \frac 45,\frac 15\right)
^{\intercal },\,e_2\right\} \right) \text{ \ \ \ and \ \ \ }\sigma _2=\text{%
pos}\left( \left\{ e_1,\,\left( \frac 45,\frac 15\right) ^{\intercal
}\right\} \right) \,. 
\]
In fact, 
\[
h_1:\mathbf{Bl}_{\text{orb}\left( \sigma _0\right) }^{\mathcal{I}}\left(
U_{\sigma _0}\right) \stackrel{}{\longrightarrow }U_{\sigma _0} 
\]
is a \textit{directed }blow-up w.r.t. $\left( 4/5,1/5\right) ^{\intercal }$.
As $\mathcal{O}_{U_{\sigma _0}}$-ideal sheaf $\mathcal{I}$ can be taken the
pullback (via $\iota $) of the ideal sheaf $\mathcal{J}\cong J^{\sim }$
supporting the $z_3$-axis of $\Bbb{C}^3$, with $J:=\left(
z_1,z_2,z_1\,z_3\right) $ (in $\Bbb{C}\left[ z_1,z_2,z_3\right] $),
because\smallskip 
\begin{eqnarray*}
\sigma _1 &=&\left\{ \mathbf{y}\in \sigma _0\ \left| \ \left\langle m_2-m_1,%
\mathbf{y}\right\rangle \geq 0\text{ \ \emph{\ }and \ \ }\left\langle m_3,%
\mathbf{y}\right\rangle \geq 0\right. \right\} , \\
&&\, \\
\ \ \ \text{ \ \ \ }\sigma _2 &=&\left\{ \mathbf{y}\in \sigma _0\ \left| \
\left\langle m_1-m_2,\mathbf{y}\right\rangle \geq 0\text{ \ \emph{\ }and \ \ 
}\left\langle m_1+m_3-m_2,\mathbf{y}\right\rangle \geq 0\right. \right\} \
.\smallskip
\end{eqnarray*}
$\bullet $ Another characterization of $\mathbf{Bl}_{\text{orb}\left( \sigma
_0\right) }^{\mathcal{I}}\left( U_{\sigma _0}\right) $ avoiding the
embedding is provided by the following commutative diagram due to the
extension of the group action on $\mathbf{Bl}_{\mathbf{0}}^{\text{red}%
}\left( \Bbb{C}^2\right) $: 
\[
\begin{array}{lll}
\mathbf{Bl}_{\mathbf{0}}^{\text{red}}\left( \Bbb{C}^2\right) & 
\twoheadrightarrow & \mathbf{Bl}_{\text{orb}\left( \sigma _0\right) }^{%
\mathcal{I}}\left( U_{\sigma _0}\right) \cong \mathbf{Bl}_{\mathbf{0}}^{%
\text{red}}\left( \Bbb{C}^2\right) /G \\ 
\downarrow &  & \downarrow \\ 
\Bbb{C}^2 & \twoheadrightarrow & U_{\sigma _0}=\Bbb{C}^2/G
\end{array}
\]
Clearly, our blown up space is isomorphic to the quotient of the blown up $%
\Bbb{C}^2$ at the origin divided by $G$. For $h_2$ we proceed analogously by
letting $\sigma _2$ play the role of $\sigma _0$. Obviously, the same is
valid for $h_3$ and $h_4$. (Comparing the above two factorizations, we see
that the first one is ``double speedy''.)

\begin{figure*}
\begin{center}
\input{fig5.pstex_t}
\makespace
 {\textbf{Fig. 5}\emph{\ }}
\end{center}
\end{figure*}

\vspace{7cm}
\noindent (iii) Finally, let us point out that also $f=f_{\mathcal{T}}$
itself can be similarly regarded as a blow-up morphism 
\[
f:\mathbf{Bl}_{\text{orb}\left( \sigma _0\right) }^{\mathcal{I}}\left(
U_{\sigma _0}\right) \longrightarrow U_{\sigma _0}\ .\smallskip 
\]
As $\mathcal{I}$ one may choose the pullback of $\mathcal{J}\cong J^{\sim }$
with \smallskip $J$ the ideal 
\[
J:=\left( z_1^4,z_1^3\,z_2,z_1^2\,z_2^3,z_1^2\,z_2\,z_3,z_1^2\,z_3^2\right)
\ \ \text{ in \ \ }\Bbb{C}\left[ z_1,z_2,z_3\right] ,
\]
because\smallskip 
\[
\begin{array}{l}
\left\{ \mathbf{y}\in \sigma _0\ \left| 
\begin{array}{l}
\left\langle m_2-m_1,\mathbf{y}\right\rangle \geq 0\text{,\smallskip \
\smallskip } \\ 
\left\langle 3\,m_2-2\,m_1,\mathbf{y}\right\rangle \geq 0,\smallskip
\smallskip  \\ 
\left\langle m_2+m_3-2\,m_1,\mathbf{y}\right\rangle \geq 0,\smallskip
\smallskip \  \\ 
\left\langle 2\,m_3-2\,m_1,\mathbf{y}\right\rangle \geq 0
\end{array}
\right. \right\} =\text{pos}\left( \left\{ \,e_2,\,\left( \frac 15,\frac
45\right) ^{\intercal }\right\} \right) , \\ 
\\ 
\left\{ \mathbf{y}\in \sigma _0\ \left| 
\begin{array}{l}
\left\langle m_1-m_2,\mathbf{y}\right\rangle \geq 0\text{,\smallskip } \\ 
\left\langle 2\,m_2-m_1,\mathbf{y}\right\rangle \geq 0,\smallskip  \\ 
\left\langle m_3-m_1,\mathbf{y}\right\rangle \geq 0,\smallskip  \\ 
\left\langle 2\,m_3-m_1-m_2,\mathbf{y}\right\rangle \geq 0
\end{array}
\right. \right\} =\text{pos}\left( \left\{ \left( \frac 15,\frac 45\right)
^{\intercal },\left( \frac 25,\frac 35\right) ^{\intercal }\right\} \right) ,
\end{array}
\]
and 
\[
\begin{array}{l}
\\ 
\left\{ \mathbf{y}\in \sigma _0\ \left| 
\begin{array}{l}
\left\langle 2\,m_1-3\,m_2,\mathbf{y}\right\rangle \geq 0,\smallskip  \\ 
\left\langle m_1-2\,m_2,\mathbf{y}\right\rangle \geq 0,\smallskip  \\ 
\left\langle m_3-2\,m_2,\mathbf{y}\right\rangle \geq 0,\smallskip  \\ 
\left\langle 2\,m_3-3\,m_2,\mathbf{y}\right\rangle \geq 0
\end{array}
\right. \right\} =\text{pos}\left( \left\{ \left( \frac 25,\frac 35\right)
^{\intercal },\left( \frac 35,\frac 25\right) ^{\intercal }\right\} \right) ,
\\ 
\\ 
\left\{ \mathbf{y}\in \sigma _0\ \left| 
\begin{array}{l}
\left\langle 2\,m_1-m_2-m_3,\mathbf{y}\right\rangle \geq 0\text{\smallskip
\smallskip } \\ 
\left\langle m_1-m_3,\mathbf{y}\right\rangle \geq 0,\smallskip  \\ 
\left\langle 2\,m_2-m_3,\mathbf{y}\right\rangle \geq 0,\smallskip  \\ 
\left\langle m_3-m_2,\mathbf{y}\right\rangle \geq 0
\end{array}
\right. \right\} =\text{pos}\left( \left\{ \left( \frac 35,\frac 25\right)
^{\intercal },\left( \frac 45,\frac 15\right) ^{\intercal }\right\} \right) ,
\end{array}
\]
and finally 
\[
\!\!\!\!\!\!\!\!\left\{ \mathbf{y}\in \sigma _0\ \left| 
\begin{array}{l}
\left\langle 2\,m_1-2\,m_3,\mathbf{y}\right\rangle \geq 0,\smallskip  \\ 
\left\langle m_1+m_2-2\,m_3,\mathbf{y}\right\rangle \geq 0,\smallskip  \\ 
\left\langle 3\,m_2-2\,m_3,\mathbf{y}\right\rangle \geq 0,\smallskip  \\ 
\left\langle m_2-m_3,\mathbf{y}\right\rangle \geq 0
\end{array}
\right. \right\} =\text{pos}\left( \left\{ \left( \frac 45,\frac 15\right)
^{\intercal },\,e_1\right\} \right) .
\]

\begin{exercise}
\emph{To represent }$f$ \emph{as the restriction (over }$U_{\sigma _0}$\emph{%
) of a proper birational morphism which comes from a single blow-up of a
torus-invariant }$0$\textit{-dimensional}\emph{\ ideal of} $\Bbb{C}\left[
z_1,z_2,z_3\right] $ \emph{\ supporting only} $\mathbf{0}\in U_{\sigma
_0}\subset \Bbb{C}^3$ \emph{(cf. thm. \ref{RESSTR}), it is enough} (\emph{%
instead of the the above }$J$\emph{)} \emph{\ to consider } 
\[
J=\left( z_1^6,\ z_1^5\,z_2,\ z_1^4\,z_2^3,\ z_1^4\,z_2\,z_3,\
z_1^4\,z_3^2,\ z_1^2\,z_3,\ z_1\,z_2\,z_3,\ z_2^3\,z_3,\ z_1\,z_3^2,\
z_3^3\right) \ .
\]
\emph{(From the point of view of toric geometry, this means that one has to
find among the defining monomials suitable monomials involving powers of the
variables }$z_1$ \emph{and} \emph{\ }$z_3$ \emph{separately,} \emph{such
that the corresponding order function becomes again linear precisely} \emph{%
on the above five maximal cones.)\smallskip }
\end{exercise}

\noindent \textsf{(b) }We now focus our attention to the new phenomena which
arise in dimension three.

\begin{theorem}[What happens in the ``intermediate'' dimension $3$ ?]
\label{WHAT3} \ \smallskip \newline
Let $\left( \Bbb{C}^3/G,\left[ \mathbf{0}\right] \right) =\left( X\left(
N_G,\Delta _G\right) ,\text{\emph{orb}}\left( \sigma _0\right) \right) $ a $3
$-dimensional Gorenstein abelian quotient msc-singularity and $\mathcal{T}$
a maximal \emph{(}and therefore basic\emph{)} triangulation of the junior
simplex $\frak{s}_G$ inducing a crepant, full resolution 
\[
f:X\left( N_G,\text{ }\widehat{\Delta _G}\left( \mathcal{T}\right) \right)
\longrightarrow X\left( N_G,\Delta _G\right) \ .\smallskip 
\]
For any $n\in $ \emph{vert}$\left( \mathcal{T}\right) $, let $D_n$ denote
the closure $V\left( \text{\emph{pos}}\left( \left\{ n\right\} \right)
\right) $.\smallskip \newline
\emph{(i) }If $n\in $ \emph{int}$\left( \frak{s}_G\right) \cap N_G$, then $%
D_n$ is a rational compact surface coming from \emph{(}usual\emph{) }%
blow-ups either of $\Bbb{P}_{\Bbb{C}}^2$ or of a Hirzebruch surface $\Bbb{F}%
_\lambda $ at finitely many $T_{N_G\left( \text{\emph{pos}}\left( \left\{
n\right\} \right) \right) }$-fixed points.\smallskip \newline
\emph{(ii) }If $\partial \frak{s}_G\cap \left( N_G\smallsetminus \left\{
e_1,e_2,e_3\right\} \right) \neq \varnothing $, and $n\in $ \emph{conv}$%
\left( e_{i_1},e_{i_2}\right) \smallsetminus \left\{ e_{i_1},e_{i_2}\right\} 
$, with $\left\{ i_1,i_2\right\} \subset \left\{ 1,2,3\right\} $, $i_1\neq
i_2$, and $\left\{ i_3\right\} =\left\{ 1,2,3\right\} \smallsetminus \left\{
i_1,i_2\right\} $, then $D_n$ is the total space of a ruled fibration over
the \emph{``}$i_3$-axis\emph{''} of $\Bbb{C}^3$. Its fibers over the
punctured $i_3$-axis are isomorphic to $\Bbb{P}_{\Bbb{C}}^1$.\smallskip 
\newline
\emph{(iii) }For three distinct vertices $n,n^{\prime },n^{\prime \prime }$
of $\mathcal{T}$, we have 
\[
\left( D_n\cdot D_{n^{\prime }}\cdot D_{n^{\prime \prime }}\right) =\left\{ 
\begin{array}{ccc}
1 & , & \text{\emph{if \ conv}}\left( \left\{ n,n^{\prime },n^{\prime \prime
}\right\} \right) \ \text{\emph{is a} }2\text{\emph{-simplex of} }\mathcal{%
T\smallskip \smallskip } \\ 
0 & , & \text{\emph{otherwise}}
\end{array}
\right. \smallskip 
\]
\emph{(iv)} If $n,n^{\prime }\in $ \emph{vert}$\left( \mathcal{T}\right) $, 
\emph{conv}$\left( \left\{ n,n^{\prime }\right\} \right) $ is a $1$-simplex
of $\mathcal{T}$, but no both $n$ and $n^{\prime }$ belong the same face of $%
\partial \frak{s}_G$, then there exist exactly two vertices $\frak{y},\frak{y%
}^{\prime }$ of $\mathcal{T}$, such that \emph{conv}$\left( \left\{
n,n^{\prime },\frak{y}\right\} \right) $, \emph{conv}$\left( \left\{
n,n^{\prime },\frak{y}^{\prime }\right\} \right) $ are $2$-simplices of $%
\mathcal{T}$ satisfying a $\Bbb{Z}$-linear dependency equation of the form 
\begin{equation}
\kappa \,n+\kappa ^{\prime }\,n^{\prime }+\frak{y}+\frak{y}^{\prime }=0\text{%
, \emph{\ for some unique} }\kappa ,\kappa ^{\prime }\in \Bbb{Z\ }\text{%
\emph{with} }\kappa +\kappa ^{\prime }=-2  \label{DD3}
\end{equation}
and 
\begin{equation}
\kappa =\left( D_n^2\cdot D_{n^{\prime }}\right) ,\ \ \kappa ^{\prime
}=\left( D_n\cdot D_{n^{\prime }}^2\right) \text{ }  \label{INTERSE}
\end{equation}
Furthermore, the normal bundle of the rational intersection curve 
\[
C=V\left( \text{\emph{pos}}\left( \left\{ n,n^{\prime }\right\} \right)
\right) 
\]
splits into the direct sum\emph{:} 
\begin{equation}
\mathcal{N}_{C\,/\,X\left( N_G,\text{ }\widehat{\Delta _G}\left( \mathcal{T}%
\right) \right) }\cong \mathcal{O}_C\left( \kappa \right) \,\,\mathbf{\oplus 
}\,\mathcal{O}_C\left( \kappa ^{\prime }\right)   \label{NORMB}
\end{equation}
\emph{(v) }If $n\in $ \emph{int}$\left( \frak{s}_G\right) \cap N_G$, then $%
D_n$ has self-intersection number 
\begin{equation}
D_n^3=12-\#\,\left( \text{\emph{Gen}}\left( \text{\emph{Star}}\left( \text{%
\emph{pos}}\left( \left\{ n\right\} \right) ,\Delta _G\right) \right)
\right)   \label{D3}
\end{equation}
\emph{(vi)} For any other maximal triangulation $\mathcal{T}^{\prime }$ of $%
\frak{s}_G$, there exists a birational morphism \emph{(}isomorphism in
codimension $1$\emph{)} 
\[
X\left( N_G,\text{ }\widehat{\Delta _G}\left( \mathcal{T}\right) \right)
\rightarrow X\left( N_G,\text{ }\widehat{\Delta _G}\left( \mathcal{T}%
^{\prime }\right) \right) 
\]
which is a composite of finitely many \emph{elementary transformations} 
\emph{(}$=$\emph{\ simple flops)} w.r.t. rational smooth curves.
\end{theorem}

\noindent \textit{Proof.} (i) follows from Oda's classification of smooth
compact toric surfaces (\cite{Oda}, thm. 1.28), (ii) is clear by
construction, and (iii) by (\ref{IN10}). The vectorial $\Bbb{Z}$-linear
dependency equation (\ref{DD3}) of (iv) with $\kappa +\kappa ^{\prime }=-2$
is obvious because $n$, $n^{\prime }$, $\frak{y}$, $\frak{y}^{\prime }$ are
``junior'' elements; (\ref{INTERSE}) follows from (\ref{INTD}) and (\ref
{NORMB}) from the splitting principle of holomorphic vector bundles over $%
C\cong \Bbb{P}_{\Bbb{C}}^1$, the normal bundle exact sequence, and the
triviality of the dualizing sheaf of $X\left( N_G,\text{ }\widehat{\Delta _G}%
\left( \mathcal{T}\right) \right) $. (\ref{D3}) in (v) can be proved by
making use of adjunction formula, combined with Noether's formula, and $\chi
\left( D_n,\mathcal{O}_{D_n}\right) =1$. For the proof of (vi) we refer to
Danilov \cite{Danilov}, prop. 2, or Oda \cite{Oda}, prop. 1.30 (ii). $_{\Box
}$

\begin{remark}
\label{NONPROJ}\emph{In dimension three one has }$\mathbf{QP}$\emph{-}$%
\mathbf{CDES}\left( X\left( N_G,\Delta _G\right) \right) \neq \varnothing $%
\emph{\ because} 
\[
\begin{array}{ccc}
\mathbf{Coh}\text{-}\mathbf{LTR}_{N_G}^{\emph{basic}}\left( \frak{s}%
_G\right)  & = & \mathbf{Coh}\text{-}\mathbf{LTR}_{N_G}^{\emph{max}}\,\left( 
\frak{s}_G\right) \neq \varnothing 
\end{array}
\]
\emph{by prop. \ref{NON-EM} and thm. \ref{TRCOR}. Nevertheless, in contrast
to what takes place in dimension }$2$\emph{, there exist lots of examples of
finite abelian subgroups }$G$\emph{\ of SL}$\left( 3,\Bbb{C}\right) $\emph{\
acting linearly on }$\Bbb{C}^3$\emph{\ whose junior simplex }$\frak{s}_G$%
\emph{\ admits basic, non-coherent triangulations }$\mathcal{T}$\emph{,
which in turn induce crepant, full, non-projective desingularizations of }$%
X\left( N_G,\Delta _G\right) $\emph{. A simple example of this kind comes
into being by taking }$G\cong \left( \Bbb{Z\,}/\,4\Bbb{Z}\right) \times
\left( \Bbb{Z\,}/\,4\Bbb{Z}\right) $\emph{\ to be defined as the abelian
subroup of SL}$\left( 3,\Bbb{C}\right) $\emph{\ generated by the diagonal
elements diag}$\left( \zeta _4,\zeta _4^3,1\right) $\emph{\ and diag}$\left(
1,\zeta _4,\zeta _4^3\right) $\emph{, and }$\mathcal{T}$\emph{\ the
triangulation of figure \textbf{6}. }$\mathcal{T}$ \emph{suffers from} \emph{%
the }``\textit{whirlpool-syndrome}\emph{'' which makes the application of
patching lemma \ref{PATCH} impossible, though strictly upper convex support
functions can be defined on each of its }$2$\emph{-simplices. The
incompatibility of these local strictly upper convex support functions along
the intersection loci of }$1$\emph{-simplices can be also explained by means
of their ``heights'' (cf. \cite{Sturm}, p. 64); the assertion of the
existence of a global upper convex support function on }$\left| \mathcal{T}%
\right| $\emph{\ would lead to a system of a finite number of inconsistent
integer linear inequalities, and hence to a contradiction.}
\end{remark}

\begin{figure*}
\begin{center}
\hspace{1.5cm}\input{fig6.pstex_t}
\end{center}
\begin{center}
 {\textbf{Fig. 6}}
\end{center}
\end{figure*}

\vspace{4cm}
\begin{proposition}[Cohomology dimensions]
\label{COH3} \ \smallskip \newline
Let $\left( \Bbb{C}^3/G,\left[ \mathbf{0}\right] \right) =\left( X\left(
N_G,\Delta _G\right) ,\text{\emph{orb}}\left( \sigma _0\right) \right) $
denote a $3$-dimensional Gorenstein cyclic quotient singularity of type $%
1/l\,\left( \alpha _1,\alpha _2,\alpha _3\right) $. Then the dimensions of
the non-trivial cohomology groups of the overlying spaces $X\left( N_G,\text{
}\widehat{\Delta _G}\right) $ of any $T_{N_G}$-equivariant crepant, full
desingularization of $X\left( N_G,\Delta _G\right) $ are given by the
formulae\emph{:\smallskip } 
\[
\emph{dim}_{\Bbb{Q}}H^{2i}\left( X\left( N_G,\text{ }\widehat{\Delta _G}%
\right) ;\Bbb{Q}\right) =\left\{ 
\begin{array}{ccc}
1 & , & \text{\emph{if \ \ \ }}i=0 \\ 
\, &  &  \\ 
\dfrac 12\,\left( l+\dsum\limits_{j=1}^3\,\text{\emph{gcd}}\left( \alpha
_j,l\right) \right) -2 & , & \text{\emph{if \ \ \ }}i=1 \\ 
\, &  &  \\ 
\dfrac 12\,\left( l-\dsum\limits_{j=1}^3\,\text{\emph{gcd}}\left( \alpha
_j,l\right) \right) +1 & , & \text{\emph{if \ \ \ }}i=2
\end{array}
\right. 
\]
\end{proposition}

\noindent \textit{Proof. }Since $\#\,\left( \partial \left( \frak{s}%
_G\right) \cap N_G\right) =\sum_{j=1}^3\,$gcd$\left( \alpha _j,l\right) -3$,
the lattice points representing the inverses of the remaining junior
group elements (w.r.t.~{\rm int}$(\frak{s}_G)$) belong
to $2\,\frak{s}_G$, and 
\[
\#\,\left( \frak{s}_G\cap N_G\right) +\#\,\left( \left( 2\,\frak{s}_G\right)
\cap \mathbf{Par}\left( \sigma _0\right) \cap N_G\right) =l-1\ , 
\]
we have 
\[
\#\,\left( \text{int}\left( \frak{s}_G\right) \cap N_G\right) =\#\,\left(
\left( 2\,\frak{s}_G\right) \cap \mathbf{Par}\left( \sigma _0\right) \cap
N_G\right) =\dfrac 12\,\left( l-\dsum\limits_{j=1}^3\,\text{gcd}\left(
\alpha _j,l\right) \right) +1 
\]
and the above formulae follow from (\ref{COH-DIM}), (\ref{CYCCOH}). $_{\Box
} $

\begin{exercise}
\emph{Generalize prop. \ref{COH3} for arbitrary abelian acting groups.
[Hint. Fix a splitting of }$G$\emph{\ into cyclic groups. Use denumerants of
weighted partitions instead of gcd's.]}
\end{exercise}

\begin{proposition}[Uniqueness criterion in dimension $3$]
\label{UNICIT}Up to isomorphism, the $3$-dimensional Gorenstein abelian
quotient msc-singularities $\left( \Bbb{C}^3/G,\left[ \mathbf{0}\right]
\right) =\left( X\left( N_G,\Delta _G\right) ,\text{\emph{orb}}\left( \sigma
_0\right) \right) $ which admit a unique \emph{(}full\emph{) }resolution,
are cyclic of type \emph{either\smallskip \smallskip }\newline
$\emph{(i)}$ $\ \ \dfrac 1l\,\left( 1,1,l-2\right) ,\ l\geq 3,$ \ \emph{%
or\smallskip \smallskip }\newline
\emph{(ii)} $\ \dfrac 17\,\left( 1,2,4\right) \smallskip \smallskip $\newline
$\bullet $ In case \emph{(i) }there are $\QTOVERD\lfloor \rfloor {l}{2}$
exceptional prime divisors 
\[
D_j:=V\left( \text{\emph{pos}}\left( \left\{ n^{(j)}\right\} \right) \right)
,\ \ \ \ n^{(j)}:=\dfrac 1l\,\left( j,j,l-2j\right) ,\ \ \ \ 1\leq j\leq
\QTOVERD\lfloor \rfloor {l}{2},
\]
all of whose are compact up to the last one for $l$ even \emph{(see fig. 
\textbf{7} and \textbf{8})}. In particular, one has\smallskip 
\[
\ \ \ \ \ \ \ \ \ \ D_j\cong \left\{ 
\begin{array}{lll}
\Bbb{F}_{l-2\,j} & \emph{,}\text{\emph{\ if}} & \ 1\leq j\leq
\QTOVERD\lfloor \rfloor {l}{2}-1,\ l\geq 4 \\ 
&  &  \\ 
\Bbb{P}_{\Bbb{C}}^2 & \emph{,}\text{\emph{\ if}} & \ j=\frac{l-1}2\text{%
\emph{,} \thinspace }l\ \text{\emph{odd}}\geq 3 \\ 
&  &  \\ 
\Bbb{P}_{\Bbb{C}}^1\times \Bbb{C} & \emph{,}\text{\emph{\ if}} & \ j=\frac
l2,\ l\text{ \thinspace \emph{even }}\geq 4
\end{array}
\right. \smallskip 
\]
The \emph{``}$3$-dimensional\emph{''} $\QOVERD\lfloor \rfloor {l}{2}\times
\QOVERD\lfloor \rfloor {l}{2}\times \QOVERD\lfloor \rfloor {l}{2}$
intersection-number-matrix is determined by\smallskip \smallskip 
\[
\ \left( D_i\cdot D_j\cdot D_k\right) =\left\{ 
\begin{array}{lll}
0 & \emph{,}\text{\emph{\ \ if}} & \ i\neq j,\ j\neq k,\smallskip  \\ 
l-2\,\left( i+1\right)  & \emph{,}\text{\emph{\ \ if}} & \ i=j=k-1\smallskip 
\\ 
0 & \emph{,}\text{\emph{\ \ if}} & \ i=j,\ k-j\geq 2\smallskip  \\ 
2\,i-l & \emph{,}\text{\emph{\ \ if}} & \ i+1=j=k\smallskip  \\ 
0 & \emph{,}\text{\emph{\ \ if}} & \ j=k,\ j-i\geq 2\smallskip  \\ 
8 & \emph{,}\text{\emph{\ \ if}} & \ i=j=k\neq \left( l-1\right) \,/\,2,\ l\ 
\text{\emph{odd}}\smallskip  \\ 
9 & \emph{,}\text{\emph{\ \ if}} & \ i=j=k=\left( l-1\right) \,/\,2,\ l\ 
\text{\emph{odd}}\smallskip  \\ 
8 & \emph{,}\text{\emph{\ \ if}} & \ i=j=k\neq l\,/\,2,\ l\ \text{\emph{even}%
}\smallskip  \\ 
-2 & \emph{,}\text{\emph{\ \ if}} & \ i=j=k=l\,/\,2,\ l\ \text{\emph{even}}
\end{array}
\right. \smallskip 
\]
for all $1\leq i\leq j\leq k\leq \QOVERD\lfloor \rfloor {l}{2}$ \emph{(}and
by the usual symmetric property of intersection numbers\emph{)}\smallskip .
Moreover, the non-trivial cohomology dimensions of the desingularizing space
equal $1,\QOVERD\lfloor \rfloor {l}{2}$ and $\QOVERD\lfloor \rfloor {l-1}{2}$%
, respectively\emph{\ (}by prop. \emph{\ref{COH3}).\smallskip }\newline
$\bullet $ In case \emph{(ii), }there are three exceptional prime divisors,
namely 
\[
D_n=V\left( \text{\emph{pos}}\left( \left\{ n\right\} \right) \right) ,\
D_{n^{\prime }}=V\left( \text{\emph{pos}}\left( \left\{ n^{\prime }\right\}
\right) \right) ,\ D_{n^{\prime \prime }}=V\left( \text{\emph{pos}}\left(
\left\{ n^{\prime \prime }\right\} \right) \right) ,
\]
with 
\[
n=\frac 17\left( 1,2,4\right) ^{\intercal },\ \ n^{\prime }=\frac 17\left(
2,4,1\right) ^{\intercal },\ \ n^{\prime \prime }=\frac 17\left(
4,1,2\right) ^{\intercal },
\]
each of which is isomorphic to $\Bbb{F}_2$. They intersect each other
paiwise along three rational curves which play interchangeably the roles of
the fibers and of the $0$-sections of the three projectivized $\Bbb{P}_{\Bbb{%
C}}^1$-bundles $\Bbb{F}_2\rightarrow \Bbb{P}_{\Bbb{C}}^1$ \emph{(see fig. 
\textbf{9}). }Obviously, in both cases\emph{\ (i) }and\emph{\ (ii) }the
desingularizing birational morphism is projective \emph{(}by prop. \emph{\ref
{NON-EM}).\smallskip }
\end{proposition}

\noindent \textit{Proof.} The uniqueness (up to automorphisms of aff$\left( 
\frak{s}_G\right) \cap N_G$) of the triangulation $\mathcal{T}$ of $\frak{s}%
_G$ (inducing a unique crepant desingularization, up to isomorphism) means
that for every $4$-tuple $\left\{ n_1,n_2,n_3,n_4\right\} $ of distinct
elements of $\frak{s}_G\cap N_G$, 
\begin{equation}
\left\{ n_1,n_2,n_3,n_4\right\} \subset \text{conv}\left( \left\{ n_{\nu
_1},n_{\nu _2},n_{\nu _3}\right\} \right) \text{, \ }  \label{QUADR}
\end{equation}
for all $\left\{ \nu _1,\nu _2,\nu _3\right\} \subset \left\{
1,2,3,4\right\} $ (so that conv$\left( \left\{ n_1,n_2,n_3,n_4\right\}
\right) $ cannot be a convex quadrilateral). Since 
\begin{equation}
\text{splcod}\left( \text{orb}\left( \sigma _0\right) ;U_{\sigma _0}\right)
=3  \label{splc3}
\end{equation}
it is easy to prove that $G$ cannot be abelian, \textit{non-cyclic}. For
cyclic $G$'s acting on $\Bbb{C}^3$ by type $1/l\,\left( \alpha _1,\alpha
_2,\alpha _3\right) $ the uniqueness condition will be examined in two
different cases.\smallskip \newline
$\bullet $ If the cardinal number of $\frak{s}_G\cap \left(
N_G\smallsetminus \left\{ e_1,e_2,e_3\right\} \right) $ is $\geq 4$, then (%
\ref{QUADR}) is equivalent to say that all points of it lie on a straight
line going through precisely one of the vertices $e_1,e_2,e_3$ of $\frak{s}_G
$ (but, of course, $\#\,\left( \partial \frak{s}_G\cap \left(
N_G\smallsetminus \left\{ e_1,e_2,e_3\right\} \right) \right) \in \left\{
0,1\right\} $, because otherwise (\ref{splc3}) would be violated). This
occurs only in the case in which at least two of the weights $\alpha
_1,\alpha _2,\alpha _3$ are equal; but then 
\[
\left( \alpha _1,\alpha _2,\alpha _3\right) \backsim \left( 1,1,l-2\right)
\,\,\ \ \text{(within\thinspace \thinspace }\Lambda \left( l;3\right) \text{)%
}
\]
(check it or see lemma \ref{wnor} below), and therefore $\left( \Bbb{C}%
^3/G,\left[ \mathbf{0}\right] \right) $ must have type of the form
(i).\smallskip \newline
$\bullet $ If $\#\,\left( \frak{s}_G\cap \left( N_G\smallsetminus \left\{
e_1,e_2,e_3\right\} \right) \right) \leq 3$, then we get the inequality 
\[
\dfrac 12\,\left( l+\dsum\limits_{j=1}^3\,\text{gcd}\left( \alpha
_i,l\right) \right) -2\leq 3\ ,
\]
which is valid only for 
\[
\frac 1l\,\left( \alpha _1,\alpha _2,\alpha _3\right) \in \left\{ 
\begin{array}{c}
1/3\,\left( 1,1,1\right) ,\,1/4\,\left( 1,1,2\right) ,\,1/5\,\left(
1,1,3\right) ,\smallskip  \\ 
1/5\,\left( 1,2,2\right) ,\,1/6\,\left( 1,1,4\right) ,\,1/7\,\left(
1,1,5\right) ,\,1/7\,\left( 1,2,4\right) 
\end{array}
\right\} \ .
\]
(It suffices to assume $\alpha _1+\alpha _2+\alpha _3=l$, because $\left(
\alpha _1,\alpha _2,\alpha _3\right) \backsim \left( l-\alpha _1,l-\alpha
_2,l-\alpha _3\right) $ within $\Lambda \left( l;3\right) $). Since $\left(
1,2,2\right) \backsim \left( 1,1,3\right) $ within $\Lambda \left(
5;3\right) $, we see that from the above $7$ possible types only (ii) $%
1/7\,\left( 1,2,4\right) $ is inequivalent (w.r.t. ``$\backsim $'') to all
those of the form (i). The non-vertex lattice points of $\frak{s}_G$ for $%
\left( \Bbb{C}^3/G,\left[ \mathbf{0}\right] \right) $ being of this ``new''
type (ii) satisfy obviously (\ref{QUADR}) and the proof is
completed.\smallskip \newline
$\bullet $ That the $\QTOVERD\lfloor \rfloor {l}{2}$ exceptional prime
divisors $D_j$ have the structure given above is an immediate consequence of
the more general thm. \ref{DESEXC} of the next section which will be proved
for all dimensions. The intersection numbers $\left( D_i\cdot D_j\cdot
D_k\right) $ are computable by (\ref{INTERSE}), (\ref{D3}).\newline
$\bullet $ Since 
\[
n+e_1=2\,n^{\prime \prime },\ \ \ n^{\prime }+e_3=2\,n,\ \ \ n^{\prime
\prime }+e_2=2\,n^{\prime },
\]
the above assertion for the structure of the exceptional prime divisors $%
D_n,D_{n^{\prime }},D_{n^{\prime \prime }}$ is obvious (e.g. by \ref{YPSILON}%
). For another proof, see Roan-Yau \cite{Roan-Yau}, pp. 272-273. $_{\Box }$
\vfill
\begin{figure*}[hb]
\begin{center}
\input{fig7.pstex_t}
\makespace
 {\textbf{Fig. 7} ($l$ odd)}
\end{center}
\end{figure*}
%
%
\begin{figure*}
\begin{center}
\input{fig8.pstex_t}
\makespace
 {\textbf{Fig. 8} ($l$ even)}
\end{center}
\end{figure*}

\begin{figure*}
\begin{center}
\input{fig9.pstex_t}
\makespace
 {\textbf{Fig. 9}}
\end{center}
\end{figure*}

\section{On the monoparametrized 
singularity-series \\ $\dfrac 1l\left( 1,\ldots ,1,l-\left( r-1\right) \right) $}
\label{MONSE}

\noindent This section contains our main results. Motivated by \ref{DIM2},
the uniqueness criterion \ref{UNICIT} (i) in dimension $3$, and Reid's
remark \cite{Reid5}, 5.4, (concerning dimension $4$), we study the
monoparametrized singularity-series of \textit{arbitrary} dimension with the
simplest possible ``lattice-geometry'', i.e. those Gorenstein cyclic
quotient singularities whose junior simplex encloses only lattice points
lying on a single straight line.

\begin{lemma}
\label{wnor}Let $\left( \Bbb{C}^r/G,\left[ \mathbf{0}\right] \right) =\left(
X\left( N_G,\Delta _G\right) ,\text{\emph{orb}}\left( \sigma _0\right)
\right) $ be the Gorenstein cyclic quotient msc-singularity of type
\smallskip $\dfrac 1l\left( \alpha _1,\ldots ,\alpha _r\right) $ \emph{(}%
with $l=\left| G\right| \geq r\geq 2$\emph{)} for which at least $r-1$ of
its defining weights are equal. Then \emph{(}in the notation of \emph{\S\ 
\ref{ABELQ}, \textsf{(c)})} 
\[
\left( \alpha _1,\ldots ,\alpha _r\right) \,\backsim \left( 1,\ldots
,1,\,l-\left( r-1\right) \right) \,\ \ \ \text{\emph{(within}\thinspace
\thinspace }\Lambda \left( l;r\right) \text{\emph{) .}}
\]
\end{lemma}

\noindent \textit{Proof. }Using a permutation sending the $r-1$ equal
defining weights of the above singularity to the first $r-1$ positions, we
have 
\[
\left( \alpha _1,\ldots ,\alpha _r\right) \backsim \left( \mu ,\ldots ,\mu
,\alpha _r\right) 
\]
for an integer $\mu $, $1\leq \mu \leq l-1$. Since $G$ contains no
pseudoreflections, gcd$\left( \mu ,l\right) =1$, i.e. there exists an
integer $\nu $, $1\leq \nu \leq l-1$, with gcd$\left( \nu ,l\right) =1$, and 
$\lambda \in \Bbb{Z}$, such that 
\begin{equation}
\nu \,\mu +\lambda \,\,l=1  \label{GGT}
\end{equation}
Moreover, since $\left( \Bbb{C}^r/G,\left[ \mathbf{0}\right] \right) $ is
Gorenstein, there must be an integer $\kappa $, $1\leq \kappa \leq r-1$,
such that 
\begin{equation}
\left( r-1\right) \,\mu +\alpha _r=\kappa \,l  \label{GBED}
\end{equation}
Equalities (\ref{GGT}) and (\ref{GBED}) imply 
\[
\nu \,\alpha _r=\left( \kappa \,\nu +\left( r-1\right) \,\lambda -1\right)
\,l+\left( l-\left( r-1\right) \right) \ . 
\]
Hence, 
\[
\left( \mu ,\ldots ,\mu ,\alpha _r\right) \backsim \left( \left[ \nu \,\mu
\right] _l,\ldots ,\left[ \nu \,\mu \right] _l,\left[ \nu \,\alpha _r\right]
_l\right) =\left( 1,\ldots ,1,\,l-\left( r-1\right) \right) 
\]
and we are done. $_{\Box }$ \medskip

\noindent By the above lemma and cor. \ref{ISOCYC} we can obviously restrict
ourselves to the study of singularities of type (\ref{montype}).

\begin{theorem}[Resolution by a unique projective, crepant morphism]
\label{Main} \ \smallskip \newline
Let $\left( \Bbb{C}^r/G,\left[ \mathbf{0}\right] \right) =\left( X\left(
N_G,\Delta _G\right) ,\text{\emph{orb}}\left( \sigma _0\right) \right) $ be
the Gorenstein cyclic quotient singularity of type\smallskip 
\begin{equation}
\fbox{$
\begin{array}{ccc}
&  &  \\ 
& \dfrac 1l\ \left( \stackunder{\left( r-1\right) \text{\emph{-times}}}{%
\underbrace{1,1,\ldots ,1,1}},\ l-\left( r-1\right) \right)  &  \\ 
&  & 
\end{array}
$}  \label{montype}
\end{equation}
with $l=\left| G\right| \geq r\geq 2.$ Then we have \emph{:\smallskip }%
\newline
\emph{(i)} This msc-singularity is isolated if and only if \emph{gcd}$\left(
l,r-1\right) =1.\smallskip $\newline
\emph{(ii)} Up to affine integral transformation, there exists a \emph{%
unique }triangulation 
\[
\mathcal{T}\in \mathbf{LTR}_{N_G}^{\emph{max}}\,\left( \frak{s}_G\right) 
\]
inducing a unique \emph{(}isomorphism class of\emph{)} crepant $T_N$%
-equivariant partial resolution-morphism 
\[
f:X\left( N_G,\widehat{\Delta }_G\left( \mathcal{T}\right) \right)
\rightarrow X\left( N_G,\Delta _G\right) 
\]
of $X\left( N_G,\Delta _G\right) $ \emph{(}with overlying space $\Bbb{Q}$%
-factorial, and maximal with respect to non-discrepancy\emph{).\smallskip }%
\newline
\emph{(iii)} $\mathcal{T}\in \mathbf{Coh}$-$\mathbf{LTR}_{N_G}^{\emph{max}%
}\,\left( \frak{s}_G\right) $, i.e., $f$ is projective.\smallskip \newline
\emph{(iv)} $\mathcal{T}\in \mathbf{Coh}$-$\mathbf{LTR}_{N_G}^{\emph{basic}%
}\,\left( \frak{s}_G\right) $ \emph{(}in other words, $f$ gives rise to a 
\emph{full }projective, crepant desingularization\emph{)} if and only
if\smallskip 
\begin{equation}
\fbox{$
\begin{array}{ccc}
&  &  \\ 
& \left[ \ l\ \right] _{\ r-1}\in \left\{ 0,1\right\}  &  \\ 
&  & 
\end{array}
$}\smallskip   \label{rescon}
\end{equation}
\emph{(i.e., iff either }$l\equiv 0$ \emph{mod}$\left( r-1\right) $ \emph{or 
}$l\equiv 1$ \emph{mod}$\left( r-1\right) .$\emph{)\smallskip }\newline
\emph{(v) }For $l$ satisfying condition \emph{(\ref{rescon})}, the
dimensions of the non-trivial cohomology groups of the resolving space $%
X\left( N_G,\widehat{\Delta }_G\right) $ are given by the formulae\emph{\ :} 
\begin{equation}
\fbox{$\text{\emph{dim}}_{\Bbb{Q}}H^{2i}\left( X\left( N_G,\widehat{\Delta }%
_G\right) ;\Bbb{Q}\right) =\left\{ 
\begin{array}{llll}
1 & , & \text{\emph{for }} & i=0 \\ 
\  &  &  &  \\ 
\QTOVERD\lfloor \rfloor {l}{r-1} & , & \text{\emph{for}} & i\in \left\{
1,2,\ldots ,r-2\right\}  \\ 
\  &  &  &  \\ 
\QTOVERD\lfloor \rfloor {l-1}{r-1} & , & \text{\emph{for}} & i=r-1
\end{array}
\right. $}  \label{cohd}
\end{equation}
\end{theorem}

\noindent \textit{Proof. }(i) This follows directly from \ref{ISOL}. Note
that if gcd$\left( l,r-1\right) \geq 2$, then 
\[
\text{Sing}\left( X\left( N_G,\Delta _G\right) \right) =U_{\sigma _0}\left( 
\text{pos}\left( \left\{ e_1,e_2,\ldots ,e_{r-2},e_{r-1}\right\} \right)
\right) , 
\]
i.e. the singular locus of $X\left( N_G,\Delta _G\right) =\Bbb{C}^r/G$ is
the entire ``$z_r$-axis'' of $\Bbb{C}^r$.\medskip \newline
(ii) Let us first introduce some notation and make certain preparatory
remarks. Define the vectors$:$ 
\[
n^{\left( j\right) }:=\left\{ 
\begin{array}{lll}
e_r & , & \text{if\ \ \ }j=0 \\ 
& \  &  \\ 
\dfrac 1l\ \left( \stackunder{\left( r-1\right) \text{-times}}{\underbrace{%
j,j,\ldots ,j,j}},\left[ j\cdot \left( l-\left( r-1\right) \right) \right]
_l\right) ^{\intercal } & , & \text{if\ \ \ }j\in \left\{ 1,2,\ldots
,l-1\right\}
\end{array}
\right. 
\]
and denote by $n_i^{\left( j\right) },1\leq i\leq r,$ the $i$-th coordinate
of each $n^{\left( j\right) },\ 0\leq j\leq l-1$, within $\left( N_G\right)
_{\Bbb{R}}\cong \Bbb{R}^r$. We have 
\[
\sum_{i=1}^r\ n_i^{\left( j+1\right) }\geq \sum_{i=1}^r\ n_i^{\left(
j\right) }+\dfrac 1l\left( r-1\right) -1,\ \ \forall j,\ \ 1\leq j\leq l-1. 
\]
If for some fixed $j=$ $j_{\bullet }\in \left\{ 1,\ldots ,l-1\right\} ,$ $%
\sum_{i=1}^r\ n_i^{\left( j_{\bullet }\right) }\geq 2$, then obviously 
\[
\sum_{i=1}^r\ n_i^{\left( k\right) }\geq 2,\ \ \forall k,\ \ 1\leq
j_{\bullet }\leq k\leq l-1. 
\]
Since 
\[
\sum_{i=1}^r\ n_i^{\left( j\right) }=\dfrac 1l\ \left( \left( r-1\right)
\cdot j+\left[ j\cdot \left( l-\left( r-1\right) \right) \right] _l\right) 
\]
we obtain the inclusion 
\[
\frak{s}_G\cap N_G\subset \left\{ e_1,\ldots ,e_r\right\} \cup \left\{
n^{\left( j\right) }\ \left| \ 1\leq j\leq \QOVERD\lfloor \rfloor
{l}{r-1}\right. \right\} \ . 
\]
And conversely, for all $j,\ 1\leq j\leq \QOVERD\lfloor \rfloor {l}{r-1},$
we have $\left[ j\cdot \left( l-\left( r-1\right) \right) \right]
_l=l-j\cdot \left( r-1\right) $, which means that $\sum_{i=1}^r\ n_i^{\left(
j\right) }=1.$ Thus, we get the equality 
\begin{equation}
\frak{s}_G\cap N_G=\left\{ e_1,\ldots ,e_r\right\} \cup \left\{ n^{\left(
j\right) }\ \left| \ 1\leq j\leq \QOVERD\lfloor \rfloor {l}{r-1}\right.
\right\} \   \label{JSEL}
\end{equation}
\smallskip \newline
$\bullet $ \textit{Construction of $\mathcal{T}$. \ }At first define 
\[
\Xi _r:=\left\{ \left( \xi _1,\xi _2,\ldots ,\xi _{r-2}\right) \in \left(
\left\{ 1,2,\ldots ,r-2,r-1\right\} \right) ^{r-2}\ \left| \ \xi _1<\xi
_2<\cdots <\xi _{r-2}\right. \right\} , 
\]
and 
\[
B\left( j;\xi _1,\xi _2,\ldots ,\xi _{r-2}\right) :=\left\{ n^{\left(
j-1\right) },n^{\left( j\right) },\ e_{\xi _1},e_{\xi _2},\ldots ,e_{\xi
_{r-2}}\right\} , 
\]
as well as 
\[
\mathbf{s}\left( j;\xi _1,\xi _2,\ldots ,\xi _{r-2}\right) :=\text{ conv}%
\left( B\left( j;\xi _1,\xi _2,\ldots ,\xi _{r-2}\right) \right) ,\smallskip
\smallskip \ 
\]
and 
\[
\sigma \left( j;\xi _1,\xi _2,\ldots ,\xi _{r-2}\right) :=\sigma _{\mathbf{s}%
\left( j;\xi _1,\xi _2,\ldots ,\xi _{r-2}\right) }=\text{ pos}\left( B\left(
j;\xi _1,\xi _2,\ldots ,\xi _{r-2}\right) \right) ,\smallskip 
\]
for all $j$, 
\[
1\leq j\leq \QOVERD\lfloor \rfloor {l}{r-1},\ \ \ \text{and all \ }r-1\text{
\ possible \ }\left( r-2\right) \text{-tuples }\ \left( \xi _1,\xi _2,\ldots
,\xi _{r-2}\right) \in \Xi _r\ .\smallskip 
\]
All lattice points $n^{\left( j\right) }$, $0\leq j\leq \QOVERD\lfloor
\rfloor {l}{r-1}$, lie on a \textit{straight line} of $\left( N_G\right) _{%
\Bbb{R}}\cong \Bbb{R}^r$, namely on\smallskip 
\[
\text{conv}\left( \left\{ n^{\left( 0\right) },\,n^{\left( \QOVERD\lfloor
\rfloor {l}{r-1}\right) }\right\} \right) 
\]
and they can be ordered canonically 
\[
n^{\left( 0\right) }\ ,\ n^{\left( 1\right) },\ \ldots \ ,\ n^{\left(
\QOVERD\lfloor \rfloor {l}{r-1}-1\right) }\ ,\ n^{\left( \QOVERD\lfloor
\rfloor {l}{r-1}\right) } 
\]
via the increasing ordering of the enumerator-superscripts $\left( j\right) $%
.\smallskip \newline
Next define the simplicial subdivision $\mathcal{T}$ (having support $\left| 
\mathcal{T}\right| \subset \frak{s}_G$) as follows :\smallskip 
\[
\mathcal{T}:=\left\{ 
\begin{array}{lll}
\frak{E\ } & , & \text{if\ \ }l\equiv 0\,\left( \text{mod }\left( r-1\right)
\right) \\ 
&  &  \\ 
\frak{E\ }\cup \ \left\{ 
\begin{array}{c}
\text{all faces of\emph{\ }} \\ 
\text{conv}\left( B_{\text{last}}\right)
\end{array}
\right\} & , & \text{otherwise\ \ }
\end{array}
\right. 
\]
where 
\[
B_{\text{last}}:=\left\{ n^{\left( \left\lfloor \frac{l-1}{r-1}\right\rfloor\right) },e_1,e_2,\ldots
,e_{r-1}\right\} , 
\]
and 
\[
\frak{E}:=\left\{ \left. 
\begin{array}{c}
\text{all faces of the simplices} \\ 
\mathbf{s}\left( j;\xi _1,\xi _2,\ldots ,\xi _{r-2}\right)
\end{array}
\ \right| \ 
\begin{array}{l}
\text{for all\ }j,\text{ }1\leq j\leq \QOVERD\lfloor \rfloor {l}{r-1},\text{
and}\ \smallskip \  \\ 
\text{all \ }\left( \xi _1,\xi _2,\ldots ,\xi _{r-2}\right) \in \Xi _r
\end{array}
\right\} \ \ .\smallskip \smallskip 
\]
$\bullet $ $\mathcal{T}$ \thinspace is a triangulation. The proof of the
assertion that the intersection of two arbitrary simplices of $\mathcal{T}$
is either a face of both or the empty set is left as an easy exercise to the
reader.\medskip \newline
$\bullet $ $\left| \mathcal{T}\right| =\frak{s}_G$, i.e., its support covers
the entire junior simplex. It is straightforward to show that 
\begin{eqnarray}
\label{martin1}
& & \text{Vol}(\text{conv}\left(\{\mathbf{0}\},\mathbf{s}\left( j;\xi _1,\xi _2,\ldots ,\xi _{r-2}\right)\right))= \nonumber \\ && \quad \frac{1}{r!}
\left| \text{det}\left( n^{\left( j-1\right)
},n^{\left( j\right) },e_{\xi _1},e_{\xi _2},\ldots ,e_{\xi _{r-2}}\right)
\right| = \nonumber \\ 
 &&\quad
\frac{1}{r!}\dfrac 1{l^2}\ \left| \left( j-1\right) \cdot \left( l-j\left( r-1\right)
\right) -j\cdot \left( l-\left( j-1\right) \left( r-1\right) \right) \right| 
 = \frac{1}{l r!}
\end{eqnarray}
and 
\begin{equation}
\label{martin2}
\text{Vol}(\text{conv}\left(\{\mathbf{0}\},B_{\text{last}}\right))=\frac{[l]_{r-1}}{lr!}.
\end{equation}
Thus the assertion  is an immediate consequence of the equality
: 
\[
\begin{array}{l}
\text{Vol}\left( \text{conv}\left( \left\{ \mathbf{0}\right\} ,\bigcup
\,\left\{ \mathbf{s\ }\left| \ \mathbf{s}\text{ simplex of }\mathcal{T}%
\right. \right\} \right) \right) = \\ 
\  \\ 
=\dfrac 1{r!}\dfrac 1l\, \left( \,\QDOVERD\lfloor \rfloor
{l}{r-1}\cdot \left( r-1\right) +\left[ \ l\ \right] _{\ r-1}\right)
 =\dfrac 1{r!}=\text{Vol}\left( \text{conv}\left( \left\{ \mathbf{0}%
\right\} ,\frak{s}_G\right) \right) \ \text{.\smallskip }
\end{array}
\]
$\bullet $ $\mathcal{T}$ \thinspace is a maximal triangulation. This is
obvious because $\frak{s}_G\cap N_G=$ vert$\left( \mathcal{T}\right) $ by
the construction of $\mathcal{T}$ and (\ref{JSEL}).\smallskip \newline
$\bullet $ $\mathcal{T}$ \thinspace is unique. Suppose $\mathbf{s}$ is an
arbitrary elementary $\left( r-1\right) $-simplex, such that vert$\left( 
\mathbf{s}\right) \subset \frak{s}_G\cap N_G$. Clearly, $\mathbf{s}%
\subsetneqq \frak{s}_G$ because $\frak{s}_G$ itself is non-elementary. This
implies 
\[
\text{vert}\left( \mathbf{s}\right) \cap \left( \frak{s}_G\cap
N_G\smallsetminus \left\{ e_1,\ldots ,e_r\right\} \right) \neq \varnothing \
. 
\]
Since $n^{\left( j\right) }$'s are collinear and $\mathbf{s}$ elementary
simplex, vert$\left( \mathbf{s}\right) $ contains at most $2$ consecutive
lattice points among them. Hence, 
\[
\#\,\left( \text{vert}\left( \mathbf{s}\right) \cap \left( \frak{s}_G\cap
N_G\smallsetminus \left\{ e_1,\ldots ,e_r\right\} \right) \right) \in
\left\{ 1,2\right\} \ . 
\]
\newline
If this number equals $2$, then $\mathbf{s}$ has to be of the form $\mathbf{s%
}\left( j;\xi _1,\ldots ,\xi _{r-2}\right) $, for some $j\in \left\{
0,\ldots ,\QOVERD\lfloor \rfloor {l}{r-1}\right\} $ and $\left( \xi
_1,\ldots ,\xi _{r-2}\right) \in \Xi _r$ (by definition). Otherwise, $%
\mathbf{s}$ possesses only one vertex belonging to the relative interior of $%
\frak{s}_G$. However, since the $n^{\left( j\right) }$'s are collinear and $n^{(0)}=e_r$ we must have  $l \not\equiv 0\,\text{mod}\,(r-1)$ and $\mathbf{s}=B_{\text{last}}$.
Hence, $\mathcal{T}$ is uniquely determined.$\medskip $\newline
(iii) The coherence of $\mathcal{T}$ will be proved by induction on the
number $\QOVERD\lfloor \rfloor {l}{r-1}$. The case in which this equals $1$
is trivial. Suppose $\QOVERD\lfloor \rfloor {l}{r-1}>1$. It is easy to check
that the simplex 
\[
\widetilde{\mathbf{s}}\mathbf{:=}\text{ conv}\left( \left\{ n_{}^{\left(
1\right) },e_1,\ldots ,e_{r-1}\right\} \right) 
\]
can be mapped by an affine integral transformation onto the junior simplex
of an equidimensional Gorenstein cyclic quotient singularity of type 
\[
\frac 1{l-\left( r-1\right) }\ \left( 1,\ldots ,1,l-2\left( r-1\right)
\right) 
\]
with 
\[
\left\lfloor \frac{l-\left( r-1\right) }{r-1}\right\rfloor =\left\lfloor
\frac l{r-1}-1\right\rfloor =\left\lfloor \frac l{r-1}\right\rfloor
-1<\QOVERD\lfloor \rfloor {l}{r-1}\ . 
\]
By induction hypothesis and by construction we may therefore assume that the
restriction $\mathcal{T}\left| _{\widetilde{\mathbf{s}}}\right. $ of $%
\mathcal{T}$ on $\widetilde{\mathbf{s}}$ is coherent. Choose a $\widetilde{%
\psi }\in $ SUCSF$_{\Bbb{R}}\left( \mathcal{T}\left| _{\widetilde{\mathbf{s}}%
}\right. \right) $ and use the abbreviation 
\[
\mathbf{s}_\iota :=\left\{ 
\begin{array}{lll}
\mathbf{s}\left( 1;\,1,2,\ldots ,\,r-3,\,r-2\right) & , & \text{for }\iota =1
\\ 
\mathbf{s}\left( 1;\,1,\ldots ,\,\iota -1,\,\iota +1,\ldots ,\,r-1\right) & ,
& \forall \iota ,\ 2\leq \iota \leq r-2 \\ 
\mathbf{s}\left( 1;\,2,\,3,\ldots ,\,r-2,\,r-1\right) & , & \text{for }\iota
=r-1
\end{array}
\right. 
\]
Note that for each $\mathbf{y}\in \mathbf{s}_\iota $ and $\iota $, $1\leq
\iota \leq r-1$, there exists an $\mathbf{y}^{\prime }$ $\in \widetilde{%
\mathbf{s}}$, such that 
\[
\mathbf{y}^{\prime }\in \text{ lin}\left( \text{conv}\left( \left\{ e_r,%
\mathbf{y}\right\} \right) \right) \cap \widetilde{\mathbf{s}}\text{ \ \ \
\& \ \ }\mathbf{y=}\left( 1-t\right) \,\mathbf{y}^{\prime }+t\,e_r\text{, \
\ for \thinspace a\thinspace \ \ }t\in \left[ 0,1\right] \,, 
\]
because of the convexity of $\widetilde{\mathbf{s}}\subset \frak{s}_G$. Now
define support functions $\psi _\iota :\left| \mathcal{T}\left| _{\mathbf{s}%
_\iota }\right. \right| \rightarrow \Bbb{R}$ by 
\[
\mathbf{s}_\iota \ni \mathbf{y}\stackrel{}{\longmapsto \psi _\iota \left( 
\mathbf{y}\right) :=}\mathbf{\,}\left( 1-t\right) \,\widetilde{\psi }\left( 
\mathbf{y}^{\prime }\right) +t\,\in \Bbb{R\,}. 
\]
One verifies easily that $\psi _\iota \in $ SUCSF$_{\Bbb{R}}\left( \mathcal{T%
}\left| _{\mathbf{s}_\iota }\right. \right) $ and 
\[
\psi _\iota \,\left| _{\widetilde{\mathbf{s\,}}\cap \,\mathbf{s}_\iota
}\right. =\widetilde{\psi }\,\left| _{\widetilde{\mathbf{s}}\,\cap \,\mathbf{%
s}_\iota }\right. \smallskip 
\]
for all $\iota $, $1\leq \iota \leq r-1$. Applying patching lemma \ref{PATCH}%
, we get SUCSF$_{\Bbb{R}}\left( \mathcal{T}\right) \neq \varnothing $ and we
are done.\smallskip \newline
$\bullet $ An alternative (but not directly constructive) method for showing
the coherence of $\mathcal{T}$, is to combine the uniqueness of $\mathcal{T}$
with prop. \ref{NON-EM}.\medskip \newline
(iv) Since $\mathcal{T}$ is uniquely determined the volume formulae in (\ref{martin1}) and (\ref{martin2}) show that  $f$ is  a full, crepant torus-equivariant resolution-morphism if and only if $[l]_{r-1}\in\{0,1\}$. Another, more direct way, to verify that condition (\ref{rescon}) is necessary is the following:
 
\noindent If we
assume that $\left[ \ l\ \right] _{\ r-1}\notin \left\{ 0,1\right\} $, the
group order $l$ can be written as 
\begin{equation}
l=\QOVERD\lfloor \rfloor {l}{r-1}\cdot \left( r-1\right) +\left[ \ l\
\right] _{\ r-1},\text{ \ with \ \thinspace }\left[ \ l\ \right] _{\ r-1}\in
\left\{ 2,\ldots ,r-2\right\} .  \label{MAGIA}
\end{equation}
Note that 
\[
n_r^{\left( j\right) }=1-\frac{j\cdot \left( r-1\right) }l\,,\ \ \forall j,\
\ 1\leq j\leq \QOVERD\lfloor \rfloor {l}{r-1}\ , 
\]
\newline
and that the decreasing sequence 
\begin{equation}
n_r^{\left( 1\right) }>n_r^{\left( 2\right) }>\cdots >n_r^{\left(
\QOVERD\lfloor \rfloor {l}{r-1}\right) }  \label{DECR}
\end{equation}
has minimum element 
\[
n_r^{\left( \QOVERD\lfloor \rfloor {l}{r-1}\right) }=\frac{\left[ \ l\
\right] _{\ r-1}}l\ \left( \geq \frac 2l\right) \ . 
\]
\newline
We shall discuss the two possible cases w.r.t. the values taken by gcd$%
\left( l,r-1\right) $ separately.\medskip \newline
I) gcd$\left( l,r-1\right) =1$. Then there exists a $\kappa \in \left\{
1,\ldots ,l-1\right\} $, such that 
\[
\kappa \,\left( l-\left( r-1\right) \right) \equiv 1\,\left( \text{mod }%
l\right) . 
\]
Clearly, 
\begin{equation}
n_r^{\left( \kappa \right) }=\frac 1l\,\left[ \kappa \,\left( l-\left(
r-1\right) \right) \right] _l=\frac 1l<n_r^{\left( \QOVERD\lfloor \rfloor
{l}{r-1}\right) }  \label{NK}
\end{equation}
and consequently, 
\[
\kappa \in \left\{ \QOVERD\lfloor \rfloor {l}{r-1}+1,\ldots ,l-1\right\}
\Rightarrow n^{\left( \kappa \right) }\notin \frak{s}_G\cap N_G\smallskip \
, 
\]
by (\ref{JSEL}) and (\ref{DECR}). If there were a basic triangulation of $%
\frak{s}_G$, then $n^{\left( \kappa \right) }$ would not belong to $\mathbf{%
Hlb}_{N_G}\left( \sigma _0\right) $ (by the equality (\ref{HILBCON}) of thm. 
\ref{KILLER}), but it would be representable as a linear combination of (at
least two) elements of the Hilbert basis $\mathbf{Hlb}_{N_G}\left( \sigma
_0\right) $ with positive integer coefficients. This would obviously
contradict to (\ref{NK}) because of (\ref{DECR}).\medskip \newline
II) gcd$\left( l,r-1\right) \geq 2$. In this case, define $\kappa :=l\,/\,$%
gcd$\left( l,r-1\right) $. Then 
\[
\kappa \,\left( l-\left( r-1\right) \right) \equiv 0\,\left( \text{mod }%
l\right) , 
\]
and repeating the same argumentation for this new $n_r^{\left( \kappa
\right) }=0$ as in I), we arrive again at a contradicting
conclusion.\medskip \newline
\medskip (v) Formula (\ref{cohd}) follows from (\ref{CYCCOH}). In
particular, 
\[
\text{dim}_{\Bbb{Q}}H^{2i}\left( X\left( N_G,\widehat{\Delta }_G\right) ;%
\Bbb{Q}\right) =\#\,\left\{ j\in \left[ 0,l\right) \cap \Bbb{Z\ }\left| \
\QTOVERD\lfloor \rfloor {\left( i-1\right) \,l}{r-1}+1\leq j\leq
\QTOVERD\lfloor \rfloor {i\,l}{r-1}\right. \right\} 
\]
and the proof of the theorem is completed. $_{\Box }\smallskip $

\begin{remark}
\emph{(i) To prove the necessity of condition (\ref{rescon}) for $\mathcal{T}
$ \thinspace to be basic, we preferred to make use of ``Hilbert-base
technology'' because it is generally applicable to }\textit{any} \emph{%
Gorenstein abelian quotient singularity. Alternative ad hoc methods (for the
above special situation) are either a suitable direct manipulation of
determinants or the use of normalized-volume-bound. According to the latter
one, violation of (\ref{rescon}) would imply for the topological
Euler-Poincar\'{e} characteristic of the overlying space:} 
\[
\begin{array}{lll}
\chi \left( X\left( N_G,\widehat{\Delta }_G\left( \mathcal{T}\right) \right)
\right) =\#\,\left( \widehat{\Delta }_G\left( \mathcal{T}\right) \left(
r\right) \right) = &  &  \\ 
\  &  &  \\ 
=\#\,\left( \left( r-1\right) \text{\emph{-simplices of }}\mathcal{T}\right)
= &  & \text{\emph{(by (\ref{Eul-Poi}))}} \\ 
&  &  \\ 
=\QOVERD\lfloor \rfloor {l}{r-1}\cdot \left( \#\,\left( \Xi _r\right)
\right) +1=\QOVERD\lfloor \rfloor {l}{r-1}\cdot \left( r-1\right) +1= &  & 
\\ 
\  &  &  \\ 
=l-\left[ \ l\ \right] _{\ r-1}+1 < l=\left| G\right|  &  & \text{%
\emph{(by (\ref{MAGIA}))}}
\end{array}
\]
\newline
\emph{which would be impossible for basic }$\mathcal{T}$\emph{\ by (\ref
{EP-DES}).\medskip }\newline
\emph{(ii) In fact, if } $\left[ \ l\ \right] _{\ r-1}\notin \left\{
0,1\right\} $\emph{\ (which is possible only for }$r\geq 4$\emph{), the only
non-basic }$\left( r-1\right) $\emph{-simplex of }$\mathcal{T}$\emph{\ is
conv}$\left( B_{\text{\emph{last}}}\right) $.
\emph{If $\text{gcd}(l,r-1)=1$ then the toric variety }$X\left( N_G,\widehat{\Delta }%
_G\left( \mathcal{T}\right) \right) $\emph{\ has a Gorenstein terminal,
isolated, cyclic quotient singularity of type} 
\[
\frac 1{\left[ \ l\ \right] _{\ r-1}}\ \left( \stackunder{\left( r-1\right) 
\text{\emph{-times}}}{\underbrace{1,1,\ldots ,1,1}},\ \left[ -\left(
r-1\right) \right] _{\ \left[ \ l\ \right] _{\ r-1}}\right) 
\]
\emph{lying on the affine piece }$U_{\text{\emph{pos}}\left( B_{\text{\emph{%
last}}}\right) }$\emph{. Otherwise, the singular locus is not a
singleton; more precisely, it is }$1$\emph{-dimensional, i.e., the
corresponding quotient singularity has splitting codimension }$r-1$\emph{,
and can be viewed as a }$1$\emph{-parameter} \emph{``Schar'' (}$\approx $ 
\emph{crowd) of terminal singularities of type } 
\[
\frac 1{\left( \frac{r-1}{\emph{gcd}\left( l,r-1\right) }\right) }\ \left( 
\stackunder{\left( r-1\right) \text{\emph{-times}}}{\underbrace{1,1,\ldots
,1,1}}\right) \smallskip 
\]
\emph{along (the sections of the normal sheaf of) the proper transform of
the ``}$z_r$\emph{-axis'' STR(}$U_{\sigma _0}\left( \text{\emph{pos}}\left(
\left\{ e_1,\ldots ,e_{r-1}\right\} \right) \right) $, $f_{\mathcal{T}}$%
\emph{)} \emph{lying on} $X\left( N_G,\widehat{\Delta }_G\left( \mathcal{T}%
\right) \right) .\smallskip $\newline
\emph{The triangulated junior terahedron of the simplest example} $%
1/5\,\left( 1,1,1,2\right) $ \emph{is drawn in fig.~\textbf{10}. The ``low''
tetrahedron induces the classical involutional terminal singularity of type }%
$1/2\,\left( 1,1,1,1\right) $.
\end{remark}
\begin{figure*}[h]
\begin{center}
\vspace{-1.3cm}
\input{fig10.pstex_t}
\makespace
 {\textbf{Fig. 10}}
\end{center}
\end{figure*}

\begin{theorem}[Exceptional prime divisors and intersection numbers]
\label{DESEXC} \ \smallskip \newline
Let $\left( \Bbb{C}^r/G,\left[ \mathbf{0}\right] \right) =\left( X\left(
N_G,\Delta _G\right) ,\text{\emph{orb}}\left( \sigma _0\right) \right) $
denote the Gorenstein cyclic quotient singularity of type \smallskip \emph{(%
\ref{montype}). }If $l$ satisfies condition \emph{(\ref{rescon}), }then 
\emph{:\smallskip \ }\newline
\emph{(i)} The exceptional locus of $f:X\left( N_G,\widehat{\Delta }_G\left( 
\mathcal{T}\right) \right) \rightarrow X\left( N_G,\Delta _G\right) $
consists of $\QOVERD\lfloor \rfloor {l}{r-1}$ prime divisors\smallskip 
\[
D_j:=D_{n^{\left( j\right) }}=V\left( \tau ^{\left( j\right) }\right)
=X\left( N_G\left( \tau ^{\left( j\right) }\right) ;\,\text{\emph{Star}}%
\left( \left( \tau ^{\left( j\right) };\widehat{\Delta }_G\left( \mathcal{T}%
\right) \right) \right) \right) ,\ 1\leq j\leq \QTOVERD\lfloor \rfloor
{l}{r-1},\smallskip \ 
\]
on $X\left( N_G,\widehat{\Delta }_G\left( \mathcal{T}\right) \right) $, with 
$\ \tau ^{\left( j\right) }:=\Bbb{R}_{\geq 0}\,n^{\left( j\right) }$, having
the following structure \emph{:\smallskip }\newline
\[
\fbox{$
\begin{array}{ccc}
&  &  \\ 
& D_j\cong \Bbb{P}\left( \mathcal{O}_{\Bbb{P}_{\Bbb{C}}^{r-2}}\oplus 
\mathcal{O}_{\Bbb{P}_{\Bbb{C}}^{r-2}}\left( l-\left( r-1\right) j\right)
\right) \text{ \ \ \ \smallskip \emph{(as}}\emph{\ }\Bbb{P}_{\Bbb{C}}^1\text{%
\emph{-bundles}}\emph{\ }\text{\emph{over}}\emph{\ }\Bbb{P}_{\Bbb{C}}^{r-2}%
\emph{)} &  \\ 
&  & 
\end{array}
$}\emph{\ \smallskip \smallskip \ }
\]
for all\emph{\ }$j,\ \ j\in \left\{ 1,2,\ldots ,\QOVERD\lfloor \rfloor
{l}{r-1}-1\right\} $, and\smallskip 
\[
\smallskip \fbox{$
\begin{array}{ccc}
&  &  \\ 
& D_{\QTOVERD\lfloor \rfloor {l}{r-1}}\cong \left\{ 
\begin{array}{lll}
\Bbb{P}_{\Bbb{C}}^{r-1} & , & \text{\emph{if}\ \ \ }l\equiv 1\ \emph{mod}%
\left( r-1\right)  \\ 
&  &  \\ 
\Bbb{P}_{\Bbb{C}}^{r-2}\times \Bbb{C} & , & \text{\emph{if} \ \ }l\equiv 0\ 
\emph{mod}\left( r-1\right) 
\end{array}
\right.  &  \\ 
&  & 
\end{array}
$}\smallskip 
\]
\emph{(ii)} For the \emph{(}highest\emph{) }intersection numbers of two
consecutive exceptional prime divisors we get \emph{:\smallskip } 
\begin{equation}
\fbox{$
\begin{array}{l}
\  \\ 
\left( D_j^{r-1}\cdot D_{j+1}\right) =\left( l-\left( r-1\right) \left(
j+1\right) \right) ^{r-2}\smallskip  \\ 
\  \\ 
\ \left( D_j\cdot D_{j+1}^{r-1}\right) =\left( \left( r-1\right) j-l\right)
^{r-2} \\ 
\ 
\end{array}
$}\smallskip   \label{INTN1}
\end{equation}
and for the self-intersection numbers \emph{:\smallskip } 
\begin{equation}
\fbox{$D_j^r=\dsum\limits_{i=0}^{r-2}\tbinom{r-1}i\ \left( -2\right)
^{r-i-1}\left( \left( l-\left( r-1\right) j\right) -r\right) ^i\left(
l-\left( r-1\right) j\right) ^{r-i-2}$}\ \smallskip   \label{INTN2}
\end{equation}
for all\emph{\ }$j,\ \ j\in \left\{ 1,2,\ldots ,\QOVERD\lfloor \rfloor
{l}{r-1}-1\right\} $, and\smallskip 
\begin{equation}
\fbox{$\left( D_{\QTOVERD\lfloor \rfloor {l}{r-1}}^{}\right) ^r=\left\{ 
\begin{array}{lll}
\left( -r\right) ^{\ r-1} & , & \text{\emph{if}\ \ \ }l\equiv 1\ \emph{mod}%
\left( r-1\right)  \\ 
&  &  \\ 
\left( -\left( r-1\right) \right) ^{r-2} & , & \text{\emph{if} \ \ }l\equiv
0\ \emph{mod}\left( r-1\right) 
\end{array}
\right. $}  \label{INTN3}
\end{equation}
\end{theorem}

\noindent \textit{Proof. }(i) We distinguish three cases depending on the
range of $j$ and the divisibility condition on $l$.\smallskip \newline
\textit{First case.} Let $j\in \left\{ 1,2,\ldots ,\QOVERD\lfloor \rfloor
{l}{r-1}-1\right\} $. Obviously, there are exactly two primitive collections
in Gen$\left( \text{Star}\left( \tau ^{\left( j\right) };\widehat{\Delta }%
_G\left( \mathcal{T}\right) \right) \right) $, namely 
\[
\left\{ e_1+\left( N_G\right) _{\tau ^{\left( j\right) }},\ldots
,e_{r-1}+\left( N_G\right) _{\tau ^{\left( j\right) }}\right\} \ \text{and\ }%
\left\{ n^{\left( j-1\right) }+\left( N_G\right) _{\tau ^{\left( j\right)
}},n^{\left( j+1\right) }+\left( N_G\right) _{\tau ^{\left( j\right)
}}\right\} 
\]
having no common elements. Furthermore, 
\[
n^{\left( j-1\right) }+n^{\left( j+1\right) }=2\,n^{\left( j\right)
}\Longrightarrow \left( n^{\left( j-1\right) }+\left( N_G\right) _{\tau
^{\left( j\right) }}\right) +\left( n^{\left( j+1\right) }+\left( N_G\right)
_{\tau ^{\left( j\right) }}\right) =\mathbf{0}_{N_G\left( \tau ^{\left(
j\right) }\right) }.
\]
Hence, each $D_j$ has to be the total space of a $\Bbb{P}_{\Bbb{C}}^1$%
-bundle over an $\left( r-2\right) $-dimensional smooth, compact toric
variety, and since Star$\left( \tau ^{\left( j\right) };\widehat{\Delta }%
_G\left( \mathcal{T}\right) \right) $ is a splitting fan, it will be, in
addition, the total space of the projectivization of a decomposable bundle
(by prop. \ref{chbudl} and thm. \ref{splfan}). On the other hand, 
\[
\#\left( \text{Gen}\left( \text{Star}\left( \tau ^{\left( j\right) };%
\widehat{\Delta }_G\left( \mathcal{T}\right) \right) \right) \right) =r+1,
\]
which means that $D_j$ has Picard number $2$, and has therefore to be
isomorphic to the total space of the projectivization of a decomposable
bundle of the form $\mathcal{O}_{\Bbb{P}_{\Bbb{C}}^{r-2}}\oplus \mathcal{O}_{%
\Bbb{P}_{\Bbb{C}}^{r-2}}\left( \lambda \right) $ over $\Bbb{P}_{\Bbb{C}%
}^{r-2}$ (by Kleinschmidt's classification theorem \ref{KLCL}). For this
reason, it suffices to determine this single twisting number $\lambda $ by
using lemma \ref{YPSILON}. It is easy to verify, that we have the $\Bbb{Z}$%
-linear dependence-relations 
\[
e_1+e_2+\cdots +e_{r-2}+e_{r-1}=\left( \left( r-1\right) \left( j+1\right)
-l\right) \ n^{\left( j\right) }+\left( l-\left( r-1\right) j\right) \
n^{\left( j+1\right) }
\]
i.e., 
\[
\sum_{\iota =1}^{r-1}\,\left( e_\iota +\left( N_G\right) _{\tau ^{\left(
j\right) }}\right) -\left( \left( l-\left( r-1\right) j\right) \ n^{\left(
j+1\right) }+\left( N_G\right) _{\tau ^{\left( j\right) }}\right) =\mathbf{0}%
_{N_G\left( \tau ^{\left( j\right) }\right) }\ .
\]
Consequently, 
\[
D_j\cong Y\left( r-1;l-\left( r-1\right) j\right) \ .\medskip 
\]
\textit{Second case. } Let $l\equiv 1\ $mod$\left( r-1\right) $ and $j=$ $%
\QOVERD\lfloor \rfloor {l}{r-1}=\frac{l-1}{r-1}\,$. Since\smallskip 
\[
\left( e_1+\left( N_G\right) _{\tau ^{\left( j\right) }}\right) +\cdots
+\left( e_{r-1}+\left( N_G\right) _{\tau ^{\left( j\right) }}\right) +\left(
n^{\left( j-1\right) }+\left( N_G\right) _{\tau ^{\left( j\right) }}\right) =%
\mathbf{0}_{N_G\left( \tau ^{\left( j\right) }\right) }\ ,
\]
we have obviously $D_j\cong \Bbb{P}_{\Bbb{C}}^{r-1}$.\medskip \newline
\textit{Third case. } Let $l\equiv 0\ $mod$\left( r-1\right) $ and $j=$ $%
\QOVERD\lfloor \rfloor {l}{r-1}=\frac l{r-1}\,$. Since $n_r^{\left( j\right)
}$ equals $0$, $n^{\left( j\right) }$ lies on the facet conv$\left( \left\{
e_1,\ldots ,e_{r-1}\right\} \right) $ of $\frak{s}_G$. This means that the
star of $\tau ^{\left( j\right) }$ in $\widehat{\Delta }_G\left( \mathcal{T}%
\right) $ can be written as a direct product of the form\smallskip 
\[
\text{Star}\left( \tau ^{\left( j\right) };\widehat{\Delta }_G\left( 
\mathcal{T}\right) \right) =\text{Star}\left( \tau ^{\left( j\right) };%
\widehat{\Delta }_G\left( \mathcal{T}\right) \left( \text{pos}\left(
e_1,..,e_{r-1}\right) \right) \right) \times \left( \text{a half-line in }%
\Bbb{R}\right) \smallskip 
\]
and $D_j$ splits into 
\[
X\left( N_G\left( \text{pos}\left( e_1,..,e_{r-1}\right) \right) ;\text{Star}%
\left( \left( \tau ^{\left( j\right) };\widehat{\Delta }_G\left( \mathcal{T}%
\right) \left( \text{pos}\left( e_1,..,e_{r-1}\right) \right) \right)
\right) \right) \times \Bbb{C\smallskip }
\]
where this first factor is isomorphic to $\Bbb{P}_{\Bbb{C}}^{r-2}$ because 
\[
e_1+\cdots +e_{r-1}=\left( r-1\right) \cdot n^{\left( j\right) }\,\,.
\]
\newline
(ii) Since 
\[
\left( D_j^{r-1}\cdot D_{j+1}\right) =D_j^{r-1}\,\left| _{D_{j+1}}\right. 
\]
equals 
\[
\left( V\left( \tau ^{\left( j\right) }+\left( N_G\right) _{\tau ^{\left(
j+1\right) }}\right) \right) ^{r-1}\,\left| _{X\left( N_G\left( \tau
^{\left( j+1\right) }\right) ;\,\text{Star}\left( \left( \tau ^{\left(
j+1\right) };\widehat{\Delta }_G\left( \mathcal{T}\right) \right) \right)
\right) }\right. ,
\]
we obtain 
\[
\left( D_j^{r-1}\cdot D_{j+1}\right) =\left( l-\left( r-1\right) \left(
j+1\right) \right) ^{r-2}
\]
by (\ref{Formula2}). Similarly for $\left( D_j\cdot D_{j+1}^{r-1}\right) $;
regarding $D_j$ as the total space of a $\Bbb{P}_{\Bbb{C}}^1$-bundle$\emph{\ 
}$over$\emph{\ }\Bbb{P}_{\Bbb{C}}^{r-2}$, and making this time use of the
``opposite'' piece of the affine covering of its typical fiber $\cong \Bbb{P}%
_{\Bbb{C}}^1$, we get the second formula of (\ref{INTN1}) by restricting $%
D_{j+1}^{r-1}$ on 
\[
D_j=X\left( N_G\left( \tau ^{\left( j\right) }\right) ;\,\text{Star}\left(
\left( \tau ^{\left( j\right) };\widehat{\Delta }_G\left( \mathcal{T}\right)
\right) \right) \right) 
\]
and by (\ref{Formula2}) after sign-change (i.e., after having identified $%
V\left( \tau ^{\left( j+1\right) }+\left( N_G\right) _{\tau ^{\left(
j\right) }}\right) $ with the divisor on $D_j$ whose associated line bundle
is $\mathcal{O}_{D_j}\left( -1\right) $).\smallskip \newline
On the other hand, $\omega _{X\left( N_G,\widehat{\Delta }_G\left( \mathcal{T%
}\right) \right) }\cong \mathcal{O}_{X\left( N_G,\widehat{\Delta }_G\left( 
\mathcal{T}\right) \right) }$, and the adjunction formula gives\smallskip 
\[
\omega _{D_j}\cong \omega _{X\left( N_G,\widehat{\Delta }_G\left( \mathcal{T}%
\right) \right) }\otimes \mathcal{O}_{X\left( N_G,\widehat{\Delta }_G\left( 
\mathcal{T}\right) \right) }\left( D_j\right) \cong \mathcal{O}_{D_j}\left(
D_j\right) \cong \mathcal{N}_{D_j/X\left( N_G,\widehat{\Delta }_G\left( 
\mathcal{T}\right) \right) }\ ,\smallskip 
\]
where $\mathcal{N}_{D_j/X\left( N_G,\widehat{\Delta }_G\left( \mathcal{T}%
\right) \right) }$ denotes the normal sheaf of $D_j$ in $X\left( N_G,%
\widehat{\Delta }_G\left( \mathcal{T}\right) \right) $. Hence, evaluating
the highest power of the first Chern class of this sheaf at the fundamental
cycle $\left[ D_j\right] $ of $D_j$, we obtain for the self-intersection
number : 
\begin{equation}
D_j^r=c_1^{r-1}\left( \mathcal{N}_{D_j/X\left( N_G,\widehat{\Delta }_G\left( 
\mathcal{T}\right) \right) }\right) \left( \left[ D_j\right] \right)
=c_1^{r-1}\left( \omega _{D_j}\right) \left( \left[ D_j\right] \right)
=K_{D_j}^{r-1}  \label{ADJF}
\end{equation}
Formula (\ref{INTN2}) for $D_j^r$ follows from (\ref{ADJF}) and (\ref
{Formula1}); formula (\ref{INTN3}) is obvious. $_{\Box }$

\section{Breaking down the desingularizing morphism\label{FACTOR}}

\noindent The unique crepant resolution-morphism of the singularities
discussed in the previous section can be factorized into (normalized)
blow-ups in several ways. We give here two canonical decompositions of $f=f_{%
\mathcal{T}}$ of this kind.\medskip \newline
\textsf{(a) }Maintaining the notation introduced in the proof of thm. \ref
{Main} for the construction of the unique, basic, coherent triangulation $%
\mathcal{T}$, besides $\mathbf{s}\left( j;\xi _1,\ldots ,\xi _{r-2}\right) $%
's (and conv$\left( B_{\text{last}}\right) $, for $l\equiv 1\ $mod$\left(
r-1\right) $), we define\smallskip 
\[
B\left( j,j^{\prime };\xi _1,\xi _2,\ldots ,\xi _{r-2}\right) :=\left\{
n^{\left( j\right) },n^{\left( j^{\prime }\right) },\ e_{\xi _1},e_{\xi
_2},\ldots ,e_{\xi _{r-2}}\right\} , 
\]
as well as 
\[
\mathbf{s}\left( j,j^{\prime };\xi _1,\xi _2,\ldots ,\xi _{r-2}\right) :=%
\text{ conv}\left( B\left( j,j^{\prime };\xi _1,\xi _2,\ldots ,\xi
_{r-2}\right) \right) ,\smallskip \smallskip \ 
\]
and 
\[
\sigma \left( j,j^{\prime };\xi _1,\xi _2,\ldots ,\xi _{r-2}\right) :=\sigma
_{\mathbf{s}\left( j,j^{\prime };\xi _1,\xi _2,\ldots ,\xi _{r-2}\right) }=%
\text{ pos}\left( B\left( j,j^{\prime };\xi _1,\xi _2,\ldots ,\xi
_{r-2}\right) \right) ,\smallskip 
\]
for all indices $j$, $j^{\prime }$ with 
\[
\ \ \ \ 0\leq j\leq j^{\prime }\leq \QOVERD\lfloor \rfloor
{l}{r-1}\,,\smallskip 
\]
and for all \ $r-1$ \ possible \ $\left( r-2\right) $-tuples $\ \left( \xi
_1,\xi _2,\ldots ,\xi _{r-2}\right) \in \Xi _r.$ Obviously, 
\[
\mathbf{s}\left( j,j^{\prime };\xi _1,\xi _2,\ldots ,\xi _{r-2}\right) =%
\mathbf{s}\left( j^{\prime };\xi _1,\xi _2,\ldots ,\xi _{r-2}\right) 
\]
whenever $j^{\prime }=j+1$. Let now $\kappa $ denote the positive
integer\smallskip 
\[
\kappa :=\left\{ 
\begin{array}{lll}
\left\lfloor \frac 12\left( \,\frac{l-1}{r-1}+1\right) \right\rfloor & ,\ 
\text{if} & l\equiv 1\ \text{mod}\left( r-1\right) \\ 
&  &  \\ 
\left\lfloor \frac 12\,\left( \frac l{r-1}\right) \right\rfloor +1 & ,\ 
\text{if} & l\equiv 0\ \text{mod}\left( r-1\right)
\end{array}
\right. \smallskip 
\]
$\bullet $ For $l\equiv 1\ $mod$\left( r-1\right) $ we introduce the
following simplicial subdivisions of the junior simplex $\frak{s}_G$
:\smallskip 
\[
\frak{T}_1:=\left\{ 
\begin{array}{c}
\mathbf{s}\left( 1;\xi _1,\xi _2,\ldots ,\xi _{r-2}\right) ,\text{ }\mathbf{s%
}\left( 1,\frac{l-1}{r-1};\xi _1,\xi _2,\ldots ,\xi _{r-2}\right) ,\text{%
\thinspace } \\ 
\text{conv}\left( B_{\text{last}}\right) =\mathbf{s}\left( \frac{l-1}{r-1},%
\frac{l-1}{r-1};\xi _1,\xi _2,\ldots ,\xi _{r-2}\right) ,\smallskip \\ 
\text{for all }\left( \xi _1,\xi _2,\ldots ,\xi _{r-2}\right) \in \Xi _r%
\text{, together with all their faces}
\end{array}
\right\} 
\]
and 
\[
\begin{array}{l}
\frak{T}_{i+1}:=\left[ \frak{T}_i\smallsetminus \left\{ 
\begin{array}{c}
\mathbf{s}\left( i,\frac{l-1}{r-1}-i+1;\xi _1,\xi _2,\ldots ,\xi
_{r-2}\right) ,\smallskip \\ 
\text{for all }\left( \xi _1,\xi _2,\ldots ,\xi _{r-2}\right) \in \Xi _r%
\text{, \smallskip } \\ 
\text{together with its faces}
\end{array}
\right\} \right] \cup \\ 
\  \\ 
\cup \ \left\{ 
\begin{array}{c}
\mathbf{s}\left( i,i+1;\xi _1,\ldots ,\xi _{r-2}\right) ,\ \mathbf{s}\left(
i+1,\frac{l-1}{r-1}-i;\xi _1,\ldots ,\xi _{r-2}\right) ,\, \\ 
\mathbf{s}\left( \frac{l-1}{r-1}-i,\frac{l-1}{r-1}-i+1;\xi _1,\ldots ,\xi
_{r-2}\right) ,\,\smallskip \\ 
\text{ for all \ }\left( \xi _1,\ldots ,\xi _{r-2}\right) \in \Xi _r\text{,\
together with all their faces}
\end{array}
\right\}
\end{array}
\smallskip 
\]
for $\kappa \geq 2$ and all $i$, $1\leq i\leq \kappa -1.\smallskip
\smallskip $\newline
$\bullet $ Analogously, for $l\equiv 0\ $mod$\left( r-1\right) $, we
define\smallskip 
\[
\frak{T}_1:=\left\{ 
\begin{array}{c}
\mathbf{s}\left( 1;\xi _1,\xi _2,\ldots ,\xi _{r-2}\right) ,\ \mathbf{s}%
\left( 1,\frac l{r-1}-1;\xi _1,\xi _2,\ldots ,\xi _{r-2}\right) ,\smallskip
\, \\ 
\smallskip \mathbf{s}\left( \frac l{r-1}-1,\frac l{r-1}-1;\xi _1,\xi
_2,\ldots ,\xi _{r-2}\right) , \\ 
\text{ for all \ }\left( \xi _1,\xi _2,\ldots ,\xi _{r-2}\right) \in \Xi _r%
\text{,\ together with all their faces}
\end{array}
\right\} 
\]
and 
\[
\begin{array}{l}
\frak{T}_{i+1}:=\left[ \frak{T}_i\smallsetminus \left\{ 
\begin{array}{c}
\mathbf{s}\left( i,\frac l{r-1}-i;\xi _1,\xi _2,\ldots ,\xi _{r-2}\right)
,\smallskip \\ 
\text{for all }\left( \xi _1,\xi _2,\ldots ,\xi _{r-2}\right) \in \Xi _r%
\text{, \smallskip } \\ 
\text{together with its faces}
\end{array}
\right\} \right] \cup \\ 
\  \\ 
\cup \ \left\{ 
\begin{array}{c}
\mathbf{s}\left( i,i+1;\xi _1,\ldots ,\xi _{r-2}\right) ,\ \mathbf{s}\left(
i+1,\frac l{r-1}-i-1;\xi _1,\ldots ,\xi _{r-2}\right) ,\smallskip \\ 
\mathbf{s}\left( \frac l{r-1}-i-1,\frac l{r-1}-i;\xi _1,\ldots ,\xi
_{r-2}\right) ,\smallskip \\ 
\text{ for all \ }\left( \xi _1,\ldots ,\xi _{r-2}\right) \in \Xi _r\text{,\
together with all their faces}
\end{array}
\right\}
\end{array}
\smallskip 
\]
for $\kappa \geq 3$ and all $i$, $1\leq i\leq \kappa -2,$ and finally 
\[
\frak{T}_\kappa :=\frak{T}_{\kappa -1}\cup \left\{ 
\begin{array}{c}
\mathbf{s}\left( \frac l{r-1};\xi _1,\xi _2,\ldots ,\xi _{r-2}\right)
,\smallskip \\ 
\text{for all }\left( \xi _1,\xi _2,\ldots ,\xi _{r-2}\right) \in \Xi _r%
\text{,\smallskip } \\ 
\text{together with its faces}
\end{array}
\right\} \ \ . 
\]
Next lemma is obvious by construction.

\begin{lemma}
\label{COVLEM}All the above defined simplicial subdivisions $\frak{T}_1$,$%
\ldots ,\frak{T}_\kappa $ are triangulations and cover the entire $\frak{s}_G
$.
\end{lemma}

\begin{proposition}[First, speedy factorization]
\label{1FACT} \ \smallskip \newline
Let $\left( \Bbb{C}^r/G,\left[ \mathbf{0}\right] \right) =\left( X\left(
N_G,\Delta _G\right) ,\text{\emph{orb}}\left( \sigma _0\right) \right) $
denote the Gorenstein cyclic quotient singularity of type \smallskip \emph{(%
\ref{montype}) }with $l$ satisfying condition \emph{(\ref{rescon}). }Then
the birational resolution-morphism $f=f_{\mathcal{T}}$ is the composite of
the $\kappa $ toric morphisms 
\[
X_0:=X\left( N_G,\Delta _G\right) \stackrel{g_1}{\longleftarrow }X_1%
\stackrel{g_2}{\longleftarrow }X_2\longleftarrow \cdots \stackrel{g_{\kappa
-1}}{\longleftarrow }X_{\kappa -1}\stackrel{g_\kappa }{\longleftarrow }%
X_\kappa :=X\left( N_G,\widehat{\Delta }_G\left( \mathcal{T}\right) \right) 
\]
with 
\[
X_i:=X\left( N_G,\widehat{\Delta }_G\left( \frak{T}_i\right) \right) ,\ \ \
\forall i,\ \ 1\leq i\leq \kappa -1,\text{ \ \emph{and\ }}\mathcal{T}=\frak{T%
}_\kappa \ .
\]
In particular, in algebraic-geometric terms, one has 
\[
X_{i+1}\cong \text{\emph{Norm}\thinspace }\left[ \mathbf{Bl}_{Z_i}^{\text{%
\emph{red}}}\left( X_i\right) \right] ,\ \ \ \forall i,\ \ 0\leq i\leq
\kappa -1\ ,
\]
with centers\smallskip 
\[
Z_i=\left\{ 
\begin{array}{lll}
\text{\emph{orb}}\left( \sigma _0\right)  & ,\ \ \text{\emph{if}} & i=0 \\ 
&  &  \\ 
V\left( \text{\emph{pos}}\left( \left\{ n^{\left( i\right) },n^{\left( \frac{%
l-1}{r-1}-i+1\right) }\right\} \right) \right)  & ,\ \ \text{\emph{if}} & 
\left\{ 
\begin{array}{l}
\kappa \geq 2,\ 1\leq i\leq \kappa -1,\ \  \\ 
\text{\emph{and} \ }l\equiv 1\ \text{\emph{mod}}\left( r-1\right) 
\end{array}
\right.  \\ 
&  &  \\ 
V\left( \text{\emph{pos}}\left( \left\{ n^{\left( i\right) },n^{\left( \frac
l{r-1}-i\right) }\right\} \right) \right)  & ,\ \ \text{\emph{if}} & \left\{ 
\begin{array}{l}
\kappa \geq 3,\ 1\leq i\leq \kappa -2,\ \  \\ 
\text{\emph{and} \ }l\equiv 0\ \text{\emph{mod}}\left( r-1\right) 
\end{array}
\right.  \\ 
&  &  \\ 
V\left( \text{\emph{pos}}\left( \left\{ e_1,e_2,\ldots ,e_{r-1}\right\}
\right) \right)  & ,\ \ \text{\emph{if}} & i=\kappa -1\ \ \text{\emph{\&} \ }%
l\equiv 0\ \text{\emph{mod}}\left( r-1\right) 
\end{array}
\right. 
\]
\end{proposition}

\noindent \textit{Sketch of proof.} Since $N_G$ is a ``skew'' lattice, it is
not so convenient to work directly with it. For this reason we consider the
linear transformation 
\[
\Phi :\Bbb{R}^r\longrightarrow \Bbb{R}^r,\ \ \mathbf{y\longmapsto \,}\Phi
\left( \mathbf{y}\right) :=\mathcal{M\cdot }\mathbf{y,}
\]
with\smallskip 
\[
\mathcal{M}:=\left( 
\begin{array}{ccccc}
l & 0 & 0 & \cdots  & 0 \\ 
\begin{array}{c}
-1 \\ 
-1 \\ 
\vdots  \\ 
-1
\end{array}
& 
\begin{array}{c}
1 \\ 
0 \\ 
\vdots  \\ 
0
\end{array}
& 
\begin{array}{c}
0 \\ 
1 \\ 
\vdots  \\ 
0
\end{array}
& \ddots  & 
\begin{array}{c}
0 \\ 
0 \\ 
\vdots  \\ 
0
\end{array}
\\ 
-\left( l-\left( r-1\right) \right)  & 0 & 0 & \cdots  & 1
\end{array}
\right) \in \text{GL}\left( r,\Bbb{Q}\right) \smallskip .
\]
The image of the lattice $N_G$ via $\Phi $ is the standard lattice $\Bbb{Z}%
^r=\sum_{i=1}^r\Bbb{Z\,}e_i$. In particular, 
\[
N_G=\Bbb{Z}^r+\Bbb{Z\ }\frac 1l\,\left( 1,1,\ldots ,1,l-\left( r-1\right)
\right) ^{\intercal }=\Bbb{Z\ }\frac 1l\,\left( 1,1,\ldots ,1,l-\left(
r-1\right) \right) ^{\intercal }+\sum_{i=2}^r\Bbb{Z\,}e_i,
\]
and $\Phi $ maps this $\Bbb{Z}$-basis of $N_G$ onto 
\[
\Phi \left( \frac 1l\,\left( 1,1,\ldots ,1,l-\left( r-1\right) \right)
^{\intercal }\right) =e_1,\ \Phi \left( e_2\right) =e_2,\ \ldots \ ,\ \Phi
\left( e_r\right) =e_r,
\]
the positive orthant $\sigma _0$ onto the cone 
\[
\overline{\sigma _0}:=\Phi \left( \sigma _0\right) =\text{ pos}\left(
\left\{ \frak{v},e_2,e_3,\ldots ,e_r\right\} \right) ,\text{ \ }
\]
with 
\[
\frak{v}:=\Phi \left( e_1\right) =\left( l,\stackunder{\left( r-2\right) 
\text{-times}}{\underbrace{-1,-1,\ldots ,-1}},-\left( l-\left( r-1\right)
\right) \right) ^{\intercal },
\]
and the Hilbert basis $\mathbf{Hlb}_N\left( \sigma _0\right) $ of $\sigma _0$
w.r.t. $N$ (cf. (\ref{HILBCON}), (\ref{JSEL})) onto\smallskip 
\begin{eqnarray*}
\Phi \left( \mathbf{Hlb}_N\left( \sigma _0\right) \right)  &=&\mathbf{Hlb}_{%
\Bbb{Z}^r}\left( \overline{\sigma _0}\smallskip \right) =\text{conv}\left(
\left\{ \frak{v},e_2,e_3,\ldots ,e_r\right\} \right) \cap \Bbb{Z}%
^r=\smallskip  \\
&=&\left\{ \frak{v},e_2,e_3,\ldots ,e_r\right\} \cup \left\{ \frak{y}%
^{\left( 1\right) },\frak{y}^{\left( 2\right) },\ldots ,\frak{y}^{\left( \nu
\right) }\right\} ,\smallskip 
\end{eqnarray*}
where $\nu :=\QOVERD\lfloor \rfloor {l}{r-1}$, and 
\[
\frak{y}^{\left( j\right) }:=\Phi \left( n^{\left( j\right) }\right) =\left(
j,\stackunder{\left( r-2\right) \text{-times}}{\underbrace{0,0,\ldots ,0,0}}%
,-j+1\right) ^{\intercal },\ \ \forall j,\ \ 1\leq j\leq \nu .
\]
(For $j=1$, $\frak{y}^{\left( 1\right) }=e_1$). The dual cone of $\overline{%
\sigma _0}$ equals 
\[
\left( \overline{\sigma _0}\right) ^{\vee }=\text{ pos}\left( \left\{
e_1^{\vee },e_1^{\vee }+l\,e_2^{\vee },\ldots ,e_1^{\vee }+l\,e_{r-1}^{\vee
},\left( l-\left( r-1\right) \right) \,e_1^{\vee }+l\,e_r^{\vee }\right\}
\right) \smallskip 
\]
(with $\left\{ e_1^{\vee },\ldots ,e_r^{\vee }\right\} $ denoting the dual
of $\left\{ e_1,\ldots ,e_r\right\} $). For every $m$ belonging to the
Hilbert basis $\mathbf{Hlb}_{\left( \Bbb{Z}^r\right) ^{\vee }}\left( \left( 
\overline{\sigma _0}\right) ^{\vee }\right) $ (w.r.t. the dual lattice $%
\left( \Bbb{Z}^r\right) ^{\vee }$) define\smallskip 
\[
\tau \left[ m\right] :=\left\{ \mathbf{y}\in \overline{\sigma _0}\ \left| \
\left\langle m\mathbf{,y}\right\rangle \leq \ \left\langle m^{\prime }%
\mathbf{,y}\right\rangle ,\ \forall m^{\prime },\ m^{\prime }\in \mathbf{Hlb}%
_{\left( \Bbb{Z}^r\right) ^{\vee }}\left( \left( \overline{\sigma _0}\right)
^{\vee }\right) \right. \right\} \ .\smallskip 
\]
$\bullet $ Suppose $l\equiv 1\ $mod$\left( r-1\right) $. At first we shall
show that for any fixed $\left( r-2\right) $-tuple $\left( \xi _1,\xi
_2,\ldots ,\xi _{r-2}\right) \in \Xi _r$, with $\left\{ q\right\} =\left\{
1,2,\ldots ,r-1\right\} \Bbb{r}\left\{ \xi _1,\xi _2,\ldots ,\xi
_{r-2}\right\} ,\smallskip $ we have 
\begin{equation}
\Phi \left( \sigma \left( 1;\xi _1,\ldots ,\xi _{r-2}\right) \right)
=\left\{ 
\begin{array}{lll}
\tau \left[ e_1^{\vee }\right]  & \ \text{if} & q=1 \\ 
&  &  \\ 
\tau \left[ e_1^{\vee }+l\,e_q^{\vee }\right]  & \ \text{if} & q\in \left\{
2,\ldots ,r-1\right\} 
\end{array}
\right.   \label{SXE1}
\end{equation}
$\bullet $ Moreover, for $l>r$, 
\begin{equation}
\Phi \!\left( \!\sigma \left( 1,\tfrac{l-1}{r-1};\xi _1,..,\xi _{r-2}\right)
\!\right) \!=\left\{ 
\begin{array}{ll}
\tau \left[ e_1^{\vee }+e_r^{\vee }\right] \text{ } & \!\!\!\text{if
\thinspace }q=1 \\ 
&  \\ 
\tau \left[ \!e_1^{\vee }+e_r^{\vee }\!+\!\left( r-1\right) \,e_q^{\vee
}\!\right] \text{ } & \!\!\!q=2,..,r\!-\!1\text{\thinspace }
\end{array}
\right.   \label{SX2}
\end{equation}
and finally 
\begin{equation}
\Phi \left( \text{pos}\left( B_{\text{last}}\right) \right) =\tau \left[
\left( l-\left( r-1\right) \right) \,e_1^{\vee }+l\,e_r^{\vee }\right] 
\label{SXE3}
\end{equation}
$\bullet $ \textit{Proof of} (\ref{SXE1}) : Suppose first $q=1.$ Obviously, $%
e_1^{\vee }\in \mathbf{Hlb}_{\left( \Bbb{Z}^r\right) ^{\vee }}\left( \left( 
\overline{\sigma _0}\right) ^{\vee }\right) $ and for all $m\in \mathbf{Hlb}%
_{\left( \Bbb{Z}^r\right) ^{\vee }}\left( \left( \overline{\sigma _0}\right)
^{\vee }\right) $ we get 
\[
\left\langle m,e_1\right\rangle \geq \left\langle e_1^{\vee
},e_1\right\rangle =1\Longrightarrow e_1\in \tau \left[ e_1^{\vee }\right] 
\]
and for every $j\in \left\{ 2,\ldots ,r\right\} $, 
\[
\left\langle m,e_1\right\rangle \geq \left\langle e_1^{\vee
},e_j\right\rangle =0\Longrightarrow e_j\in \tau \left[ e_1^{\vee }\right] \
.
\]
Conversely, let $\mathbf{y}$ denote an arbitrary element of $\tau \left[
e_1^{\vee }\right] $. Write 
\[
\mathbf{y=}\sum_{i=1}^r\,\mu _i\,e_i
\]
as an $\Bbb{R}$-linear combination w.r.t. the basis $\left\{ e_1,\ldots
,e_r\right\} $. Since $e_1^{\vee }$, $e_1^{\vee }+l\,e_j^{\vee }$, $2\leq
j\leq r-1$, and $e_1^{\vee }+e_r^{\vee }$ belong to $\left( \overline{\sigma
_0}\right) ^{\vee }$ we obtain 
\[
\left\langle e_1^{\vee },\mathbf{y}\right\rangle \geq 0,\ \ \ \left\langle
e_1^{\vee }+l\,e_j^{\vee },\mathbf{y}\right\rangle \geq \left\langle
e_1^{\vee },\mathbf{y}\right\rangle ,\ \ \ \left\langle e_1^{\vee
}+e_r^{\vee },\mathbf{y}\right\rangle \geq \left\langle e_1^{\vee },\mathbf{y%
}\right\rangle ,\smallskip 
\]
i.e., $\mu _1,..,\mu _r\in \Bbb{R}_{\geq 0}$ and therefore $\mathbf{y}\in $ $%
\Phi \left( \sigma \left( 1;2,3,\ldots ,r-1\right) \right) =$ pos$\left(
\left\{ e_1,..,e_r\right\} \right) .\medskip $ \newline
$\bullet $ Suppose now $q\in \left\{ 2,\ldots ,r-1\right\} .$ Since 
\[
\left\langle e_1^{\vee }+l\,e_q^{\vee },\frak{v}\right\rangle =\left\langle
e_1^{\vee }+l\,e_q^{\vee },e_p\right\rangle =0\text{,\ for\ }p\in \left\{
2,..,q-2,q-1,q+1,q+2,..,r\right\} ,
\]
and 
\[
\left\langle m,e_1\right\rangle \geq \left\langle e_1^{\vee }+l\,e_q^{\vee
},e_1\right\rangle =1\text{, \ for all \ }m\in \mathbf{Hlb}_{\left( \Bbb{Z}%
^r\right) ^{\vee }}\left( \left( \overline{\sigma _0}\right) ^{\vee }\right)
,
\]
the inclusion ``$\subset $'' is obvious. Conversely, let $\mathbf{y}$ denote
an arbitrary element of $\tau \left[ e_1^{\vee }+l\,e_q^{\vee }\right] $ and
write it as $\Bbb{R}$-linear combination 
\[
\mathbf{y=}\sum_{i\in \left\{ 1,2,\ldots ,q-2,q-1,q+1,q+2,\ldots ,r\right\}
}\,\mu _i\,e_i+\mu _q\,\frak{v}\,.
\]
Since 
\[
\mu _1=\left\langle e_1^{\vee }+l\,e_q^{\vee },\mathbf{y}\right\rangle \geq
0,\ \ \ \ \left\langle e_1^{\vee }+l\,e_q^{\vee },\mathbf{y}\right\rangle
\leq \left\langle e_1^{\vee },\mathbf{y}\right\rangle \Longleftrightarrow
\mu _q\geq 0,
\]
and 
\[
\left\langle e_1^{\vee }+l\,e_q^{\vee },\mathbf{y}\right\rangle \leq
\left\langle e_1^{\vee }+l\,e_p^{\vee },\mathbf{y}\right\rangle
\Longleftrightarrow \mu _p\geq 0,\text{ \ }\smallskip 
\]
for \ $p\in \left\{ 2,\ldots ,q-2,q-1,q+1,q+2,\ldots ,r\right\} $, we have 
\[
\mathbf{y}\in \Phi \left( \sigma \left( 1;2,..,q-2,q-1,q+1,q+2,..,r\right)
\right) =\text{pos}\left( \left\{ \frak{v},e_1,..,e_{q-1},e_{q+1},..,e_r%
\right\} \right) .\smallskip 
\]
\newline
$\bullet $ \textit{Proof of }(\ref{SX2}) : Suppose first $q=1.$ The
inclusion ``$\subset $'' can be easily checked as before. Let $\mathbf{y}$
be an element of $\tau \left[ e_1^{\vee }+e_r^{\vee }\right] $ and write it
as linear combination 
\[
\mathbf{y=}\sum_{i=1}^{r-1}\,\mu _i\,e_i+\mu _r\,\frak{y}^{\left( \nu
\right) }\,.
\]
We have 
\[
\left\langle e_1^{\vee }+e_r^{\vee },\mathbf{y}\right\rangle \leq
\left\langle e_1^{\vee },\mathbf{y}\right\rangle \Longleftrightarrow -\mu
_r\nu +\mu _r\leq 0\Longleftrightarrow \mu _r\geq 0\ \ \ \left( \text{by }%
\nu =\tfrac{l-1}{r-1}>1\right) 
\]
and 
\[
\left\langle e_1^{\vee }+e_r^{\vee },\mathbf{y}\right\rangle \leq
\left\langle e_1^{\vee }+e_r^{\vee }+\left( r-1\right) \,e_j^{\vee },\mathbf{%
y}\right\rangle \Longleftrightarrow \mu _j\geq 0,\ \forall j,\ 2\leq j\leq
r-1,
\]
as well as 
\[
\left\langle e_1^{\vee }+e_r^{\vee },\mathbf{y}\right\rangle \leq
\left\langle \left( l-\left( r-1\right) \right) \,e_1^{\vee }+l\,e_r^{\vee },%
\mathbf{y}\right\rangle 
\]
which is equivalent to 
\[
0\leq \left\langle \left( l-r\right) e_1^{\vee }+\left( l-1\right) e_r^{\vee
},\mathbf{y}\right\rangle =\left( l-r\right) \left( \mu _1+\mu _r\nu \right)
+\left( l-1\right) \left( \mu _r-\mu _r\nu \right) =\left( l-r\right) \mu
_1,\smallskip 
\]
i.e., $\mu _1\geq 0$. Hence, $\mathbf{y}\in \Phi \left( \sigma \left( 1,%
\frac{l-1}{r-1};2,3,..,r-1\right) \right) =$ pos$\left( \left\{ \frak{y}%
^{\left( \nu \right) },e_1,e_2,..,e_{r-1}\right\} \right) .\medskip $\newline
$\bullet $ For $q\in \left\{ 2,\ldots ,r-1\right\} $ the proof of the
inclusion ``$\subset $'' is again easy. To prove ``$\supset $'' it is enough
to consider an element $\mathbf{y}\in \tau \left[ e_1^{\vee }+e_r^{\vee
}+\left( r-1\right) \,e_q^{\vee }\right] $ and write it as linear
combination 
\[
\mathbf{y=}\sum_{i\in \left\{ 1,2,\ldots ,q-2,q-1,q+1,q+2,\ldots
,r-1\right\} }\,\mu _i\,e_i+\mu _q\,\frak{v}+\mu _r\,\frak{y}^{\left( \nu
\right) }\,.
\]
By definition, $\mathbf{y}$ satisfies the three inequalities\smallskip 
\[
\begin{array}{l}
\left\langle e_1^{\vee }+e_r^{\vee }+\left( r-1\right) \,e_q^{\vee },\mathbf{%
y}\right\rangle \leq \left\langle e_1^{\vee }+e_r^{\vee },\mathbf{y}%
\right\rangle , \\ 
\, \\ 
\left\langle e_1^{\vee }+e_r^{\vee }+\left( r-1\right) \,e_q^{\vee },\mathbf{%
y}\right\rangle \leq \left\langle \left( l-\left( r-1\right) \right)
\,e_1^{\vee }+l\,e_r^{\vee },\mathbf{y}\right\rangle ,
\end{array}
\]
and
\[
\left\langle e_1^{\vee }+e_r^{\vee }+\left( r-1\right) \,e_q^{\vee },\mathbf{%
y}\right\rangle \leq \smallskip \left\langle e_1^{\vee }+e_r^{\vee }+\left(
r-1\right) \,e_p^{\vee },\mathbf{y}\right\rangle ,
\]
for all \ $p\in \left\{ 2,\ldots ,q-2,q-1,q+1,q+2,\ldots ,r\right\} $.
Direct evaluation combined with 
\[
\left\langle \left( l-\left( r-1\right) \right) \,e_1^{\vee }+l\,e_r^{\vee },%
\frak{v}\right\rangle =0,\ \left\langle \left( l-\left( r-1\right) \right)
\,e_1^{\vee }+l\,e_r^{\vee },\frak{y}^{\left( \nu \right) }\right\rangle
=l-\nu \,\left( r-1\right) =1,
\]
gives $\mu _i\in \Bbb{R}_{\geq 0}$, for all $i$, $1\leq i\leq r,$
i.e.,  $\mathbf{y}$ belongs to 
\[
\Phi \left( \sigma \left( 1,\tfrac{l-1}{r-1};2,..,q-1,q+1,..,r\right)
\right) =\text{pos}\left( \left\{ \frak{v},\frak{y}^{\left( \nu \right)
},e_1,..,e_{q-1},e_{q+1},..,e_{r-1}\right\} \right) .\smallskip 
\]
$\bullet $ \textit{Proof of }(\ref{SXE3}) : Since 
\[
\left\langle \left( l-\left( r-1\right) \right) \,e_1^{\vee }+l\,e_r^{\vee },%
\frak{v}\right\rangle =\left\langle \left( l-\left( r-1\right) \right)
\,e_1^{\vee }+l\,e_r^{\vee },e_j\right\rangle =0,\ \forall j,\ 2\leq j\leq
r-1,\smallskip 
\]
we have $\frak{v},e_2,\ldots ,e_{r-1}\in \tau \left[ \left( l-\left(
r-1\right) \right) \,e_1^{\vee }+l\,e_r^{\vee }\right] $. On the other hand,
by the definition of $\left( \overline{\sigma _0}\right) ^{\vee }$%
,\smallskip 
\[
\left\langle \left( l-\left( r-1\right) \right) \,e_1^{\vee }+l\,e_r^{\vee },%
\frak{y}^{\left( \nu \right) }\right\rangle =1\leq \left\langle m,\frak{y}%
^{\left( \nu \right) }\right\rangle ,\ \forall m,\ m\in \left( \overline{%
\sigma _0}\right) ^{\vee }\cap \left( \Bbb{Z}^r\Bbb{r}\left\{ \mathbf{0}%
\right\} \right) .
\]
\newpage

Consequently, 
\[
\Phi \left( \text{pos}\left( B_{\text{last}}\right) \right) =\text{ pos}%
\left( \left\{ \frak{v},\frak{y}^{\left( \nu \right) },e_2,e_3,\ldots
,e_{r-1}\right\} \right) \subset \tau \left[ \left( l-\left( r-1\right)
\right) \,e_1^{\vee }+l\,e_r^{\vee }\right] \ .\smallskip 
\]
To show the converse inclusion take again a $\mathbf{y}\in \tau \left[
\left( l-\left( r-1\right) \right) \,e_1^{\vee }+l\,e_r^{\vee }\right] $,
write it as linear combination 
\[
\mathbf{y}=\mu _1\,\frak{v}+\sum_{j=2}^{r-1}\,\mu _j\,e_j+\mu _r\,\frak{y}%
^{\left( \nu \right) },
\]
and use the inequalities 
\begin{eqnarray*}
0&\leq&\left\langle \left( l-\left( r-1\right) \right) \,e_1^{\vee
}+l\,e_r^{\vee },\mathbf{y}\right\rangle =\mu _r, 
\\
0&\leq&\left\langle \left( l-\left( r-1\right) \right) \,e_1^{\vee
}+l\,e_r^{\vee },\mathbf{y}\right\rangle =\mu _r\leq \left\langle e_1^{\vee
}+e_r^{\vee },\mathbf{y}\right\rangle =\left( r-1\right) \,\mu _1+\mu _r,
\end{eqnarray*}
and for all $j,\ 2\leq j\leq r-1$, the inequalities 
\[
\left\langle \left( l-\left( r-1\right) \right) \,e_1^{\vee }+l\,e_r^{\vee },%
\mathbf{y}\right\rangle =\mu _r\leq \left\langle e_1^{\vee }+e_r^{\vee
}+\left( r-1\right) \,e_j^{\vee },\mathbf{y}\right\rangle =\left( r-1\right)
\,\mu _j+\mu _r.\smallskip 
\]
$\bullet $ By (\ref{SXE1}), (\ref{SX2}), (\ref{SXE3}), and lemma \ref{COVLEM}
we obtain 
\[
\frak{T}_1=\bigcup_{m\in \Phi \left( \mathbf{Hlb}_N\left( \sigma _0\right)
\right) }\ \left\{ 
\begin{array}{c}
\Phi ^{-1}\left( \tau \left[ m\right] \right) \text{ together } \\ 
\text{with all their faces}
\end{array}
\right\} \ .
\]
This means that $\widehat{\Delta }_G\left( \frak{T}_1\right) =\left( \Delta
_G\right) _{\mathbf{bl}}\left[ \text{orb}\left( \sigma _0\right) \right] $
by proposition \ref{USNBL}, and therefore $g_1$ is indeed the proper
birational morphism corresponding to the normalized, usual blow-up of $%
X\left( N_G,\Delta _G\right) $ at the closed point orb$\left( \sigma
_0\right) $.\smallskip \newline
$\bullet $ If $\kappa \geq 2$, in the second step we blow up
(simultaneously) the $\left( r-2\right) $-dimensional common singular locus
of all affine charts $U_{\sigma \left( 1,\frac{l-1}{r-1};\xi _1,\xi
_2,\ldots ,\xi _{r-2}\right) }$. Note that locally a neighbourhood of such a
singular point within $X\left( N_G,\widehat{\Delta }_G\left( \frak{T}%
_1\right) \right) $ can be viewed like a 
\[
\left( 2\text{-dimensional }A_{\frac{l-1}{r-1}-2}\text{-singularity}\right)
\times \Bbb{C}^{r-2}\ .
\]
To prove that the above defined triangulation $\frak{T}_2$ induces the
normalization of this blow-up, one applies theorem \ref{TORNBL} and
techniques similar to those used for $\frak{T}_1$. The details are left as
an exercise to the reader. Repeating the described procedure altogether $%
\kappa -1$ times we arrive at the entire basic triangulation $\mathcal{T}=%
\frak{T}_\kappa $.\smallskip \newline
$\bullet $ The proof in the case in which $l\equiv 0\ $mod$\left( r-1\right) 
$ can be done analogously and will be omitted. The only difference is that
in the last step one blows up (once) the remaining $1$-dimensional singular
locus $V\left( \text{pos}\left( \left\{ e_1,e_2,\ldots ,e_{r-1}\right\}
\right) \right) $ inherited from the single non-basic facet of $\frak{s}_G$. 
$_{\Box }\medskip $\newline
Figures \textbf{11} (a) and (b) show this speedy factorization of $f$ for
the singularities of type $1/10\,\left( 1,1,8\right) $ and $1/11\,\left(
1,1,9\right) $, respectively. \newline 

\begin{figure*}
\begin{center}
\input{fig11.pstex_t}
\makespace
 {Fig. \textbf{11 }}
\end{center}
\end{figure*}

\noindent
\textsf{(b) }A second canonical factorization of $f=f_{\mathcal{T}}$ is
constructed by means of the following $\nu =$ $\QOVERD\lfloor \rfloor
{l}{r-1}$ triangulations of the junior simplex:\smallskip 
\[
\frak{T}_1^{\,\prime }:=\left\{ 
\begin{array}{c}
\mathbf{s}\left( 1;\xi _1,\xi _2,\ldots ,\xi _{r-2}\right) ,\text{ conv}%
\left( \left\{ n^{\left( 1\right) },e_{\xi _1},\ldots ,e_{\xi
_{r-2}}\right\} \right) ,\smallskip \smallskip \\ 
\text{for all }\left( \xi _1,\xi _2,\ldots ,\xi _{r-2}\right) \in \Xi _r%
\text{, together with all their faces}
\end{array}
\right\} 
\]
and 
\[
\begin{array}{l}
\frak{T}_{i+1}^{\prime }:=\left[ \frak{T}_i^{\prime }\smallsetminus \left\{ 
\begin{array}{c}
\text{conv}\left( \left\{ n^{\left( i\right) },e_{\xi _1},\ldots ,e_{\xi
_{r-2}}\right\} \right) ,\smallskip \\ 
\text{for all }\left( \xi _1,\xi _2,\ldots ,\xi _{r-2}\right) \in \Xi _r%
\text{,\smallskip } \\ 
\text{together with its faces}
\end{array}
\right\} \right] \cup \\ 
\  \\ 
\cup \ \left\{ 
\begin{array}{c}
\mathbf{s}\left( i+1;\xi _1,\ldots ,\xi _{r-2}\right) ,\ \text{conv}\left(
\left\{ n^{\left( i+1\right) },e_{\xi _1},\ldots ,e_{\xi _{r-2}}\right\}
\right) ,\smallskip \\ 
\text{ for all \ }\left( \xi _1,\ldots ,\xi _{r-2}\right) \in \Xi _r\text{,\
together with all their faces}
\end{array}
\right\}
\end{array}
\]
for all $i$, $1\leq i\leq \nu -1.$

\begin{proposition}[Second factorization]
\label{2FACT} \ \smallskip \newline
Let $\left( \Bbb{C}^r/G,\left[ \mathbf{0}\right] \right) =\left( X\left(
N_G,\Delta _G\right) ,\text{\emph{orb}}\left( \sigma _0\right) \right) $ be
the Gorenstein cyclic quotient singularity of type \smallskip \emph{(\ref
{montype}) }with $l$ satisfying condition \emph{(\ref{rescon}). }Then the
birational resolution-morphism $f=f_{\mathcal{T}}$ can be expressed also as
the composite of $\nu =$ $\QOVERD\lfloor \rfloor {l}{r-1}$ toric morphisms 
\[
X_0:=X\left( N_G,\Delta _G\right) \stackrel{h_1}{\longleftarrow }X_1%
\stackrel{h_2}{\longleftarrow }X_2\longleftarrow \cdots \stackrel{h_{\nu -1}%
}{\longleftarrow }X_{\nu -1}\stackrel{h_\nu }{\longleftarrow }X_\nu
:=X\left( N_G,\widehat{\Delta }_G\left( \mathcal{T}\right) \right) 
\]
with 
\[
X_i:=X\left( N_G,\widehat{\Delta }_G\left( \frak{T}_i^{\,\prime }\right)
\right) ,\forall i,\ 0\leq i\leq \nu -1,\text{\emph{and\ }}\frak{T}%
_0^{\,\prime }=\text{\emph{trival triangulation,} }\frak{T}_\nu ^{\,\prime
}\ =\mathcal{T}.
\]
In particular, 
\[
X_{i+1}\cong \text{\emph{Norm}\thinspace }\left[ \mathbf{Bl}_{Z_i}^{\mathcal{%
I}_i}\left( X_i\right) \right] ,\ \ \ \forall i,\ \ 0\leq i\leq \nu -1\ ,
\]
for appropriate $\mathcal{O}_{X_i}$-ideal sheaves $\mathcal{I}_i$, such that 
$Z_i=$ \emph{supp}$\left( \mathcal{O}_{X_i}\,/\,\mathcal{I}_i\right) $,
where\smallskip 
\[
Z_i=\left\{ 
\begin{array}{lll}
\text{\emph{orb}}\left( \sigma _0\right)  & ,\ \ \text{\emph{if}} & \ \ i=0
\\ 
&  &  \\ 
\text{\emph{orb}}\left( \text{\emph{pos}}\left( \left\{ n^{\left( i\right)
},e_1,\ldots ,e_{r-1}\right\} \right) \right)  & ,\ \ \text{\emph{if}} & \
\left\{ 
\begin{array}{l}
\text{\emph{either }}\nu \geq 2,1\leq i\leq \nu -1,\smallskip  \\ 
\text{\emph{and }}l\equiv 1\ \text{\emph{mod}}\left( r-1\right)  \\ 
\\ 
\text{\emph{or }}\nu \geq 3,1\leq i\leq \nu -2,\smallskip  \\ 
\text{\emph{and }}l\equiv 0\ \text{\emph{mod}}\left( r-1\right) 
\end{array}
\right.  \\ 
&  &  \\ 
V\left( \text{\emph{pos}}\left( \left\{ e_1,e_2,\ldots ,e_{r-1}\right\}
\right) \right)  & ,\ \ \text{\emph{if}} & \ \ l\equiv 0\ \text{\emph{mod}}%
\left( r-1\right) \ \ \text{\emph{\&} \ }i=\nu -1
\end{array}
\right. \smallskip 
\]
\emph{(Up to the above last case and the case in which }$\nu =1$\emph{, all }%
$Z_i$\emph{'s are endowed with a non-reduced scheme structure).\smallskip }%
\newline
Moreover, one has 
\[
X_{i+1}\cong \text{\emph{Norm}\thinspace }\left[ \mathbf{Bl}_{\text{\emph{orb%
}}\left( \theta _i\right) }^{\text{\emph{red}}}\left( X\left( N_{\left[
i\right] }\emph{;\,}\text{\emph{Star}}\left( \theta _i;\widehat{\Delta }%
_G\left( \frak{T}_i^{\,\prime }\right) \right) \right) \right)
\,/\,G_i\right] ,\ \ \ 
\]
with respect to the different lattice 
\[
N_{\left[ i\right] }:=\Bbb{Z\,}n^{\left( i\right) }\oplus \Bbb{Z\,}e_1\oplus
\cdots \oplus \Bbb{Z\,}e_{r-1}\ ,
\]
and the star of the cone 
\[
\theta _i:=\text{\emph{pos}}\left( \left\{ n^{\left( i\right) },e_1,\ldots
,e_{r-1}\right\} \right) \text{\emph{, \ for all} }i\ \text{\emph{with}
\thinspace }\ \ \left\{ 
\begin{array}{l}
\text{\emph{either\ }}0\leq i\leq \nu -1,\smallskip  \\ 
\text{\emph{and }}l\equiv 1\ \text{\emph{mod}}\left( r-1\right)  \\ 
\\ 
\text{\emph{or \thinspace }}\nu \geq 2,\ 0\leq i\leq \nu -2,\smallskip  \\ 
\text{\emph{and }}l\equiv 0\ \text{\emph{mod}}\left( r-1\right) 
\end{array}
\right. \smallskip 
\]
where $G_0=G$ and $G_i=\left\langle g_i\right\rangle $ denotes the cyclic
group of analytic automorphisms of 
\[
X\left( N_{\left[ i\right] }\emph{;\,}\text{\emph{Star}}\left( \theta _i;%
\widehat{\Delta }_G\left( \frak{T}_i^{\,\prime }\right) \right) \right)
\cong \Bbb{C}^r
\]
generated by 
\[
g_i:\Bbb{C}^r\ni \left( z_1,..,z_r\right) \longmapsto \left( \zeta
_{l-i\left( r-1\right) }\cdot z_1,..,\zeta _{l-i\left( r-1\right) }\cdot
z_{r-1},\zeta _{l-i\left( r-1\right) }^{l-\left( i+1\right) \left(
r-1\right) }\cdot z_r\right) \in \Bbb{C}^r\,.\smallskip 
\]
\end{proposition}

\noindent The proof of \ref{2FACT} is an immediate generalization of that of
the case in which $r=2$ (see \ref{DIFFBL} (ii)), relies on a successive
application of \ref{TORNBL}, and is left as an exercise to the reader. (The
only difference is that whenever $r\geq 3$ and $l\equiv 0\ $mod$\left(
r-1\right) $, we also blow up the remaining $1$-dimensional singular locus
in the last step). Figures \textbf{12} (a) and (b) illustrate the
triangulations inducing the factorization of $f$ for the singularities of
type $1/6\,(1,1,4)$ and $1/7\,(1,1,5)$, respectively.

\begin{figure*}
\begin{center}
\input{fig12.pstex_t}
\makespace
 {Fig. \textbf{12}}
\end{center}
\end{figure*}

\vspace{8cm}
\begin{remark}
\emph{(i) Combining propositions \ref{1FACT} and \ref{2FACT} with corollary 
\ref{MBLU} one may obtain alternative proofs of the projectivity of }$f=f_{%
\mathcal{T}}$\emph{\ .\smallskip \newline
(ii) The proofs of \ref{1FACT} and \ref{2FACT} work after minor
modifications even if one omits the assumption for }$\mathcal{T}$ \emph{to
be basic. }
\end{remark}

\begin{exercise}
\emph{Under the assumption of prop. \ref{1FACT}, determine a} \textit{single}
$\mathcal{O}_{X\left( N_G,\Delta _G\right) }$\emph{-ideal sheaf }$\mathcal{I}
$\emph{, such that Sing}$\left( X\left( N_G,\Delta _G\right) \right) =$\emph{%
\ supp}$\left( \mathcal{O}_{X\left( N_G,\Delta _G\right) }/\mathcal{I}%
\right) $\emph{\ and }$f=f_{\mathcal{T}}$ \emph{\ itself is nothing but the
normalized blow-up:} 
\[
f:\text{\emph{Norm\thinspace }}\left[ \mathbf{Bl}_{\text{\emph{Sing}}\left(
X\left( N_G,\Delta _G\right) \right) }^{\mathcal{I}} \left(X\left( N_G,%
\Delta_G\left( \mathcal{T}\right) \right) \right) \right] \longrightarrow X\left(
N_G,\Delta _G\right) =\Bbb{C}^r/G
\]
\emph{of }$X\left( N_G,\Delta _G\right) $\emph{\ (cf. thm. \ref{RESSTR} and
rem. \ref{DIFFBL} (iii)).}
\end{exercise}

\section{Further remarks and a conjecture}

\noindent As we already saw in \ref{UNICIT}\thinspace (ii), in dimension
three, besides $\frac 1l\left( 1,1,l-2\right) $'s there is also another
``new'' Gorenstein, cyclic quotient singularity having a unique, projective,
crepant resolution, namely $\frac 17\left( 1,2,4\right) $. This can be
generalized in arbitrary dimensions too!

\begin{theorem}
\label{222}The cyclic Gorenstein quotient singularity of type 
\[
\fbox{$
\begin{array}{lll}
&  &  \\ 
& \dfrac 1{2^r-1}\ \,\left( 1,2,2^2,2^3,\ldots ,2^{r-2},2^{r-1}\right)  & 
\\ 
&  & 
\end{array}
$}
\]
can be fully resolved by a torus-equivariant projective crepant morphism in 
\textbf{all }dimensions $r\geq 2$. Moreover, up to isomorphism, this
resolution is \textbf{unique}.
\end{theorem}

\noindent The proof of theorem \ref{222} will be given in \cite{DHZ2}. As
you guess, the required triangulation $\mathcal{T}$ will be the
high-dimensional analogue of that of figure \textbf{9}. The details of the
proof of the uniqueness of $\mathcal{T}$ and of the fact that it is indeed
basic are somewhat lengthy, and involve binary representations, explicit
Hilbert-basis determination and some tricks with determinants. The coherence
of $\mathcal{T}$, on the other hand, can be shown directly by using tools
from the theory of polytopes, i.e., by avoiding both patching lemma and
factorization arguments.\smallskip \newline
Forgetting completely the uniqueness-property, we believe that this single
singularity is again nothing but ``the first member'' of an infinite family
of Gorenstein cyclic quotient singularities (\ref{kkk}), called for
simplicity $r$-dimensional \textit{geometric progress singularity-series} of
ratio $k$ (in notation: GPSS$\left( r;k\right) $), all of whose members
admit the desired resolutions.

\begin{conjecture}[GPSS$\left( r;k\right) $-Conjecture]
\label{GPCON}All cyclic Gorenstein quotient singularities of type 
\begin{equation}
\fbox{$
\begin{array}{lll}
&  &  \\ 
& \dfrac 1{\left( \dfrac{k^r-1}{k-1}\right) }\ \,\left( 1,k,k^2,k^3,\ldots
,k^{r-2},k^{r-1}\right)  &  \\ 
&  & 
\end{array}
$}  \label{kkk}
\end{equation}
admit torus-equivariant projective, crepant, full resolutions for \textbf{%
all }$r\geq 4$ and \textbf{all} $k\geq 2.$
\end{conjecture}

\begin{exercise}
\emph{As a first approach to \ref{GPCON} (e.g. to GPSS}$\left( 4;k\right) $%
\emph{-conjecture), consider the example }$\frac 1{40}\left( 1,3,9,27\right) 
$\emph{\ (}$\emph{with}$ $k=3$\emph{) and normalize the blow-up of }$X\left(
N_G,\Delta _G\right) $\emph{\ at orb}$\left( \sigma _0\right) $\emph{\
(equipped with the reduced structure). What kind of triangulation of the
junior tetrahedron }$\frak{s}_G$\emph{\ will be induced by this procedure ?
What would you expect as ``next step'' ? [Hint. Relate what you ``see'' with
the singularities being studied in \cite{DHZ1}.]\medskip }
\end{exercise}

\noindent \textit{Acknowledgements.} The first author would like to express
his thanks to I. Nakamura and M. Reid for having pointed out that the
existence of projective crepant resolutions of $\Bbb{C}^r/G$ in dimensions $%
r\geq 4$ does not necessarily imply the smoothness of the corresponding
Hilbert scheme of $G$-orbits, to V. Batyrev for a discussion about the
structure of the Chow ring $A^{\bullet }\left( \Bbb{P}\left( \mathcal{E}%
\right) \right) $, to G. M. Ziegler for an e-mail about the counterexample 
\ref{Hibi}\thinspace (ii), to DFG for the support by an
one-year-research-fellowship, to Mathematics Institute of Bonn
University for hospitality until the end of March, 1997, and to Mathematics Department of T\"ubingen University, where this work was completed.\smallskip \newline
The second author thanks R. Firla for explaining to him his counterexample 
\ref{FIRLA}\thinspace (i), DFG for the support by the Leibniz-Preis awarded
to M.Gr\"{o}tschel, and Konrad-Zuse-Zentrum in Berlin for ideal working
conditions during the writing of this paper.

\vspace{2cm}

\hbox{
\vbox{\noindent Dimitrios I.~Dais \\
\noindent Mathematisches Institut \\ 
\noindent Universit{\"a}t T{\"u}bingen \\
\noindent Auf der Morgenstelle 10 \\
\noindent 72076 T{\"u}bingen, Germany \\
\noindent dais@wolga.mathematik.uni-tuebingen.de}
\hspace{-5cm}
\vbox{ \noindent Martin Henk \\
\noindent Konrad-Zuse-Zentrum Berlin \\
\noindent Takustra{\ss}e 7 \\
\noindent 14195 Berlin, Germany \\
\noindent henk@zib.de \\
\noindent \hphantom{space}}
}


\begin{thebibliography}{99}
\bibitem{Bat}  \textsc{Batyrev V.V.}: \textit{On the classification of
smooth projective toric varieties}, T\^{o}hoku Math. Jour. \textbf{43},
(1991), 569-585.


\bibitem{Bat2}  \textsc{Batyrev V.V.}: \textit{Non-Archimedian integrals and stringy Euler numbers of
 log-terminal pairs}, preprint, alg-geom / 9803071.

\bibitem{BD}  \textsc{Batyrev V.V., Dais D.I.}: \textit{Strong McKay
correspondence, string theoretic Hodge numbers and mirror symmetry},
Topology \textbf{35}, (1996), 901-929.


\bibitem{BM1}  \textsc{Bierstone E., Milman P.D.}: \textit{A simple
constructive proof of canonical resolution of singularities.} In:
``Effective Methods in Algebraic Geometry'', (edited by T.Mora, C.Traverso),
Progress in Math., Vol. \textbf{94}, Birkh\"{a}user, (1991), 11-30.


\bibitem{BM2}  \textsc{Bierstone E., Milman P.D.}: \textit{Canonical
desingularization in characteristic zero by blowing up the maximum strata of
a local invariant}, Inventiones Math. \textbf{128}, (1997), 207-302.


\bibitem{BFS}  \textsc{Billera L.J., Filliman P., Sturmfels B.}: \textit{%
Construction and complexity of secondary polytopes,} Adv. in Math. \textbf{83%
}, (1990), 155-179.


\bibitem{Br1}  \textsc{Brieskorn E.}: \textit{Holomorphe }$\Bbb{P}^n$\textit{%
-B\"{u}ndel \"{u}ber }$\Bbb{P}^1$, Math. Ann. \textbf{157}, (1965), 343-357.


\bibitem{Br2}  \textsc{Brieskorn E.}: \textit{Rationale Singularit\"{a}ten
komplexer Fl\"{a}chen}, Inventiones Math. \textbf{4}, (1968), 336-358.


\bibitem{BGT}  \textsc{Bruns W., Gubeladge J., Trung N.V.}: \textit{Normal
polytopes, triangulations and Koszul algebras}, Jour. f\"{u}r die reine und
ang. Math. \textbf{485}, (1997), 123-160.


\bibitem{Bryl}  \textsc{Brylinski J.L.}: \textit{Eventails et vari\'{e}t\'{e}%
s toriques}. In : ``S\'{e}minaire sur les Singularit\'{e}s des Surfaces'',
(edited by M.Demazure, H.Pinkham \& B.Teissier), Lecture Notes in Math.,
Vol. \textbf{777}, Springer-Verlag, (1980), pp. 247-288.


\bibitem{DHZ1}  \textsc{Dais D.I., Henk M., Ziegler G.M.}: \textit{All
abelian quotient c.i.-singularities admit projective crepant resolutions in
all dimensions}, preprint, alg-geom / 9704007.


\bibitem{DHZ2}  \textsc{Dais D.I., Henk M., Ziegler G.M.}: \textit{On the
existence of crepant resolutions of Gorenstein abelian quotient
singularities in dimensions }$\geq 4$, in preparation.


\bibitem{Danilov}  \textsc{Danilov V.I.}: \textit{The birational geometry of
toric }$3$\textit{-folds}, Math. USSR Izvestiya \textbf{21}, (1983), 269-280.


\bibitem{DuVal1}  \textsc{Du Val P.}: \textit{On the singularities which do
not affect the condition of adjunction I, II, III,} Proc. Camb. Phil. Soc. 
\textbf{30}, (1934), 453-459 \& 483-491.


\bibitem{DuVal2}  \textsc{Du Val P.}: \textit{Homographies, Quaternions and
Rotations}, Oxford at the Clarendon Press, (1964).


\bibitem{Ehlers}  \textsc{Ehlers F.}: \textit{Eine Klasse komplexer
Mannigfaltigkeiten und die Aufl\"{o}sung einiger isolierter Singularit\"{a}%
ten}, Math. Ann. \textbf{218}, (1975), 127-156.


\bibitem{Ewald}  \textsc{Ewald G.}: \textit{Combinatorial Convexity and
Algebraic Geometry}, Graduate Texts in Mathematics, Vol. \textbf{168},
Springer-Verlag, (1996).\textit{\ }


\bibitem{Ew-Sch}  \textsc{Ewald G., Schmeink A.}: \textit{Representation of
the Hirzebruch-Kleinschmidt varieties by quadrics}, Beitr\"{a}ge zur Algebra
und Geometrie, Vol. \textbf{34}, (1993), 151-156.


\bibitem{Firla}  \textsc{Firla R.T.}: \textit{Hilbert-Cover- \ und
Hilbert-Partitions-Probleme}, Diplomarbeit, TU-Berlin, (1997).

\bibitem{Fir-Zi}  \textsc{Firla R.T., Ziegler G.M.}: \textit{Hilbert bases,
unimodular triangulations, and binary covers of rational polyhedral cones},
preprint, (1997); to appear in Discrete \ \& Comp. Geom.


\bibitem{Fujiki}  \textsc{Fujiki A.}: \textit{On resolutions of cyclic
quotient singularities}, Publ. RIMS \textbf{10}, (1974), 293-328.


\bibitem{Ful}  \textsc{Fulton W.}: \textit{Intersection Theory, }Ergebnisse
der Mathematik und ihrer Grenzgebiete, 3. Folge, Bd. \textbf{2},
Springer-Verlag, (1984). \textit{\ }


\bibitem{Fulton}  \textsc{Fulton W.}: \textit{Introduction to Toric Varieties%
}, Annals of Math. Studies, Vol. \textbf{131}, Princeton University Press,
(1993).


\bibitem{GKZ}  \textsc{Gelfand I.M., Kapranov M.M., Zelevinsky A.V.}: 
\textit{Discriminants, Resultants and Multidimensional Determinants. }%
Mathematics: Theory \& Applications. (Ed. V.Kadison, I.M.Singer),
Birkh\"{a}user, (1994).


\bibitem{GSV}  \textsc{Gonzalez-Sprinberg G., Verdier J.L.}: \textit{%
Construction g\'{e}om\'{e}trique de la correspondance de McKay}, Ann. Sc.
E.N.S. \textbf{16}, (1983), 409-449.


\bibitem{Gordan}  \textsc{Gordan P.}: \textit{\"{U}ber die Aufl\"{o}sung
linearer Gleichungen mit reellen Coefficienten}, Math. Ann.\textit{\ }%
\textbf{6}, (1873), 23-28.


\bibitem{Hart}  \textsc{Hartshorne R.}: \textit{Algebraic Geometry},
Graduate Texts in Mathematics, Vol. \textbf{52}, Springer-Verlag, (1977).


\bibitem{HEW}  \textsc{Henk M., Weismantel R.}: \textit{On Hilbert bases of
polyhedral cones}, preprint, (Konrad-Zuse-Zentrum, Berlin), SC 96-12,
(1996); available from:\smallskip \newline
http://www.zib-berlin.de/paperweb/paperweb?query$=$Henk\&mode$=$all


\bibitem{HO}  \textsc{Hibi T., Ohsugi H.}: \textit{A normal }$\left(
0,1\right) $\textit{-polytope none of whose regular triangulations is
unimodular}, preprint, (1997); to appear in Discrete \ \& Comp. Geom.


\bibitem{Hilbert}  \textsc{Hilbert D.}: \textit{\"{U}ber die Theorie der
algebraischen Formen}, Math. Ann. \textbf{36}, (1890), 473-531. [See also:
``Gesammelte Abhandlungen'', Band \textbf{II}, Springer-Verlag, (1933),
119-257.]


\bibitem{Hiro}  \textsc{Hironaka H.}: \textit{Resolution of singularities of
an algebraic variety over a field of characteristic zero I, II}, Annals of
Math. \textbf{79}, (1964), 109-326.


\bibitem{Hirz1}  \textsc{Hirzebruch F.}: \textit{\"{U}ber eine Klasse von
einfach zusammenh\"{a}ngenden komplexen Mannigfaltigkeiten}, Math. Ann. 
\textbf{124}, (1951), 77-86. [See also: ``Gesammelte Abhandlungen'', Band 
\textbf{I}, Springer-Verlag, (1987), pp. 1-10.]


\bibitem{Hirz2}  \textsc{Hirzebruch F.}: \textit{\"{U}ber vierdimensionale
Riemannsche Fl\"{a}chen mehrdeutiger analytischer Funktionen von zwei
komplexen Ver\"{a}nderlichen}, Math. Ann. \textbf{126}, (1953), 1-22. [See
also: ``Gesammelte Abhandlungen'', Band \textbf{I}, Springer-Verlag, (1987),
pp. 11-32.]


\bibitem{HH}  \textsc{Hirzebruch F., H}\"{o}\textsc{fer T.}: \textit{On the
Euler number of an orbifold}, Math. Ann.\textbf{\ 286}, (1990), 255-260.


\bibitem{Ito1}  \textsc{Ito Y.}: \textit{Crepant resolution of trihedral
singularities}, Proc. of Japan Acad., Ser. A, Vol. \textbf{70}, (1994),
131-136.


\bibitem{Ito2}  \textsc{Ito Y.}: \textit{Crepant resolution of trihedral
singularities and the orbifold Euler characteristic}, Intern. Jour. of Math. 
\textbf{6}, (1995), 33-43.


\bibitem{Ito3}  \textsc{Ito Y.}: \textit{Gorenstein quotient singularities
of monomial type in dimension three}, Jour. of Math. Sciences, University of
Tokyo, Vol. \textbf{2}, Nr. 2, (1995), 419-440.


\bibitem{Ito-Nak1}  \textsc{Ito Y., Nakamura I.}: \textit{McKay
correspondence and Hilbert schemes}, Proc. Japan Acad. \textbf{72}, (1996),
135-138.


\bibitem{Ito-Nak2}  \textsc{Ito Y., Nakamura I.}: \textit{Hilbert schemes
and simple singularities }$A_n$ \textit{and }$D_n$, Hokkaido Univ. preprint
series, Nr. \textbf{348}, (1996).


\bibitem{Ito-Nak3}  \textsc{Ito Y., Nakamura I.}: \textit{Hilbert schemes
and simple singularities}, preprint, (1997).


\bibitem{Ito-Reid}  \textsc{Ito Y., Reid M.}: \textit{The McKay
correspondence for finite subgroups of }SL$\left( 3,\Bbb{C}\right) $. In :
``Higher Dimensional Complex Varieties'', Proceedings of the International
Conference held in Trento, Italy, June 15-24, 1994, (edited by M.Andreatta,
Th.Peternell), Walter de Gruyter, (1996), 221-240.


\bibitem{KKMS}  \textsc{Kempf G., Knudsen F., Mumford D., Saint-Donat D.}: 
\textit{Toroidal Embeddings I}, Lecture Notes in Mathematics, Vol. \textbf{%
339}, Springer-Verlag, (1973).


\bibitem{Kl}  \textsc{Klein F}.: \textit{Vorlesungen \"{u}ber das Ikosaeder}%
. Erste Ausg. Teubner Verlag, (1884). Zweite Ausg. von Teubner und
Birkh\"{a}user mit einer Einf\"{u}hrung und mit Kommentaren von P. Slodowy,
(1993).


\bibitem{Klein}  \textsc{Kleinschmidt P.}: \textit{A classification of toric
varieties with few generators}, Aequationes Mathematicae \textbf{35},
(1988), 254-266.


\bibitem{Kn}  \textsc{Kn}\"{o}\textsc{rrer H.}:\textit{\ Group
representations and the resolution of rational double points}, Contemp.
Math. \textbf{45}, A.M.S., (1985), 175-222.


\bibitem{Lamotke}  \textsc{Lamotke K.}: \textit{Regular Solids and Isolated
Singularities}, Vieweg Adv. Lectures in Math., (1986).


\bibitem{Mark}  \textsc{Markushevich D.G.}: \textit{Resolution of
singularities }(\textit{Toric Method}). Appendix of the article:\smallskip 
\newline
\textsc{Markushevich D.G., Olshanetsky M.A., Perelomov A.M.}: \textit{%
Description of a class of superstring compactifications related to
semi-simple Lie algebras}, Commun. in Math. Phys. \textbf{111}, (1987),
247-274.


\bibitem{Markush}  \textsc{Markushevich D.G.}: \textit{Resolution of }$\Bbb{C%
}^3/H_{168}$, Math. Ann. \textbf{308}, (1997), 279-289.


\bibitem{McK}  \textsc{McKay J.}: \textit{Graphs, singularities and finite
groups}. In : ``The Santa-Cruz Conference of Finite Groups'', Proc. of Symp.
in Pure Math., A.M.S., Vol. \textbf{37}, (1980), 183-186.


\bibitem{Mo-Ste}  \textsc{Morrison D.R., Stevens G.}: \textit{Terminal
quotient singularities in dimension three and four}, Proc. A.M.S.\textbf{\ 90%
}, (1984), 15-20.


\bibitem{Nak1}  \textsc{Nakamura I.}: \textit{Simple singularities, McKay
correspondence and Hilbert-schemes of }$G$\textit{-orbits}, preprint, (1996).


\bibitem{Nak2}  \textsc{Nakamura I.}: \textit{Hilbert schemes and simple
singularities }$E_6,E_7$ \textit{and }$E_8$, Hokkaido Univ. preprint series,
Nr. \textbf{362}, (1996).


\bibitem{Nak3}  \textsc{Nakamura I.}: \textit{Hilbert schemes of }$G$\textit{%
-orbits for abelian }$G$\textit{, }(in preparation).


\bibitem{Oda}  \textsc{Oda T.}:\textit{\ Convex Bodies and Algebraic
Geometry. An Introduction to the theory of toric varieties}, Ergebnisse der
Mathematik und ihrer Grenzgebiete, 3. Folge, Bd. \textbf{15},
Springer-Verlag, (1988).


\bibitem{OP}  \textsc{Oda T., Park H.S.}: \textit{Linear Gale Transforms and
GKZ-Decompositions}, T\^{o}hoku Math. Jour. \textbf{43}, (1991), 375-399.


\bibitem{Pottier1}  \textsc{Pottier L.}: \textit{Subgroups of }$\Bbb{Z}^n$%
\textit{, standard bases, and linear diophantine systems}, Research Report 
\textbf{1510}, Universit\'{e} de Nice, Sophia Antipolis, Valbonne, France,
(1992).


\bibitem{Pottier2}  \textsc{Pottier L.}: \textit{Euclid's algorithm in
dimension }$n$, preprint, Universit\'{e} de Nice, Sophia Antipolis,
Valbonne, France, (1996).


\bibitem{Prill}  \textsc{Prill D.}: \textit{Local classification of
quotients of complex manifolds by discontinuous groups}, Duke Math. Jour. 
\textbf{34}, (1967), 375-386.


\bibitem{Reid1}  \textsc{Reid M.}: \textit{Canonical threefolds}, Journ\'{e}%
e de G\'{e}om\'{e}trie Alg\'{e}brique d'Angers, A. Beauville ed., Sijthoff
and Noordhoff, Alphen aan den Rijn, (1980), 273-310.


\bibitem{Reid2}  \textsc{Reid M.}: \textit{Decompositions of toric morphisms}%
. In : ``Arithmetic and Geometry II'', (edited by M.Artin and J.Tate),
Progress in Math. \textbf{36}, Birkh\"{a}user, (1983), 395-418.


\bibitem{Reid3}  \textsc{Reid M.}: \textit{Young person's guide to canonical
singularities.} In : ``Algebraic Geometry, Bowdoin 1985'', (edited by
S.J.Bloch), Proc. of Symp. in Pure Math., A.M.S., Vol. \textbf{46}, Part I,
(1987), 345-416.


\bibitem{Reid4}  \textsc{Reid M.}: \textit{The McKay correspondence and the
physicist's Euler number}. Notes from lectures given at the Utah University
(Sept. 1992) and at MSRI (Nov. 1992).


\bibitem{Reid5}  \textsc{Reid M.}: \textit{\ McKay correspondence},
preprint, alg-geom / 9702016.


\bibitem{Reid6}  \textsc{Reid M.}: \textit{\ Rational Scrolls}, Ch. 2 of
author's article : \textit{Chapters on Algebraic Surfaces}. In : ``Complex
Algebraic Geometry'' (J.Koll\'{a}r, ed.), IAS/Park City Mathematics Series,
Volume \textbf{3}, AMS, (1997), pp. 5-159.


\bibitem{Roan1}  \textsc{Roan S.-S.}:\textit{\ On the generalization of
Kummer surfaces}, Jour. of Differential Geometry \textbf{30}, (1989),
523-537.


\bibitem{Roan2}  \textsc{Roan S.-S.}:\textit{\ On }$c_1=0$ \textit{%
resolution of quotient singularities}, Intern. Jour. Math. \textbf{5},
(1994), 523-536.


\bibitem{Roan3}  \textsc{Roan S.-S.}:\textit{\ Minimal resolutions of
Gorenstein orbifolds in dimension three}, Topology \textbf{35}, (1996),
489-508.


\bibitem{Roan4}  \textsc{Roan S.-S.}:\textit{\ Orbifold Euler characteristic}%
. In : ``Mirror Symmetry II'' (B.Greene \& S.-T.Yau, eds.), AMS /\thinspace
I\thinspace P Studies in Adv. Math. \textbf{1, }(1997), pp. 129-140.


\bibitem{Roan-Yau}  \textsc{Roan S.-S., Yau S.-T.}:\textit{\ On Ricci flat
3-folds}, Acta Mathematica Sinica, New Series, (Chinese Jour. of Math.),
Vol. \textbf{3}, Nr. 3, (1987), 256-288.


\bibitem{Schrijver}  \textsc{Schrijver A.}: \textit{Theory of Linear and
Integer Programming,} Wiley-Intersience Series in Discrete Mathematics and
Optimization, (Adv.Eds. R.L.Graham, J.K.Lenstra), John Wiley \& Sons,
sec.cor. repr., (1989).


\bibitem{Segre}  \textsc{Segre C.}: \textit{Recherches g\'{e}n\'{e}rales sur
les courbes et les surfaces r\'{e}gl\'{e}es alg\'{e}briques}, Math. Ann. 
\textbf{30,} (1887), 203-226; ibid. \textbf{34}, (1889), 1-25.


\bibitem{Slodowy}  \textsc{Slodowy P.}: \textit{Algebraic groups and
resolutions of Kleinian singularities}, R.I.M.S. Series, Nr. \textbf{1086},
Kyoto University, (1996), 1-80.


\bibitem{Sturm}  \textsc{Sturmfels B.}: \textit{Gr\"{o}bner Bases and Convex
Polytopes}, University Lecture Series, Vol. \textbf{8}, AMS, (1996).


\bibitem{Ueno1}  \textsc{Ueno K.}: \textit{On fibre spaces of normally
polarized abelian varieties of dimension }$2$, I. (\textit{Singular fibers
of the first type}); Jour. Fac. Sc. Univ. Tokyo, Sect. IA, Vol. \textbf{18},
(1971), 37-95.


\bibitem{Ueno2}  \textsc{Ueno K.}: \textit{The canonical resolution of
certain quotient singularities. }Appendix of the author's paper:\textit{\
Classification theory of algebraic varieties}\textbf{\ }I, Compositio Math. 
\textbf{27}, (1973), 277-342.


\bibitem{Ueno3}  \textsc{Ueno K.}: \textit{Classification Theory of
Algebraic Varieties and Compact Complex Spaces,} Lecture Notes in
Mathematics, Vol. \textbf{\ 439}, Springer-Verlag, (1975).


\bibitem{vanderCorput1}  \textsc{van der Corput J.G.}: \textit{\"{U}ber
Systeme von linear-homogenen Gleichungen und Ungleichungen}, Koninklijke
Akademie van Wetenschappen te Amsterdam \textbf{34}, (1931), 368-371.


\bibitem{vanderCorput2}  \textsc{van der Corput J.G.}: \textit{\"{U}ber
Diophantische Systeme von linear-homogenen Gleichungen und Ungleichungen},
Koninklijke Akademie van Wetenschappen te Amsterdam \textbf{34}, (1931),
372-382.


\bibitem{Wat1}  \textsc{Watanabe K.}: \textit{Certain invariant subrings are
Gorenstein I,II}, Osaka Jour. Math. \textbf{11}, (1974), 1-8 and 379-388.


\bibitem{Wat2}  \textsc{Watanabe K.}: \textit{Invariant subrings which are
complete intersections, I, (Invariant subrings of finite Abelian subgroups),}
Nagoya Math. Jour. \textbf{77}, (1980), 89-98.


\bibitem{Wessels}  \textsc{Wessels U.}: \textit{Kombinatorisch-geometrische
Kennzeichnung und Berechnung der Schnittzahlen von Cartierdivisoren in
kompakten torischen Variet\"{a}ten}, Dissertation, Ruhr-Universit\"{a}t
Bochum, (1993).
\end{thebibliography}
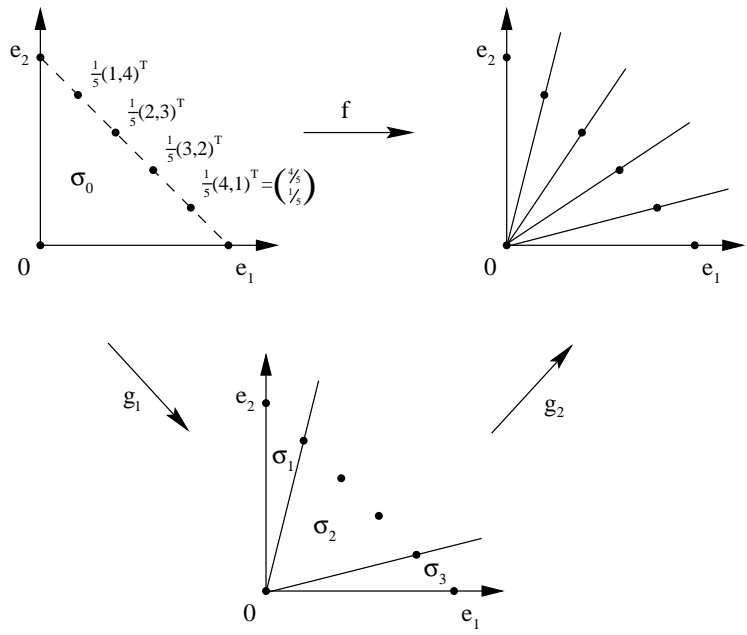
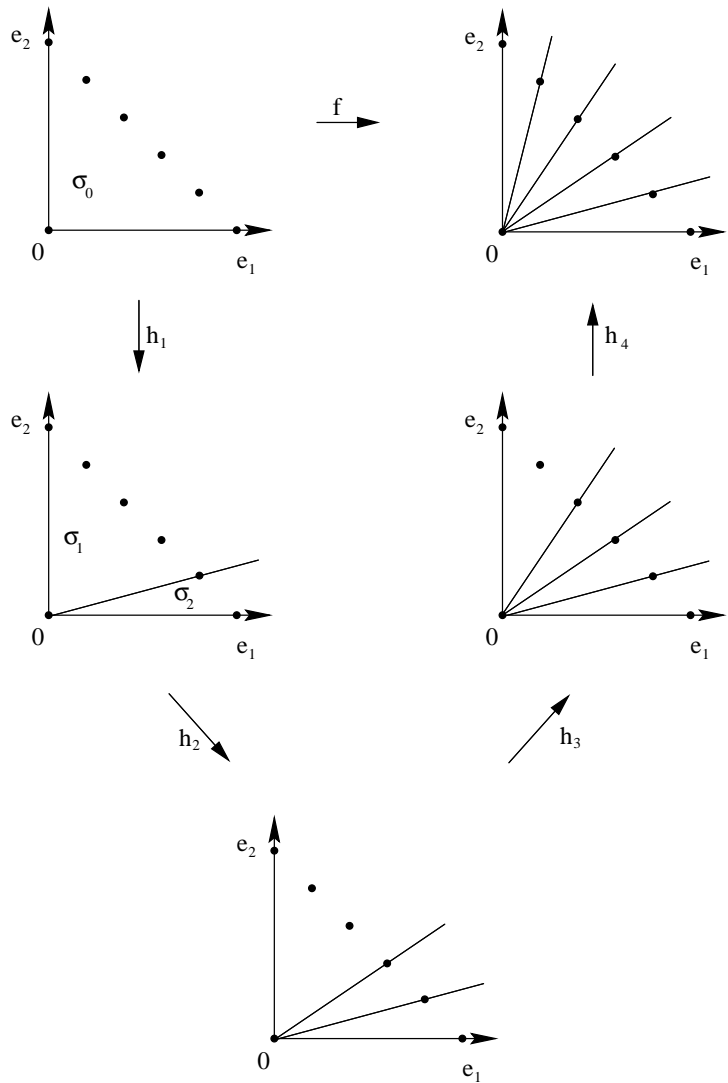
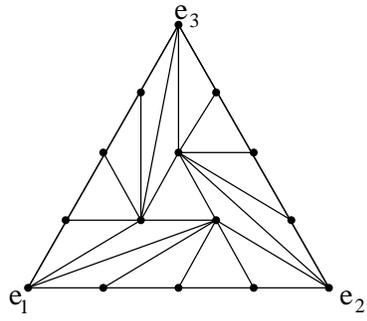
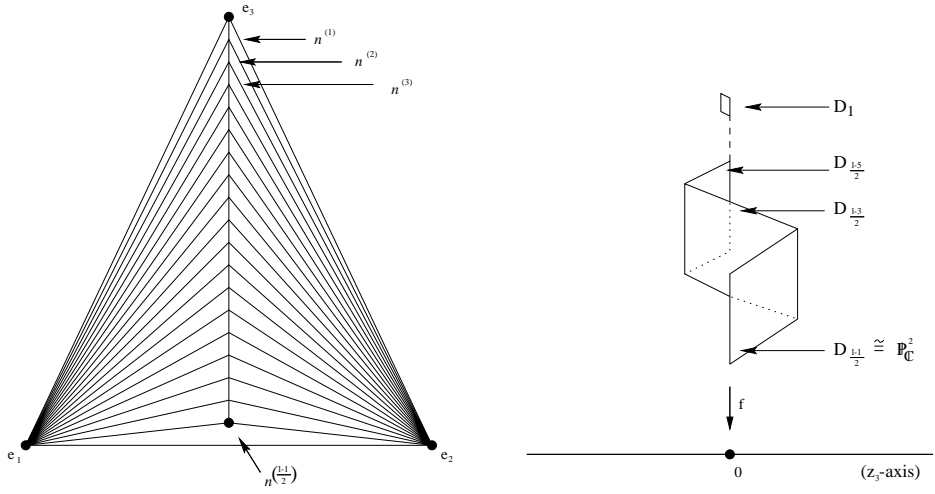
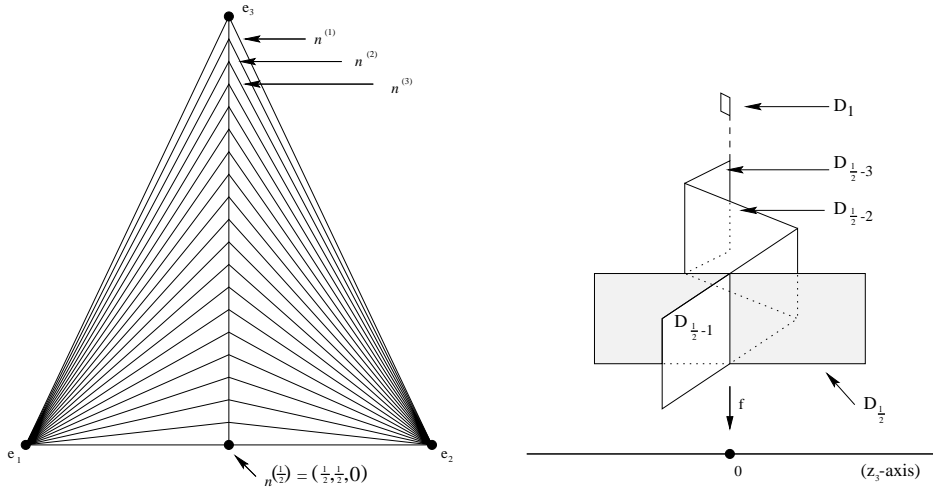
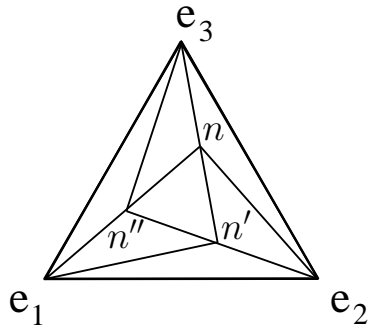
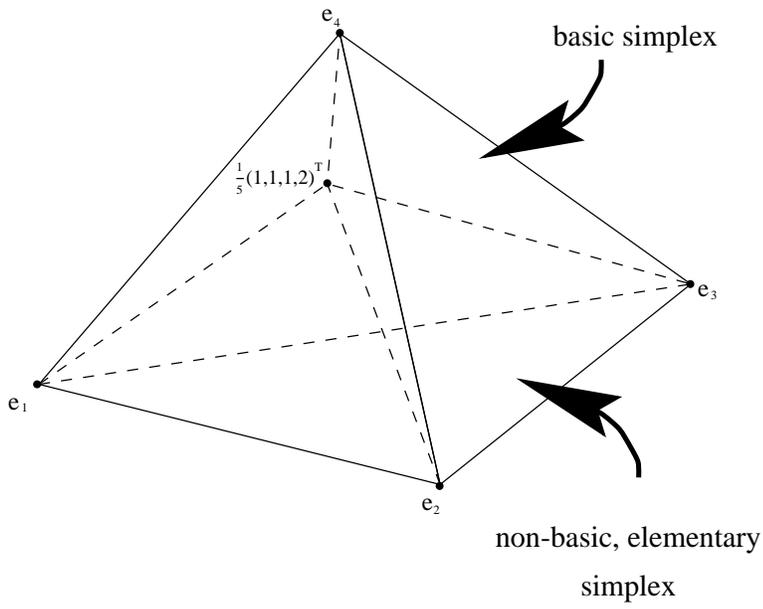
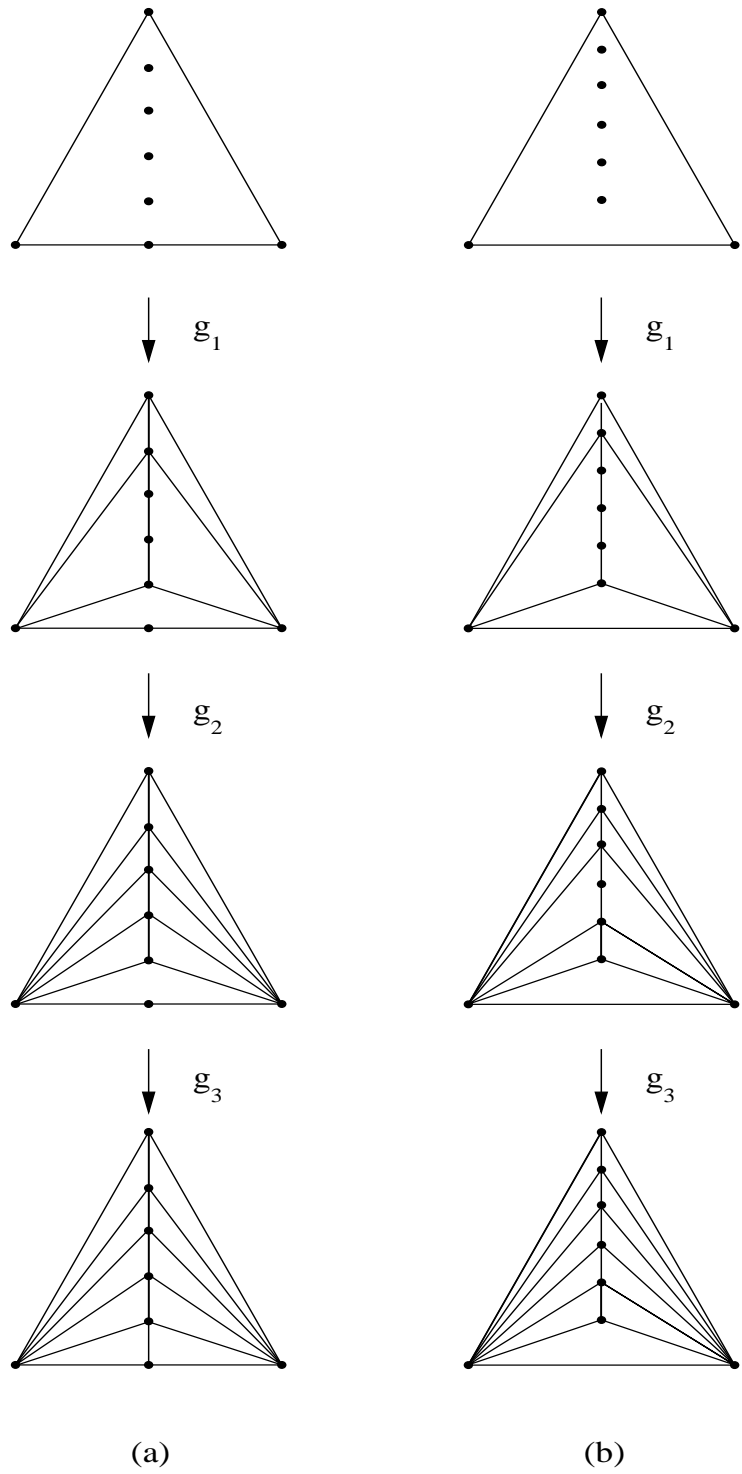
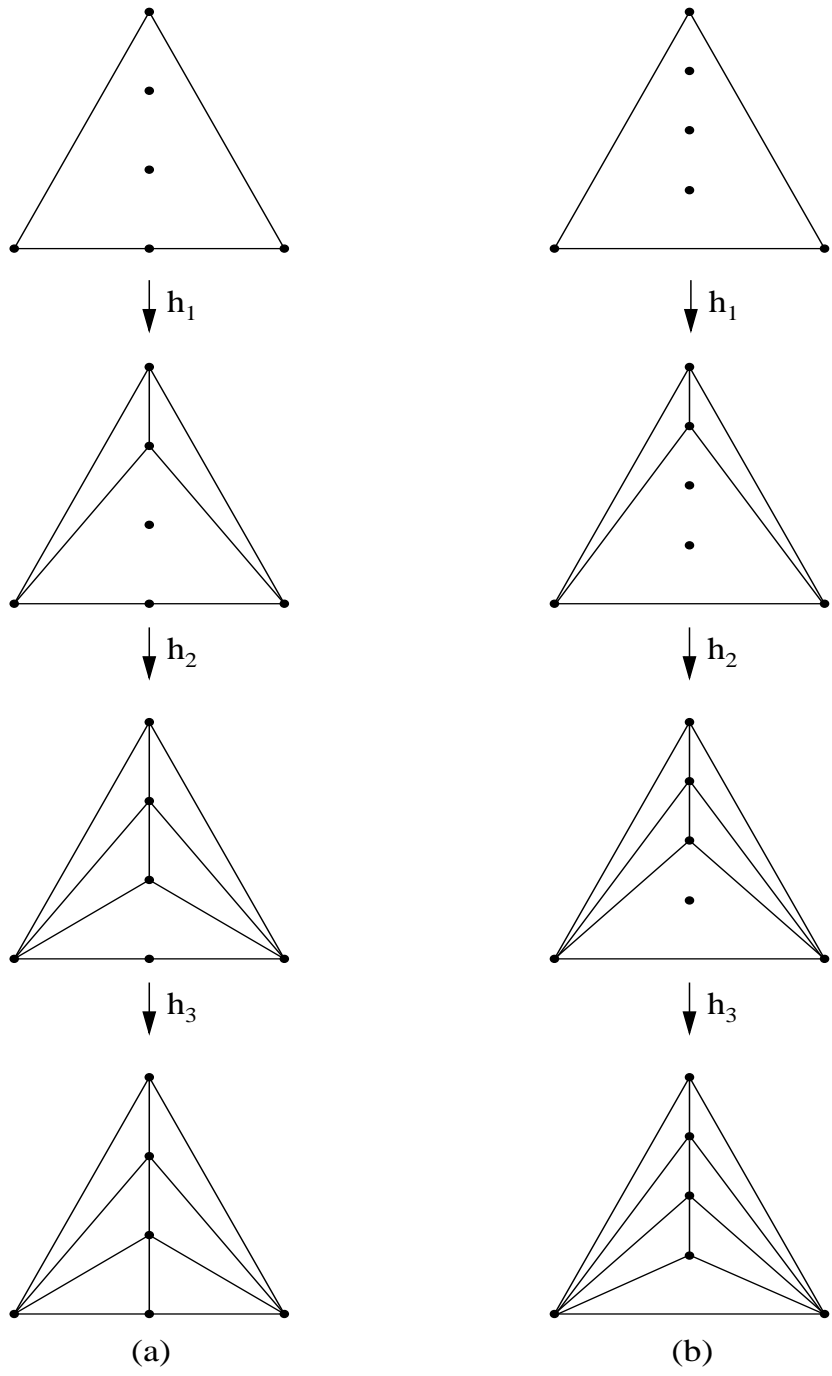
\end{document}